%% file: Main.tex
\pdfoutput=1
\documentclass[reqno]{gokova}

\usepackage{amsmath}
\usepackage{amssymb}
\usepackage{MnSymbol}
\usepackage{exercise}

\usepackage{amsfonts}
\input xy
\xyoption{all}

\usepackage{graphicx}           
\usepackage{boxedminipage}
\usepackage{color}
\usepackage{epstopdf}

\DeclareMathOperator{\trop}{trop}
\newcommand{\lt}{\lim\limits_{t\to+\infty}\hspace{-5pt}\trop\hspace{2pt}}

\DeclareMathOperator{\Log}{Log}
\DeclareMathOperator{\codim}{codim}

\begin{document}
\input goksty.tex      
\setcounter{page}{1}
\volume{13}

\input Define.tex

\title{Brief introduction to tropical geometry} 

\thanks{
Research is supported in part by
the
FRG Collaborative Research grant 
DMS-1265228 of the U.S. National Science Foundation (I.I.);
the project TROPGEO of the European Research Council, 
the grants 140666 and 141329 of the Swiss National Science Foundation,
the National Center of Competence in Research SwissMAP of the Swiss National Science Foundation (G.M.);
and an  Alexander von Humboldt Foundation Postdoctoral Research Fellowship (K.S.).}

\author{Erwan Brugall\'e} 
\address{Erwan Brugall\'e, \'Ecole Polytechnique, Centre de Math\'ematiques Laurent Schwartz, 91 128 Palaiseau Cedex, France}  
\email{erwan.brugalle@math.cnrs.fr} 

\author{Ilia Itenberg}
\address{Ilia Itenberg, Institut de Math\'ematiques de Jussieu - Paris Rive Gauche,
Universit\'e Pierre et Marie Curie, 
4 place Jussieu,
75252 Paris Cedex 5,
France \hskip5pt and \hskip5pt 
D\'epartement de Math\'emati\-ques et Applications, 
Ecole Normale Sup\'erieure, 
45 rue d'Ulm, 75230 Paris Cedex 5, France}
\email{ilia.itenberg@imj-prg.fr} 

\author{Grigory Mikhalkin}
\address{Grigory Mikhalkin,
Universit\'e de Gen\`eve, Section de Math\'ematiques, Villa Battelle,
7 route de Drize, 1227 Carouge, Switzerland.}
\email{grigory.mikhalkin@unige.ch}

\author{Kristin Shaw}
\address{Kristin Shaw, 
Technische Universit\"at
Berlin, MA 6-2, 10623 Berlin, Germany.}
\email{shaw@math.tu-berlin.de} 


\begin{abstract}
The paper consists of lecture notes for a mini-course given by the authors
at the G\"okova Geometry \& Topology conference in May 2014.
We start the exposition with tropical curves in the plane and their applications
to problems in classical enumerative geometry, and continue with a  look at more general
tropical varieties and their homology theories. 
\end{abstract}

\maketitle

\tableofcontents

The goal of these lectures is to give a basic
introduction to tropical
geometry
focusing
on some of its particularly simple and visual aspects.
The first section is devoted to tropical arithmetic
and its relations to classical arithmetic.
The second section reviews tropical curves in $\R^2$.
The content of these two sections is quite standard,
and we refer to \cite{Br11,Br29,BrSh13} for their 
extended versions.
Section 3 contains a tropical version of the combinatorial patchworking
construction
for plane curves,
 as well as a tropical reformulation of Haas' theorem. 
Section 4 presents some enumerative applications of tropical geometry, as
well as the floor diagram technique.
Section 5 looks at 
general tropical subvarieties of $\RR^n$ and their approximation by complex
algebraic varieties. Section 6 is devoted to a basic study
of tropical curves inside non-singular affine tropical surfaces.
Finally, in Section 7 we 
define abstract tropical manifolds and review their homology theories.

For other introductions to tropical geometry 
one can look, 
for
example, at 
 \cite{BIT,St2,Mik-It-11,V11,V12,Gath1,St7} and references therein.
A more advanced reader may refer to \cite{Mik3,Mik8,Mik9}.

\bigskip
{\bf Acknowledgements: }The authors are grateful to Matteo Ruggiero for providing helpful comments on a preliminary version of these notes.

\input{Erwan.tex}

\input{Ilia.tex}

\input{Kristin.tex}

\input{Kristin2.tex}

\input{Grisha.tex}

\small

\def\rightmark{\em Bibliography}

\bibliographystyle{alpha}

\bibliography{Biblio.bib}

\end{document}

%% file: goksty.tex
\def\E{\ifmmode{\mathbb E}\else{$\mathbb E$}\fi} 
\def\N{\ifmmode{\mathbb N}\else{$\mathbb N$}\fi} 
\def\R{\ifmmode{\mathbb R}\else{$\mathbb R$}\fi} 
\def\Q{\ifmmode{\mathbb Q}\else{$\mathbb Q$}\fi} 
\def\C{\ifmmode{\mathbb C}\else{$\mathbb C$}\fi} 
\def\H{\ifmmode{\mathbb H}\else{$\mathbb H$}\fi} 
\def\Z{\ifmmode{\mathbb Z}\else{$\mathbb Z$}\fi} 
\def\P{\ifmmode{\mathbb P}\else{$\mathbb P$}\fi} 
\def\T{\ifmmode{\mathbb T}\else{$\mathbb T$}\fi} 
\def\SS{\ifmmode{\mathbb S}\else{$\mathbb S$}\fi} 
\def\DD{\ifmmode{\mathbb D}\else{$\mathbb D$}\fi} 

\renewcommand{\a}{\alpha}
\renewcommand{\b}{\beta}
\renewcommand{\d}{\delta}
\newcommand{\D}{\Delta}
\newcommand{\e}{\varepsilon}
\newcommand{\g}{\gamma}
\newcommand{\G}{\Gamma}
\newcommand{\la}{\lambda}
\newcommand{\La}{\Lambda}
\newcommand{\n}{\nabla}
\newcommand{\var}{\varphi}
\newcommand{\s}{\sigma}
\newcommand{\Sig}{\Sigma}
\renewcommand{\t}{\tau}
\renewcommand{\th}{\theta}
\renewcommand{\O}{\Omega}
\renewcommand{\o}{\omega}
\newcommand{\z}{\zeta}

\newcommand{\ben}{\begin{enumerate}}
\newcommand{\een}{\end{enumerate}}
\newcommand{\be}{\begin{equation}}
\newcommand{\ee}{\end{equation}}
\newcommand{\bea}{\begin{eqnarray}}
\newcommand{\eea}{\end{eqnarray}}
\newcommand{\bc}{\begin{center}}
\newcommand{\ec}{\end{center}}

\newtheorem{thm}{Theorem}[section]
\newtheorem{cor}[thm]{Corollary}
\newtheorem{lem}[thm]{Lemma}
\newtheorem{prop}[thm]{Proposition}
\newtheorem{ax}{Axiom}
\newtheorem{conj}[thm]{Conjecture}

\theoremstyle{definition}
\newtheorem{defn}[thm]{Definition}

\theoremstyle{remark}
\newtheorem{rem}[thm]{\rm\bfseries{Remark}}
\newtheorem*{notation}{Notation}

\newtheorem{ques}[thm]{\rm\bfseries{Question}}
\newtheorem{cons}[thm]{\rm\bfseries{Construction}}
\newtheorem{exm}[thm]{\rm\bfseries{Example}}


%% file: Define.tex
\newcommand{\ZZ}{{\mathbb Z}}
\renewcommand{\Z}{{\mathbb Z}}
\newcommand{\NN}{{\mathbb N}}
\newcommand{\PP}{{\mathbb P}}
\newcommand{\KK}{{\mathbb K}}
\newcommand{\CC}{{\mathbb C}}
\newcommand{\RR}{{\mathbb R}}
\renewcommand{\R}{{\mathbb R}}
\newcommand{\QQ}{{\mathbb Q}}
\newcommand{\FF}{{\mathbb F}}
\newcommand{\TT}{{\mathbb T}}
\renewcommand{\T}{{\mathbb T}}
\newcommand{\TP}{{\mathshbb{CP}}}
\newcommand{\m}{{\mathfrak M}}
\newcommand{\A}{{\mathcal A}}
\newcommand{\B}{{\mathcal B}}
\renewcommand{\D}{{\mathcal D}}
\renewcommand{\C}{{\mathcal C}}
\renewcommand{\E}{{\mathcal E}}
\newcommand{\F}{{\mathcal F}}
\newcommand{\X}{{\mathcal X}}
\newcommand{\M}{{\mathcal M}}
\renewcommand{\N}{{\mathcal N}}
\newcommand{\I}{{\mathfrak I}}
\newcommand{\U}{{\mathcal U}}
\renewcommand{\H}{{\mathcal H}}
\renewcommand{\S}{{\mathcal S}}
\renewcommand{\P}{{\mathcal P}}
\renewcommand{\L}{{\mathcal L}}
\newcommand{\V}{{\mathcal V}}
\renewcommand{\Im}{\text{Im\ }}
\newcommand{\Ve}{\text{Vert}}
\newcommand{\Ed}{\text{Edge}}
\newcommand{\Hod}{\text{Hod}}
\renewcommand{\div}{\text{div}}
\newcommand{\val}{\text{val}}
\newcommand{\trp}{\text{Trop}}
\newcommand{\oC}{\overset{\circ}C}
\renewcommand{\Q}{{\mathcal Q}}
\newcommand{\Sh}{\text{Sh}}
\newcommand{\f}{{\mathcal f}}
\newcommand{\fb}{{\mathcal f}}
\newcommand{\Vol}{\text{Vol}}
\newcommand{\tlim}{\lim_{t\to +\infty} \Log_t}
\newcommand{\tP}{\Pi}
\renewcommand{\epsilon}{\varepsilon}

\newcommand{\td}{{\text{''}}}
\newcommand{\tg}{{\text{``}}}
\newcommand\Rwan[1]{\textbf{[R1: #1]}}
\newcommand{\mnote}{\marginpar}

 \theoremstyle{definition}
\newtheorem{defi}[thm]{Definition}
 \newtheorem{proposition}[thm]{Proposition}
\newtheorem{example}[thm]{Example}
 \newtheorem{exa}[thm]{Example}
\newtheorem{theorem}[thm]{Theorem}
 \newtheorem{exo}[]{Exercises}
 \newtheorem{definition}[thm]{Definition}
 \newtheorem{construction}{Construction}[section]
 \newtheorem{corollary}[thm]{Corollary}

%% file: Erwan.tex
\section{Tropical algebra}\label{sec:trop alg}

\subsection{Tropical semi-field}
The set of \emph{tropical numbers} is defined as $\TT=\RR\cup
\{-\infty\}$. We endow $\TT$
  with the following
operations, called \emph{tropical addition} and
\emph{multiplication}: 
$$``x+y" = \max \{x, y \} \qquad ``x \times y" = x+y$$
with the usual conventions:
$$ \forall x \in \T, \quad ``x + (-\infty)" = \max (x, -\infty) = x \quad \text{and} \quad ``x \times (-\infty)" = x + (-\infty) = -\infty. $$

In the entire text, tropical operations will be placed under quotation
marks. 
Just as in classical algebra we often abbreviate $``x \times y"$ to $``xy"$.  
The tropical numbers
along with these two operations
 form a semi-field, {\it i.e.}, they satisfy all the axioms
of a field except the existence of an inverse for the law $\tg + \td$. 

To familiarize ourselves with these two operations, let us 
do some simple calculations: 

$$\tg 1+1\td=1, \ \ \tg 1+2\td=2, \ \ 
\tg 1+2+3\td=3, \ \  \tg 1\times 2\td=3,  \ \   \tg 1\times
(2+(-1))\td=3, $$
$$   \tg
1\times (-2)\td=-1, \ \ \tg (5+3)^2 \td= 10, \ \ \tg (x+y)^n \td= \tg x^n+y^n\td.
$$

Be 
careful when writing tropical formulas!  As, 
$``2x" \neq ``x+x"$ but $ ``2x" = x+2$, similarly $``1x" \neq x$ but $``1x" = x+1$, 
and again $``0x" = x$ and $``(-1)x" = x-1$. 

\medskip
A very important feature of the tropical 
semi-field 
is that it is 
\emph{idempotent}, 
which means that $``x+x" = x$ for all $x$ in
$\T$. This implies, in particular, that one cannot solve the problem of
non-existence of tropical substraction 
by adding more elements to $\T$ 
(see Exercise 1(1)).

\subsection{Maslov Dequantization} 
Let us explain how the tropical semi-field arises naturally as the 
limit of  some classical semi-fields.  This procedure, studied by 
V.~Maslov  
and his collaborators, 
is known as
\emph{dequantization of  
positive real numbers}.

The non-negative real numbers form a semi-field $\R_{\geq 0}$
under 
the usual addition and multiplication. 
If $t$ is a  real number greater than 1, then the logarithm of base $t$ 
provides a bijection between the sets $\R_{\geq 0}$ and $\T$. This bijection 
induces a semi-field structure on $\T$ with the operations denoted by 
$\tg +_t\td$ and $\tg\times_t\td$: 
$$\tg x +_t  y\td= \log_t(t^x +t^y) \  \ \ \text{ and } \ \ \  \tg x \times_t  y\td=
\log_t(t^x t^y) = x+y.$$ 

The equation on the right-hand side already shows classical addition arising from the multiplication on $(\R_{\geq 0}, +, \times)$. 
Notice that by construction, all 
semi-fields 
$(\T,\tg +_t\td,\tg
\times_t\td)$ are isomorphic to $(\R_{\geq 0},+,\times)$. 
The 
inequalities 
$\max(x,y)\le x+y\le 2\max(x,y)$ on $\R_{\geq 0}$ together with the fact that the logarithm 
of base $t > 1$ 
is 
an increasing function gives us the following bounds for $\tg +_t\td$:
$$\forall t>1,  \ \  \max(x,y) \le \tg x +_t  y\td \le \max(x,y)
+\log_t 2.$$
If we let $t$ tend to infinity,  $\log_t 2$ tends to $0$, 
and the operation
$\tg +_t\td$ therefore tends to 
the tropical addition $\tg +\td$! 
Hence
 the tropical semi-field comes naturally from degenerating the classical semi-field 
$(\R_{\geq 0},+,\times)$.
From an alternative perspective, we can view the classical semi-field 
$(\R_{\geq 0},+,\times)$ as a deformation of the tropical semi-field. This explains the use of the 
term 
{\it dequantization}. 

\subsection{Tropical polynomials}
As in classical algebra,
a tropical polynomial expression  $P(x)=\tg\sum_{i = 0}^d a_ix^i
\td $ induces  a tropical polynomial function, still denoted by $P$, on
$\TT$:
$$
\begin{array}{cccc}
P : & \TT &\longrightarrow & \TT
\\ & x &\longmapsto & \tg\sum_{i = 0}^d a_ix^i\td= \max_{i=1}^d(a_i + ix ) 
\end{array}.$$

Note that the map which associates a tropical polynomial function to a
tropical polynomial is surjective, by definition, but is \emph{not}
injective.
In the whole text, tropical polynomials have to be understood as
\emph{tropical polynomial functions}.

Let us look at some examples of tropical polynomials:
$$\tg x\td = x, \ \ \  \tg 1+ x\td = \max(1,x), \ \ \  
\tg 1+ x +3x^2\td =\max(1,x,2x+3), $$
$$\tg 1+ x +3x^2+(-2)x^3\td =\max(1,x,2x+3, 3x-2).$$

\begin{figure}[h]
\begin{center}
\begin{tabular}{ccc}
\includegraphics[width=4cm, angle=0]{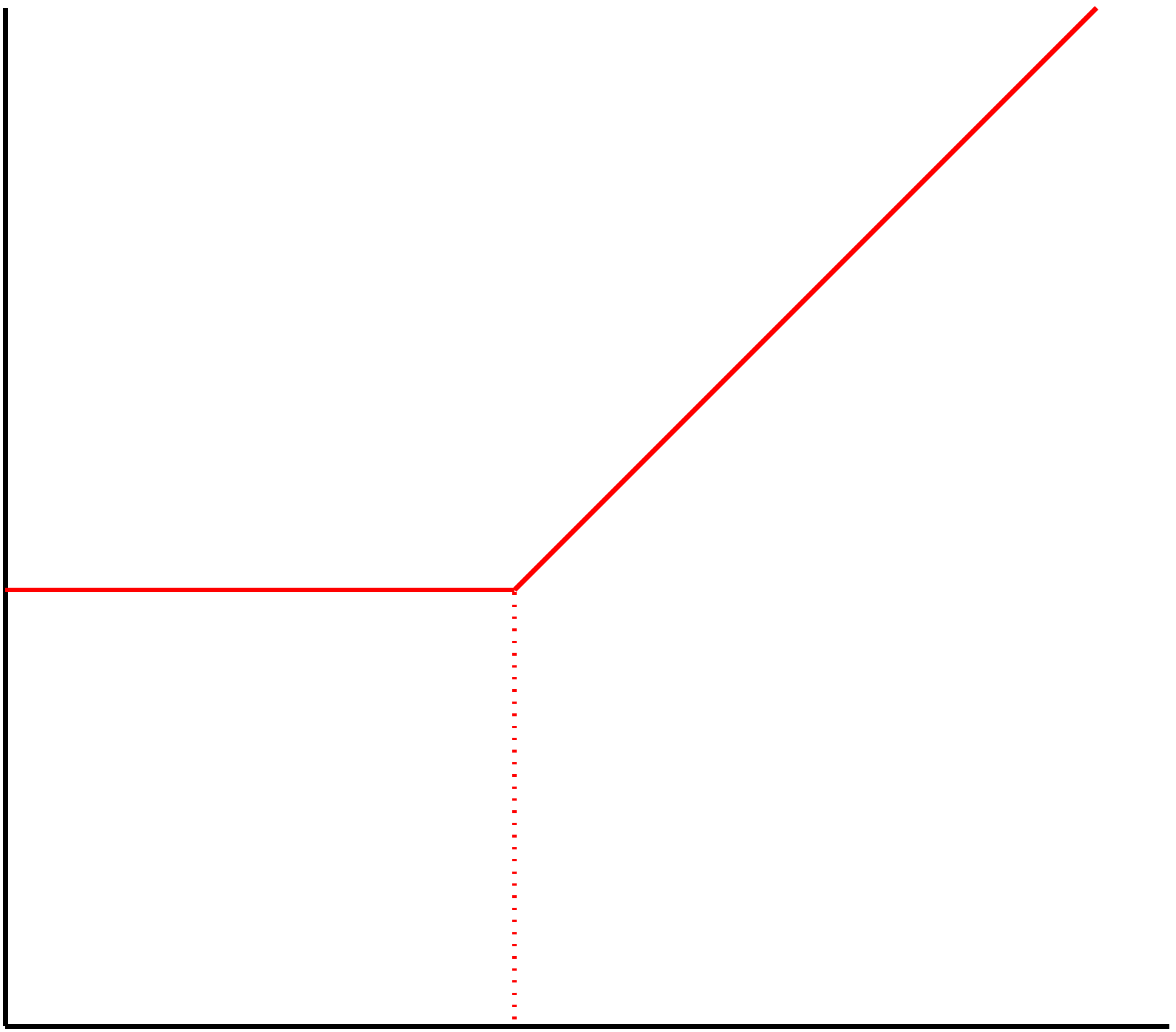}
\put(-120, 40){\small{$0$}}
\put(-140, -10){\small{$(-\infty,-\infty)$}}
\put(-66, -10){\small{$0$}}
&\hspace{20ex} &
\includegraphics[width=4cm, angle=0]{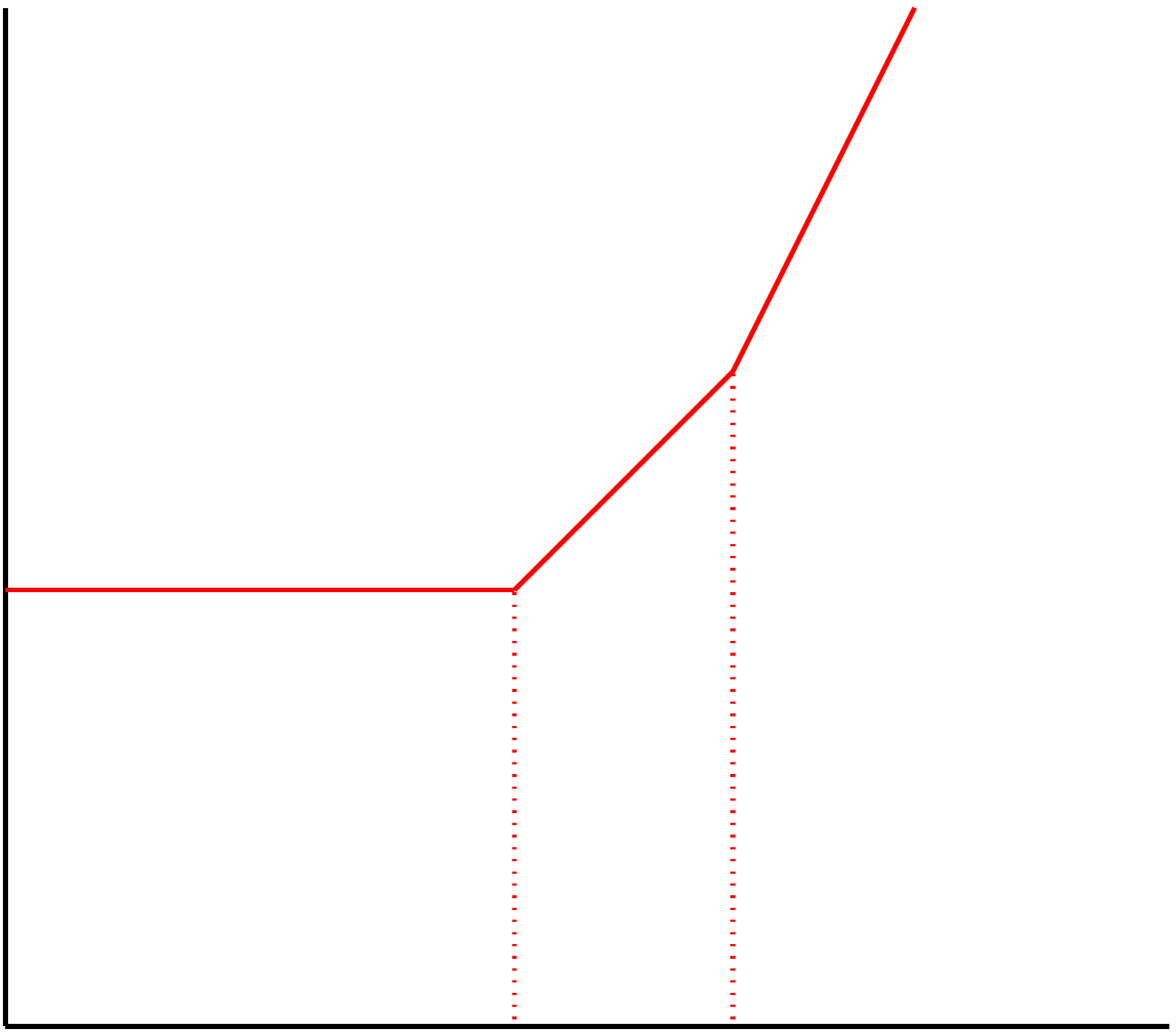}
\put(-120, 40){\small{$0$}}
\put(-140, -10){\small{$(-\infty,-\infty)$}}
\put(-66, -10){\small{$0$}}
\put(-46, -10){\small{$1$}}
\\ \\a) $P(x)=\tg 0+x \td$ && b) $P(x)=\tg 0+x + (-1)x^2\td$
\end{tabular}

\bigskip

\begin{tabular}{c}
\includegraphics[width=4cm, angle=0]{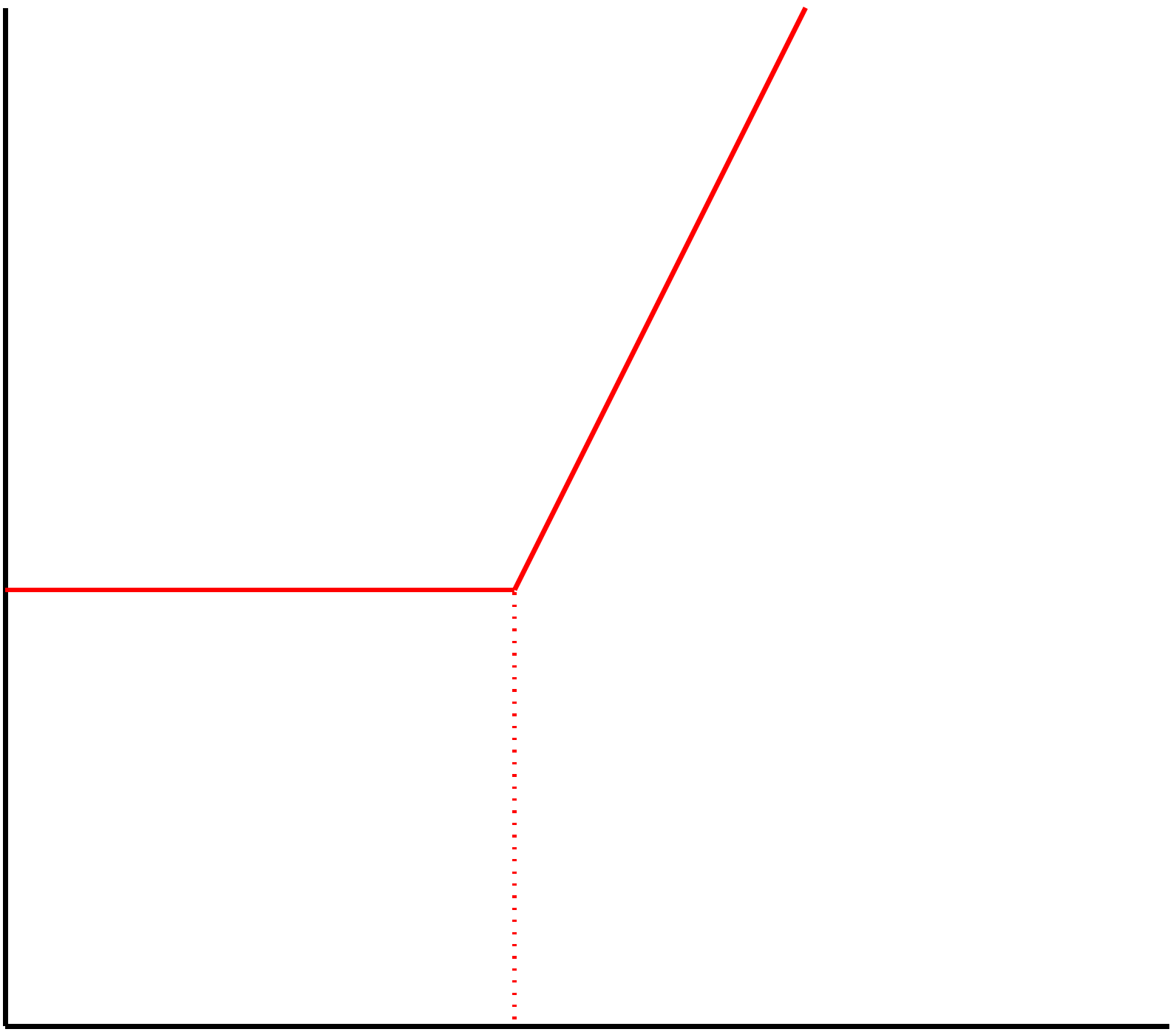}
\put(-120, 40){\small{$0$}}
\put(-140, -10){\small{$(-\infty,-\infty)$}}
\put(-66, -10){\small{$0$}}
\\ \\ c) $P(x)=\tg 0+x^2\td$
\end{tabular}
\end{center}
\caption{The graphs of some tropical polynomials}
\label{graphes}
\end{figure}
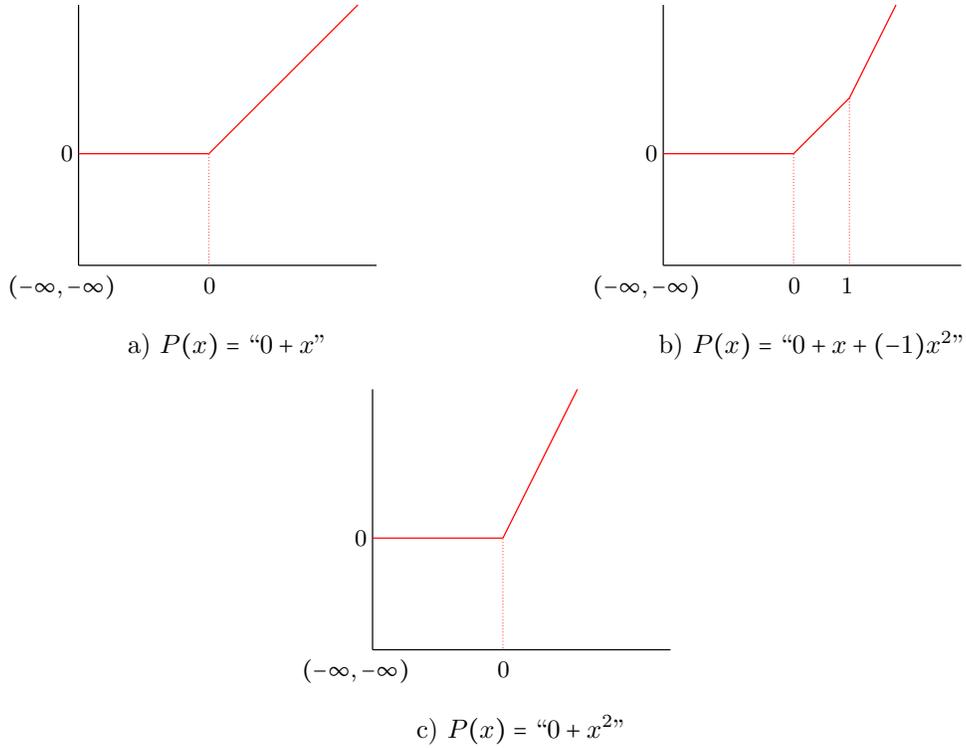

Now 
define 
the roots of a tropical polynomial. For this, let us 
 take a geometric point of view of the
problem.  A tropical polynomial is a convex piecewise affine function
and each piece has an integer slope (see Figure \ref{graphes}). We call
\emph{tropical roots} of the polynomial $P(x)$ all points $x_0$ of
$\T$ for which the graph of $P(x)$ has a corner at $x_0$. 
Notice, this is equivalent to $P(x_0)$ being equal to the value of at least two of its monomials evaluated at $x_0$. 
   Moreover,
the difference in the slopes of the two pieces adjacent to a corner
gives the \emph{order} of  the corresponding root. 

\begin{defi}\label{defi:roots}
The roots of a tropical polynomial 
 $P(x)=\tg \sum_{i=0}^d a_ix^i\td$
 are  the tropical numbers $x_0$ for which 
 either $P(x_0) = -\infty$, or there exists 
 a pair $i \neq j$ such that 
$P(x_0)=\tg a_ix_0^i  \td=\tg a_jx_0^j\td$. 

The order of a root  $x_0$ is the maximum of $|i-j|$ for all 
possible such pairs $i$, $j$ if $x_0\ne -\infty$, and is the minimal
$i$ such that $a_i\ne -\infty$ if $x_0=-\infty$. 
\end{defi}
Thus, the
polynomial $``0+x"$ has a simple root at $x_0 = 0$, the polynomial
$``0 + x + (-1)x^2\td $ has simple roots $0$ and $1$, and the polynomial
$``0+x^2\td $ has a double root at $0$. 

\begin{proposition}\label{prop:alg closed}
The tropical semi-field is algebraically closed. In other words, every tropical polynomial of 
degree $d > 0$ has exactly $d$ roots when counted with multiplicities. 
\end{proposition} 

\begin{exa}
We have the following factorizations:
$$\tg 0+x + (-1)x^2\td = \tg (-1)(x+0)(x + 1)\td \ \ \ \mbox{and} \ \ \ \tg 0 +
x^2\td = \tg (x+0)^2\td.$$ 
Once again the 
equalities hold in terms of polynomial functions and not on the level of the 
polynomial expressions.
For example, $\tg 0+x^2\td$ and $ \tg (0+x)^2\td$ are equal as polynomial functions but not as polynomials.
\end{exa}

\subsection{Relation to classical algebra}\label{sec:relationtoclassical}
Let $P_t(z)=\sum \alpha_i(t) z^i$ be a family of complex polynomials 
parameterized by
$t$ which we assume to be a sufficiently large 
positive 
number. We make the assumption that 
$$\forall i, \quad \exists a_i\in \TT, \quad \exists \beta_i\in \CC,\quad \alpha_i(t)\sim_{t\to+\infty} 
\beta_i t^{a_i}.$$ 

Then, we define the tropical polynomial, called the
\emph{tropical limit} of the family $P_t$, by 
$$P_{trop}(x)= \tg \sum a_i x^i \td.$$
We also define the map
$$\begin{array}{cccc}
\Log_t : & \CC &\longrightarrow & \TT
\\ & z &\longmapsto & \log_t(|z|)
\end{array}.
$$

The following theorem can be seen as a dual version of Newton-Puiseux
method.

\begin{thm}\label{NP1} 
One has
$$\Log_t \left( \{\text{roots of }P_t \}\right) \xrightarrow[t\to+\infty]{}
 \{\text{roots of } P_{trop}\}.$$
 Moreover, the order of 
any 
tropical root 
$x_0$ of $P_{trop}$ is exactly the number of roots of $P_t$
whose logarithms converge to $x_0$. 
\end{thm} 

\begin{exo} 

\begin{enumerate}
\item Why does the idempotent property of tropical addition prevent the
existence of additive inverses?

\item Draw the graphs of the tropical polynomials  $P(x)=\tg x^3+2x^2+3x
  +(-1)\td$ and $Q(x)=\tg x^3+(-2)x^2+2x+(-1)\td$, and determine their tropical 
  roots. 
\item Prove that the only root of the tropical polynomial $P(x)=\tg x\td$
  is $-\infty$. 
\item Prove that  $x_0$ is a root of order $k$ of a tropical
polynomial $P(x)$ if and only if there exists a 
tropical polynomial $Q(x)$ such that $P(x) =\tg (x+x_0)^kQ(x) \td$ and
$x_0$ is not a root of $Q(x)$. 
(Note that a factor $x-x_0$ in classical algebra gets transformed to the 
factor  $\tg x+x_0\td$, since the root of the polynomial
$\tg x+x_0\td$ is
$x_0$ and not  $-x_0$.)

\item Prove Proposition \ref{prop:alg closed}.
\item Let $a  \in \R$ and  $b, c , d\in \T$. Determine the roots of the 
polynomials  $ \tg ax^2+bx+c\td$ and $ \tg ax^3+bx^3+cx + d\td$. 
\end{enumerate}
\end{exo}

\section{Tropical curves in $\RR^2$}\label{sec:trop curves}
Let us now extend the preceding notions to the case of
 tropical polynomials in two variables.
Since this makes all definitions,  statements and  drawings simpler,
we restrict ourselves 
to tropical 
curves in $\RR^2$ instead of $\TT^2$. 

\subsection{Definition}\label{sec:deftropcurves}
 A {\it tropical polynomial in two variables} is 
$$P(x,y)=\tg \sum_{(i,j) \in A}a_{i,j}x^iy^j\td = 
\max_{(i,j) \in A}(a_{i,j}+ix+jy),$$ 
where $A$ is a finite subset of $(\ZZ_{\geq 0})^2$.  
Thus,
a tropical polynomial is  
a convex piecewise affine function,
and we denote by $\widetilde V(P)$  the corner locus of
this function. That is to say,
$$\widetilde V(P)=\left\{(x_0,y_0)\in\RR^2\ | \ 
\exists (i,j)\ne(k,l),\quad   P(x_0,y_0)= \tg a_{i,j}x_0^iy_0^j \td
=\tg a_{k,l}x_0^ky_0^l \td\right\}. $$

\begin{figure}[h]
\begin{center}
\begin{tabular}{ccc}
\includegraphics[width=3.5cm, angle=0]{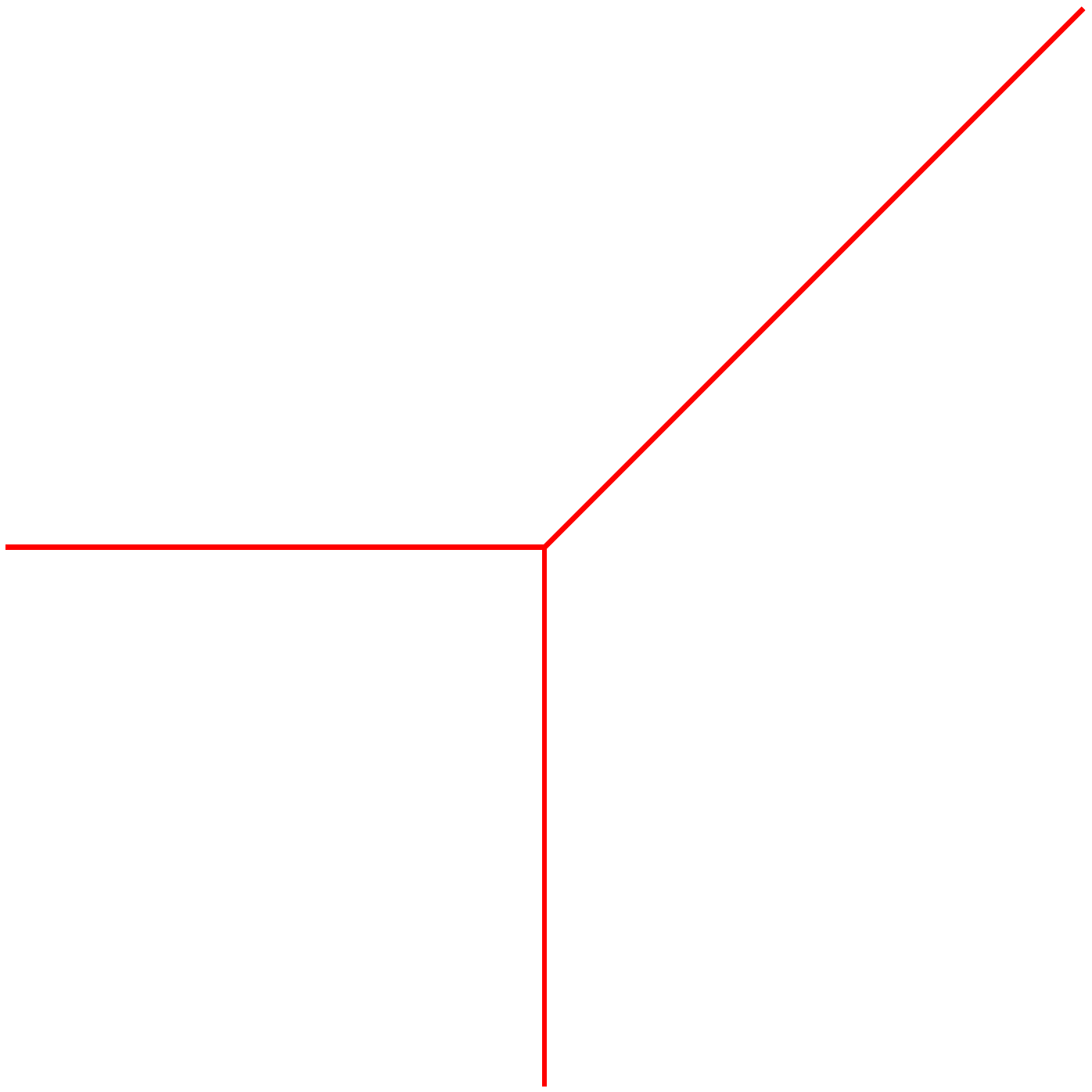}&
\includegraphics[width=3.5cm, angle=0]{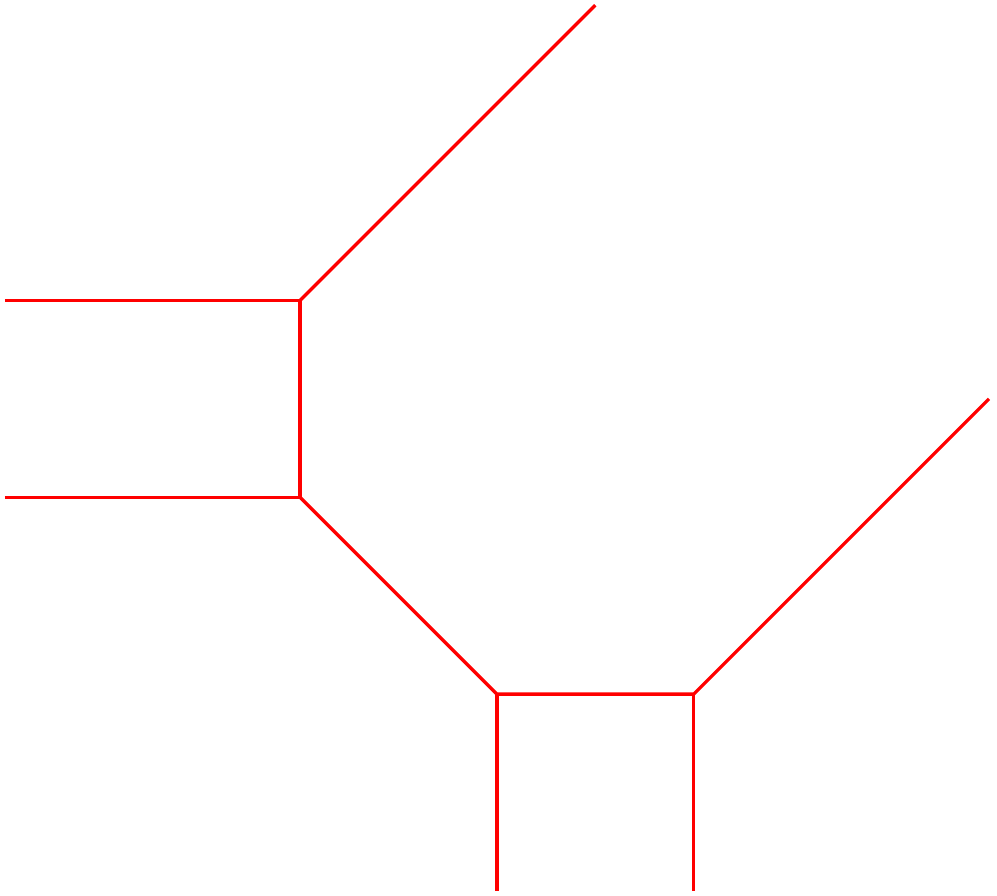}&
\includegraphics[width=3.5cm, angle=0]{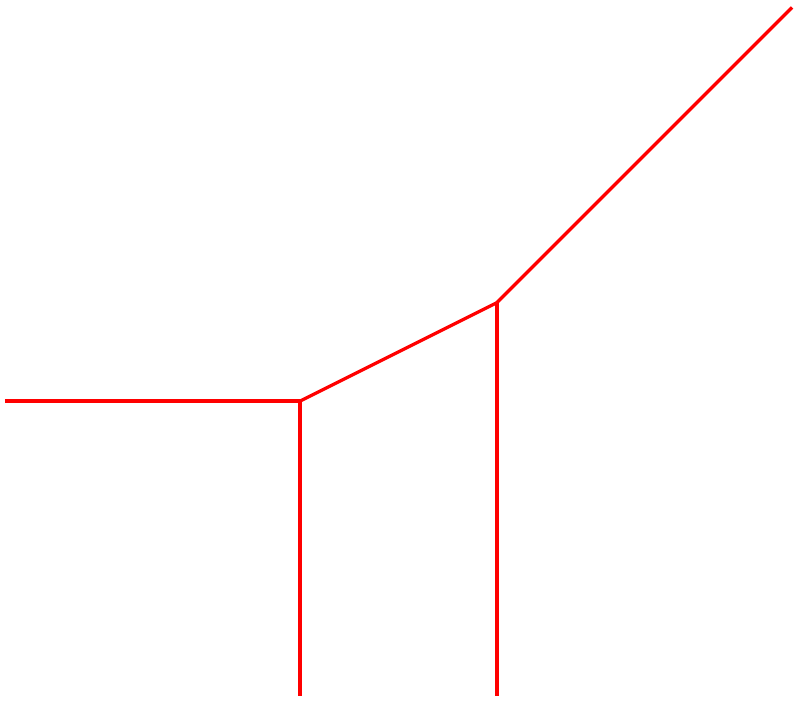}
\put(-100, 40){\small{$2$}}
\put(-30, 68){\small{$2$}}
\\ \\a) $\tg x+y+0\td$ & b)  $\tg 3+ 2x + 2y + 3xy+y^2+x^2\td$ 
 & c) $\tg 0+ x +y^2+(-1)x^2\td$

\end{tabular}
\end{center}
\caption{A tropical line and two tropical conics.}
\label{droite}
\end{figure}

\begin{exa}\label{ref:line}
Let us look
at the tropical line defined by the polynomial
$P(x,y)=\tg x+y+0 \td$. 
We must find the points 
$(x_0,y_0)$ in $\R^2$ that satisfy one of the following three 
conditions: 
$$x_0=0\ge y_0, \ \ \ \ \ \ \ \ \ y_0=0\ge x_0,
\ \ \ \ \ \ \ \ \ x_0=y_0\ge 0 $$ 

We see that the set $\widetilde V(P)$ is made 
of three standard half-lines  (see  Figure \ref{droite}a): 
$$\{(0,y) \in \R^2 \ | \ y\le 0 \},  \  \{(x,0) \in \R^2\ | \ x\le 0 
\}, \text{ and } \{(x,x) \in \R^2 \ | 
\ x\ge  0\}.$$
\end{exa}

The set $\widetilde V(P)$ is a piecewise linear graph in $\RR^2$ :
it is a finite union of possibly infinite straight edges in $\R^2$.
As in the case of
polynomials in one variable, 
we 
take into account the
difference in the slopes of $P(x, y)$  
on the two sides of  
an edge. 

\begin{defi}\label{def:trop curve}
The weight of an edge of $\widetilde V(P)$
 is defined as  the maximum of the 
greatest common divisor {\rm (}gcd{\rm )\/} of the numbers $|i-k|$ and 
$|j-l|$ for all pairs $(i,j)$ and $(k,l)$ such that the value of
$P(x,y)$ on this edge is given by the corresponding monomials. 

The tropical curve  
defined by  $P(x,y)$  is the 
graph $\widetilde V(P)$ equipped with
this weight function on 
the edges. 

\end{defi}

\begin{exa}
Figures \ref{droite}a,b,c depict tropical curves in $\R^2$. 
The weight of an edge is assumed to be $1$, unless indicated otherwise.
For example, in the case of the tropical 
line, all edges are of weight 1. 
\end{exa}

\subsection{Dual subdivision}\label{ref:dualSub} 
A polynomial $P(x,y)$ over any field or semi-field
always comes with a Newton polygon. 
Let 
$P(x,y)=\tg \sum_{i,j} a_{i,j}x^iy^j\td$ be a tropical polynomial. 
 The  \emph{Newton polygon} of $P(x,y)$, denoted by  $\Delta(P)$, is
 defined by
$$\Delta(P)=Conv\{(i,j) \in (\Z_{\geq 0})^2 \ | \ a_{i,j}\ne - \infty \} \subset \R^2 . $$ 
In classical algebra, one just replaces
   $-\infty$ by $0$ in the definition of $\Delta(P)$. 

A tropical polynomial
also determines 
a subdivision of 
$\Delta(P)$, called its \emph{dual subdivision}.
Given $(x_0,y_0)\in\RR^2$, let 
$$\Delta_{(x_0,y_0)}=Conv\{(i,j) \in (\Z_{\geq 0})^2 \ | \ P(x_0,y_0)=\tg a_{i,j}x_0^iy_0^j\td \} \subset \Delta(P). $$  
The tropical curve $C$ defined by $P(x,y)$ induces a polyhedral
decomposition of $\RR^2$,
and 
the polygon $\Delta_{(x_0,y_0)}$ only
depends on the cell $F \ni (x_0, y_0)$ of the decomposition given by $C$. 
Thus, we define $\Delta_F 
= \Delta_{(x_0,y_0)}$ for $(x_0, y_0) \in F$. 

\begin{figure}[h]
\begin{center}
\begin{tabular}{ccccc}
\includegraphics[width=1cm, angle=0]{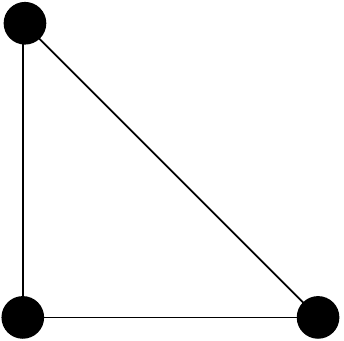}&\hspace{6ex} &
\includegraphics[width=2cm, angle=0]{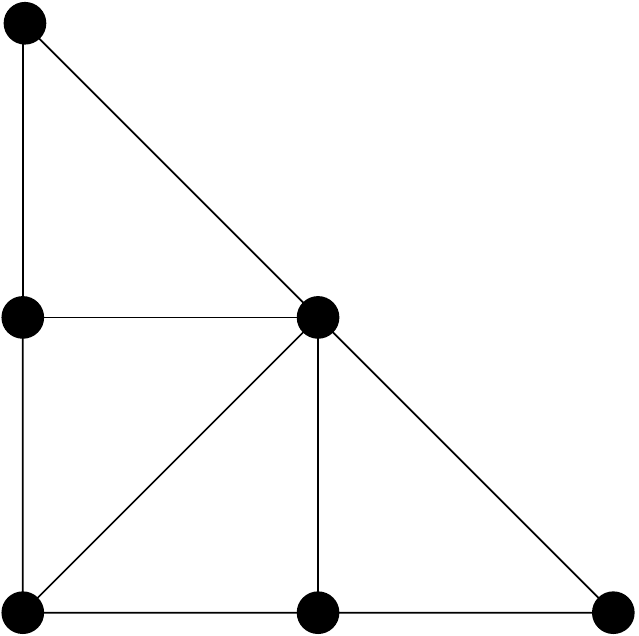}&\hspace{6ex} &
\includegraphics[width=2cm, angle=0]{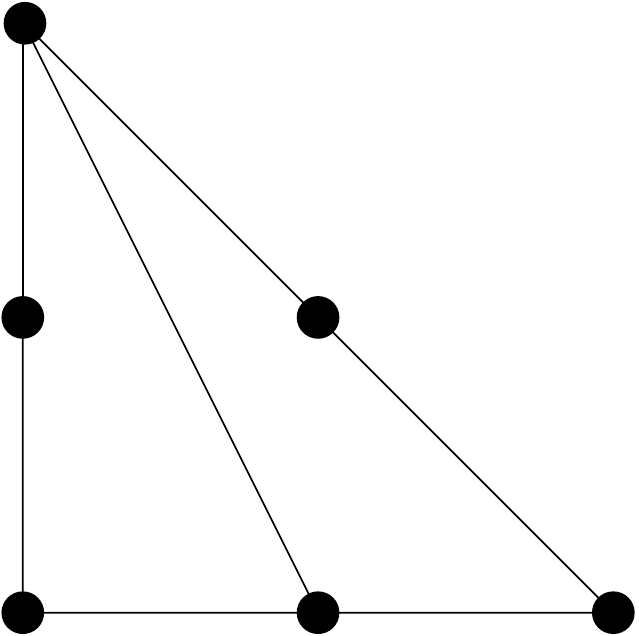} 
\\ \\a) && b) &&c)

\end{tabular}
\end{center}
\caption{Subdivisions dual to the tropical curves depicted in Figure \ref{droite}}
\label{subd}
\end{figure}

\begin{exa}\label{ex:celldecomp}
Let us go back to the tropical line $L$ defined by the 
polynomial  
$P(x,y)=\tg 
x+y+0\td$ (see  
Figure \ref{droite}a). 
On the 2-cell $F_1=\{\max(x,y)<0\}$,  
the value of $P(x,y)$ is given by the monomial $0$, and so
$\Delta_{F_1}=\{(0,0)\}$. Similarly, we have 
 $\Delta_{F_2}=\{(1,0)\}$ and $\Delta_{F_3}=\{(0,1)\}$ for the cells
$F_2=\{x>\max(y,0)\}$ and $F_3=\{y>\max(x,0)\}$. 

Along the horizontal edge $e_1$ 
of $L$ the value
of  $P(x, y)$ is given by the 
monomials $0$ and $y$, and so $\Delta_{e_1}$ is the 
 the vertical edge of $\Delta(P)$.  
In the same way, $\Delta_{e_2}$ is the 
 the horizontal edge of $\Delta(P)$ for the vertical edge $e_2$ of
 $L$, and $\Delta_{e_3}$ is the 
 the  edge of $\Delta(P)$ with endpoints $(1,0)$ and $(0,1)$
 for the diagonal edge $e_3$ 
 of $L$. 

The point
 $(0,0)$ is the vertex $v$ of the line $C$. This
is where the three monomials  
 $0$, 
$x$ and
$y$ take the same value, and so  $\Delta_v=\Delta(P)$, 
(see Figure \ref{subd}a). 
\end{exa} 

All polyhedra $\Delta_F$ form a subdivision of the
Newton polygon $\Delta(P)$. 
This subdivision 
is 
dual
to the tropical curve defined by $P(x,y)$
in the following sense.

\begin{prop}\label{prop:dual subd}
One has
\begin{itemize}
\item $\Delta(P)=\bigcup_F \Delta_F$, where 
the union is taken over all cells $F$ of
the polyhedral subdivision of 
$\RR^2$ induced by the tropical curve defined by $P(x,y)$;
\item $\dim F=\codim \Delta_F$;
\item $\Delta_F$ and $F$ are orthogonal;
\item $\Delta_F\subset \Delta_{F'}$ 
if and only if 
$F'\subset F$;
furthermore, in this case $\Delta_F$ is a face of $\Delta_{F'}$; 
\item $\Delta_F\subset \partial\Delta(P)$
if and only if 
$F$ is unbounded.
\end{itemize}
\end{prop} 

\begin{exa}
The dual subdivisions of the tropical curves in 
Figure \ref{droite} are drawn in 
 Figure \ref{subd} (the black points represent the points of 
$\Delta(P)$
which have integer
coordinates;
note that 
these points 
are not necessarily 
vertices of the dual
subdivision). 
\end{exa} 

The weight of an edge of a tropical curve 
can be read off from
the dual subdivision.
\begin{prop}\label{prop:weight}
 An edge $e$ of a tropical curve  has weight $w$ if and only if the \emph{integer
   length} of $\Delta_e$ is $w$, {\it i.e.} $Card(\Delta_e\cap \ZZ^2) - 1=w$.
\end{prop}

\subsection{Balanced graphs and tropical curves.}\label{sec:balancing}
\begin{figure}[h]
\begin{center}
\begin{tabular}{ccc}
\includegraphics[width=4cm,
  angle=0]{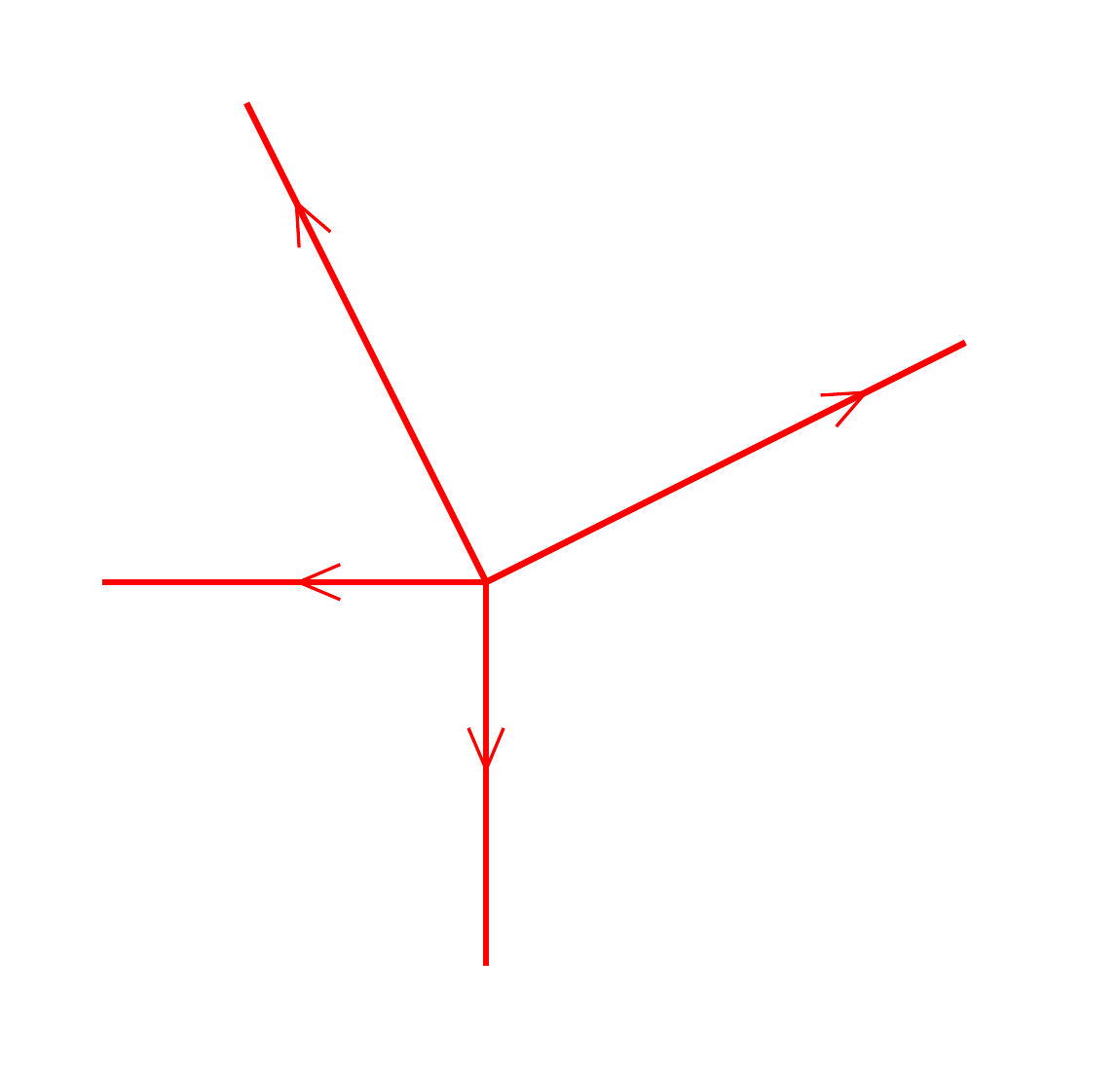}& \hspace{4ex} &
\put(-100, 55){\small{$v$}}
\put(-105, 20){\small{$3$}}
\includegraphics[width=4cm, angle=0]{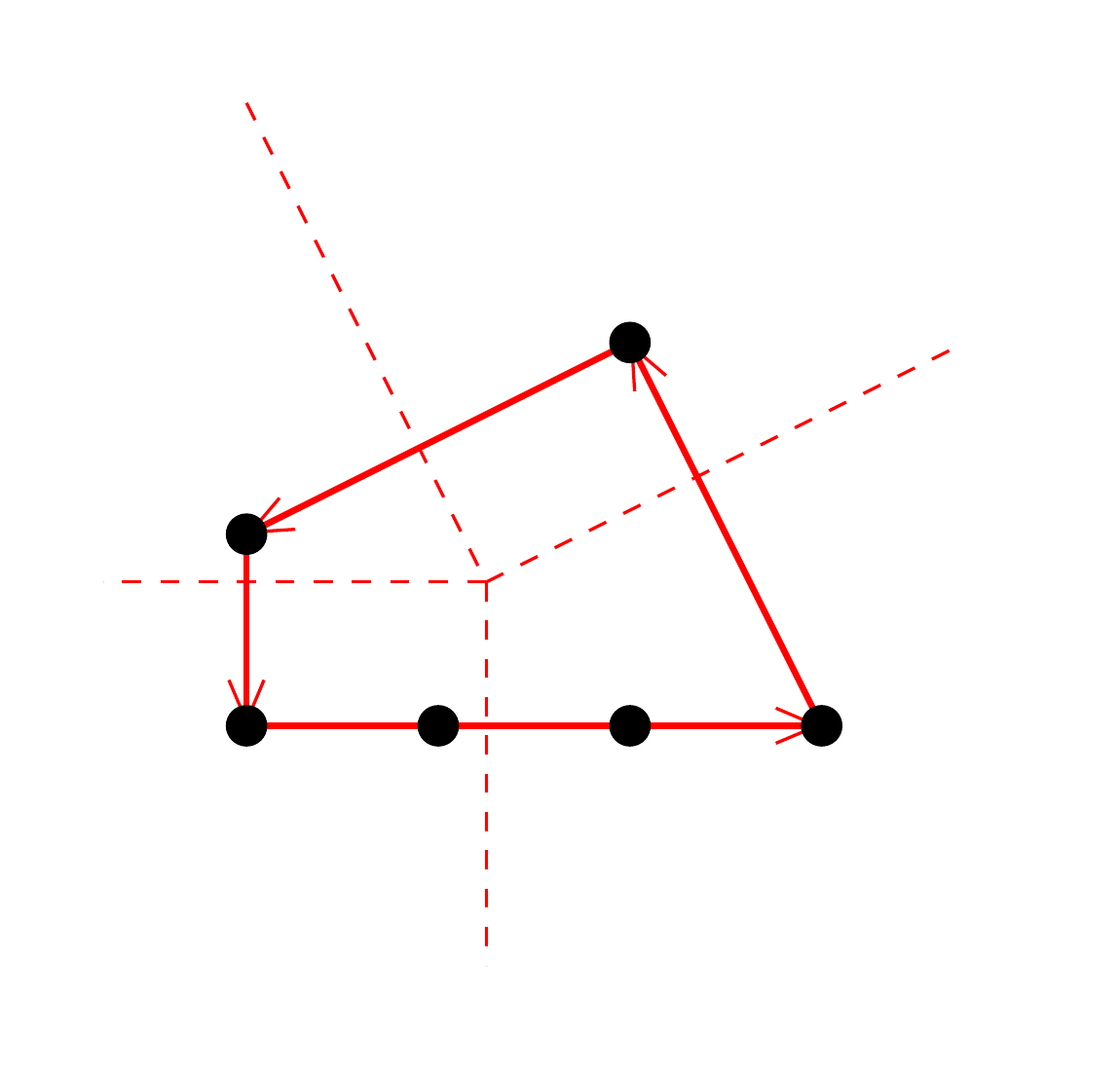}
\put(-30, 40){\small{$\Delta_v$}}
\\ a) && b)
\end{tabular}
\end{center}
\caption{Balancing condition.}
\label{equ}
\end{figure}

Let $v$ be a vertex of a tropical curve $C$, and let 
$e_1, \dots, e_k$ be the edges adjacent to $v$. 
Denote by 
$w_1, \dots , w_k$
the weights of $e_1, \dots, e_k$.
Let $v_i$, $i = 1$, $\ldots$, $k$, 
be the primitive integer vector ({\it i.e.}, having mutually prime $\Z$-coordinates) 
in the direction of $e_i$
and pointing outward from $v$, 
see Figure  \ref{equ}a. 
The vectors $w_1 v_1, \dots , w_kv_k$ are obtained from the sides of the closed polygon $\Delta_v$
(oriented counter-clockwise) via rotation by $\pi/2$, 
(see Figure
\ref{equ}b). Hence, 
the tropical curve $C$ satisfies the so 
called \emph{balancing condition} at each of its vertices $v$.

\begin{proposition}[Balancing condition]
One has
$$
\sum_{i=1}^k w_i v_i=0. 
$$
\end{proposition} 

A rectilinear graph $\Gamma \subset \R^2$ whose edges have rational slopes and 
are equipped with positive 
integer weights is called a \emph{balanced graph} if $\Gamma$ satisfies the
balancing condition 
at each 
vertex.
We have just 
seen that every tropical curve is a balanced graph. 
The converse is 
also true. 

\begin{thm} [\cite{Mik12}]\label{prop:trop balanced} 
Any balanced graph 
in $\R^2$  
is a tropical curve.
\end{thm} 

\begin{figure}[h]
\begin{center}
\begin{tabular}{ccc}
\includegraphics[width=5cm, angle=0]{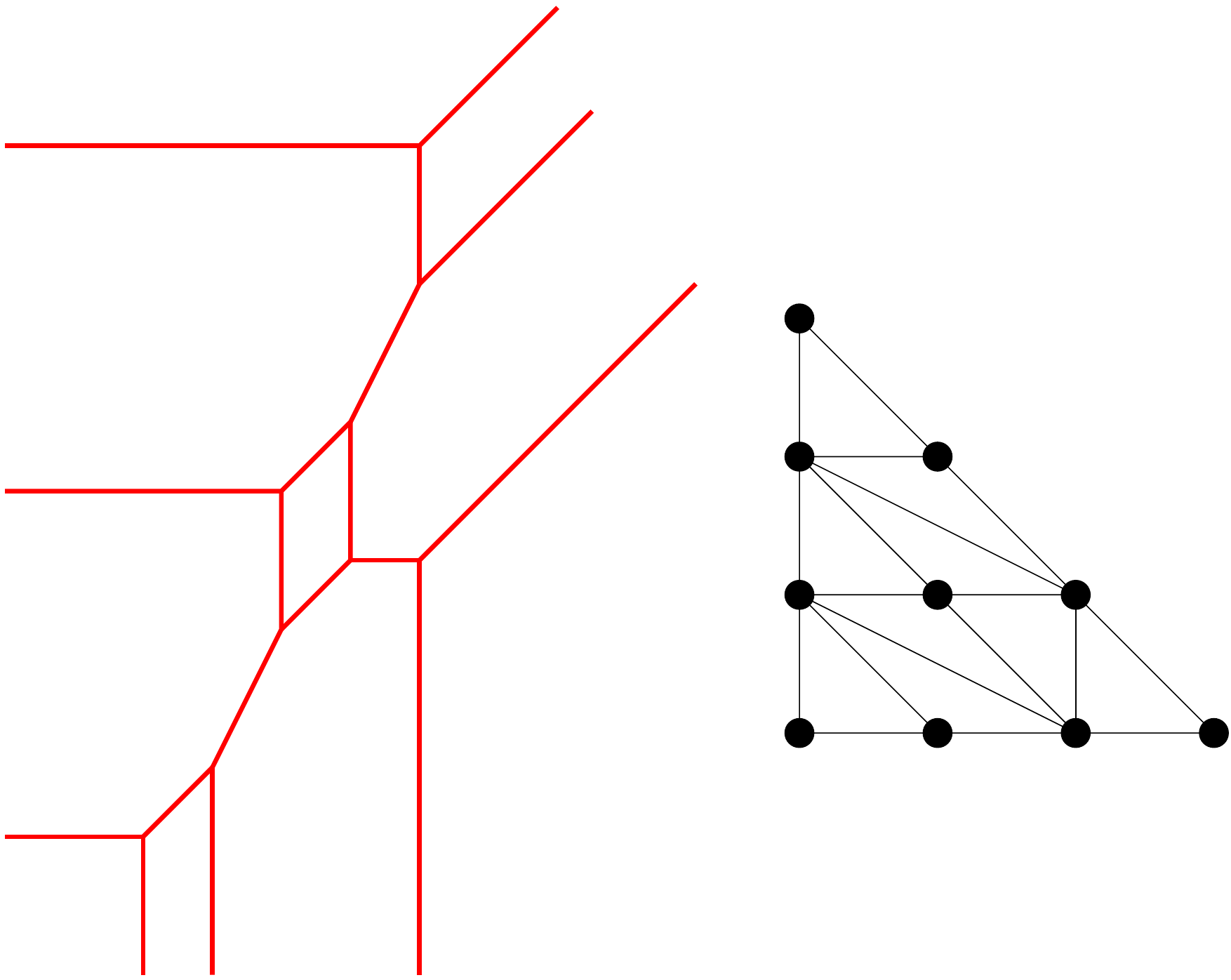}&\hspace{10ex} &
\includegraphics[width=5cm, angle=0]{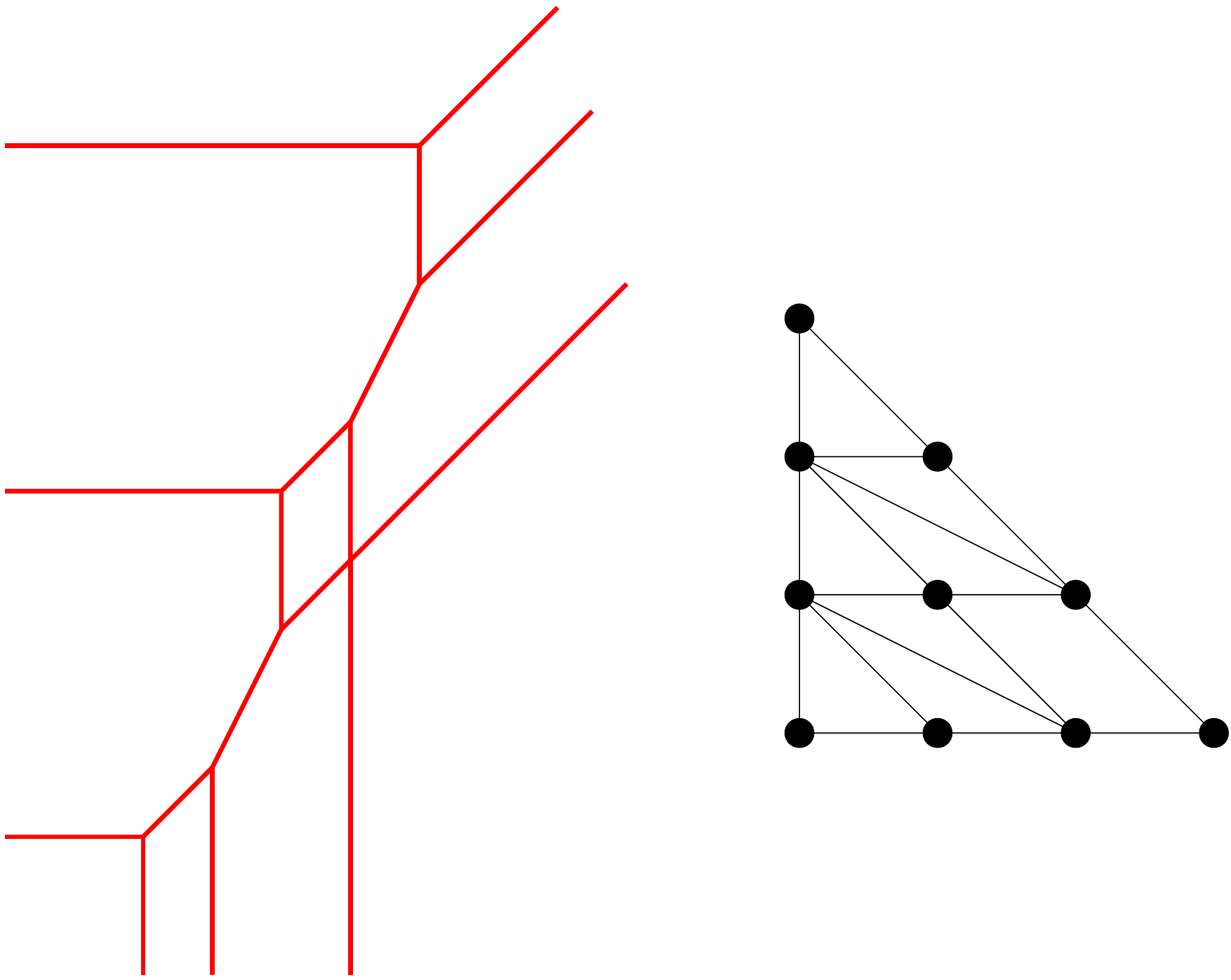}
\\ \\a) && b)
\end{tabular}

\bigskip
\begin{tabular}{c}
\includegraphics[width=5cm, angle=0]{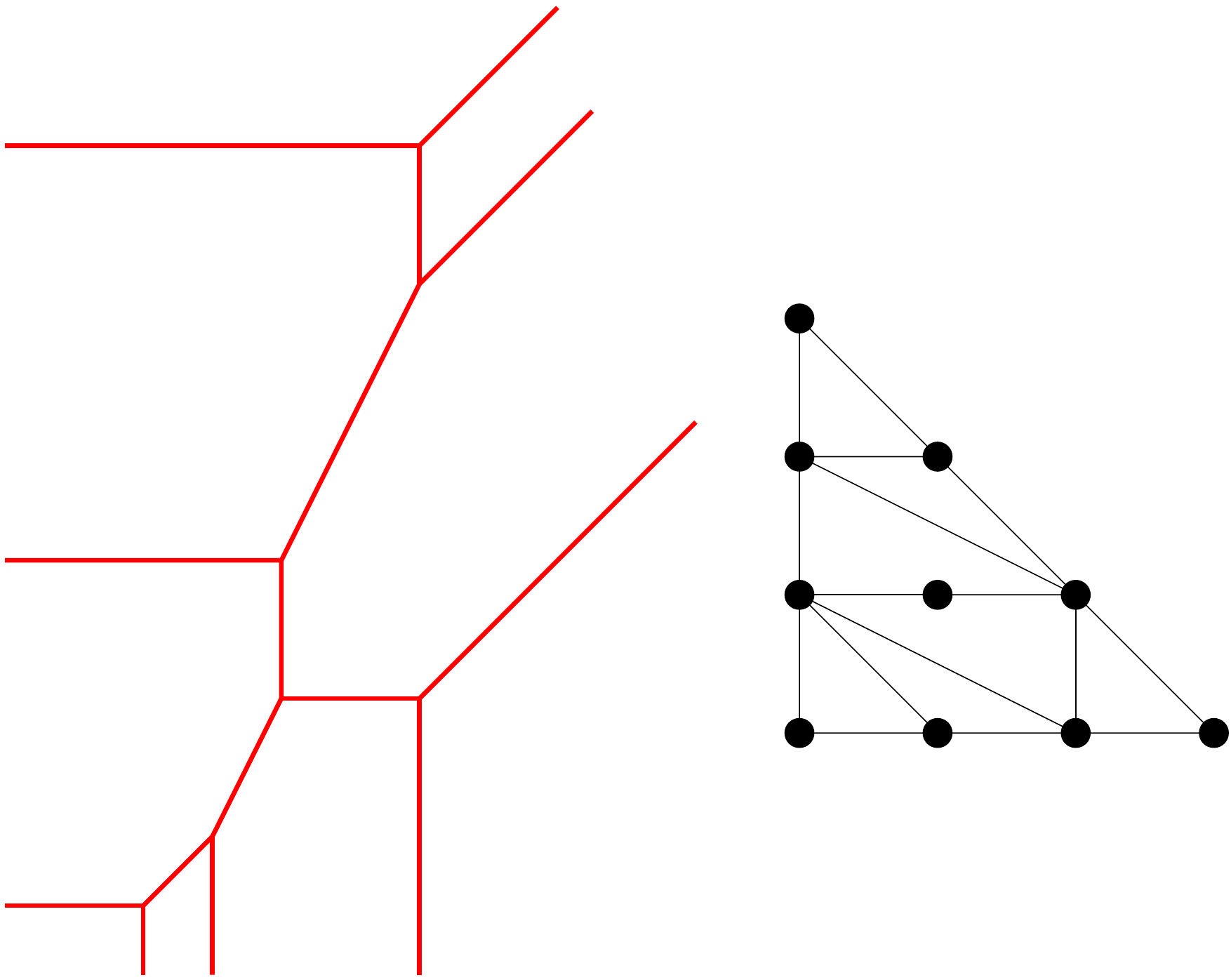}
\put(-105, 38){\small{$2$}}
\\ \\c)

\end{tabular}
\end{center}
\caption{Some tropical cubics and their dual subdivisions} 
\label{equil}
\end{figure}

For example, 
there exist tropical polynomials of
degree $3$ whose  
tropical curves are the weighted graphs depicted in Figure \ref{equil}. 
The figure also contains the dual subdivisions of these curves. 

\subsection{Tropical curves as limits of amoebas}\label{sec:limit amoeba}
As in the case of polynomials in one variable, tropical curves can be
approximated, via the logarithm map, by algebraic curves in $(\CC^\times)^2$.
For
this, we need the following map
(where $t > 1$): 
$$\begin{array}{cccc}
\text{Log}_t: & (\CC^\times)^2 &\longrightarrow & \RR^2
\\ & (z,w) &\longmapsto& (\log_t |z|,\log_t|w|)
\end{array}.$$ 
\begin{defi}[Gelfand-Kapranov-Zelevinsky \cite{GKZ}]
The \emph{amoeba} {\rm (}in base $t${\rm )\/} of  
$V \subset (\CC^\times)^2$ is 
$\text{Log}_t(V)$. 
\end{defi} 

\begin{figure}[h]
\centering
\begin{tabular}{ccc}
\includegraphics[width=4cm, angle=0]{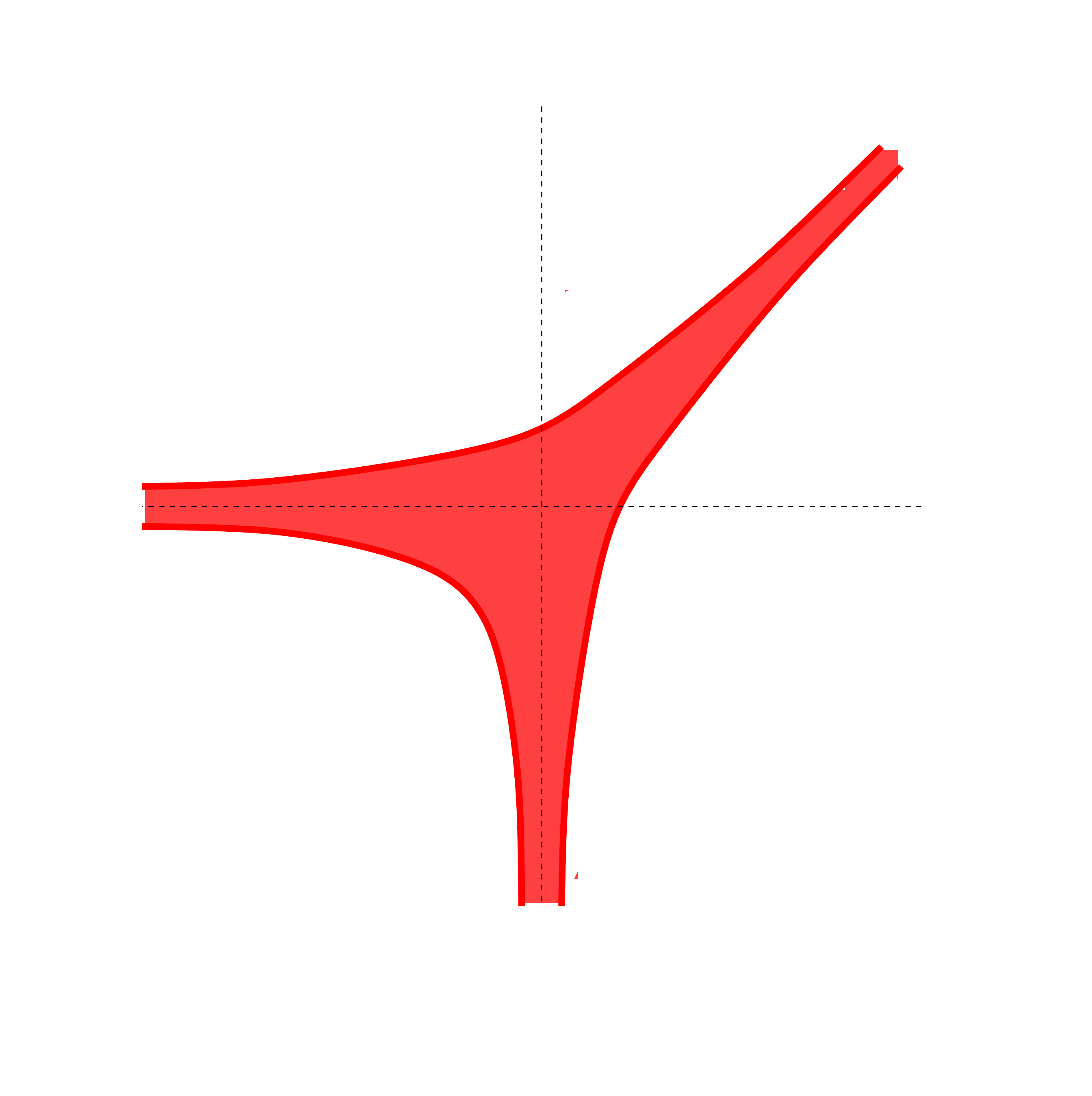} & \hspace{10ex}
&\includegraphics[width=4cm, angle=0]{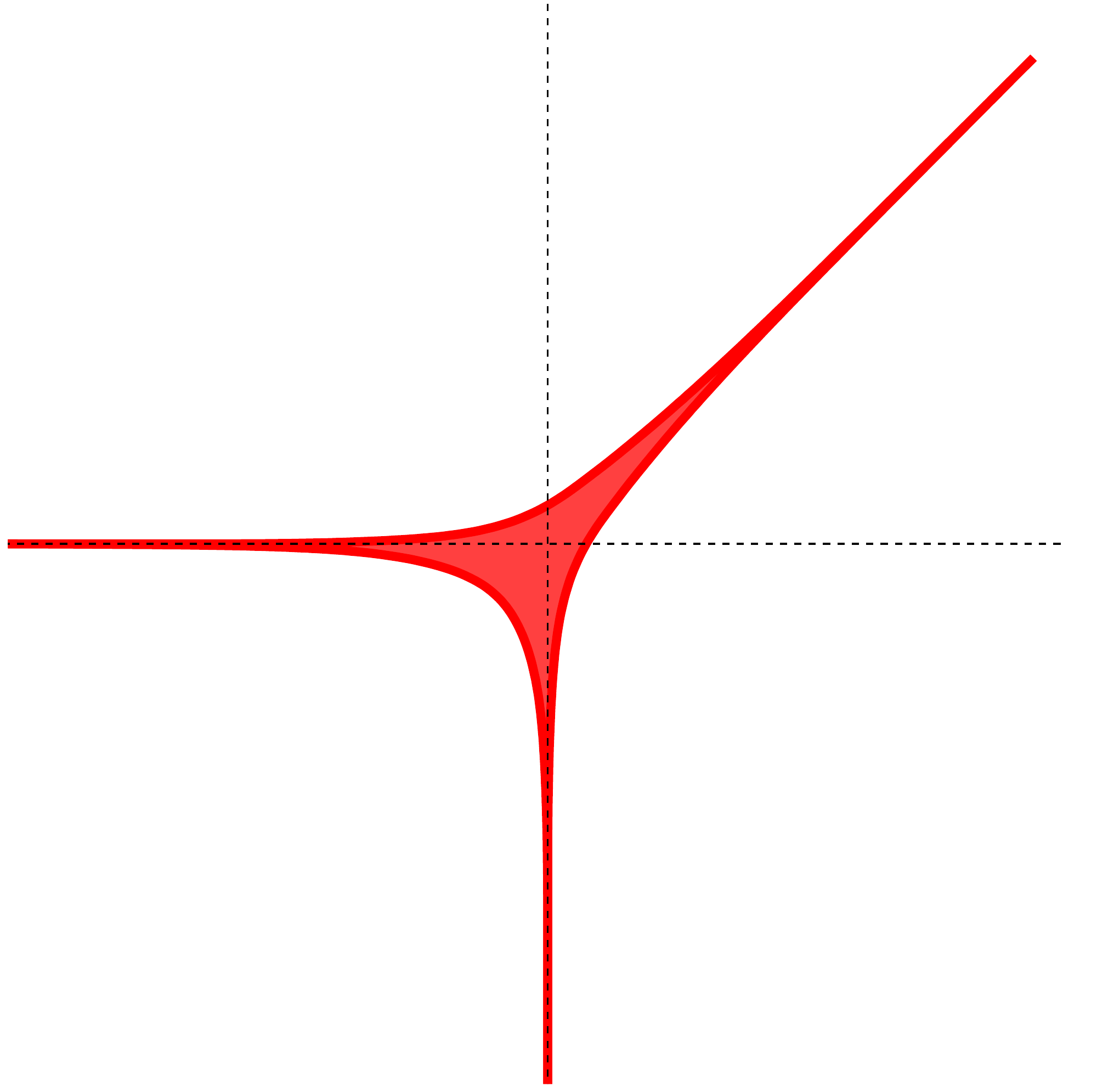} 
\\
\\ a) $\text{Log}(\L)$ && b) $\text{Log}_{t_1}(\L)$ \end{tabular}

\bigskip
\begin{tabular}{ccc}
\includegraphics[width=4cm, angle=0]{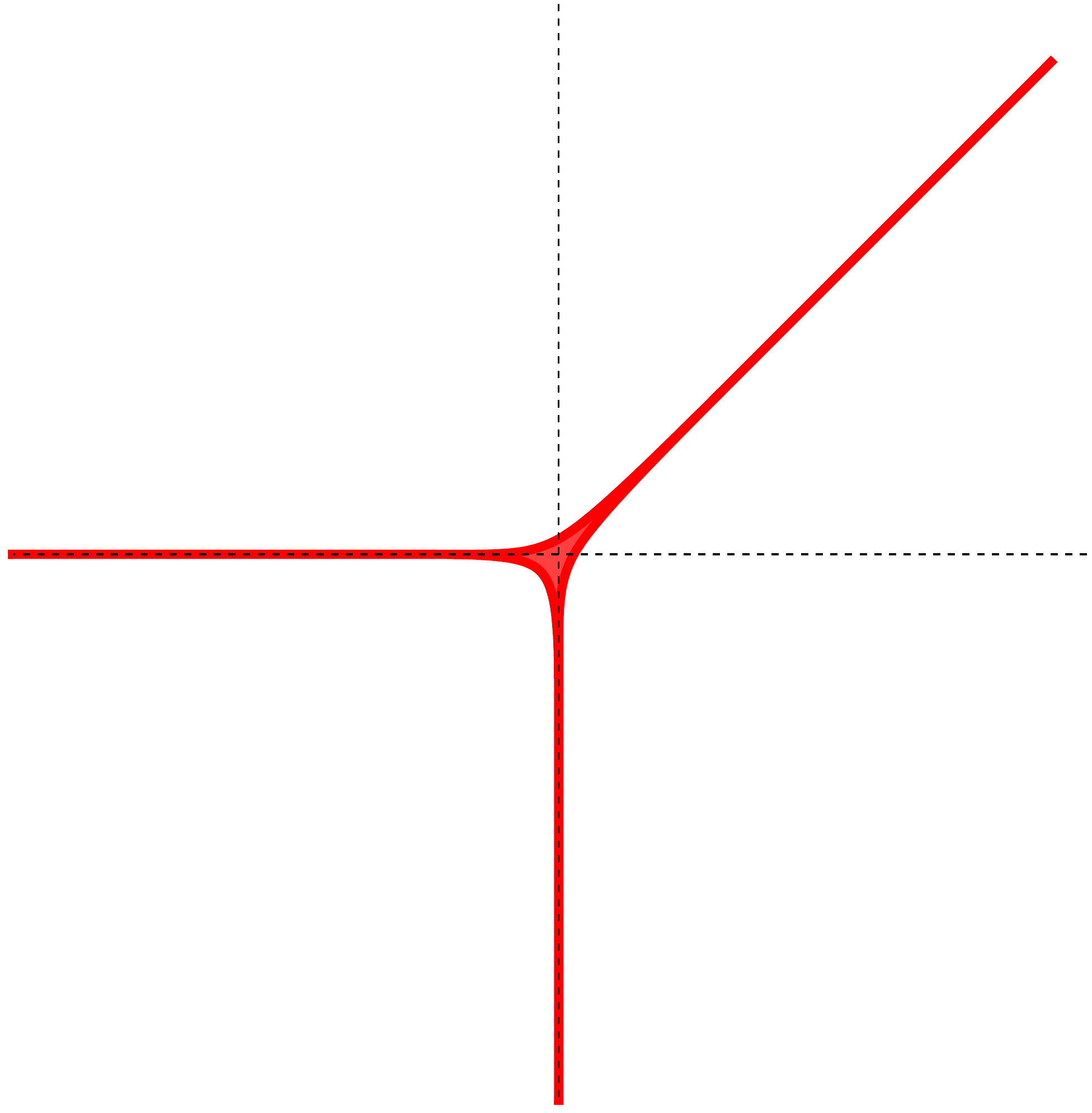} & \hspace{10ex}
&\includegraphics[width=4cm, angle=0]{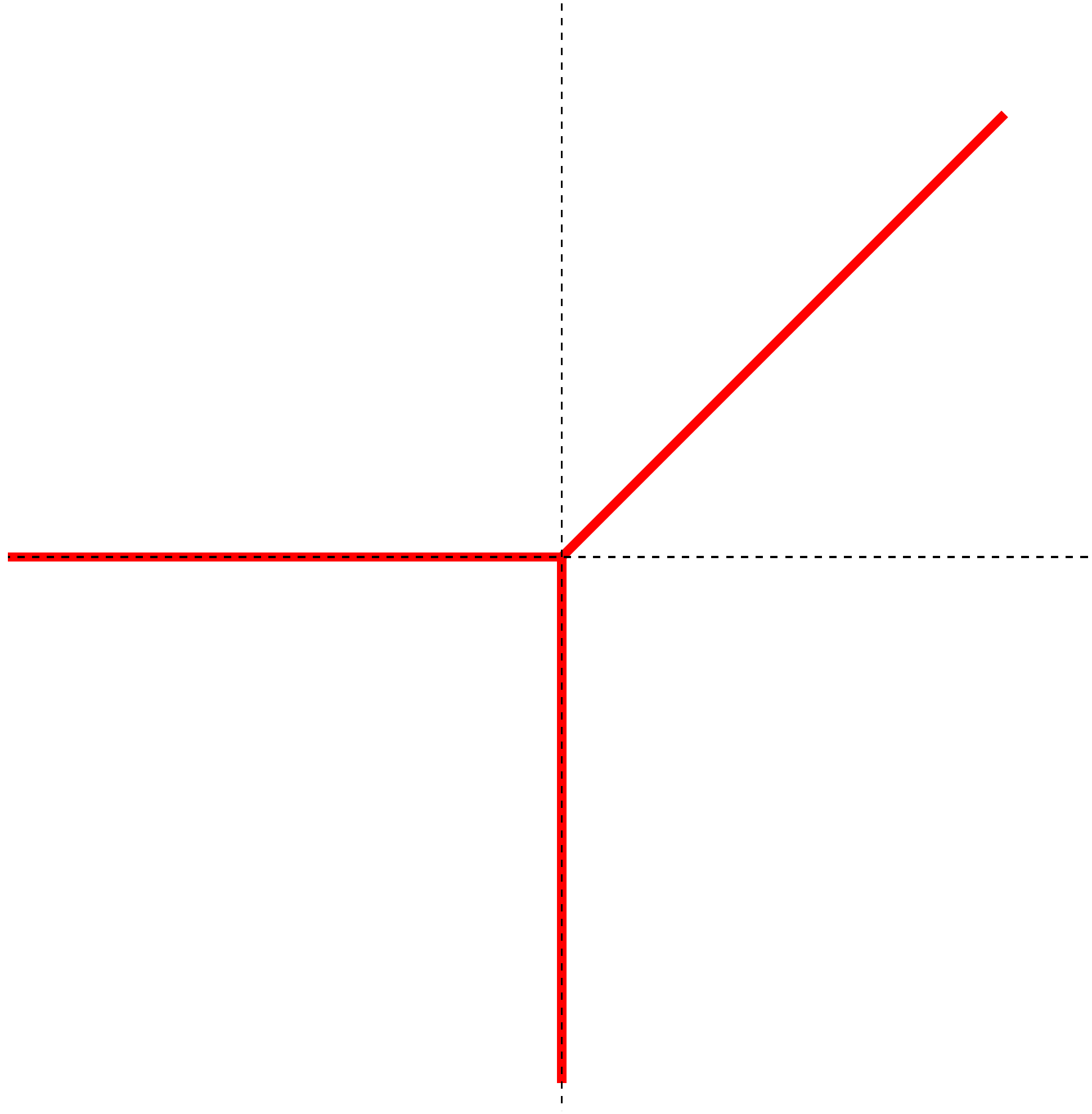}

\\
\\ c) $\text{Log}_{t_2}(\L)$ &&
d) $\lim_{t\to\infty}\text{Log}_t(\L)$ 
\end{tabular}
\caption{Dequantization of a line ($e<t_1<t_2$)}\label{amoeba line}
\end{figure}

For example, the amoeba of the line 
$\L$ 
defined by $z+w+1=0$  in
$(\CC^\times)^2$ 
is depicted in Figure \ref{amoeba line}a. 
This amoeba 
has three asymptotic directions: $(-1,0)$, $(0,-1)$, and $(1,1)$. 

The amoeba of $\L$ in base $t$ is a 
contraction by a factor $\log t$ of the amoeba of $\L$ in base
$e$ (see Figures \ref{amoeba line}b, c). 
Hence, when $t$ goes to $+\infty$, the 
amoeba is contracted
to the origin, only the three asymptotic
directions remain.
In other words, what we see at the limit in Figure \ref{amoeba line}d 
is a tropical line! 

Of course, the same strategy applied to any 
complex curve in $(\CC^\times)^2$ 
produces a 
similar picture at the limit: 
a collection of rays emerging from the origin
in the asymptotic directions of the amoeba. 
To get a more interesting limit, one
can
look 
at 
the family of amoebas, 
$\text{Log}_t(\C_t)$,
 where $(\C_t)_{t\in\R_{>1}}$ 
 is 
a non-trivial family of complex curves. 

\begin{exa}\label{ex:AmoebaFamily}
Figure \ref{amoeba
  conic} depicts 
  the 
amoeba of the complex curve 
given by the equation $1-z-w+t^{-2} z^2   - t^{-1}zw+t^{-2}w^2=0$
for $t$ sufficiently large, 
and 
the limiting object which is... a tropical conic. 
\end{exa}
\begin{figure}[h]
\centering
\begin{tabular}{ccc}
\includegraphics[width=5cm, angle=0]{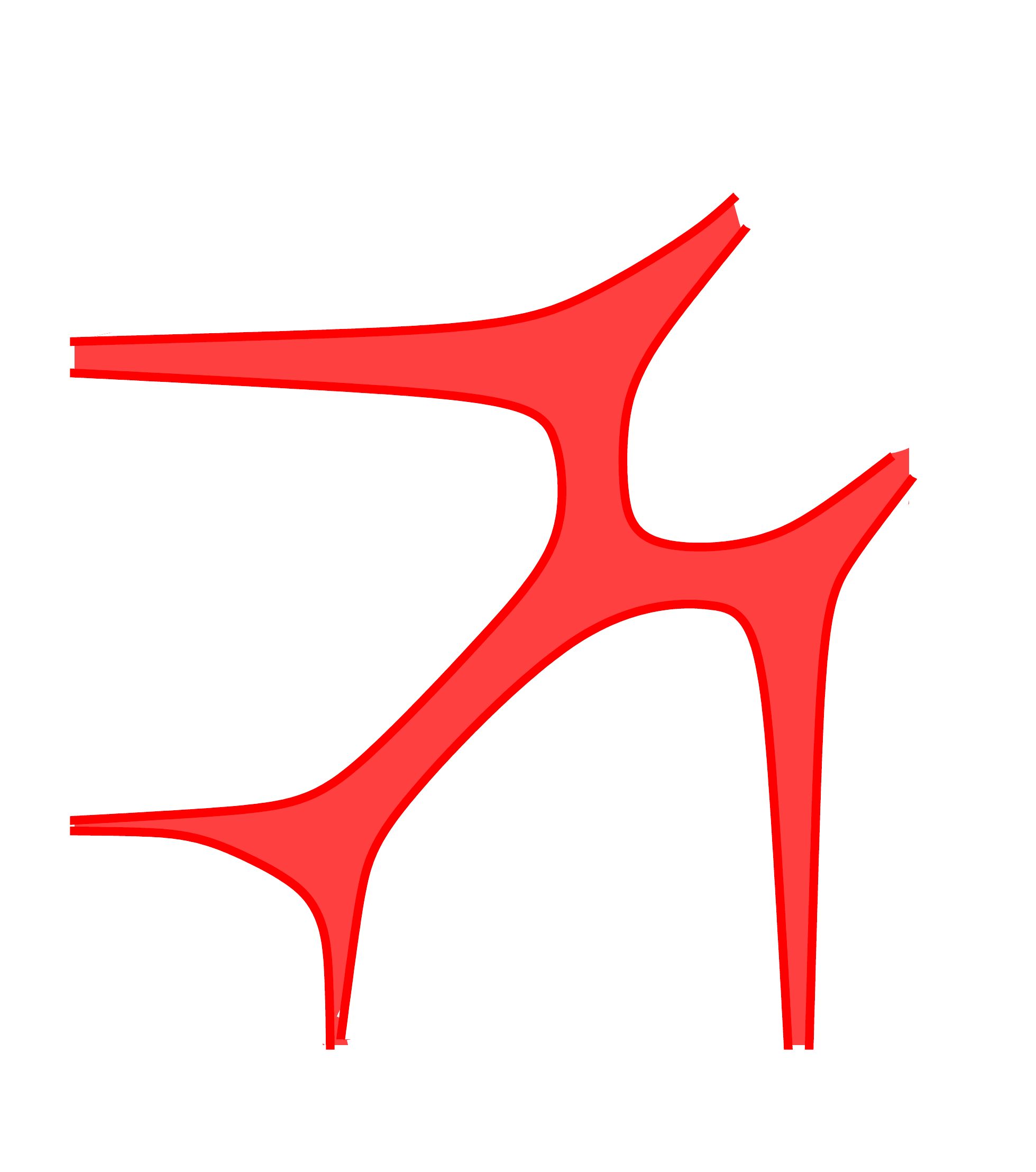}& \hspace{3ex} 
&\includegraphics[width=5cm, angle=0]{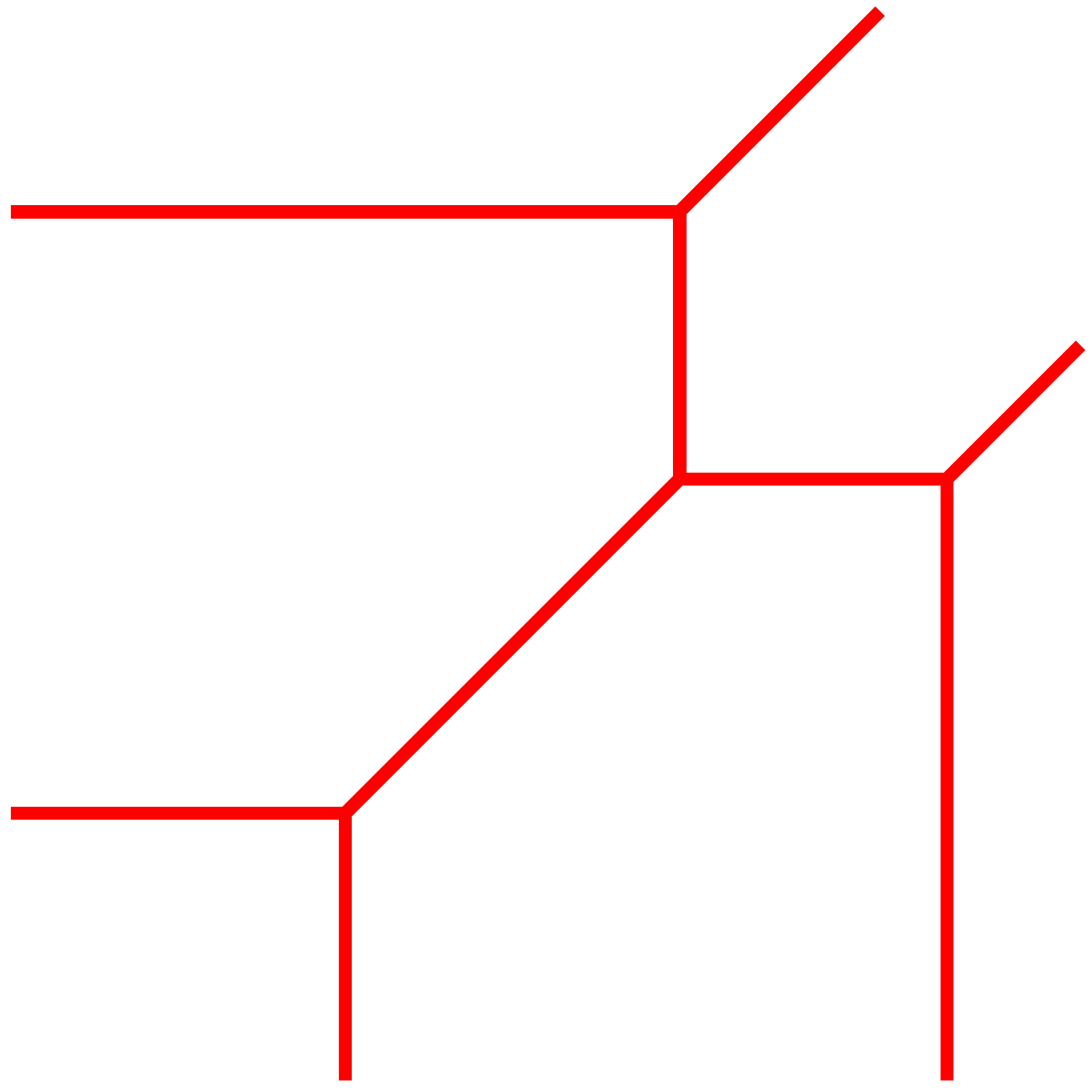}

\\
\\   $\text{Log}_t(\C_t)$ && $\lim_{t \to +\infty}\text{Log}_t(\C_t)$ 
\end{tabular}
\caption{$\C_t: 1-z-w+t^{-2} z^2   - t^{-1}zw+t^{-2}w^2=0$}\label{amoeba conic}
\end{figure}

The following statement is the two-dimensional counterpart of Theorem \ref{NP1}.

\begin{thm}\label{approx} 
{\rm (}cf. \cite{Kap1,Mik12}{\rm )}  
Let $P_t(z,w)=\sum_{i,j}\alpha_{i,j}(t)z^iw^j$ be a polynomial whose
coefficients are functions $\alpha_{i,j}:\RR \to \CC$, and suppose
that $\alpha_{i,j}(t)\sim \gamma_{i,j}t^{a_{i,j}}$ when
$t$ goes to $+\infty$ with $\gamma_{i,j}\in\CC^\times$ and
$a_{i,j}\in\TT$. 

If $\C_t$ denotes 
the curve in $(\CC^\times)^2$ defined by the polynomial $P_t(z,w)$, 
then  
the amoebas $\text{Log}_t(\C_t)$ converge to the tropical curve defined
by the tropical polynomial $P_{trop}(x,y)=\tg \sum_{i,j}a_{i,j}x^iy^j \td$.
\end{thm} 

It remains to explain the relation between amoebas and  weights
of a tropical curve.
Let $P_t(z,w)$ and $P'_t(z,w)$ be two families of complex
polynomials, defining two families of complex algebraic curves
$(\C_t)_{t\in\R_{>1}}$
and $(\C'_t)_{t\in\R_{>1}}$,
respectively. 
As in Theorem \ref{approx}, these two families of polynomials 
induce  two tropical
polynomials $P_{trop}(x,y)$ and $P'_{trop}(x,y)$, which in turn
define
two tropical curves
 $C$ and $C'$. 
 
\begin{prop}[cf. \cite{Mik8, PR}]\label{limit inter}
Let $p\in C\cap C'$ be a point which is not a vertex of $C$ or $C'$. 
Assume that the edges of $C$ and $C'$
which contain $p$ are not parallel. 
Then, 
the number of intersection points  of $\C_t$
and $\C'_t$ whose image under $\Log_t$ converges to $p$ is 
equal to 
the Euclidean area of the polygon $\Delta_p$ dual to
$p$ in the subdivision dual to $C\cup C'$. 
\end{prop} 
The above number is denoted by $(C \cdot C')_p$ and is called the
\emph{multiplicity} of the 
intersection point $p$ of $C$ and $C'$. 
It is worth 
noting 
that the number of intersection points
which converge to $p$ 
depends only on $C$ and $C'$,
that is, 
only on the order at infinity of the coefficients of $P_t(z,w)$
and $P'_t(z,w)$. 

If $e$ is an edge of a tropical curve $C$, and $p$ is an interior point of $e$, 
then the weight of $e$ is equal to the minimal multiplicity
$(C \cdot C')_p$ for all possible tropical curves $C' \ni p$
such that $(C \cdot C')_p$ is defined. 

\begin{figure}[h]
\centering
\begin{tabular}{ccc}
\includegraphics[height=3cm, angle=0]{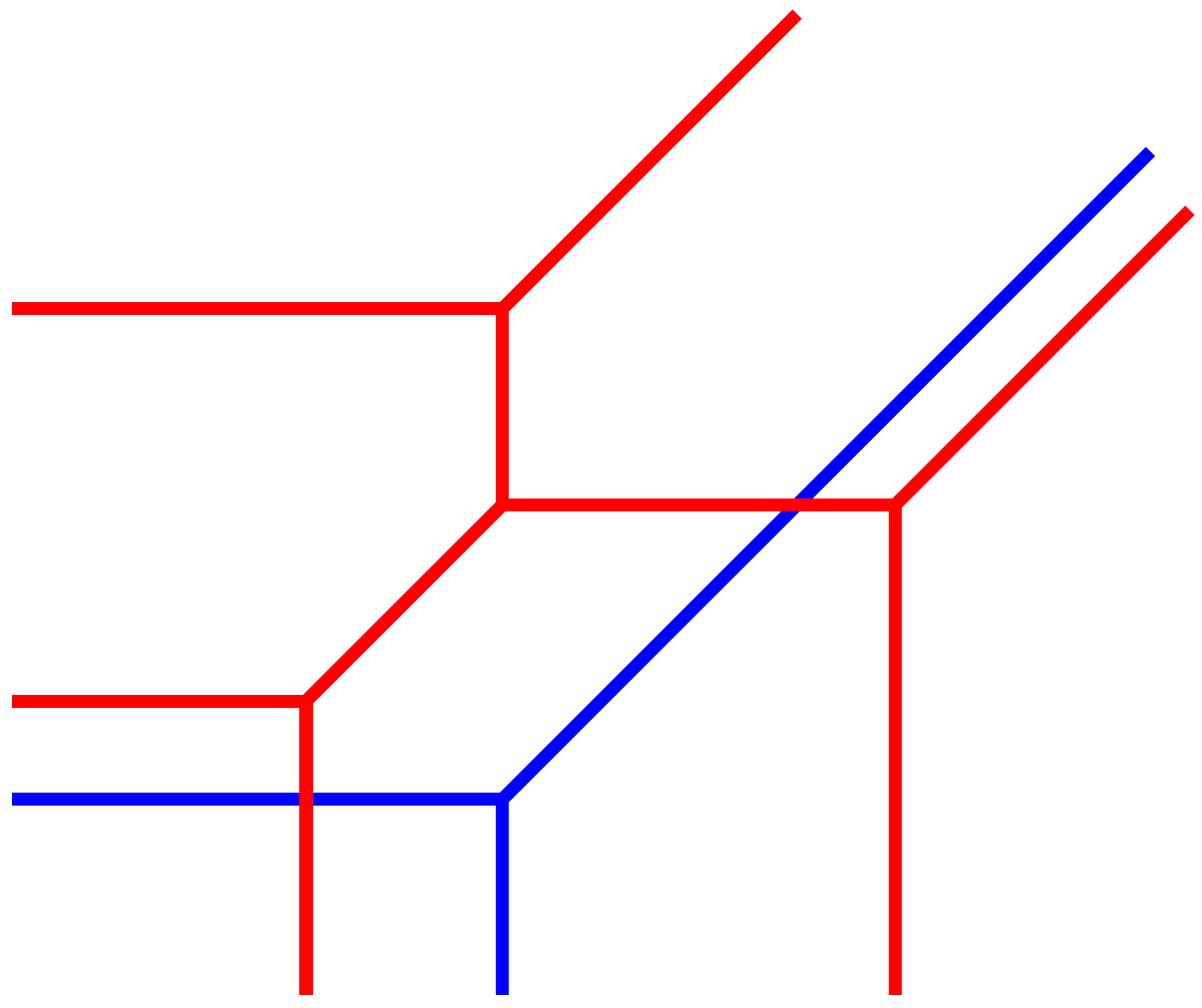}&\hspace{10ex} &
\includegraphics[height=3cm, angle=0]{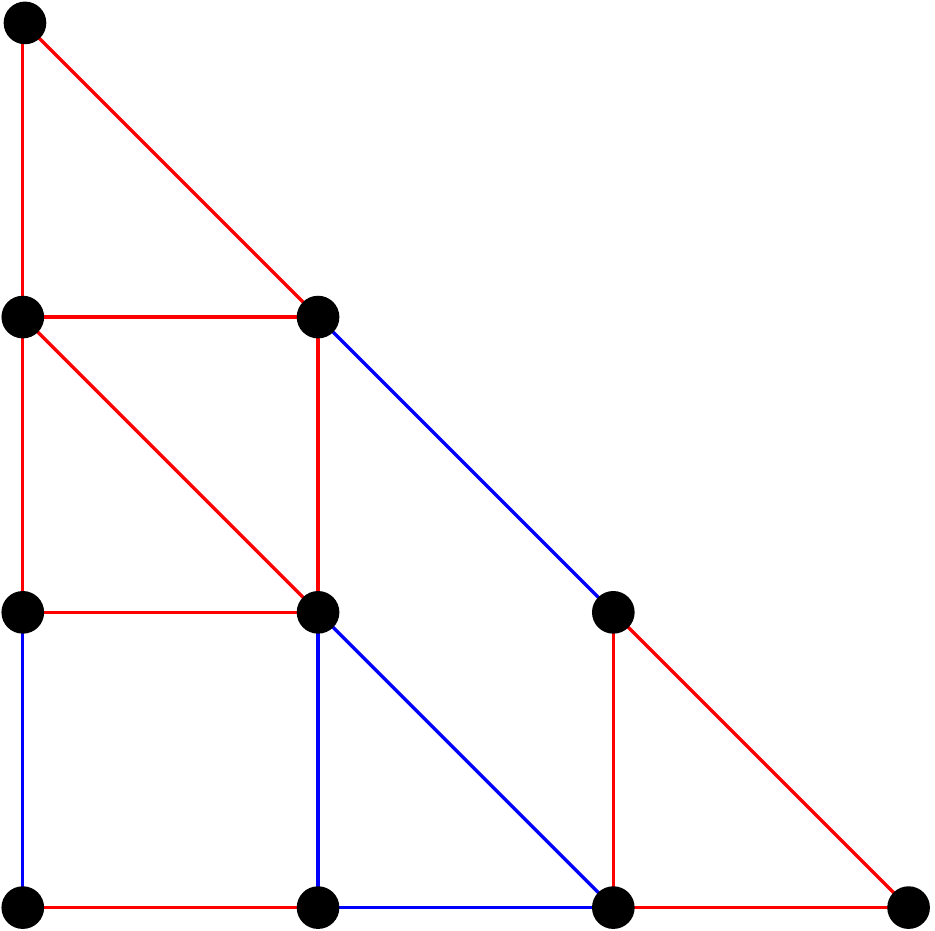}
\\ \\ a) && b)\end{tabular}

\bigskip
\begin{tabular}{ccc}
\includegraphics[height=3cm, angle=0]{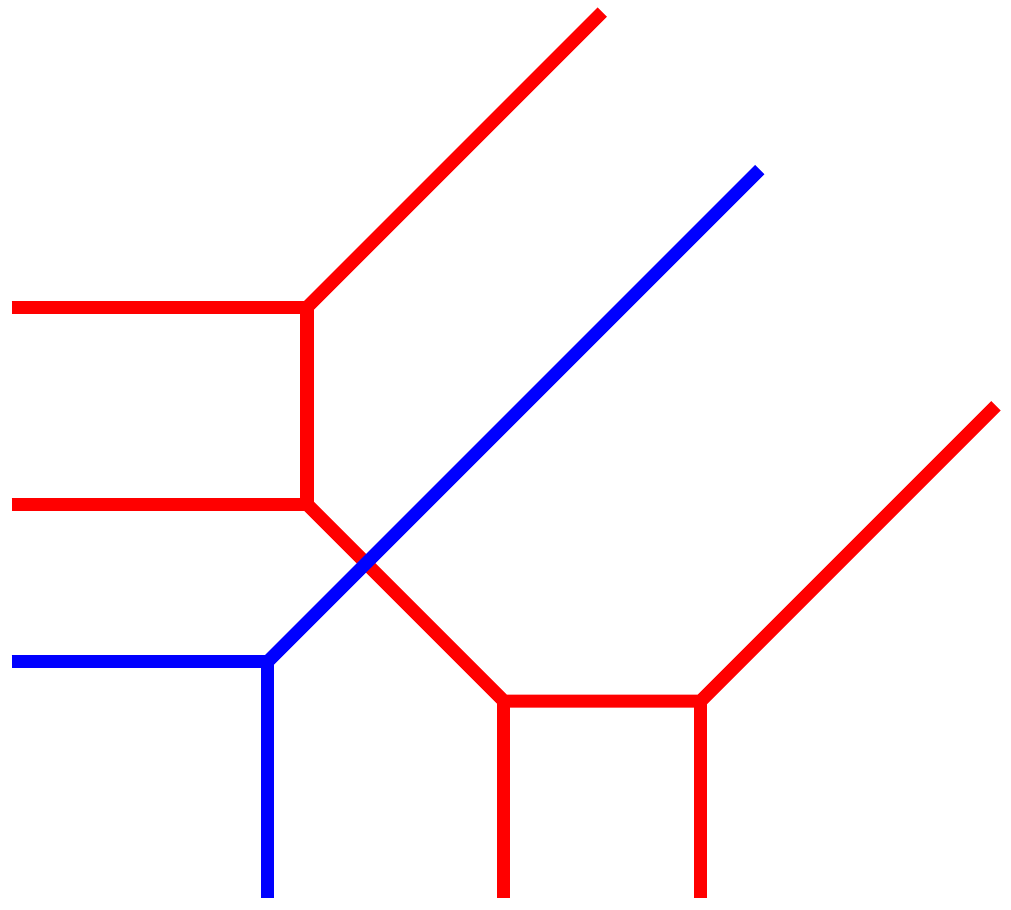}&\hspace{10ex} &
\includegraphics[height=3cm, angle=0]{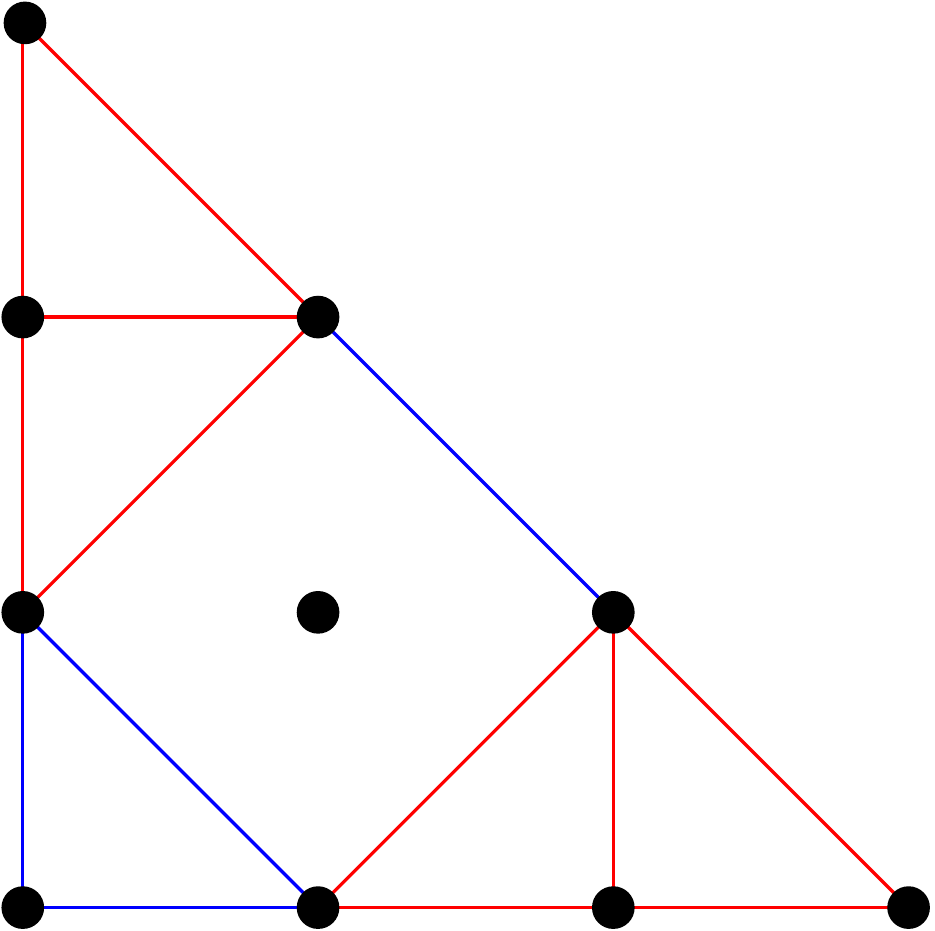}
\\ \\c) && d)
\end{tabular}
\caption{Tropical intersections}\label{inters amoeb}
\end{figure}

\begin{exa}
Figures \ref{inters amoeb}a,c
depict
different mutual
positions of a tropical line and a tropical conic. The corresponding
dual subdivisions of the union of the two curves 
are depicted in Figures
\ref{inters amoeb}b,d. 

In Figure \ref{inters amoeb}a (respectively, Figure \ref{inters amoeb}c),
the tropical line intersects the
tropical conic in two points of multiplicity $1$
(respectively, in one point of
multiplicity 2). 

\end{exa}

The combination of Theorem \ref{approx} and Proposition \ref{limit
  inter}
allows one, for example, to deduce the Bernstein Theorem \cite{Bern} in classical algebraic geometry
from the tropical Bernstein Theorem (see Exercise $2(4)$). 

\begin{exo}
\

\begin{enumerate}
\item Draw the tropical curves defined by the tropical polynomials 
$P(x,y)=\tg 5 + 5x + 5y + 4xy+1y^2+x^2\td$ and
$Q(x,y)=\tg  7+ 4x + y + 4xy+3y^2+(-3)x^2\td$, as well as 
the dual subdivisions of these tropical curves.

\item Show that a tropical curve defined by a polynomial of degree $d$ has at most $d^2$ vertices. 

\item Find an equation for each of the tropical curves in Figure \ref{equil}. 
The following reminder might be helpful: if $v$ is a vertex of a tropical curve defined by a 
tropical polynomial $P(x, y)$, then the value of $P(x, y)$ in a neighborhood of $v$ is given
uniquely by the monomials corresponding to the polygon dual to $v$. 

\item Prove the tropical Bernstein Theorem: let $C$ and $C'$ be two
  tropical curves such that $C\cap C'$ does not contain any vertex
  of
  $C$ 
  or $C'$;
  then, the sum of the 
  multiplicities 
  of all
  intersection points of $C$ and $C'$ is equal to 
the mixed volume of
  $\Delta(C_1)$ and $\Delta(C_2)$, {\it i.e.} to
  $$Area(\Delta(C\cup C')) - Area(\Delta(C)) - Area(\Delta(C')) .$$
Here, $Area(\Delta(C))$ is the Euclidean area of the Newton polygon
of $C$. 
Deduce the classical
Bernstein Theorem from its tropical counterpart.

\end{enumerate}
\end{exo}

\section{Patchworking}\label{sec:patchwork}
In this section, we present a first application of the material discussed above to real
algebraic geometry. The patchworking technique invented by O. Viro at
the end of the 1970's constitutes one the roots of tropical geometry. At
that time the formalism of tropical geometry did not exist yet, and the original
formulation of patchworking is dual to the presentation we give here.
Tropical formulation of patchworking became natural
with the introduction of amoebas 
due to I.M.~Gelfand, M.M.~Kapranov and A.V.~Zelevinsky \cite{GKZ}, 
further relations between amoebas and real algebraic curves 
were found by O.~Viro and G.~Mikhalkin (see \cite{V9} and \cite{Mik11}).
We discuss here only
a particular case of the general patchworking theorem.
This case, called {\it combinatorial patchworking}, 
turned out to be a powerful tool to construct plane real algebraic
curves (and, more generally, real algebraic hypersurfaces of toric
varieties). 

\medskip
A \emph{real algebraic curve  in $(\CC^\times)^2$} is an
algebraic curve defined by a polynomial with real coefficients.
Given such a real algebraic curve  $\C$, we denote by $\RR\C$ the set
of real points of $\C$, 
{\it i.e.}  $\RR\C=\C\cap (\RR^\times)^2$.

\subsection{Patchworking of a line}

\begin{figure}[h]
\centering
\begin{tabular}{ccc}
\includegraphics[height=4.5cm, angle=0]{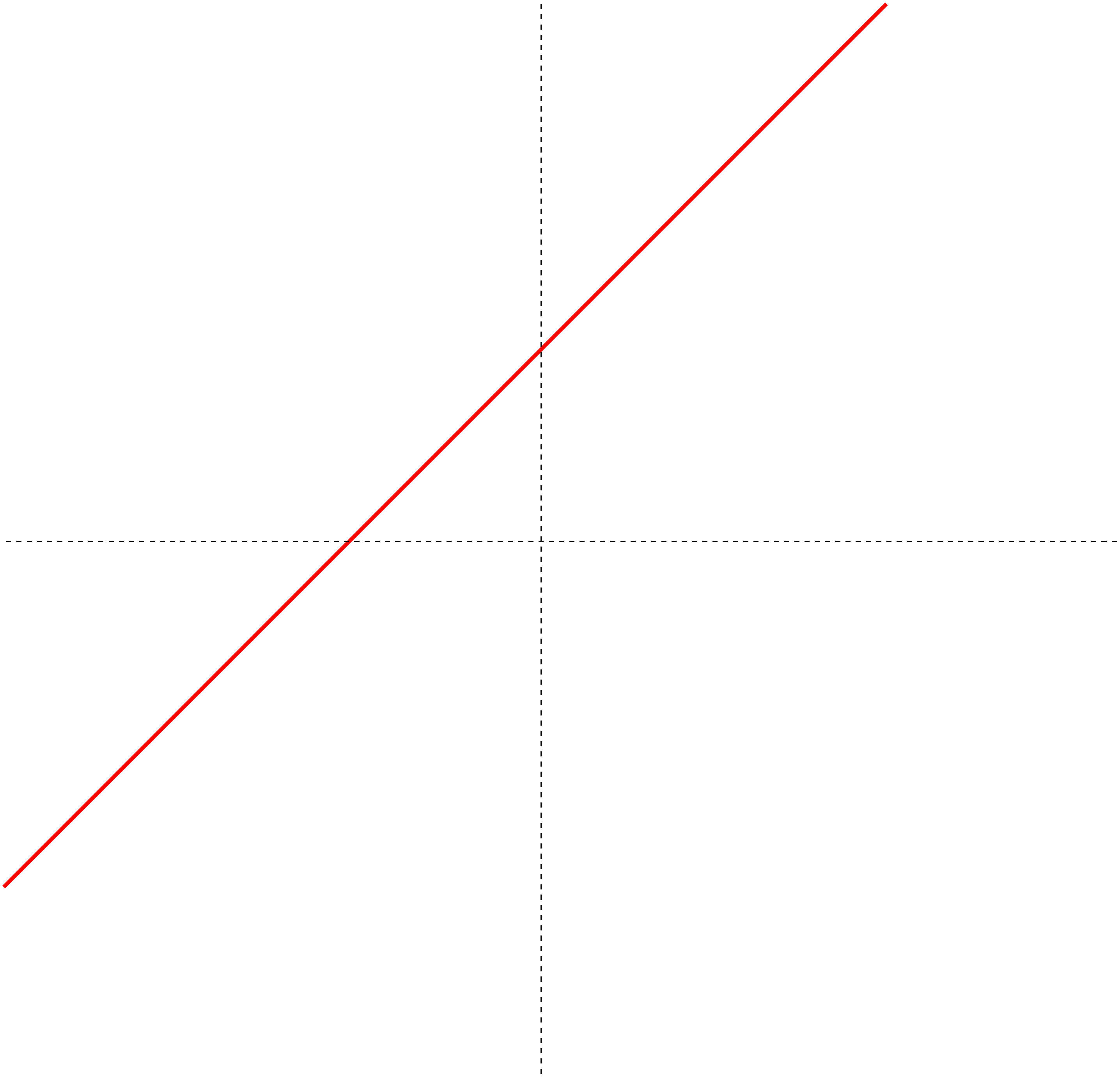}& 
&
\includegraphics[height=4.5cm, angle=0]{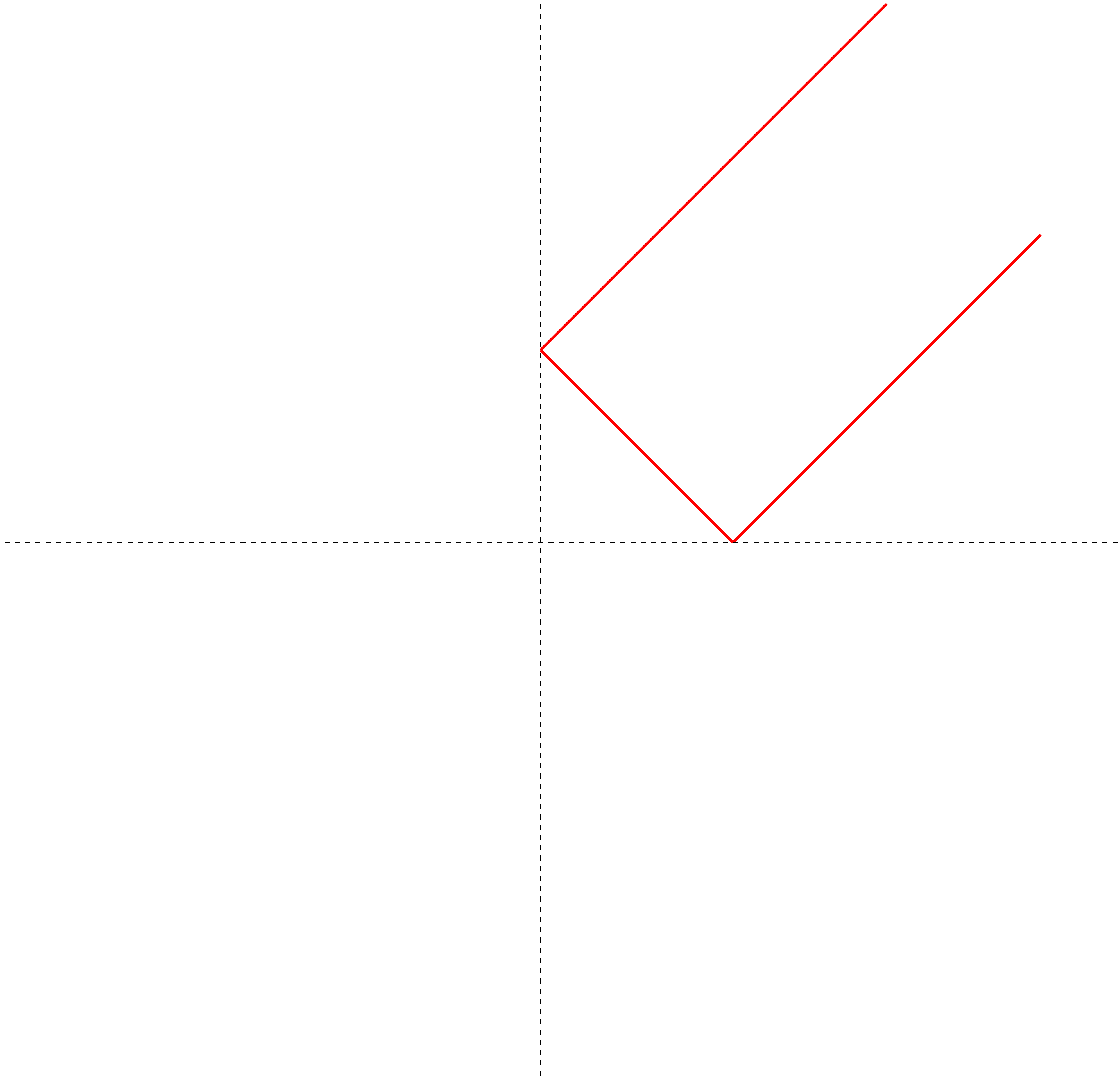}
\\  a) $\L: z-w+1=0$&&
b)
\\ \\ 
\includegraphics[height=4.5cm,
  angle=0]{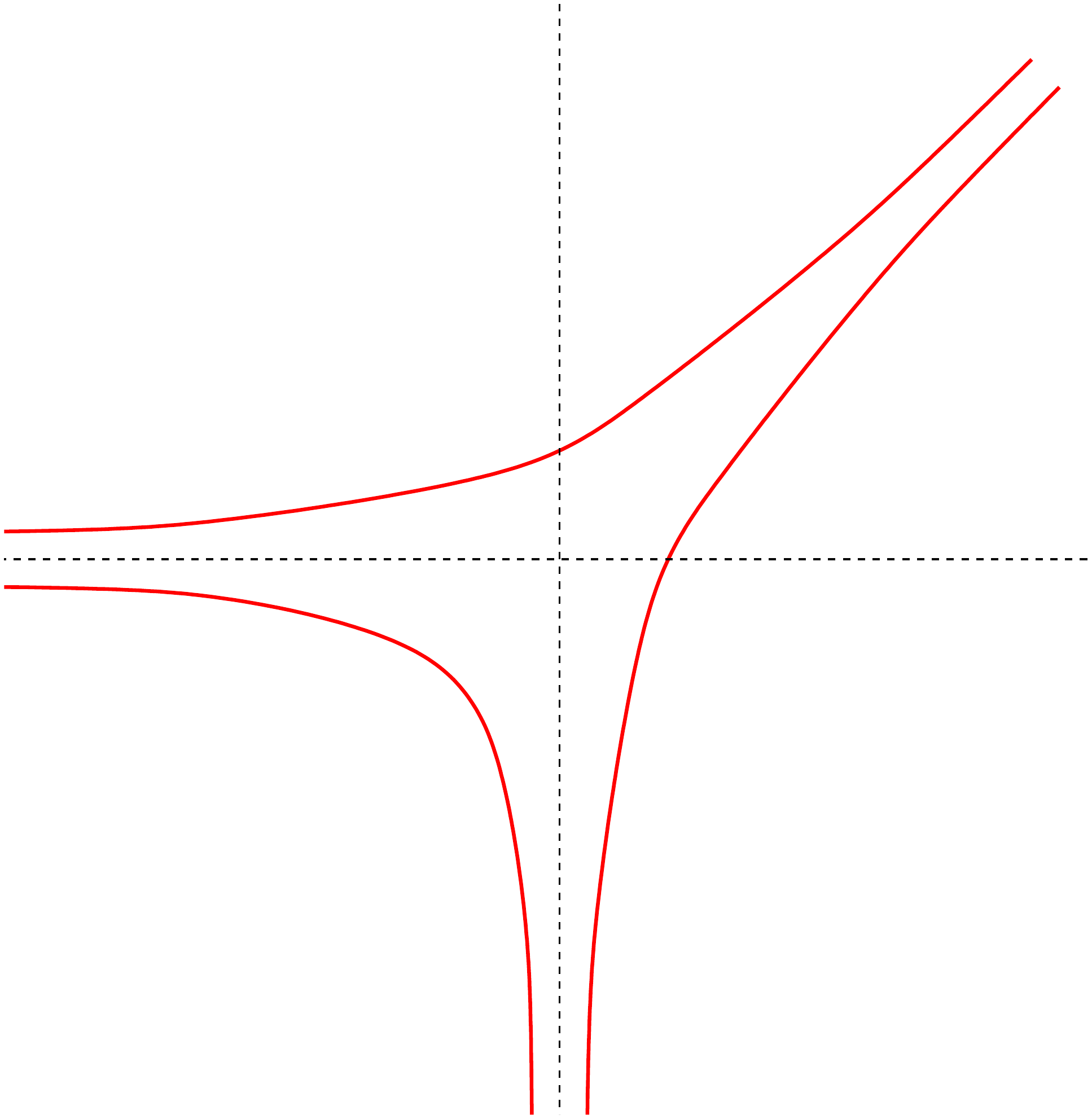}& \hspace{5ex} &
\includegraphics[height=4.5cm, angle=0]{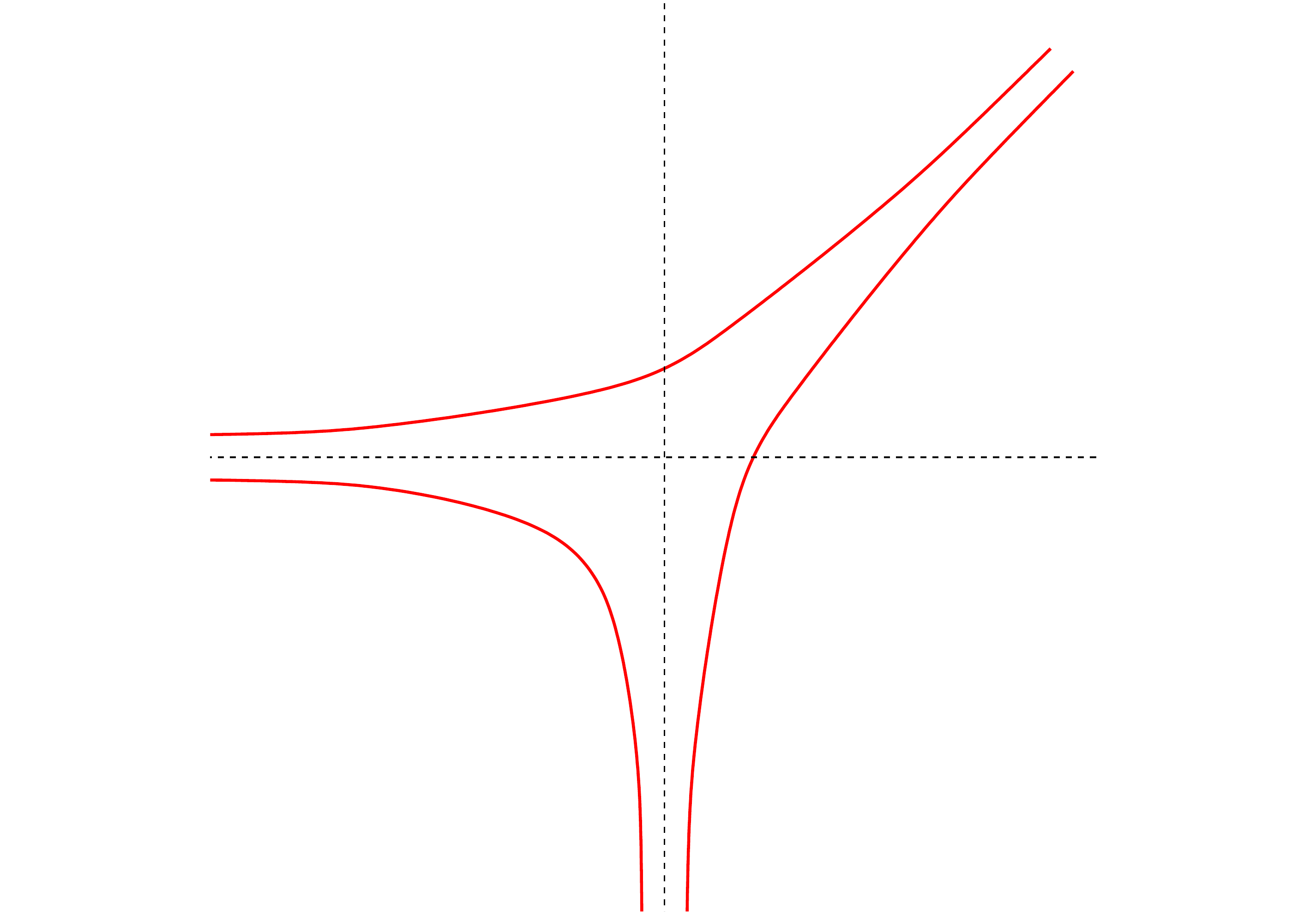}
\put(-165, 25){\small{$(-\epsilon(ac),-\epsilon(bc))$}}
\put(-80, 30){\small{$(-\epsilon(ac),\epsilon(bc))$}}
\put(-143, 85){\small{$(\epsilon(ac),-\epsilon(bc))$}}
\\ 
c) $\A(\RR\L)$&&
d) $\A(\RR\L)$ with  signs
\end{tabular}
\caption{Amoeba of a real line}\label{fig:real amoeba}
\end{figure}

Let us start by looking more closely at the amoeba of the real
algebraic line $\L\subset(\R^\times)^2$ given by the equation $az+bw+c=0$
with $a,b,c\in\RR^\times$. 
The whole amoeba 
$\A(\L)$ is depicted in Figure 
 \ref{amoeba line}a, and the amoeba of $\RR\L$ is depicted in Figure
 \ref{fig:real amoeba}c. Note that 
$\A(\L)$ does not depend on $a,b,$ and $c$ up to 
translation in
 $\RR^2$, and that
$\partial \A(\L)=\A(\RR \L)$.

We 
can label each arc of $\A(\RR \L)$ by the pair of signs
corresponding to the quadrant of $(\RR^\times)^2$ through which the
corresponding arc of 
$\RR\L$
passes (see Figure 
 \ref{amoeba line}d,
where $\epsilon(x)$ denotes the sign of $x$). 
This labeling only depends on the signs of
 $a,b,$ and $c$. Moreover, if two arcs of $\A(\RR \L)$
have an asymptotic direction $(u,v)$ in common, then these pairs
of signs differ by a factor $((-1)^u,(-1)^v)$. 

\begin{figure}[h]
\centering
\begin{tabular}{ccccc}
\includegraphics[height=3.5cm, angle=0]{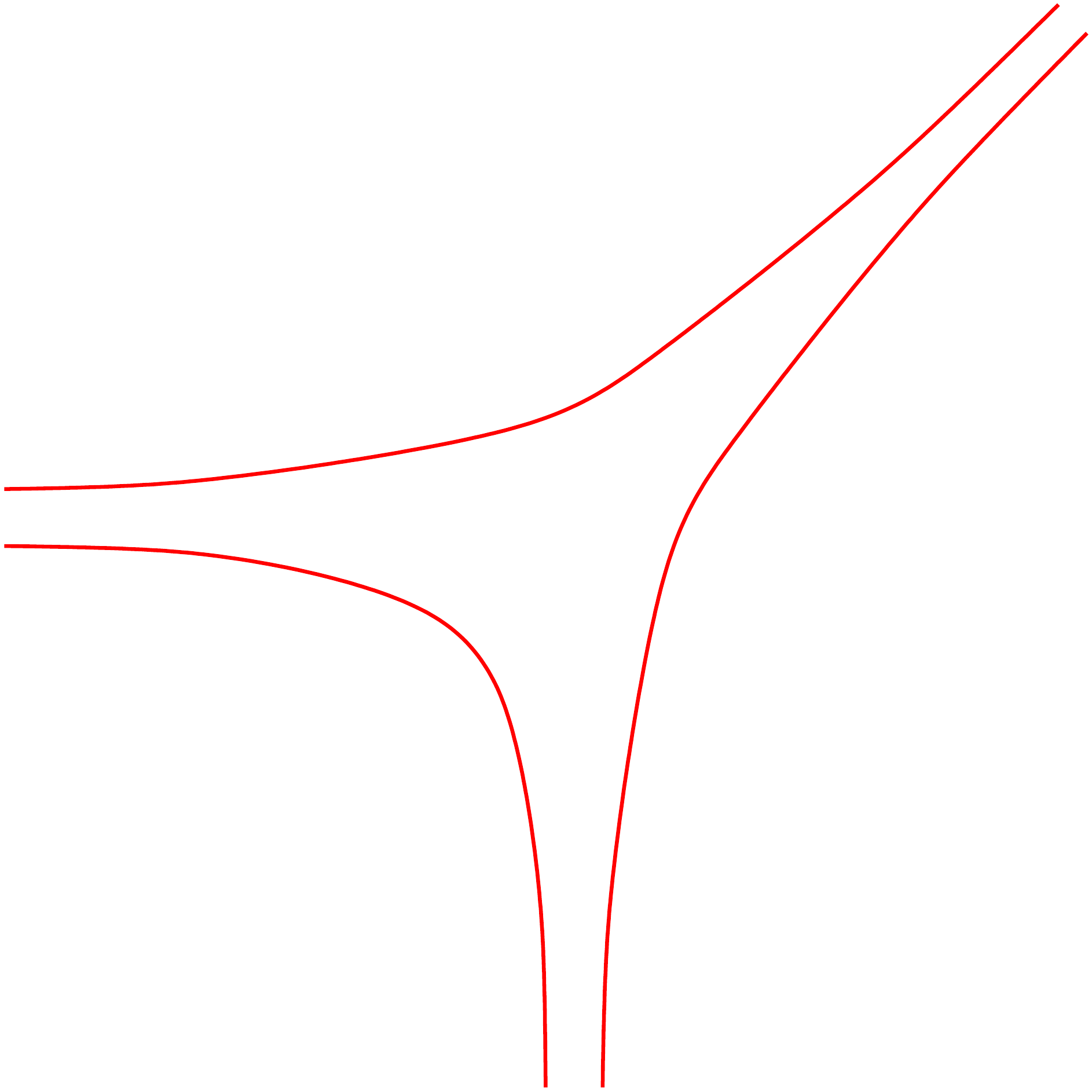}
\put(-80, 65){\small{$(+,+)$}}& \hspace{3ex} &
\includegraphics[height=3.5cm,angle=0]{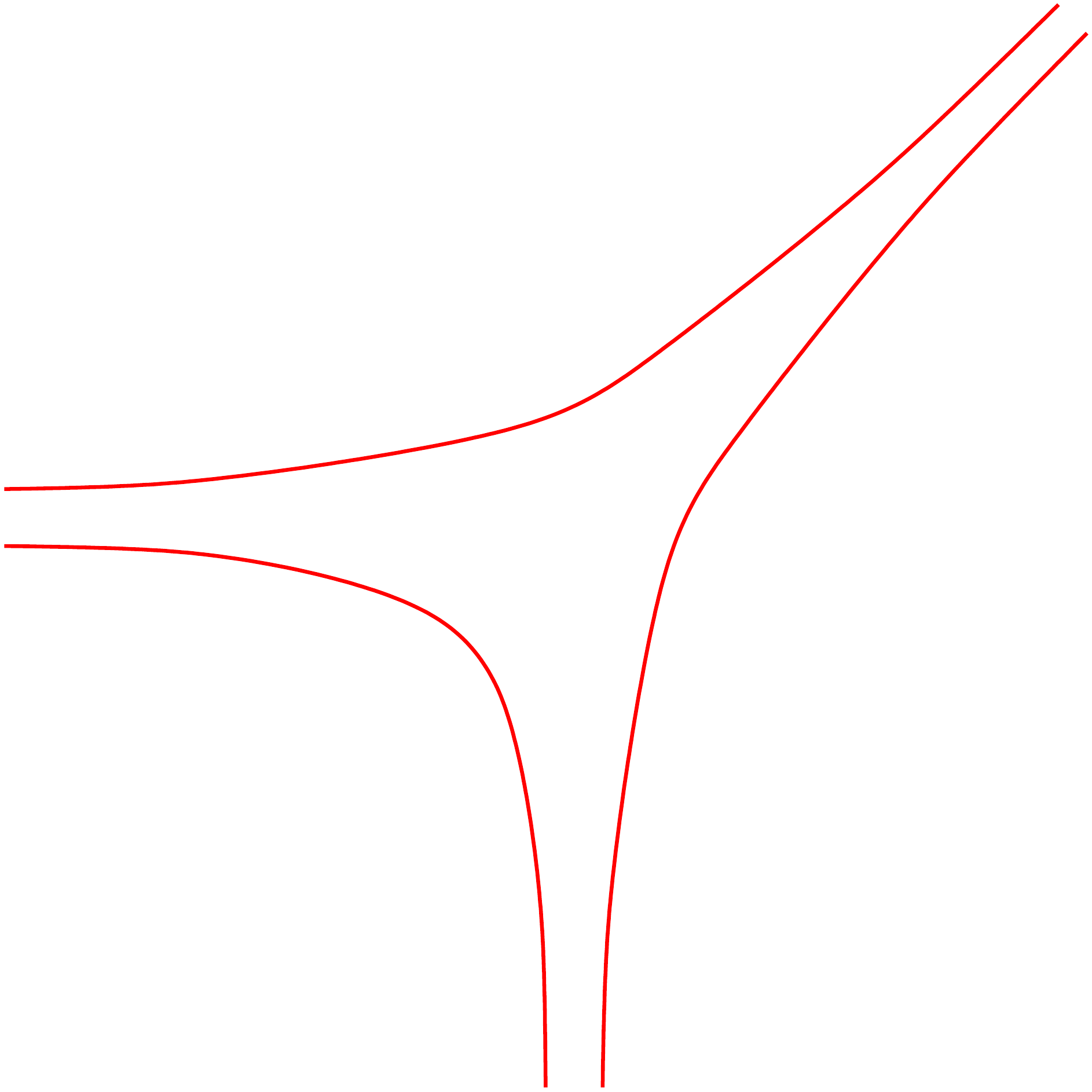}
\put(-80, 65){\small{$(+,+)$}}
\put(-80, 25){\small{$(-,+)$}}
\put(-40, 25){\small{$(-,-)$}}&  \hspace{3ex} &
\includegraphics[height=3.5cm, angle=0]{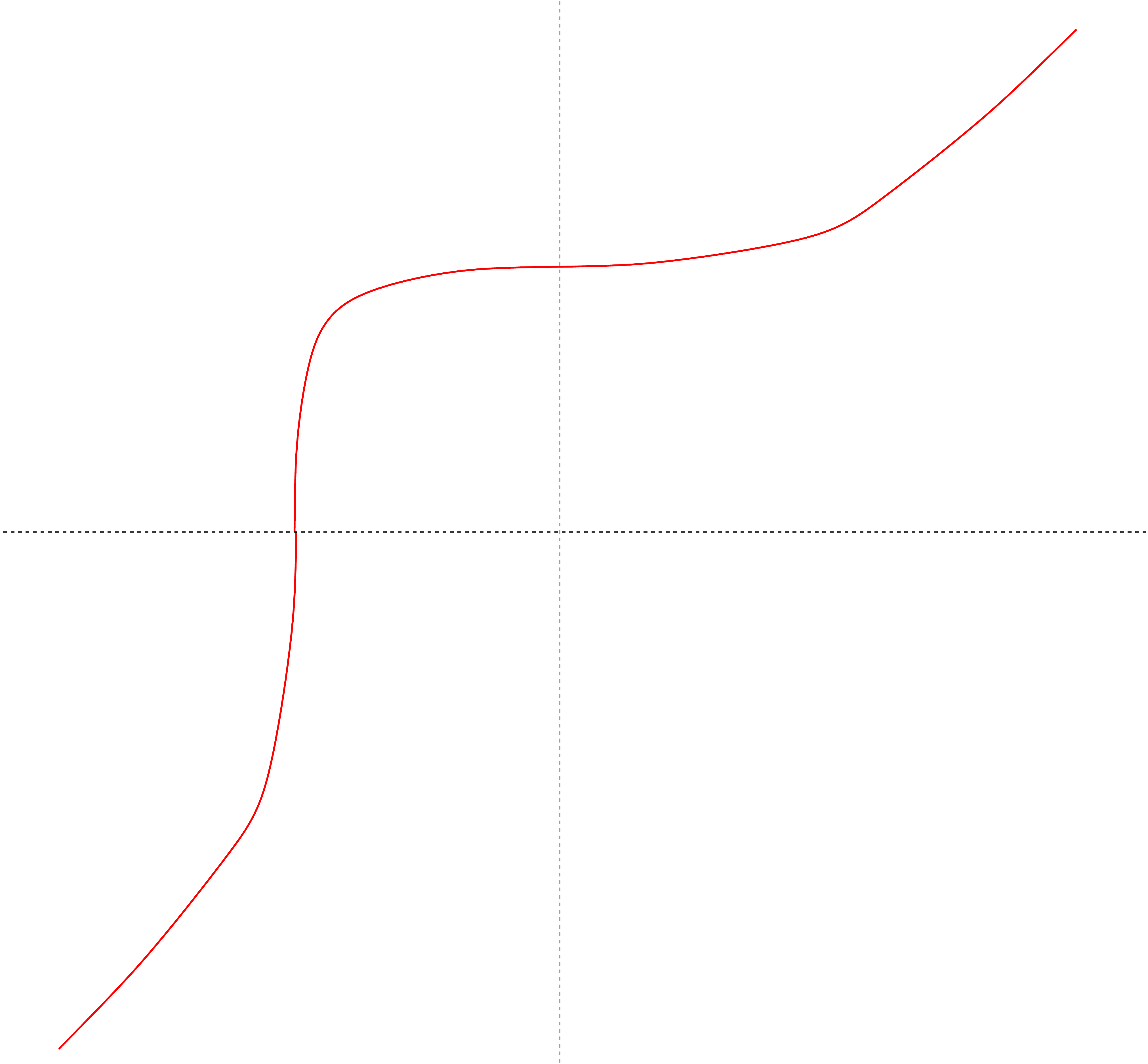}
\\ a) && b) && c)
\end{tabular}
\caption{Patchworking of a line}\label{fig:patch line}
\end{figure}

We may recover from $\A(\RR \L)$ the isotopy class (up to axial symmetries) of $\RR\L$ in $(\RR^\times)^2$. 
To do this, we assign a pair of signs to some
arc of $\A(\RR \L)$ (Figure \ref{fig:patch line}a).
As we have seen, this determines a pair of signs for the two other
arcs  of 
$\A(\RR \L)$ (Figure \ref{fig:patch line}b).
For an arc $A\subset\A(\RR \L)$ labeled by $((-1)^{\epsilon_1},(-1)^{\epsilon_2})$
we draw its image under the map $(x,y)\mapsto
((-1)^{\epsilon_1}e^x,(-1)^{\epsilon_2}e^y)$ in $(\RR^\times)^2$. 
Clearly, the image of $A$ is contained in the $((-1)^{\epsilon_1},(-1)^{\epsilon_2})$-quadrant.
The union of the images of the three arcs (shown on Figure \ref{fig:patch line}c)
is isotopic to a straight line
(shown on Figure \ref{fig:real amoeba}a).

\subsection{Patchworking of non-singular tropical curves}
Viro's patchworking theorem is a generalization of the previous 
observations, 
in the case of an approximation of a \emph{non-singular tropical
  curve} by a family of amoebas of real algebraic curves.

\begin{defi}
A tropical curve in $\RR^2$ is \emph{non-singular} if its dual
subdivision 
is formed by 
triangles of Euclidean area $\frac{1}{2}$.
\end{defi} 
Equivalently, a tropical curve is non-singular if and only if it
has exactly 
$2\operatorname{Area}(\Delta(C))$ 
vertices. 
Recall that a triangle  
with vertices in $\ZZ^2$
and having Euclidean area $\frac{1}{2}$ can be mapped, via the composition of a 
 translation and an element of $SL_2(\Z)$, to the triangle with
 vertices $(0,0)$, $(1,0)$, $(1,1)$. In other words, an algebraic
 curve in $(\CC^\times)^2$ with Newton polygon 
of Euclidean area $\frac{1}{2}$
is nothing else but
 a line in suitable coordinates. 
We use this observation in the following construction.

\medskip

\begin{figure}[h]
\centering
\begin{tabular}{ccc}
\includegraphics[height=2cm,
  angle=0]{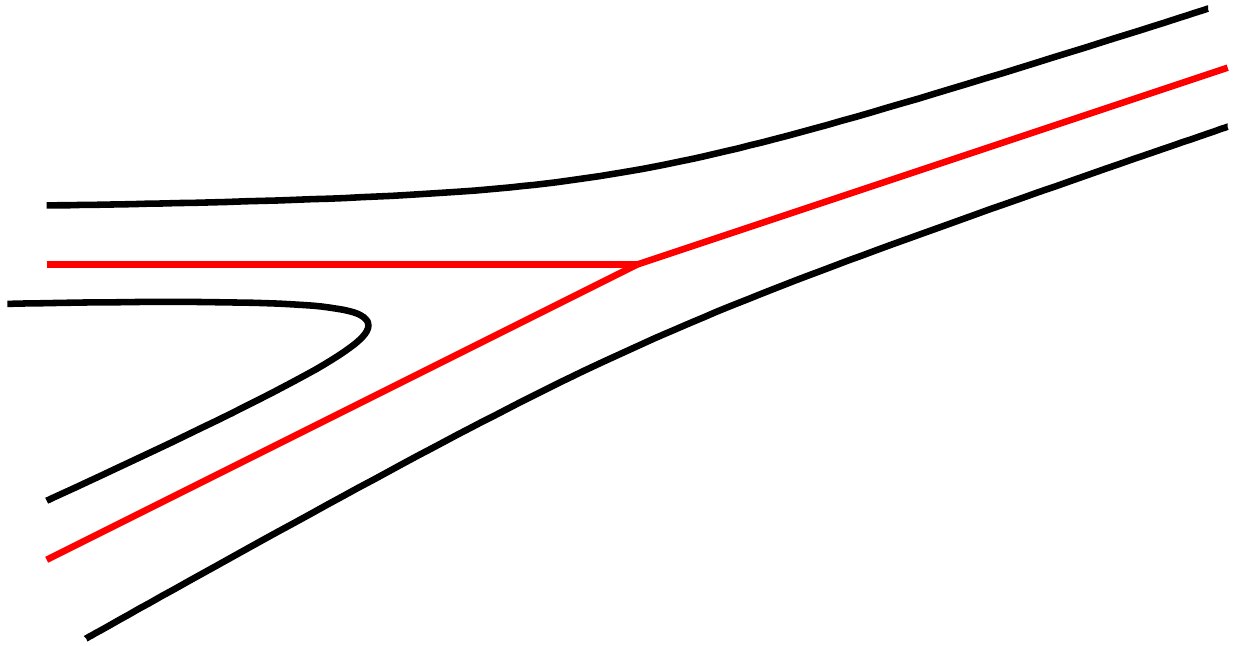}& \hspace{10ex} &
\includegraphics[height=3cm, angle=0]{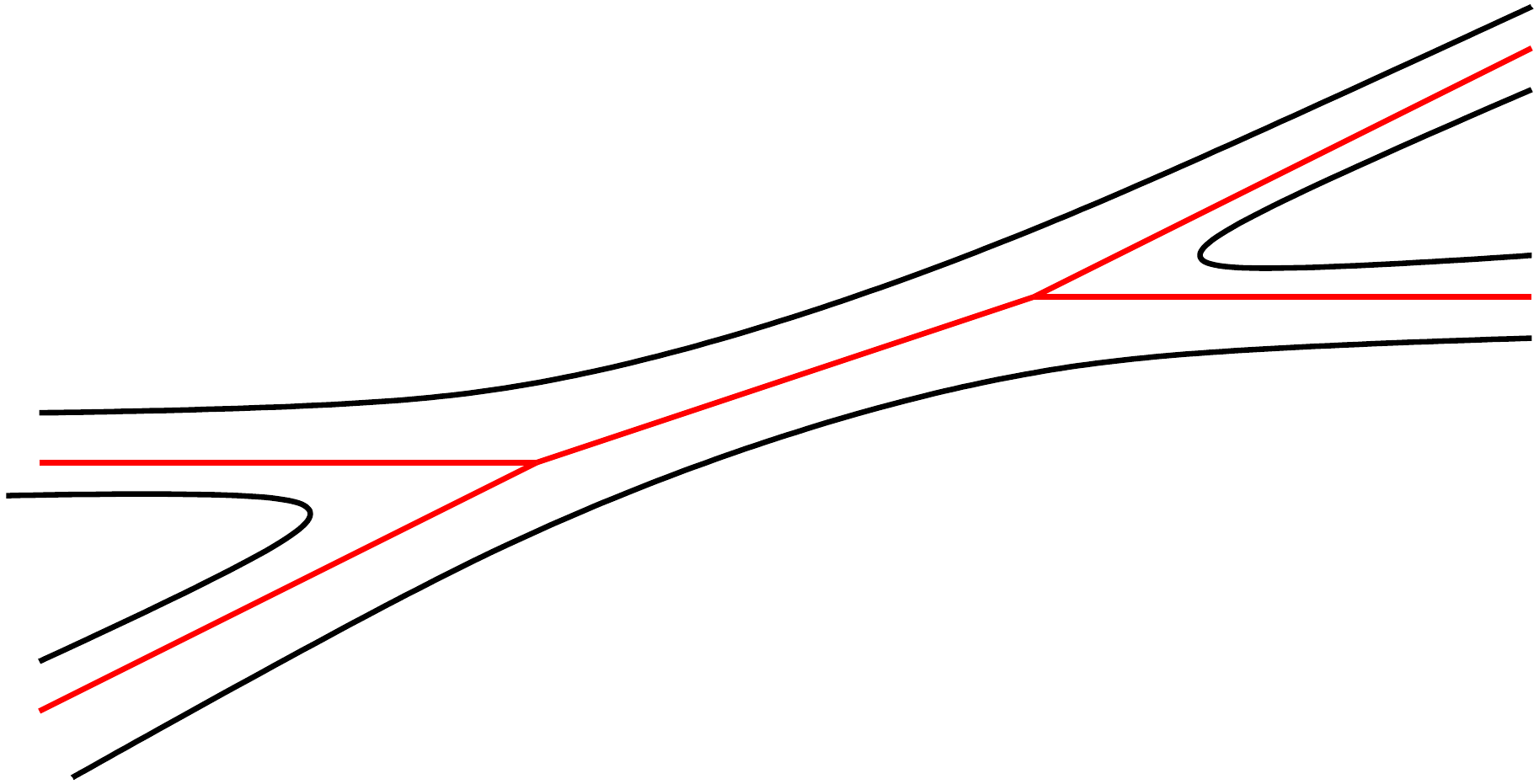}
\put(-150, 55){\small{$(\epsilon_1,\epsilon_2)$}}
 \put(-130, 5){\small{$((-1)^u\epsilon_1,(-1)^v\epsilon_2)$}}
\\  a)&& b)
\end{tabular}
\bigskip

\begin{tabular}{c}
\includegraphics[height=3cm, angle=0]{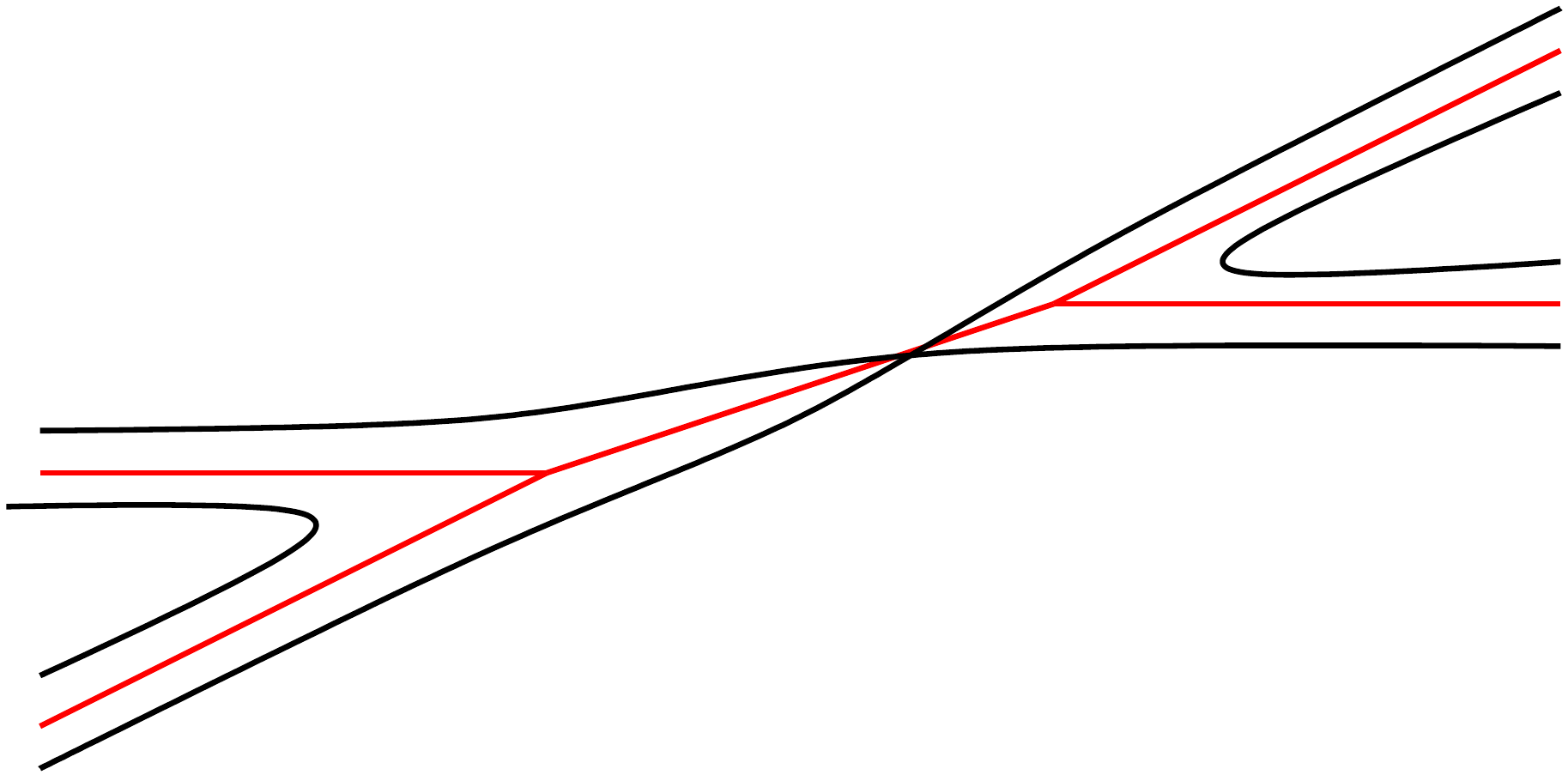}
\put(-150, 55){\small{$(\epsilon_1,\epsilon_2)$}}
 \put(-130, 5){\small{$((-1)^u\epsilon_1,(-1)^v\epsilon_2)$}}
\\ c) 
\end{tabular}
\caption{$\A_t(\RR \C_t)\cap U$ for large $t$}\label{fig:patch gen}
\end{figure}

Let $C$ be a non-singular tropical curve in $\RR^2$. 
Let $(\C_t)_{t\in\R_{>1}}$ be a family of real algebraic curves whose amoebas approximate $C$
in the sense of Theorem \ref{approx}.
Then, one can show that when $t$ is sufficiently large, 
the following hold:
\begin{itemize} 

\item for any vertex $a$ of $C$, in a small neighborhood $U_a$ of $a$ 
the amoeba $\A_t(\RR \C_t)\cap U_a$ is made of 
three 
arcs
as depicted in Figure \ref{fig:patch gen}a, corresponding to 
three 
arcs on 
$\RR \C_t$; 

\item for each bounded edge $e$ of $C$, 
in a small neighborhood $U_e$ of $e$ 
the amoeba $\A_t(\RR \C_t)\cap U_e$ is made of four arcs, corresponding to four arcs on 
$\RR \C_t$. The position of the arcs with respect to the edge is 
as depicted in either Figure \ref{fig:patch gen}b or c. Moreover, if $e$ has primitive integer direction $(u,v)$, then 
the two arcs of $\A_t(\RR \C_t)\cap U_e$ converging to $e$ correspond to 
arcs
of $\RR \C_t$ contained in quadrants of $(\RR^\times)^2$ whose corresponding
pairs
of signs 
differ by a factor $((-1)^u,(-1)^v)$.
\end{itemize}

The above 
two properties 
imply that 
the position of $\RR\C_t$ in $(\RR^\times)^2$
up to the action of
$(\ZZ/2\ZZ)^2$
by axial symmetries 
$(z,w)\mapsto (\pm z,\pm w)$
is entirely determined by the
partition of the edges of $C$ between the two types of
 edges depicted in Figures
 \ref{fig:patch gen}b, c.

\begin{defi}
An edge of $C$ as in Figure \ref{fig:patch gen}c is said to be \emph{twisted}.
\end{defi}

Not any subset of the edges of $C$  
may arise
as the 
set of twisted edges.
Nevertheless, the possible distributions of twists  are easy to describe.

\begin{defi}
A subset $T$ of the set of bounded edges of $C$ is called {\em twist-admissible}
if they satisfy 
the following
condition:

\noindent for any cycle
$\gamma$ of $C$, if $e_1,\ldots, e_k$  are the edges in $\gamma\cap
T$, and if $(u_i,v_i)$ is a primitive integer vector in the direction of $e_i$, then 
\begin{equation}\label{Men}
\sum_{i=1}^k (u_i,v_i)=0\quad \text{mod }2.
\end{equation}
\end{defi}

The Viro patchworking theorem \cite{V9}
may be reformulated in terms of twist-admissible
sets as follows.
\begin{thm}
\label{thm:viro}
For any twist-admissible set $T$ in a non-singular tropical curve $C$ in $\RR^2$,
there exists 
a family of non-singular real algebraic curves
$(\C_t)_{t\in\R_{>1}}$ in $(\CC^\times)^2$ which converges to $C$ in the sense of Theorem 
\ref{approx}  
and such that the corresponding set of twisted edges  is $T$. 
\end{thm} 

\begin{exa}
One may choose $T$ to be empty as 
the empty set clearly satisfies 
\eqref{Men}. 
The resulting curve corresponds to the
construction of  \emph{simple Harnack curves} described in \cite{Mik11}
via Harnack distribution of signs, see \cite{IV2}.
\end{exa}

\begin{exa}
If $C$ is a tree
then any set of edges is twist-admissible.
\end{exa}

\begin{figure}[h]
\centering
\begin{tabular}{ccc} 
\includegraphics[width=4cm, angle=0]{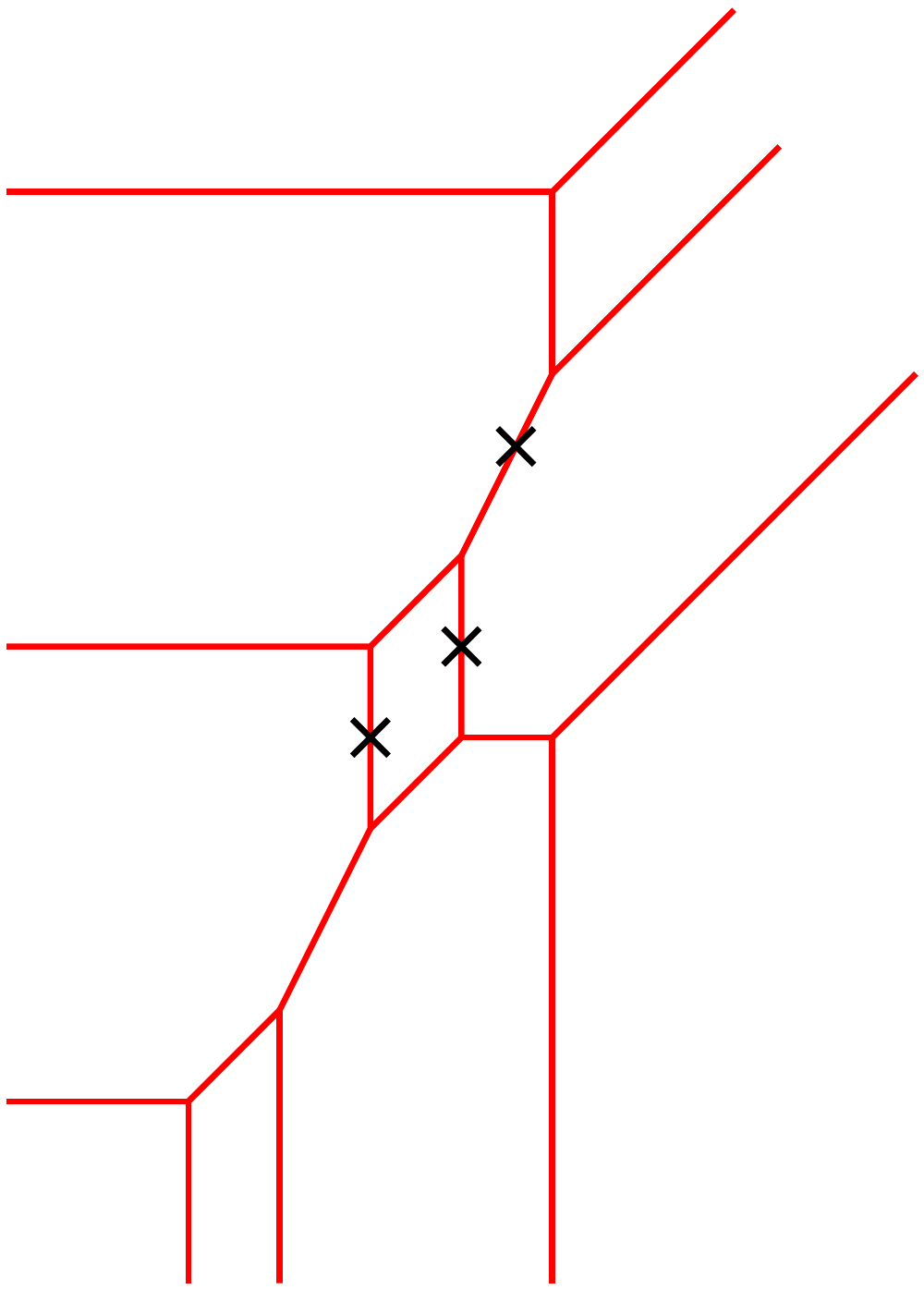}& 
\hspace{4ex} &
\includegraphics[width=4cm, angle=0]{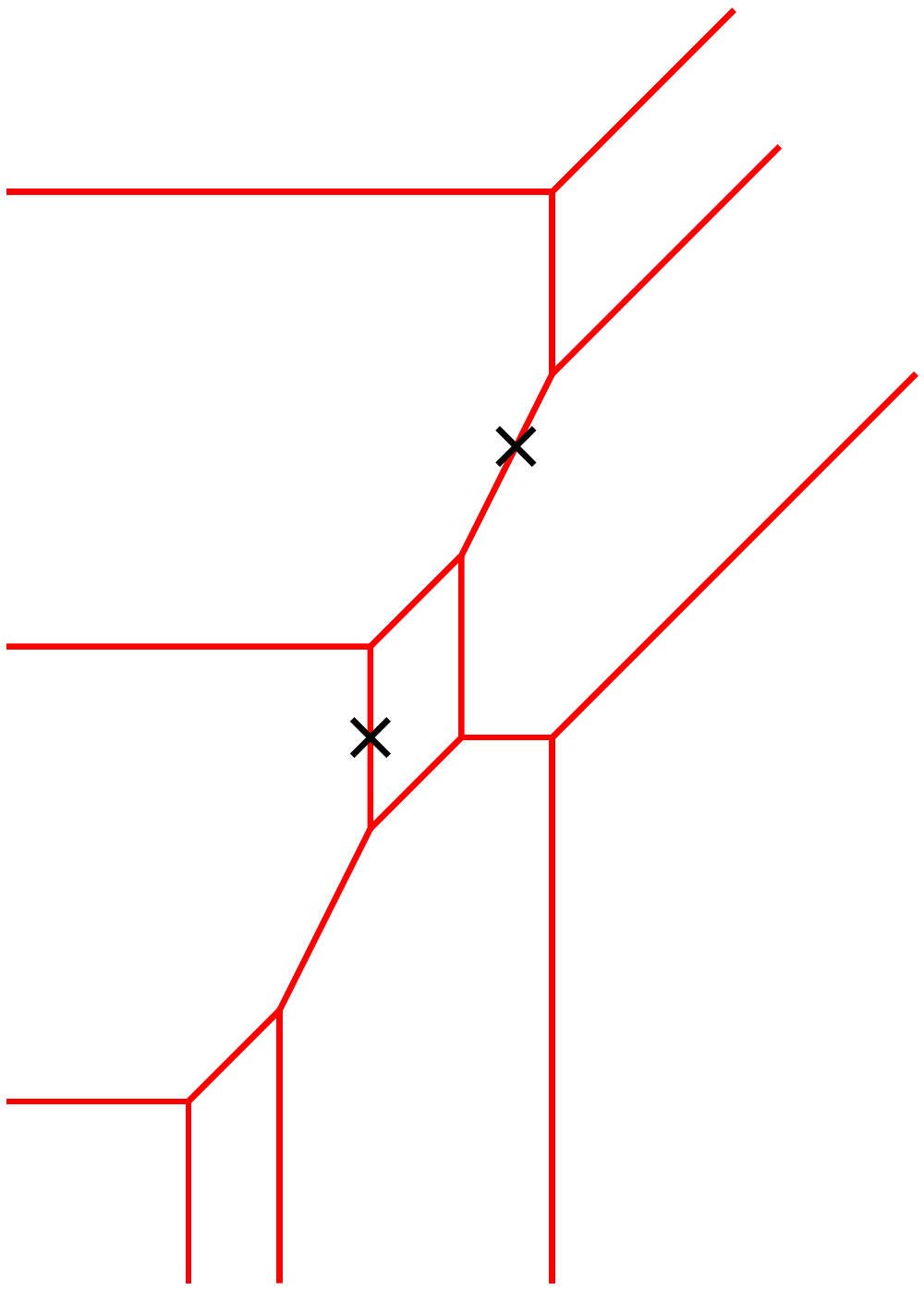}

\\  a) Twist-admissible set&&
b)   Not a twist-admissible set 
\end{tabular}
\caption{}\label{fig:patch ex}
\end{figure}

\begin{exa}
Consider the tropical cubic depicted in Figure \ref{equil}a, and 
choose
two subsets $T$ of the set of edges (marked by a cross) of $C$ 
as depicted in Figures \ref{fig:patch ex}a, b. The first one
is twist-admissible,  
while the second is not. 
\end{exa}

Below is a summary
of the procedure recovering the isotopy type of $\RR\C_t\subset (\RR^\times)^2$
(up to axial symmetries)  
from a smooth tropical curve $C\subset\R$ and a twist-admissible
set $T$ of edges in $C$.

\begin{itemize}
\item At each vertex of $C$, we draw three arcs as depicted in Figure
  \ref{fig:patch gen}a.

\item For each bounded edge $e$ adjacent to the vertices $v$ and $v'$
we join the two corresponding arcs at $v$  to the
  corresponding ones for $v'$ in the following way: if $e\notin T$, then join these arcs
 as depicted in Figure
  \ref{fig:patch gen}b; if $e\in T$, then join these arcs
 as depicted in Figure
  \ref{fig:patch gen}c; 
denote by $\P$ the curve obtained.

\item We choose arbitrarily a connected component of $\P$
and a pair of signs for it. (Different choices produce
the same isotopy type in $(\RR^\times)^2$ up to axial symmetries.)

\item We associate pairs of signs to all connected components of $\P$ using the
  following rule.  Given an edge 
$e$ with primitive integer direction $(u,v)$, 
the pairs
of signs of the two connected components of $\P$ corresponding to $e$ 
differ by a factor $((-1)^u,(-1)^v)$.
(The compatibility condition
\eqref{Men} ensures that this rule is consistent.)

\item We map each connected component $A$ of $\P$ to $(\RR^\times)^2$
by $(x,y)\mapsto (\epsilon_1 e^x, \epsilon_2 e^y)$, where 
$(\epsilon_1, \epsilon_2)$ is the pair of signs  
associated to $A$.
The resulting curve 
is the union of these images over all connected components
of $\P$.
\end{itemize}

\begin{figure}[h]
\centering
\begin{tabular}{cccccc}
\includegraphics[height=3.5cm, angle=0]{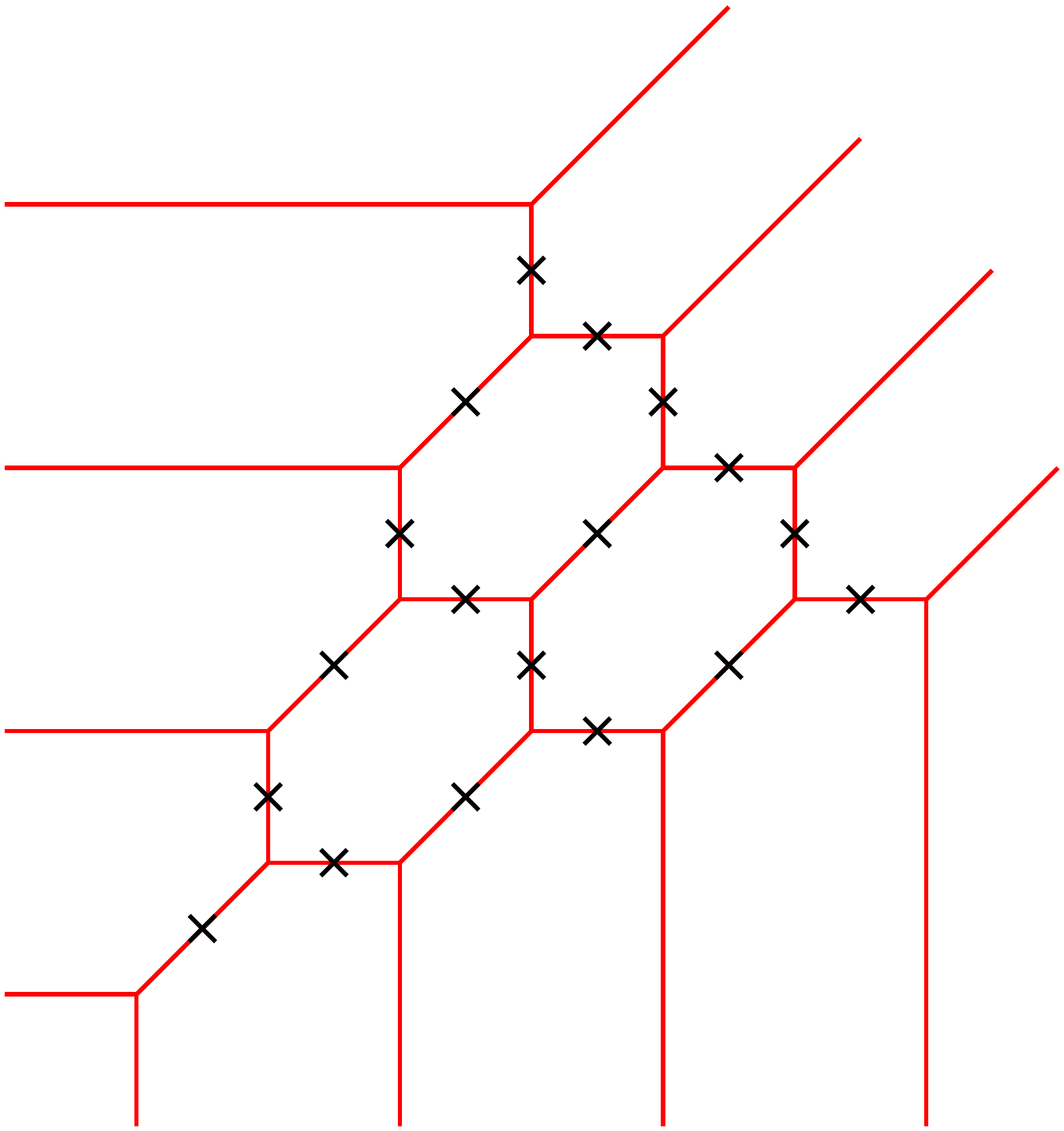}& 
\includegraphics[height=3.5cm, angle=0]{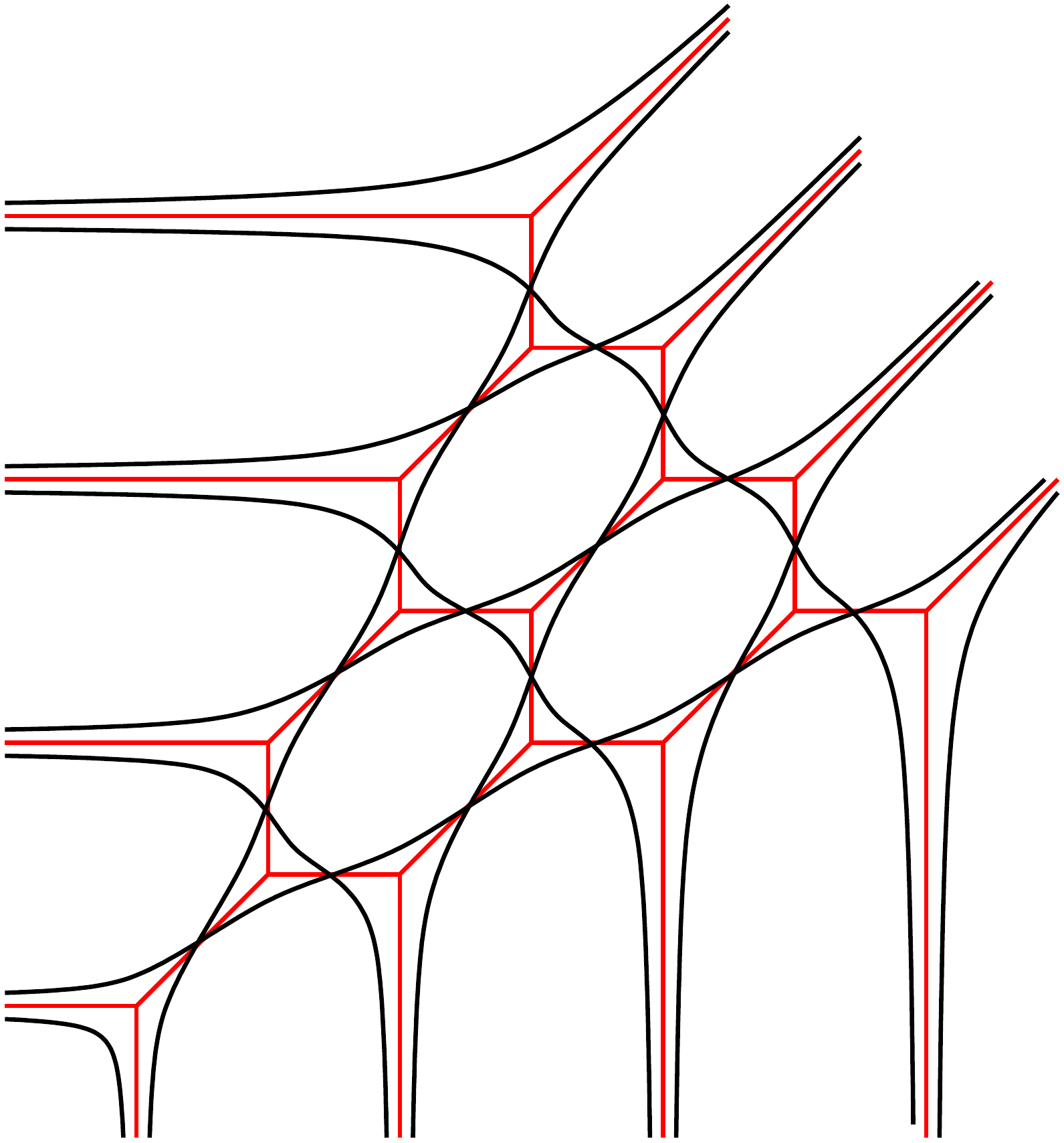}& 
\includegraphics[height=3.5cm, angle=0]{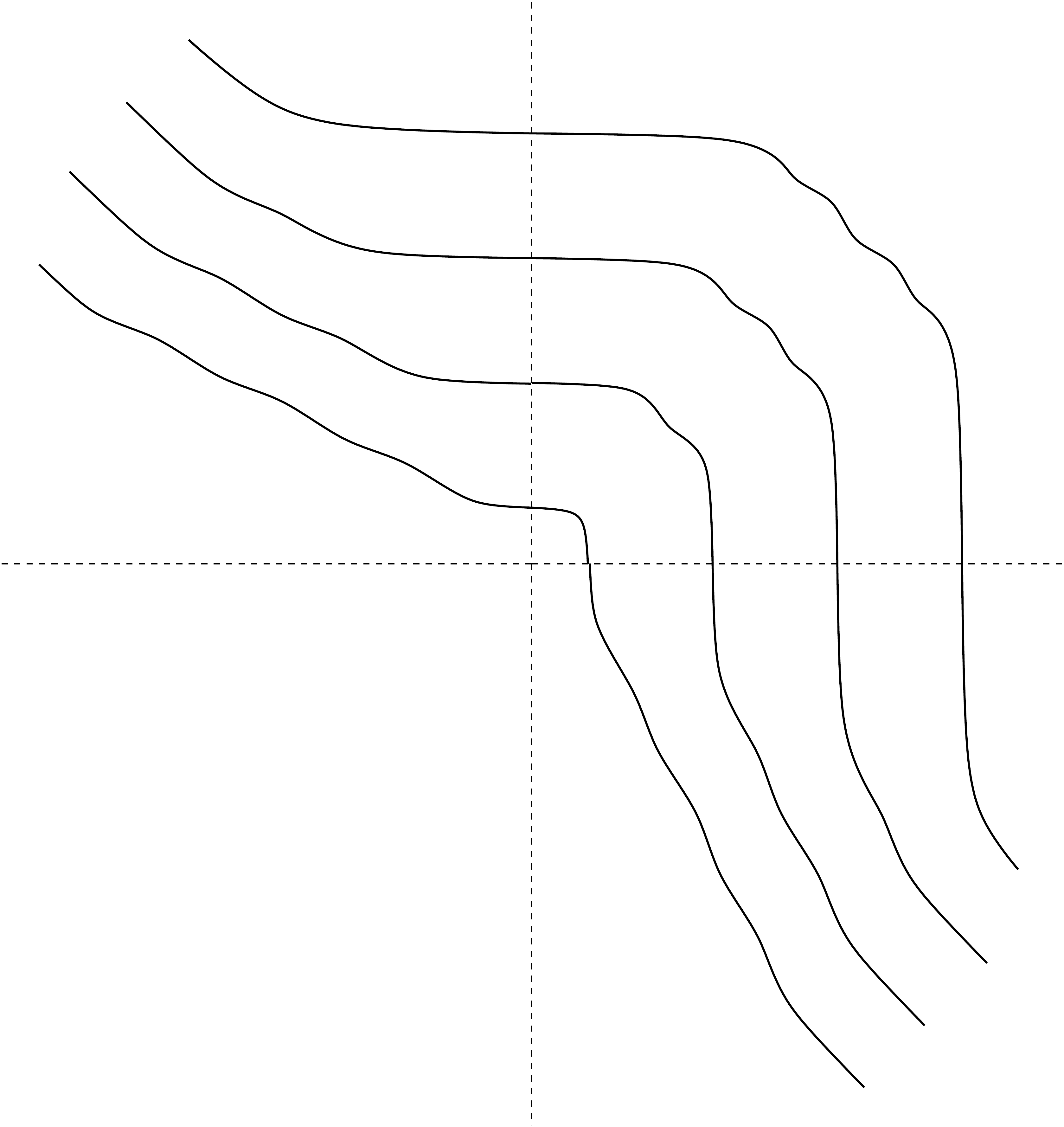}& &
\includegraphics[height=3.5cm, angle=0]{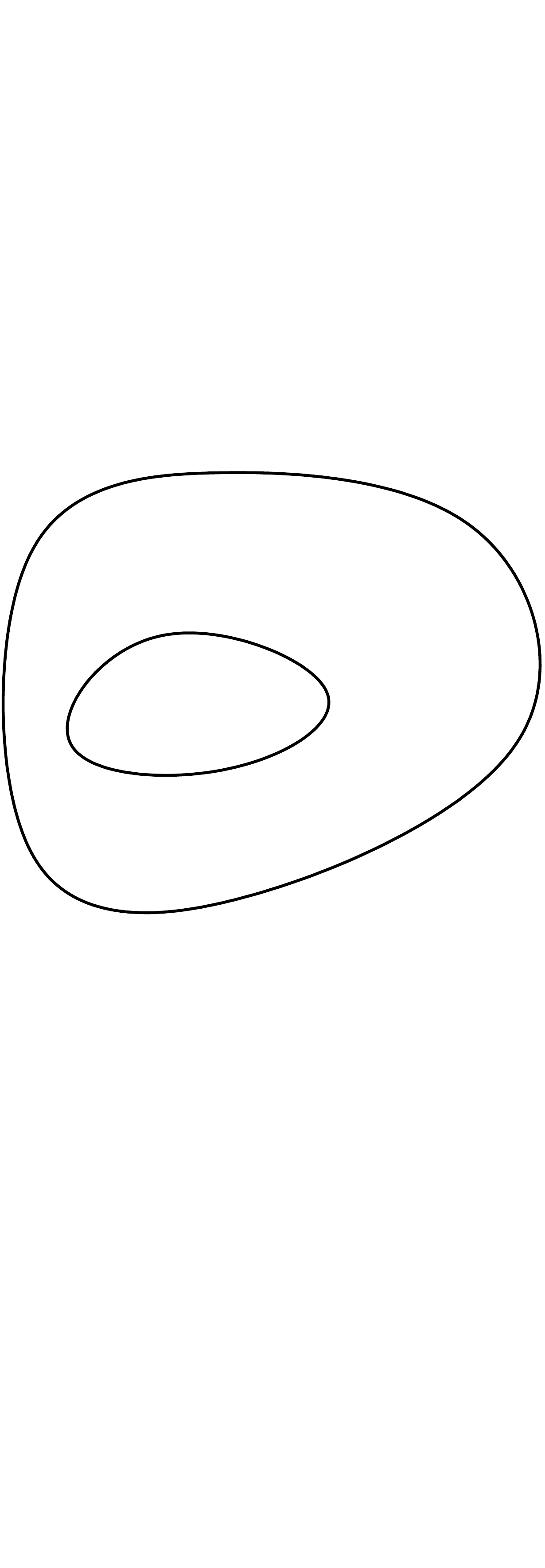}

\\  a)&
b) &
c) & \hspace{3ex} & d)
\end{tabular}
\caption{A hyperbolic quartic}\label{fig:patch hyperquartic}
\end{figure}

\begin{figure}[h]
\centering
\begin{tabular}{ccc}
\includegraphics[height=5cm, angle=0]{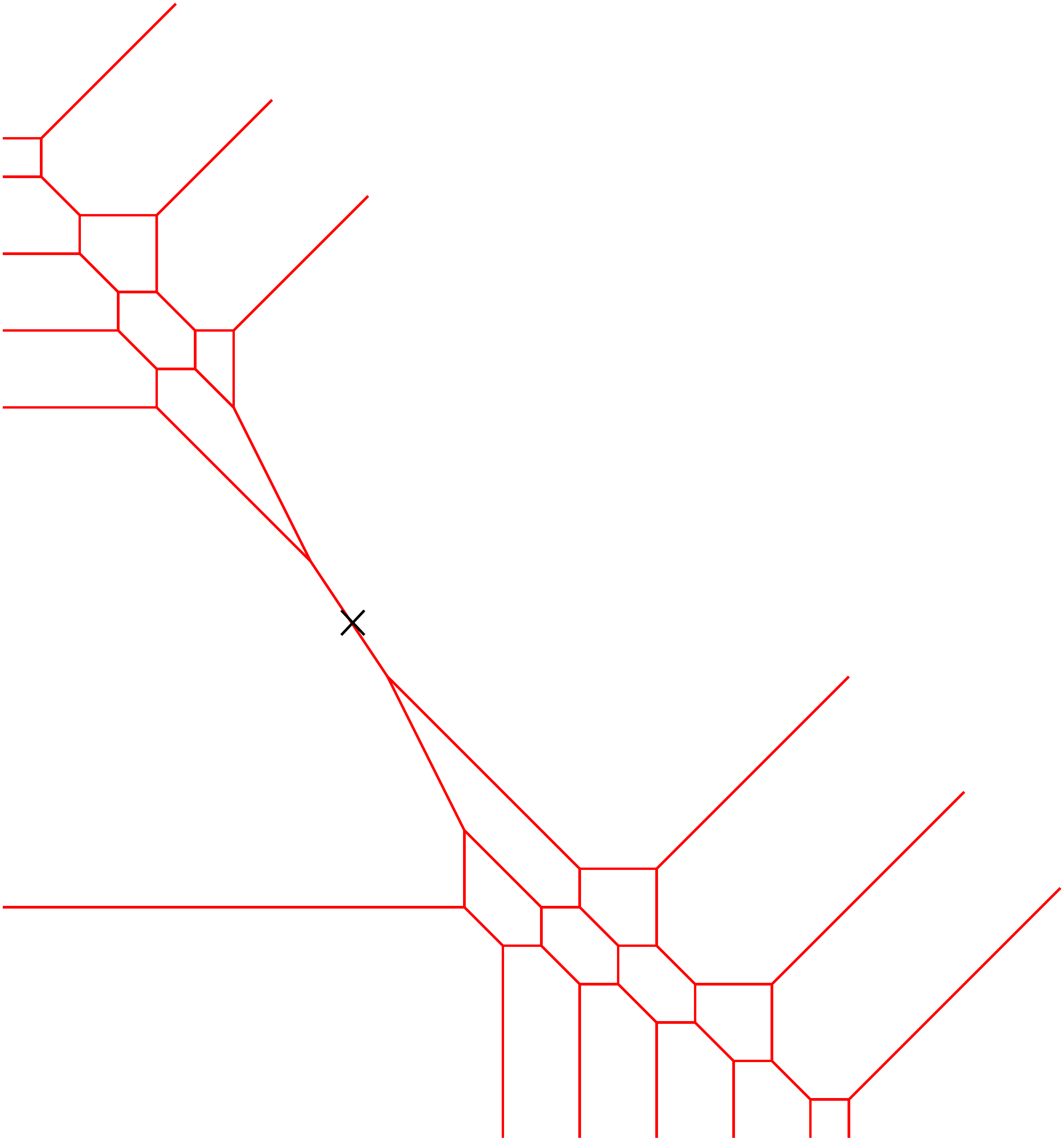}& \hspace{10ex}  &
\includegraphics[height=5cm, angle=0]{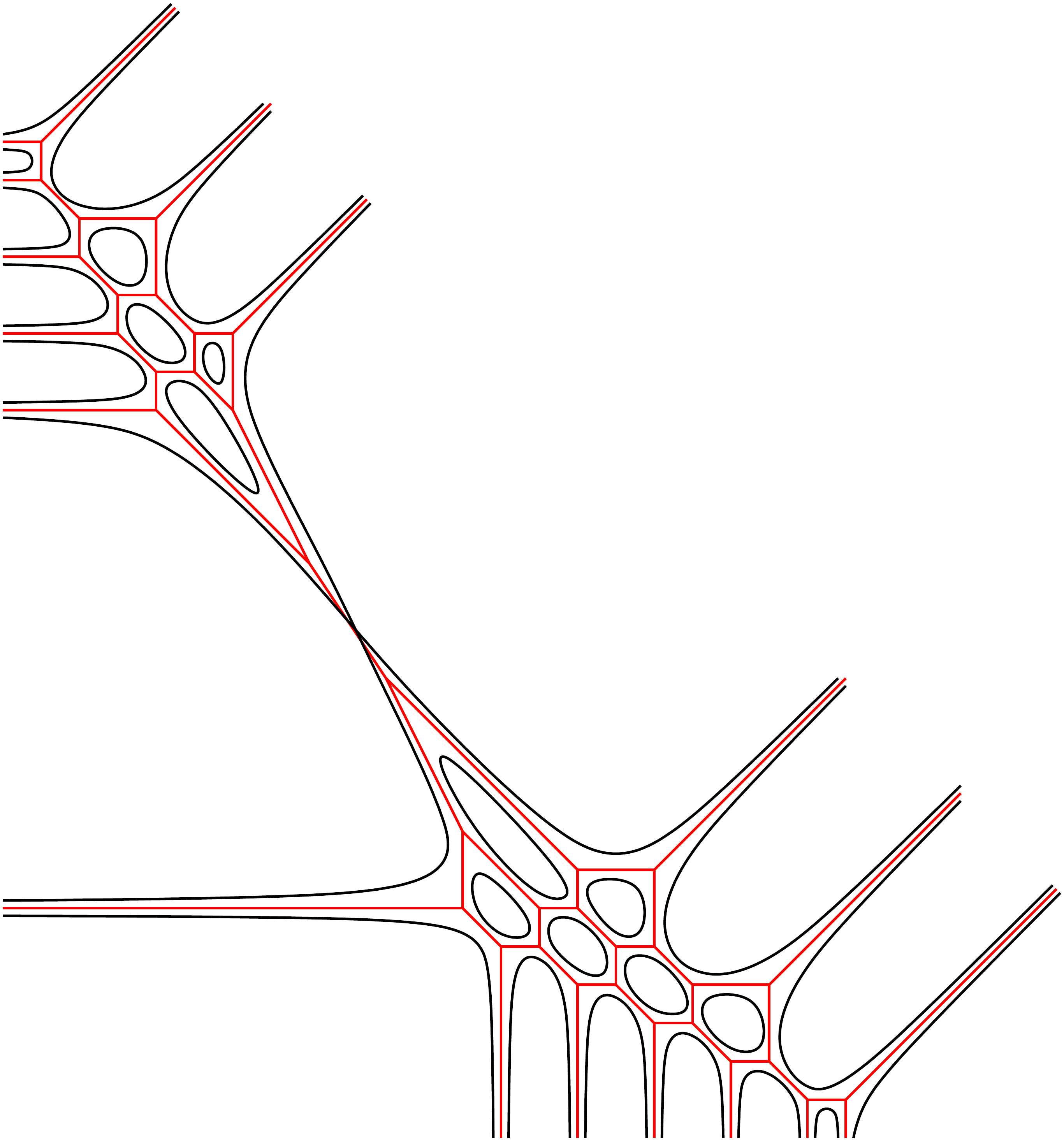}
\\  a)&&
b) 
\\ \\
\includegraphics[height=5cm, angle=0]{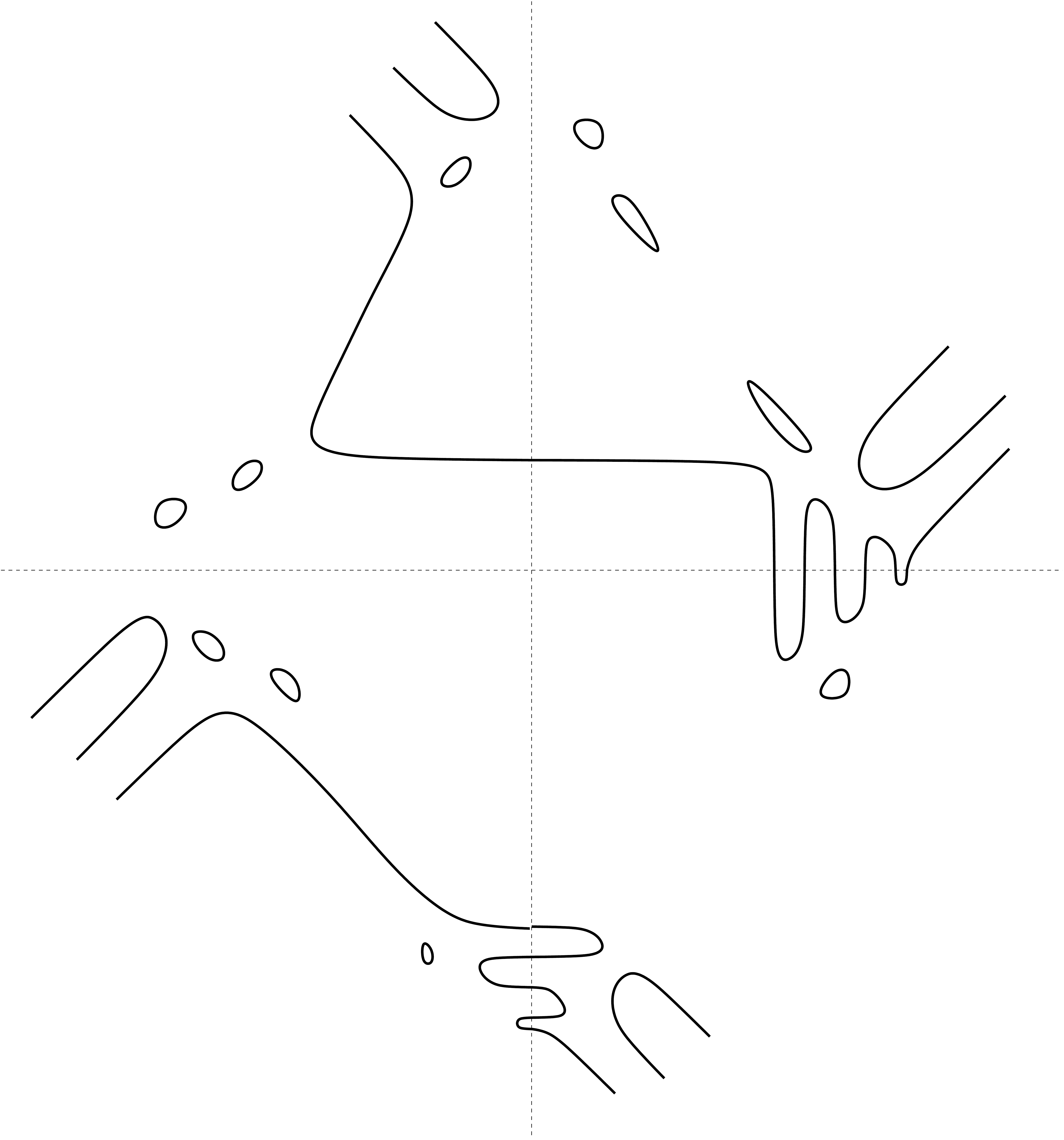}& &
\includegraphics[height=5cm, angle=0]{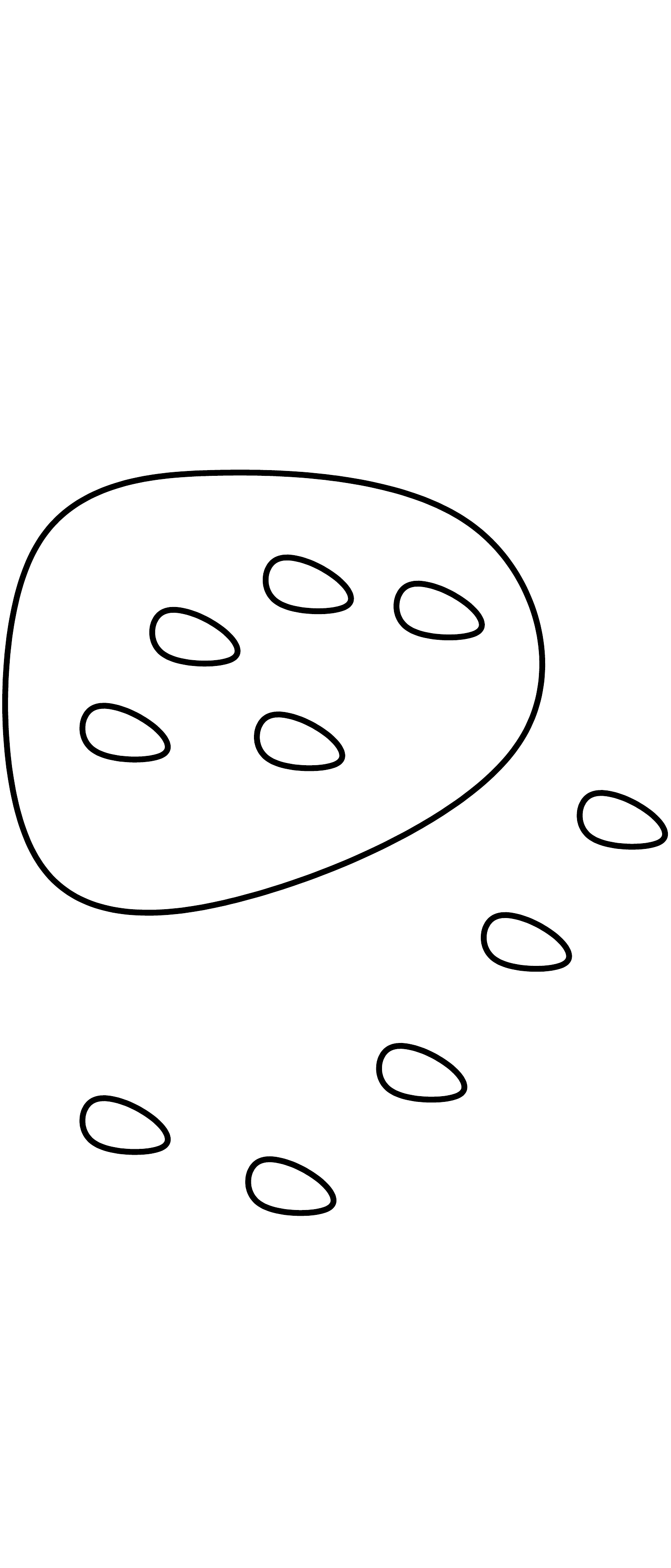}

\\ c) & \hspace{3ex} & d)
\end{tabular}
\caption{Gudkov's sextic}\label{fig:patch gudkov}
\end{figure}

\begin{exa}\label{ex:patch}
Figures
\ref{fig:patch cubic1}, \ref{fig:patch cubic2}, 
\ref{fig:patch hyperquartic}, and \ref{fig:patch gudkov}
depict tropical curves
of degree 3, 4, and 6
 enhanced with twist-admissible collections of edges,
and the isotopy types of the corresponding real algebraic curves 
in both $(\RR^\times)^2$ and $\RR P^2$.

\begin{figure}[h]
\centering
\begin{tabular}{cccccc}
\includegraphics[height=4cm, angle=0]{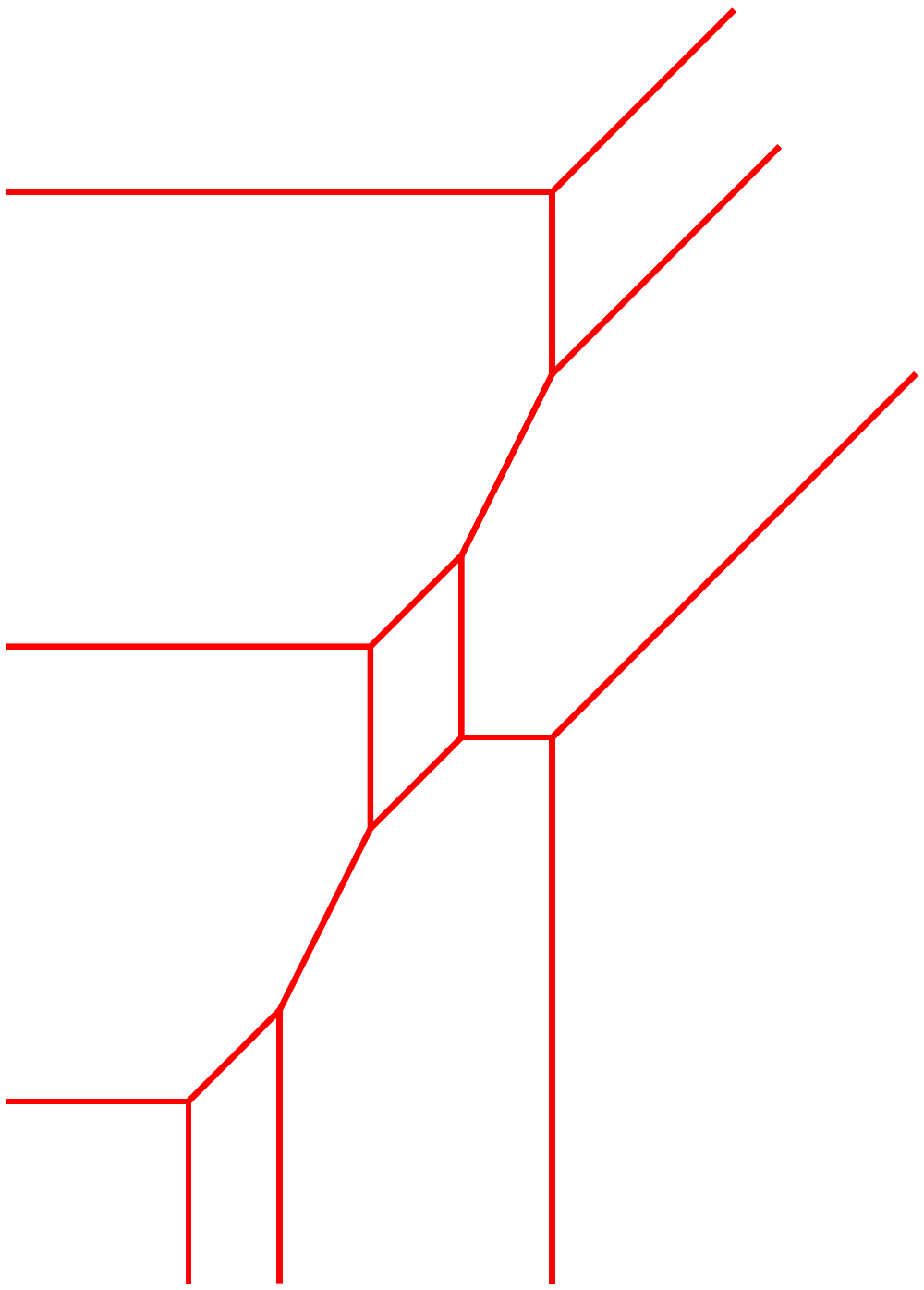}& 
\includegraphics[height=4cm, angle=0]{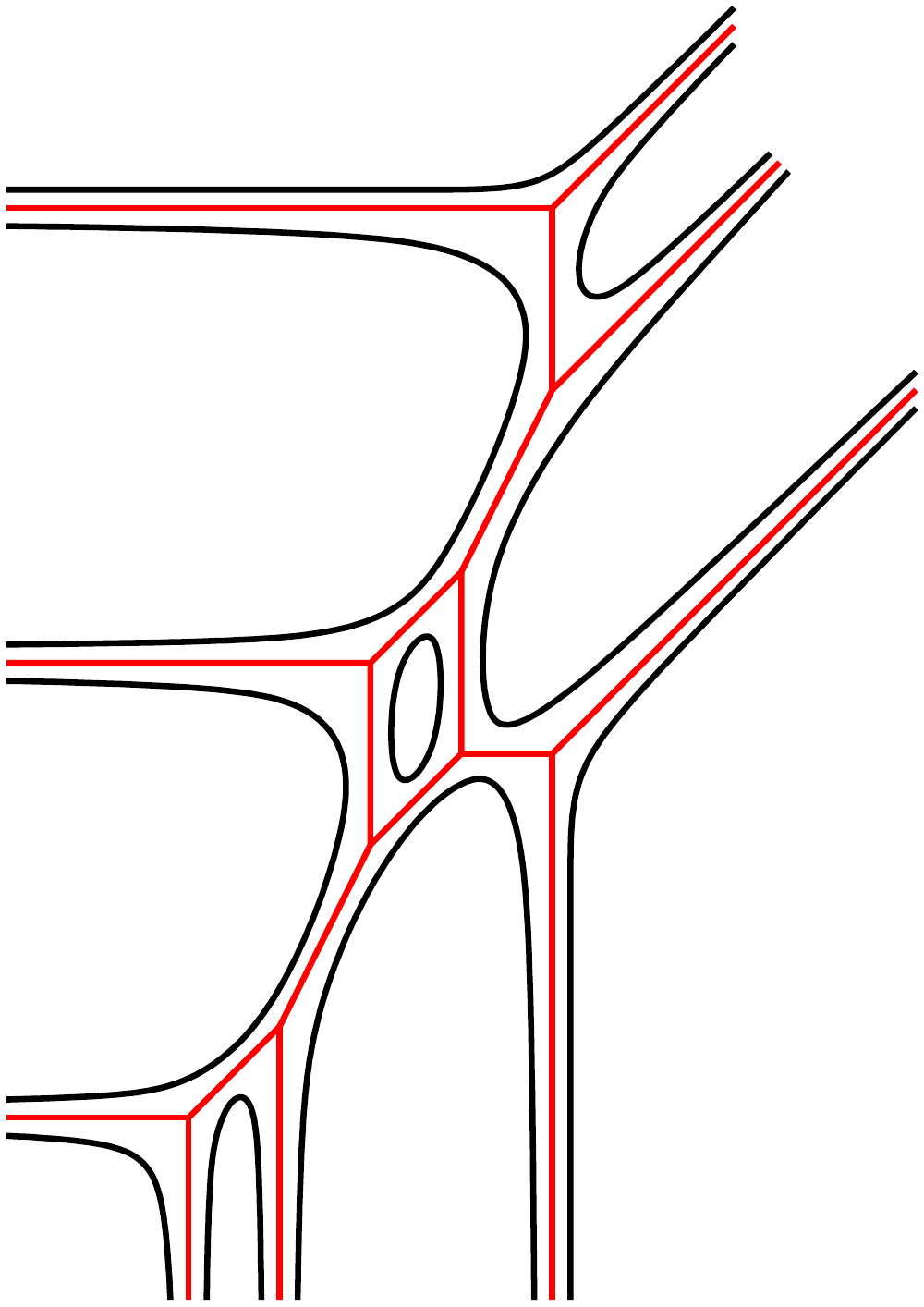}& 
\includegraphics[height=4cm, angle=0]{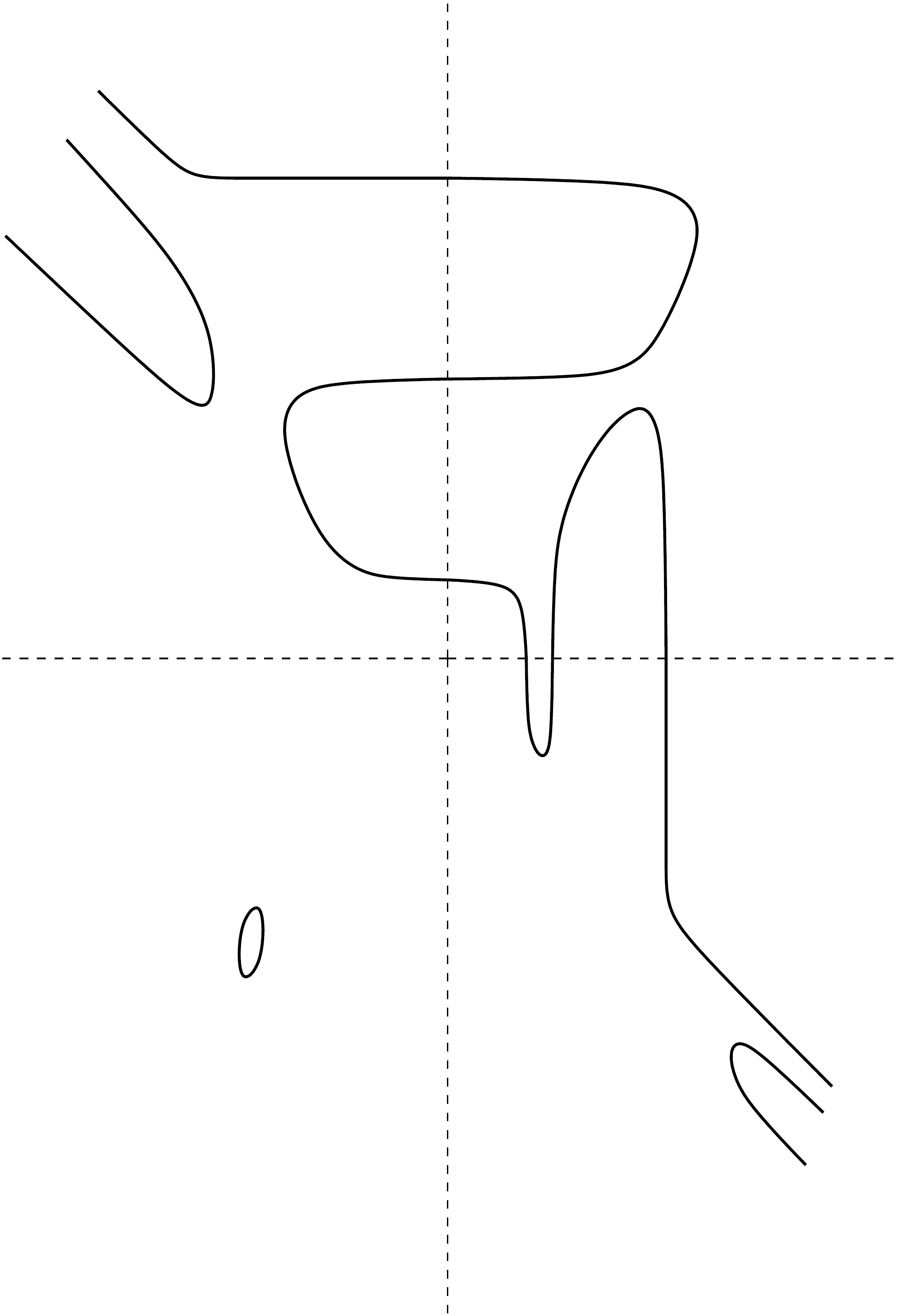}& &
\includegraphics[height=4cm, angle=0]{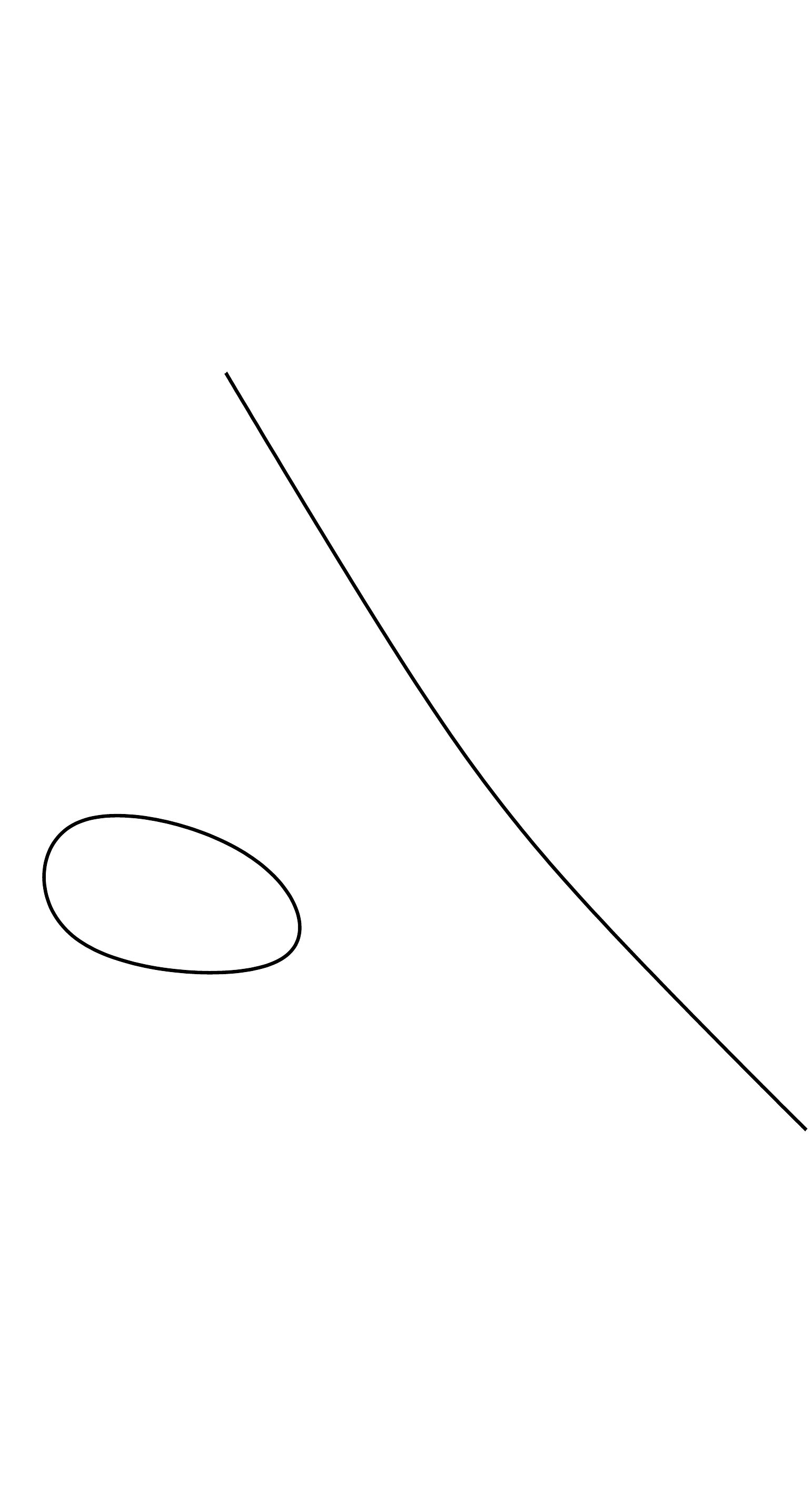}

\\  a)&
b) &
c) & \hspace{3ex} & d)
\end{tabular}
\caption{A Harnack cubic}\label{fig:patch cubic1}
\end{figure}

\begin{figure}[h]
\centering
\begin{tabular}{cccccc}
\includegraphics[height=4cm, angle=0]{Figures/PatchCub1.pdf}& 
\includegraphics[height=4cm, angle=0]{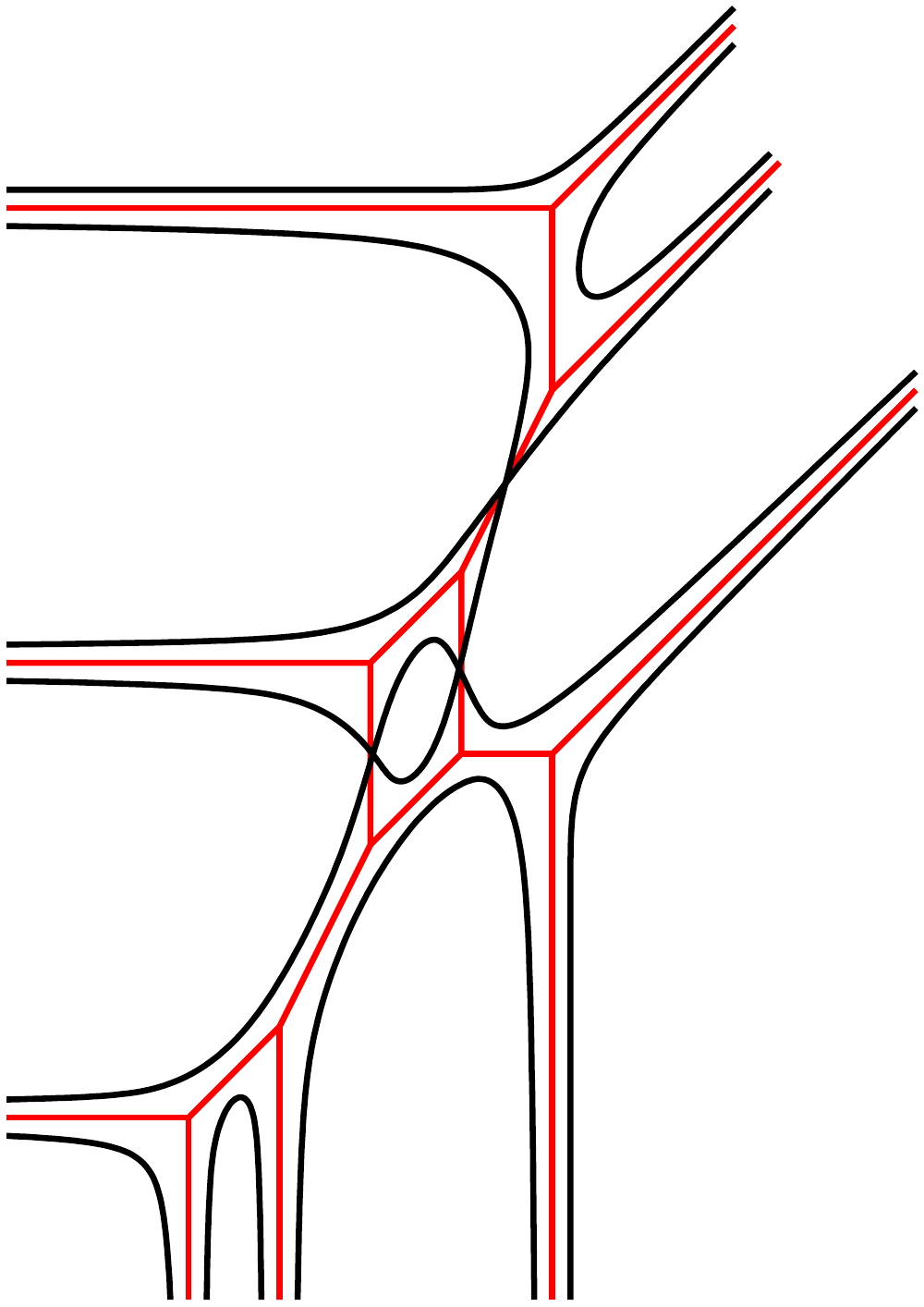}& 
\includegraphics[height=4cm, angle=0]{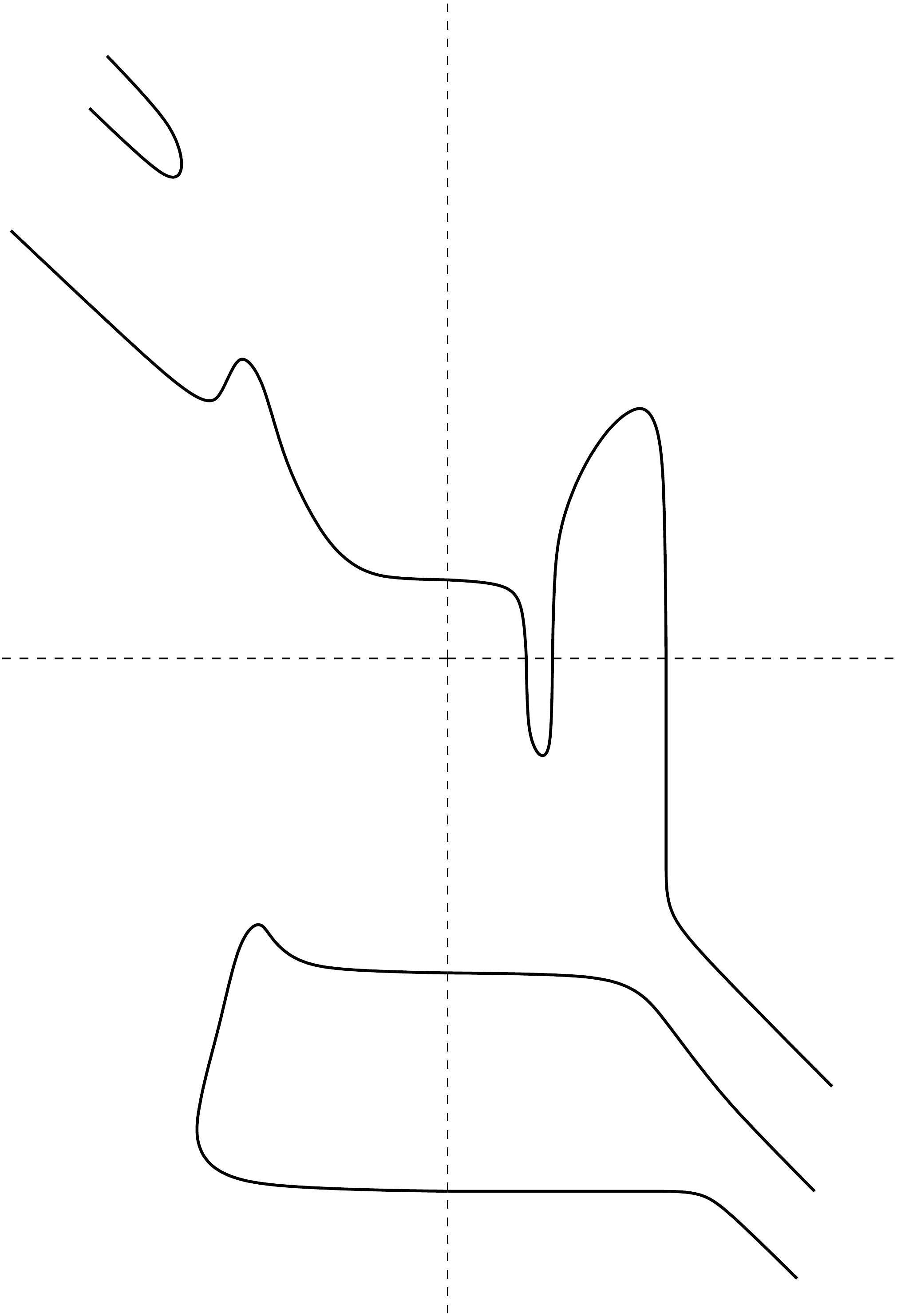}& &
\includegraphics[height=4cm, angle=0]{Figures/PatchCub34.pdf}

\\  a)&
b) &
c) & \hspace{3ex} & d)
\end{tabular}
\caption{Another patchworking of a cubic}\label{fig:patch cubic2}
\end{figure}

A real algebraic curve of degree 6 which realizes the  
isotopy type in $\RR P^2$
depicted in Figure \ref{fig:patch gudkov}d
 was first constructed by D.~A.~Gudkov in 1960's
by a different technique. 
This answered one of the questions posed by D. Hilbert in 1900. 
\end{exa}

\begin{figure}[h]
\centering
\begin{tabular}{ccc}
\includegraphics[scale=0.85]{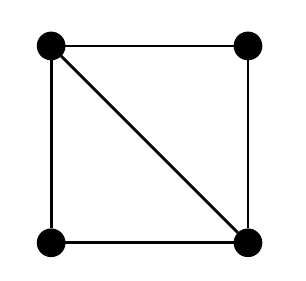}
\put(-10, -5){$p_1^e$}
\put(-75, 65){$p_2^e$}
\put(-75, -5){$ p_3^e$}
\put(-10, 65){$ p_4^e$}
& \hspace{3cm} &
\includegraphics[scale=0.85]{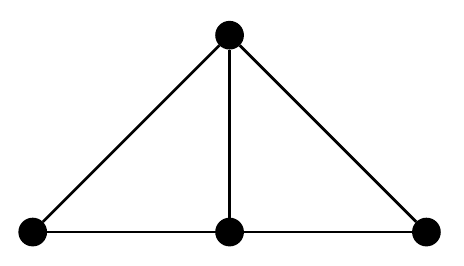}
\put(-7, -5){$p_3^e$}
\put(-63, 67){$p_2^e$}
\put(-63, -5){$ p_1^e$}
\put(-115, -5){$ p_4^e$}
\\ \\ a) && b)
\end{tabular}
\caption{The dual subdivision surrounding $\Delta_e$ from Remark
  \ref{rem:equation}.}
\label{fig:signTwists}
\end{figure}

\begin{rem}\label{rem:equation}
We may find an explicit equation  
for a family $(\C_t)_{t\in\R_{>1}}$ from Theorem \ref{thm:viro} as follows.
Suppose that the
 tropical curve $C$ is given by the tropical polynomial
 $P_{trop}(x,y)=\tg \sum_{i,j} a_{i,j}x^iy^j\td$.
We define the family
$(\C_t)_{t\in\R_{>1}}$ 
by 
a
family of real polynomials 
$P_t(z,w)=\sum_{i,j} \gamma_{i,j}t^{a_{i,j}}z^iw^j$, 
with
$\gamma_{i,j}\in\RR^\times$. 
For any choice of $\gamma_{i,j}$ the resulting family will converge to 
$C$ in the sense of Theorem \ref{approx}.
In turn the set of twisted edges only depends on
the signs of $\gamma_{i,j}$, as described below.

For each bounded edge $e$ of $C$, denote by
 $p^e_1$ and $p^e_2$ the two vertices of the segment $\Delta_e$ dual
 to $e$. The segment $\Delta_e$ is adjacent to exactly two other triangles 
 of the dual subdivision. Let $p^e_3$ and $p^e_4$ denote the vertices of these 
 two triangles not equal to $p^e_1, p^e_2$. 
 Then, the set of twisted edges of the family 
 $(\C_t)_{t\in\R_{>1}}$ is exactly $T$ if and only if for each bounded edge $e$ of the dual subdivision, the
 following holds:
\begin{itemize}
\item if the coordinates modulo 2 of  $p^e_3$ and 
  $p^e_4$ are distinct, (see 
Figure \ref{fig:signTwists}a), then 
$e$ is twisted if and only if
  $\gamma_{p^e_1}\gamma_{p^e_2}\gamma_{p^e_3}\gamma_{p^e_4}>0$;

\item if the coordinates modulo 2 of  $p^e_3$ and
  $p^e_4$ coincide (see 
Figure \ref{fig:signTwists}b), then 
$e$ is twisted if and only if
  $\gamma_{p^e_3}\gamma_{p^e_4}<0$.
\end{itemize}
\end{rem}

\begin{rem}
We may also determine the topological type of the surface $\C_t/conj$ for sufficiently large
$t$, where $conj$ is the restriction on $\C_t$ of the complex
conjugation 
in $(\CC^\times)^2$.
Start with a small tubular neighborhood
$S$ of
$C$ in $\RR^2$. For each twisted edge $e$ of $C$, we cut $S$ along a fiber
at some point inside $e$ and glue it back with a half-twist.
In other words, we start with
the amoeba $\A_t(\C_t)$ (for large $t$)  and add a half-twist 
wherever we see a double point of 
$\A_t( \RR\C_t)$. 

The result is a surface with boundary (and punctures)
diffeomorphic to $\C_t/conj$ for $t$ large enough. For example, the surface 
$\C_t/conj$ corresponding to the patchworking depicted in Figure
\ref{fig:patch ex2}a
is depicted in Figure \ref{fig:patch ex2}c
 (compare with Figure \ref{fig:patch ex2}b).
\begin{figure}[h]
\centering
\begin{tabular}{ccc}
\includegraphics[width=4cm, angle=0]{Figures/PatchCub1.pdf}& 
\includegraphics[width=4cm, angle=0]{Figures/PatchCub12.pdf}
& 
\includegraphics[width=4cm, angle=0]{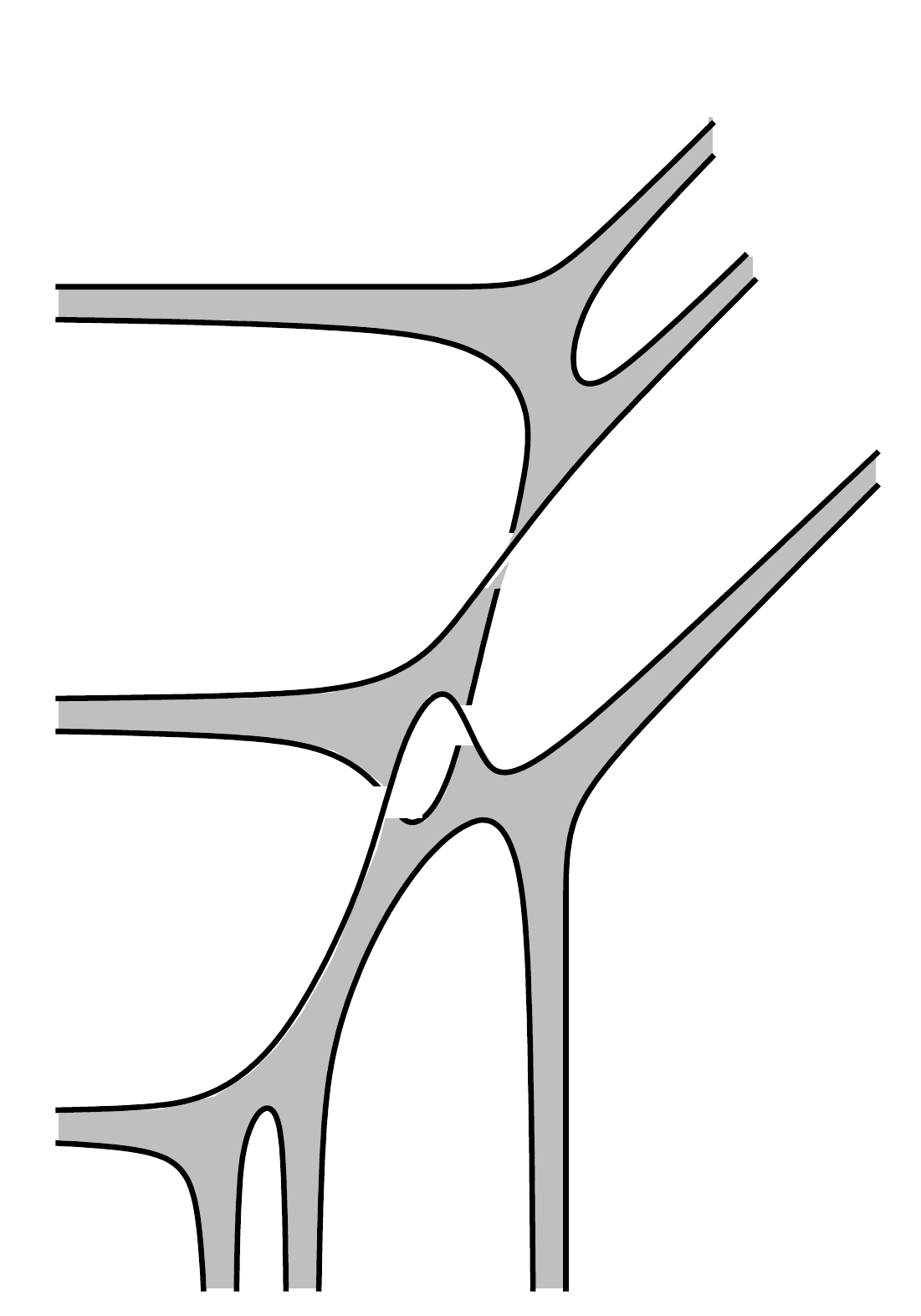}

\\  a) &
b)    & c) The surface  $\C_t/conj$
\end{tabular}
\caption{}\label{fig:patch ex2}
\end{figure}
\end{rem}

Recall that a non-singular real curve $\R\C_t$ is 
said to be 
{\em of type I} 
if $\C_t/conj$ is orientable. 
\begin{prop}[\cite{Haa2}]
Let $C\subset\R^2$ be a smooth tropical curve,
and let $(\C_t)_{t\in\R_{>1}}$ 
be  
a family of non-singular real curves converging to $C$ in the sense of Theorem \ref{approx},
so that the corresponding 
set of twisted edges is $T$. 
A curve $\R\C_t$ is of type I for sufficiently large $t$ if and only if
each cycle in $C$ contains an even number of edges from $T$.  
\end{prop} 
Accordingly, we say that a twist-admissible set $T$ is {\em of type I} 
if each cycle in $C$ contains an even number of edges from $T$.

\subsection{Haas Theorem}

Let $\C$ be a 
non-singular 
real algebraic curve in $(\CC^\times)^2$. 
Denote by $\Delta(\C)$ the Newton polygon of $\C$,
that is, the Newton polygon of a polynomial defining $\C$;
the polygon $\Delta(\C)$ is defined up to translation by a vector with integer coordinates. 
Topologically, the curve $\C$ may be obtained from a closed Riemann
surface $\bar\C$ by removing a finite number of points 
(corresponding to the ends of $\C$).
We denote 
by $\R\bar\C$ the topological closure of $\R \C$ in $\bar\C$.
By Harnack's inequality, the maximal number of connected components 
of $\R\bar\C$ is equal to $g(\bar\C) + 1$, where $g(\bar\C)$
is the genus of $\bar\C$; this genus 
in turn equals to
$Card(\ZZ^2\cap Int(\Delta(\C)))$. 
We say that $\RR C$ is an
\emph{$M$-curve}, or a \emph{maximal curve}, if the number of connected components of $\R\bar\C$
is equal to $g(\bar\C) + 1$. 
The maximal curves constitute extremal objects in real algebraic geometry;
their study goes back to  A.~Harnack and F.~Klein in the XIXth century. We
refer to \cite{VirIntroduction} for an introduction to the subject. 

B. 
Haas in \cite{Haa2} 
found a nice criterion for a curve obtained by
the combinatorial patchworking to be an $M$-curve. 

\begin{defi}
Let $C$  be a 
non-singular 
tropical curve in $\RR^2$, and let $T$
be a twist-admissible set of edges in $C$. 
We say that $T$ is \emph{maximal} if it is of type I, and for any
edge $e\in T$, either $C\setminus e$ is disconnected or there exists
an edge $e'\in T$ such that $C\setminus e$ and $C\setminus e'$ are
connected, but
 $C\setminus (e\cup e')$ is disconnected.
\end{defi} 

For example,  an empty collection $T$ 
is not only twist-admissible, but also maximal. 

\begin{thm}[\cite{Haa2}]
Let $C$  be a 
non-singular 
tropical curve in $\RR^2$, and let $T$ be a
twist-admissible 
set of edges in $C$.
Suppose that $(\C_t)_{t\in\R_{>1}}$ is a family of non-singular real algebraic
curves in $(\CC^\times)^2$ such that
$(\C_t)_{t\in\R_{>1}}$ converges to $C$ in the sense of Theorem \ref{approx} and 
$T$ is the corresponding set of twisted edges. 
Then, 
the real curve $\C_t$ is maximal
for sufficiently large $t$ if and only if
$T$ is maximal.
\end{thm}

\begin{exa} 
If $T$ is a twist-admissible set of edges in a non-singular tropical curve $C$, 
and $T$ contains an edge which is not adjacent
to an unbounded connected component of $\RR^2\setminus C$, 
then $T$ is 
not 
maximal (cf. Figure \ref{fig:patch hyperquartic}).
\end{exa}

\begin{exo}
\

\begin{enumerate}
\item Show that the first Betti number of a non-singular tropical
  curve in $\RR^2$ is equal to the number of integer points contained
  in the interior of the Newton polygon of this curve. 

\item Let $C$ be a non-singular tropical curve, and let $(\C_t)_{t\in\R_{>1}}$ be a
family of non-singular real algebraic 
curves converging to $C$ so that there are no twisted edges. 
Prove that, for sufficiently large values of $t$, the isotopy type of $\C_t$ in $(\RR^\times)^2$
is determined 
{\rm (}up to the action of $(\ZZ/2\ZZ)^2$
by axial symmetries{\rm )}
by the Newton polygon $\Delta(C)$ of $C$ 
and does not depend on the choice of 
a particular tropical curve with given Newton polygon. 

\item Show that any twist-admissible set of edges in the tropical cubic depicted at 
Figure \ref{equil}a 
produces an $M$-cubic. Find
another tropical cubic together with a twist-admissible set for it 
such that they produce a connected real cubic in $\RR P^2$. 

\item 
Following Remark \ref{rem:equation}, write equations
of families of non-singular real algebraic curves $(\C_t)_{t\in\R_{>1}}$ corresponding to
tropical curves and twisted edges in Example \ref{ex:patch}.
\end{enumerate}
\end{exo}

%% file: Ilia.tex
\section{Applications in enumerative geometry}\label{sec:enumerative} 

\subsection{Complex and real enumerative problems}\label{enum_problems} 
Tropical geometry has various applications in enumerative geometry.
We restrict our attention here to one of the classical enumerative problems:
enumeration of curves of given degree $d>0$ and (geometric) genus $g\ge 0$ 
that pass through the appropriate number (equal to $3d-1+g$)
of points in general position in the projective plane. 
This problem can be considered over different fields (and even semi-fields). 
The easiest and the most studied framework is the one of complex geometry.
Consider a collection $\omega$ of $3d-1+g$ points in general position in the complex projective plane $\CC P^2$.
The number $N_{d, g}$ of irreducible 
complex algebraic
curves of degree $d$ and genus $g$ which pass through the points of $\omega$
depends only on $d$ and $g$ and not on the choice of points as long as this choice is generic:
as we will see later, this number can be interpreted as the degree of certain algebraic variety.

A non-singular algebraic curve in $\CC P^2$
of degree $d$ has genus $\frac{(d-1)(d-2)}{2}$, so we obtain
\[ 
N_{d, g}=0\quad\mbox{if}\ g>\frac{(d-1)(d-2)}{2}.
\]

The case $g=\frac{(d-1)(d-2)}{2}$ is also easy. 
If $d = 1$ and $g = 0$, we are counting straight lines which
 pass through two points, 
 so $N_{1, 0}$ is equal to $1$. 
 In the case $d = 2$ and $g = 0$, we are counting conics which pass
 through five points 
 in general position, so $N_{2, 0}$ is also equal to $1$.
More generally, for any choice of
$\frac{d(d+3)}{2}$ points in general position
in $\CC P^2$, 
there exists exactly one curve of degree~$d$ in~$\CC P^2$ 
and genus $g=\frac{(d-1)(d-2)}{2}$
which passes through these points.
Indeed,
the space $\CC C_d$ of all curves of degree $d$ in $\CC P^2$ 
can be identified with a
projective space $\CC P^N$ of dimension
$N=\frac{d(d+3)}{2}$: 
the coefficients of a polynomial 
$\sum \alpha_{i,j}z^iw^ju^{d-i-j}$
defining a given curve
can be taken for homogeneous coordinates 
$[\alpha_{0,0}:\alpha_{0,1}:\alpha_{1,0}:\ldots:\alpha_{d,0} ]$
of the corresponding point
in $\CC C_d$. The condition to pass through a given point 
$[z_0:w_0:u_0]$
in $\CC P^2$ gives rise to 
the
linear condition 
$$\sum  \alpha_{i,j}z_0^iw_0^ju_0^{d-i-j}=0 $$
on the coefficients 
$\alpha_{i,j}$
of a polynomial
defining the curve, and thus
defines a hyperplane in $\CC C_d$. If the 
collection 
of the $\frac{d(d+3)}{2}$
chosen points is sufficiently generic, the corresponding
$\frac{d(d+3)}{2}$ hyperplanes in $\CC C_d$ have exactly one
common point, and this point corresponds
to a non-singular curve. 
Hence, we proved the following statement. 
\begin{prop}
For any positive integer $d$, we have
$$N_{d, \frac{(d - 1)(d - 2)}{2}} = 1.$$
\end{prop}

A more interesting situation arises in the cases $g < \frac{(d - 1)(d - 2)}{2}$, as illustrated by the following example.
\begin{exa}\label{exa:rat cubic}
The number of rational 
cubic curves in $\CC P^2$ 
which pass 
through a collection $\omega$ of
$8$ points in general position  is equal to $12$, {\it i.e.} $N_{3,0}=12$. 
Indeed, the collection $\omega$ determines a straight line ${\mathcal P}$ 
in $\CC C_3$, that is, a pencil ${\mathcal P}$ 
of  cubics. Since this  pencil is generated by any two of its
elements, the intersections of all cubics of ${\mathcal P}$ consists
of $\omega$ together with a ninth point.
Let
$\widetilde{\CC P^2}$ be the  projective plane $\CC P^2$ blown up at
these $9$ points.
The Euler characteristic $\chi(\widetilde{\CC P^2})$
of $\widetilde{\CC P^2}$ is then  equal
to $3 + 9 = 12$. On the other hand, the pencil ${\mathcal P}$ induces
a projection $\widetilde{\CC P^2}\to {\mathcal P}$, which has two
types of fibers: either a smooth elliptic curve, or a nodal
cubic. Since the former has   Euler characteristic
0, and the latter has   Euler characteristic $1$, we obtain that 
 $\chi(\widetilde{\CC P^2}) = N_{3,0}$. 
\end{exa}

Example \ref{exa:rat cubic} generalizes to 
any degree.
\begin{prop}
For any positive integer $d$, we have
$$N_{d, \frac{(d - 1)(d - 2)}{2}-1} =3(d- 1)^2.$$
\end{prop}
We may give another geometric interpretation of the numbers $N_{d, \frac{(d - 1)(d - 2)}{2}-1}$.
Consider the hypersurface $D \subset \CC C_d$
formed by the points corresponding to singular curves of degree $d$ in $\CC P^2$.
This hypersurface is called the {\it discriminant}
of $\CC C_d$. 
The smooth part of~$D$ is formed by the points
corresponding to curves
whose only singular point is a non-degenerate double point.
A generic collection of $\frac{d(d+3)}{2}-1$ points determines 
a straight line in $\CC C_d$, 
and moreover, this line intersects
the discriminant only in its smooth part and transversally. 
Thus, the number $N_{d, \frac{(d - 1)(d - 2)}{2} - 1}$ coincides with the degree of~$D$. 

To reformulate in a similar way the general problem formulated in the beginning of the section,
choose a collection $\omega$ of $\frac{d(d + 3)}{2} - \delta = 3d - 1 + g$ points in
general position in $\CC P^2$,
where  $\delta = \frac{(d - 1)(d - 2)}{2} - g$.
The expression ``{\it in a general position}'' can be made precise
in the following way.
Denote by $S_d(\delta)$ the subset of $\CC C_d$
formed by the points
corresponding to irreducible curves of degree~$d$
satisfying the following property:
each of these curves has $\delta$ non-degenerate double points 
and no other singularities
(such curves are called {\it nodal}). 
The {\it Severi variety} $\overline{S}_d(\delta)$
is the closure of $S_d(\delta)$ in $\CC C_d$.
It is an algebraic variety of codimension~$\delta$
in $\CC C_d$. 
The smooth part of $\overline{S}_d(\delta)$ 
contains $S_d(\delta)$ (see, for example, \cite{Zar}).
We say that the points of $\omega$ are {\it in a general position}
(or that $\omega$ is {\it generic}) if the following hold, 
\begin{itemize}
\item the dimension of the projective subspace $\Pi(\omega) \subset \CC C_d$
defined by the points of~$\omega$
is equal to~$\delta$;
\item the intersection $\Pi(\omega) \cap \overline{S}_d(\delta)$
is contained in
$S_d(\delta)$;
\item and the above intersection
is transverse. 
\end{itemize}
It can be 
proved (see, for example, \cite{Kleiman-Piene}) that 
the generic collections form an open dense subset 
in the space of all collections of $\frac{d(d+3)}{2} - \delta$
points in $\CC P^2$. 
If~$\omega$ is generic, 
any irreducible curve of degree~$d$ and genus $g$ in $\CC P^2$
which passes through the points of~$\omega$
corresponds to a point of $\Pi(\omega) \cap \overline{S}_d(\delta)$.
Thus, the number $N_{d, g}$ coincides with the degree of
the Severi variety $\overline{S}_d(\delta)$.

\begin{rem}\label{torus}
Since we consider generic collections $\omega \subset \CC P^2$,
we can restrict ourselves to the situation where all points of $\omega$
are contained in the complex torus $(\CC^\times)^2 \subset \CC P^2$;
then, $N_{d, g}$ becomes the number of irreducible nodal curves in $(\CC^\times)^2$ 
which pass through the points of $\omega$, are defined by polynomials of degree $d$ in two variables
and have $\delta = \frac{(d - 1)(d - 2)}{2} - g$ double points. 
\end{rem} 

The numbers $N_{d, g}$
are {\it Gromov-Witten invariants} of $\CC P^2$.
The number $N_d = N_{d, 0}$
is the number of rational curves of degree~$d$
which pass through a generic collection
of $3d - 1$
points in $\CC P^2$.
A recursive formula for the numbers
$N_d$, 
was found
by M.~Kontsevich (see~\cite{KonMan1}).
A recursive formula that allows one to calculate all numbers
$N_{d, g}$ was obtained by
L.~Caporaso and J.~Harris~\cite{CapHar1}. 

\vskip10pt

The enumerative problem discussed above can as well be considered over the real numbers.
Consider a collection $\omega$ of $3d-1+g$ points in general position in the real projective plane $\RR P^2$.
This time, the number $R_{d, g}(\omega)$ of irreducible {\bf real} curves of degree $d$ and genus $g$ 
which pass through the points of $\omega$,
in general, depends on $\omega$. For example, for $d = 3$ and $g = 0$,
the number $R_{d, g}(\omega)$ can take values $8$, $10$, and $12$ (see
\cite{DK},
or Example \ref{exa:rat real cubic}).

\begin{figure}[h]
\begin{center}
\begin{tabular}{ccc}
\includegraphics[height=2cm, angle=0]{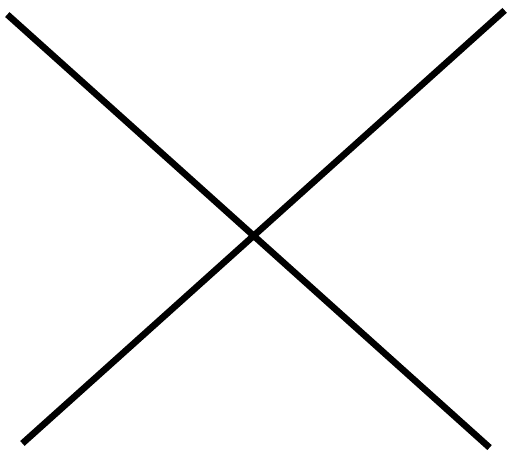}&
\hspace{15ex}  &
\includegraphics[height=2cm, angle=0]{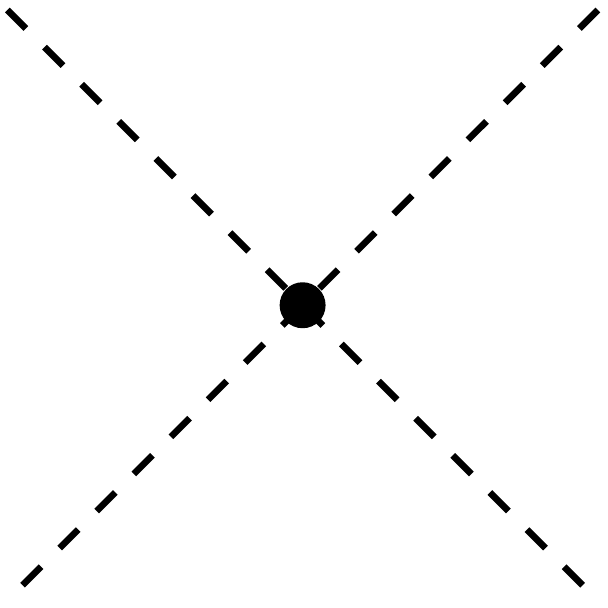}

\\ a)  hyperbolic && b) elliptic
\end{tabular}
\end{center}
\caption{Two types of real nodes}
\label{fig:real nodes}
\end{figure}

J.-Y.~Welschinger 
suggested to treat real curves differently for enumeration, 
so that some real curves are counted with multiplicity $+1$ and some with multiplicity $-1$. 
He proved that the result is invariant on the choice of points in general position in the case $g = 0$
(see \cite{Wel1}).  
To define the Welschinger signs, recall that all the curves under enumeration 
are nodal. 
A real non-degenerate double point of a nodal real curve $\C$ can be 
\begin{itemize}
\item
{\it hyperbolic} ({\it i.e.}, intersection of two real branches
of the curve, see Figure \ref{fig:real nodes}a),
\item
or {\it elliptic} ({\it i.e.}, intersection of two imaginary conjugated branches, see Figure \ref{fig:real nodes}b).
\end{itemize} 
Denote by $s(\C)$ the number of elliptic double points of $\C$. 

\begin{theorem}[Welschinger, \cite{Wel1}]\label{real_enumeration} 
Let $\omega$ be a collection of $3d - 1$ points in general position in $\RR P^2$. 
The number 
$$
W_{d}(\omega) = \sum_\C (-1)^{s(\C)},
$$
where the sum runs over all real 
rational curves of degree $d$ in $\RR P^2$ 
which pass through the points of $\omega$, 
does not depend on the choice of {\rm (}generic{\rm )\/} collection $\omega$. 
\end{theorem} 

The number $W_d(\omega)$ is denoted by $W_d$, and is called {\it Welschinger invariant}.
Notice that the absolute value of $W_d$ provides a lower bound for the
numbers $R_{d, 0}(\omega)$.

\begin{rem}
More generally we
could have chosen a real 
collection 
of points in $\CC P^2$, {\it i.e.} a 
collection 
made of real points and pairs of complex
conjugated points. In this more general situation, Welschinger proved
in \cite{Wel1} that the signed enumeration of real rational curves given
in Theorem \ref{real_enumeration} only depends on $d$ and on the
number of real points in $\omega$.
For the sake of simplicity, we 
consider here only the case 
$\omega\subset \RR P^2$. 
\end{rem}

\begin{exa}
Since rational curves of degree 1 and 2 in $\CC P^2$ are non-singular, we 
 have $W_1 = W_2 = 1$.
\end{exa}
\begin{exa}\label{exa:rat real cubic}
Adapting to the real setting the 
calculation 
performed in Example
\ref{exa:rat cubic}, we can compute 
$W_3$. 
A collection~$\omega$ of $8$ points in general position in $\RR P^2$ 
defines a pencil $\RR {\mathcal P}$ 
of real cubics which pass through all the points of~$\omega$.
In particular, the $9$ intersection points of all cubics in ${\mathcal P}$ are real. 
Let
$\widetilde{\RR P^2}$ be the real projective plane $\RR P^2$ blown up
at these~$9$ points.
The Euler characteristic $\chi(\widetilde{\RR P^2})$
of $\widetilde{\RR P^2}$ is equal
to $1 - 9 = -8$. On the other hand, the calculation of the Euler
characteristic of $\widetilde{\RR P^2}$ via the pencil $\RR {\mathcal P}$ 
gives $\chi(\widetilde{\RR P^2}) = -W_3$. 
Thus $W_3 = 8$. 

The lower bound $8$ for the number of real rational cubics
passing through $8$ points in general position in $\RR P^2$
is sharp and
was proved by V.~Kharlamov before
the discovery of the Welschinger invariants
(see, for example, \cite{DK}).
Using floor diagrams (see Section \ref{sec:floor_diagrams}), E. Rey showed that
the lower 
bound 
provided by the Welschinger invariant $W_4$ is also sharp.
It is not known whether the lower bounds
provided by $W_d$, $d \geq 5$ are sharp.
\end{exa} 

The  values of $W_d$ with $d\ge 4$ are more difficult to calculate; the first calculation
of the numbers $W_d$ for $d \geq 4$ was obtained 
via tropical geometry.  

\medskip
Note that
the enumeration of real curves with Welschinger signs does not give
rise to an invariant count 
if $g > 0$ (see \cite{Wel1,IKS1}).

\subsection{Enumeration of tropical curves}\label{tropical_enumeration} 
The above enumerative problem can be also considered over the tropical semi-field. 
Since we did not yet introduce the notion of {\it tropical projective plane}, 
as well as the notion of {\it genus} for tropical curves
(these notions are 
respectively 
introduced 
 in  Section \ref{proj spaces}
and 
in Definition \ref{def:genus}),
we reformulate our problem similarly to Remark \ref{torus}
and consider tropical curves in the tropical torus $\RR^2 = (\TT^\times)^2$. 

A tropical curve $C$ in $\R^2$ is said to be \emph{irreducible} if it is not a
union of two tropical curves, both different from $C$.
A tropical curve $C$ 
in $\RR^2$ is called {\it nodal} if 
each unbounded edge of $C$ is of weight $1$ and 
each polygon of the subdivision dual to $C$
is either a triangle or a parallelogram. 
We define the {\it number of double points} of a nodal tropical curve $C$ of degree $d$ in $\RR^2$ 
to be equal to $I(C) + P(C)$, where $I(C)$ is the number of integer points
of the triangle 
$Conv\{(0, 0), (d, 0), (0, d)\} \subset \RR^2$
which are not vertices of the subdivision dual to $C$, and $P(C)$  
is the number of parallelograms
in the dual subdivision.
For any irreducible nodal tropical curve $C$ having exactly $\delta$ double points in $\RR^2$, 
we can define the genus $g$ of $C$ putting 
$g = \frac{(d - 1)(d - 2)}{2} - \delta$ (cf. the definition of genus
of a tropical curve 
at the end of Section \ref{sec:manifold}).
In the case of genus $0$, we speak about {\it rational} tropical curves. 

\begin{exa}
A non-singular tropical curve 
whose Newton polygon is the triangle with vertices 
$(0, 0)$, $(d, 0)$, and $(0, d)$ has genus 
$\frac{(d - 1)(d - 2)}{2}$ (compare with Exercise 3(1)).
\end{exa}

\begin{exa}
The tropical curves depicted in Figures \ref{droite}a, b are irreducible, nodal,
and rational.
This is also the case for the tropical curves in Figures \ref{equil}b,c.
The tropical curve depicted in Figure \ref{equil}a is irreducible, nodal,
and has genus $1$. 
\end{exa} 

For different choices of $3d - 1 + g = \frac{d(d + 3)}{2} - \delta$ points in general position
in $\RR^2$, 
the numbers of irreducible nodal tropical curves of degree $d$ 
which pass through these points and have $\delta$ double points can be different. 
Nevertheless, these tropical curves may 
be prescribed multiplicities in such a way that the resulting numbers are invariant.
Consider a collection $\omega$ of $3d -1 + g = \frac{d(d + 3)}{2} - \delta$ points in general position in $\RR^2$.
Let $C$ be an irreducible nodal tropical curve of degree $d$ such that $C$ passes through the points of $\omega$ 
and has $\delta$ double points ({\it i.e.}, has genus $g$).
Each vertex of $C$ is either trivalent (dual to a triangle)
or four-valent (dual to a parallelogram).  
To each trivalent vertex $v$ of such a tropical curve $C$ we associate
two numbers:
\begin{itemize}
\item $m_{\CC}(v)$ equal to twice the Euclidean area of the dual
triangle; 
\item $m_{\RR}(v)$ equal to $0$ if $m_{\CC}(v)$ is even and $(-1)^{i(v)}$ if $m_{\CC}(v)$ is odd,
where $i(v)$ is the number of integer points in the interior of the triangle dual to $v$. 
\end{itemize}
Put 
$$
\displaylines{
m_{\CC}(C) = \prod_v m_{\CC}(v), \cr
m_{\RR}(C) = \prod_v m_{\RR}(v),
}
$$
where the products are taken over all trivalent vertices of $C$.

\begin{thm}[Mikhalkin's correspondence theorem, \cite{Mik1}]\label{correspondence} 
Let $\omega$ be a collection of $3d - 1 + g$ points in general
position in $\RR^2$. Then, 
\begin{enumerate}
\item
the number
of 
irreducible
nodal tropical curves $C$ of degree $d$ and genus $g$ in $\RR^2$,
counted with multiplicities $m_{\CC}(C)$, which pass through the points of $\omega$
is equal to $N_{d, g}$; 
\item
if $g=0$, 
then the number
of irreducible nodal rational tropical curves $C$ of degree $d$ in $\RR^2$,
counted with multiplicities $m_{\RR}(C)$, which pass through the points of $\omega$ is equal to $W_d$. 
\end{enumerate} 
\end{thm} 
\begin{exa}
For each integer $1 \leq d \leq 3$, 
we depicted in Figure \ref{fig:enum smooth} 
a generic collection of $\frac{d(d + 3)}{2}$ points in $\RR^2$
and the 
unique nodal tropical curve $C$ of degree
$d$ and genus $\frac{(d-2)(d-1)}{2}$ 
which passes through the points of the chosen 
collection. 
In each case 
we have $m_{\CC}(C) =m_{\RR}(C) = 1$, and Theorem \ref{correspondence}
gives $N_{1,0}=N_{2,0}=N_{3,1}=W_1=W_2=1$.
\begin{figure}[h]
\begin{center}
\begin{tabular}{ccccc}
\includegraphics[height=2cm, angle=0]{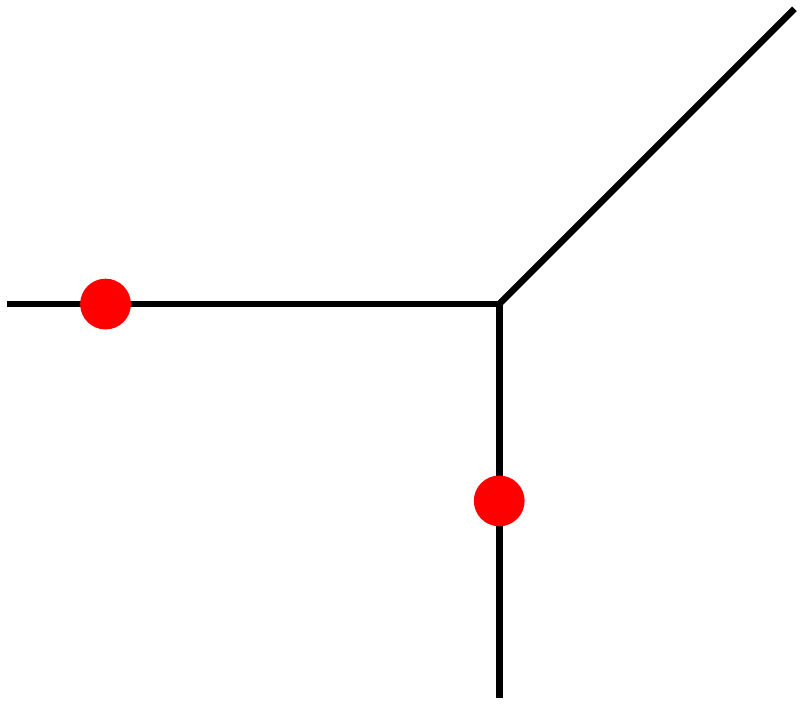}&
\hspace{5ex}  &
\includegraphics[height=3cm, angle=0]{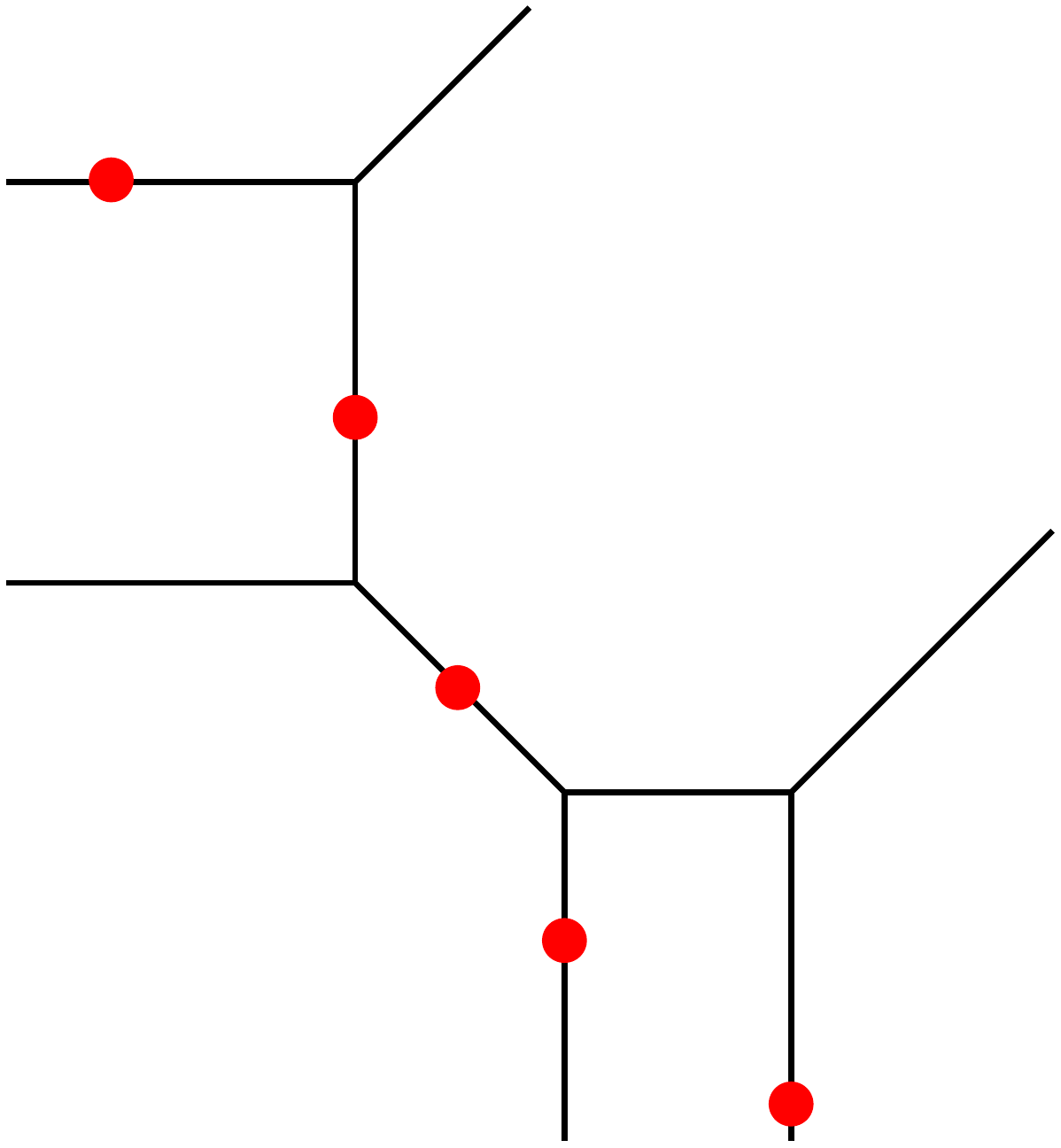}&
\hspace{5ex}  &
\includegraphics[height=4cm, angle=0]{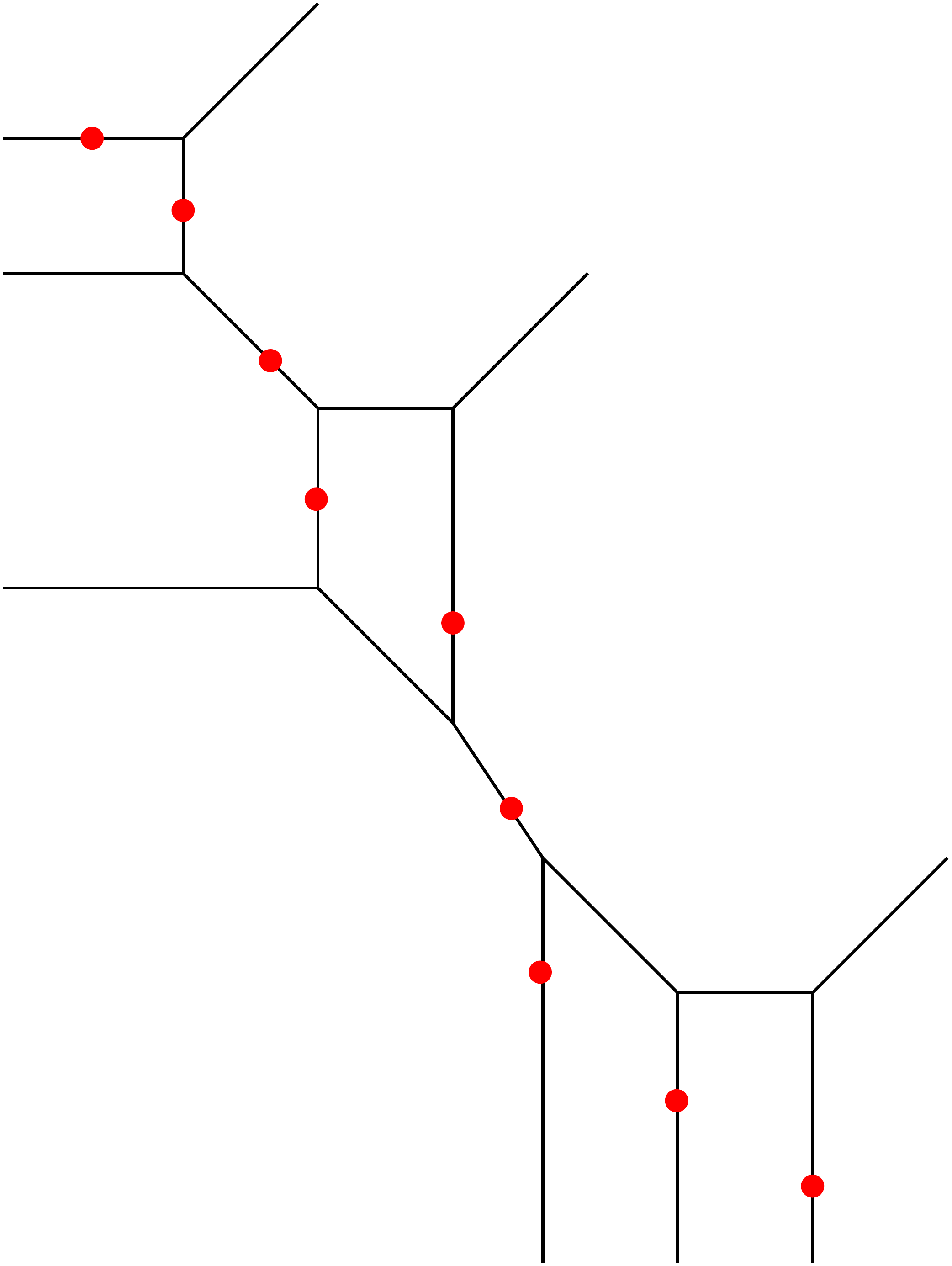}

\\  $d=1$ &&  $d=2$ &&  $d=3$
\end{tabular}
\end{center}
\caption{}
\label{fig:enum smooth}
\end{figure}
\end{exa}
\begin{exa}
We depicted in Figure \ref{fig:enum 30} the rational nodal tropical cubics 
that pass through the points of 
a given generic 
collection 
of $8$ points in $\RR^2$. For each curve, we precise
its real and complex multiplicities (as well as the multiplicity $G$  
which will be
defined in Section \ref{quantum}). 
In particular, using Theorem \ref{correspondence}, 
we 
obtain 
again
$N_{3,0}=12$ and $W_3=8$. 
\begin{figure}[h]
\centering
\begin{tabular}{ccccc}
\includegraphics[width=2cm, angle=0]{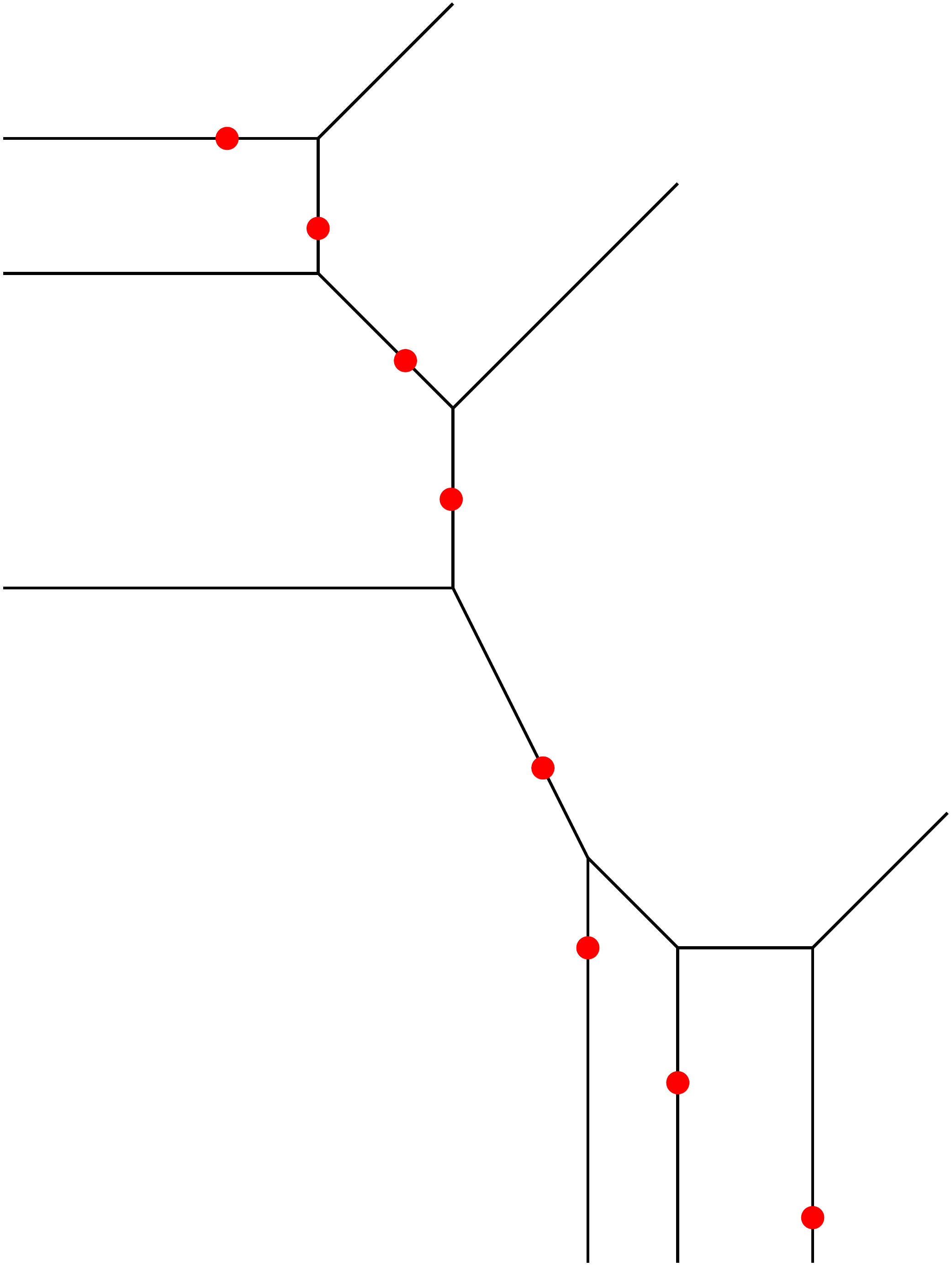}
\put(-28,43){\small{$2$}}& \hspace{4ex} &
\includegraphics[width=2cm, angle=0]{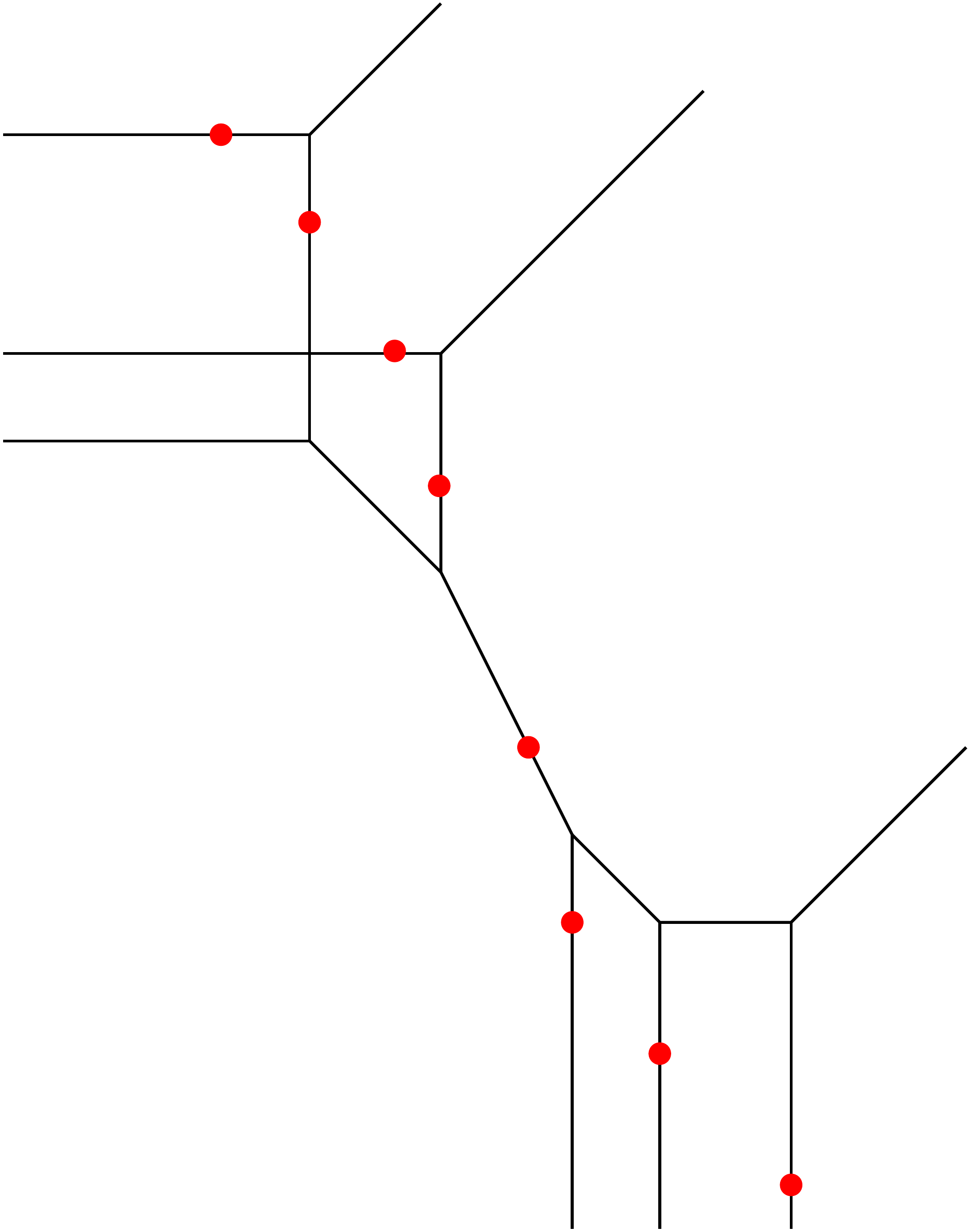}&  \hspace{4ex} &
\includegraphics[width=2cm, angle=0]{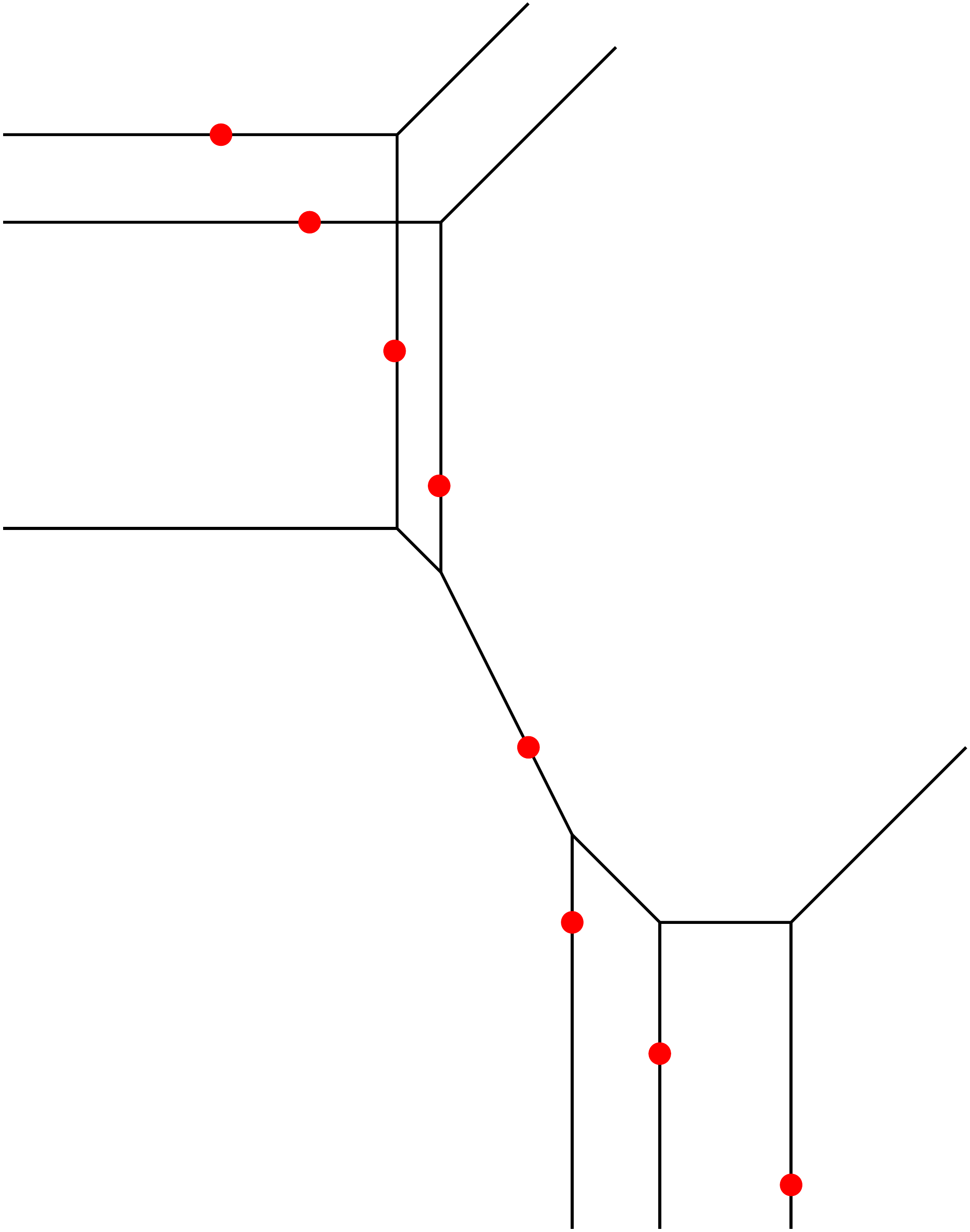}
\\\begin{tabular}{l}
$m_{\CC} = 4$
\\ $m_{\RR}=0$
\\ $G= q^{-1}+2 +q$
\end{tabular}&&
$m_{\CC} = m_{\RR}=G= 1$&&$m_{\CC} = m_{\RR}=G= 1$

\\ \\ 
\includegraphics[width=2cm, angle=0]{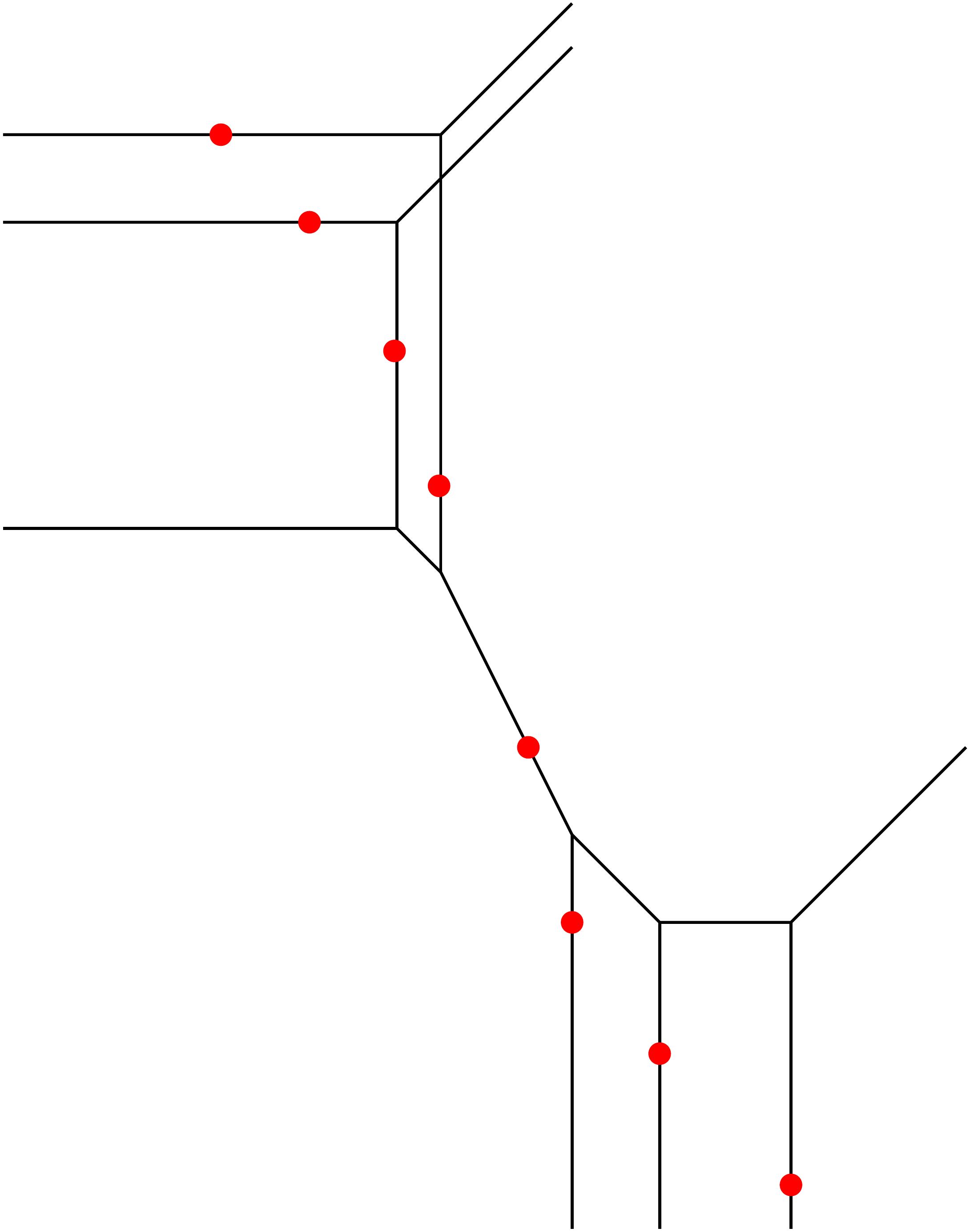} &&
\includegraphics[width=2cm, angle=0]{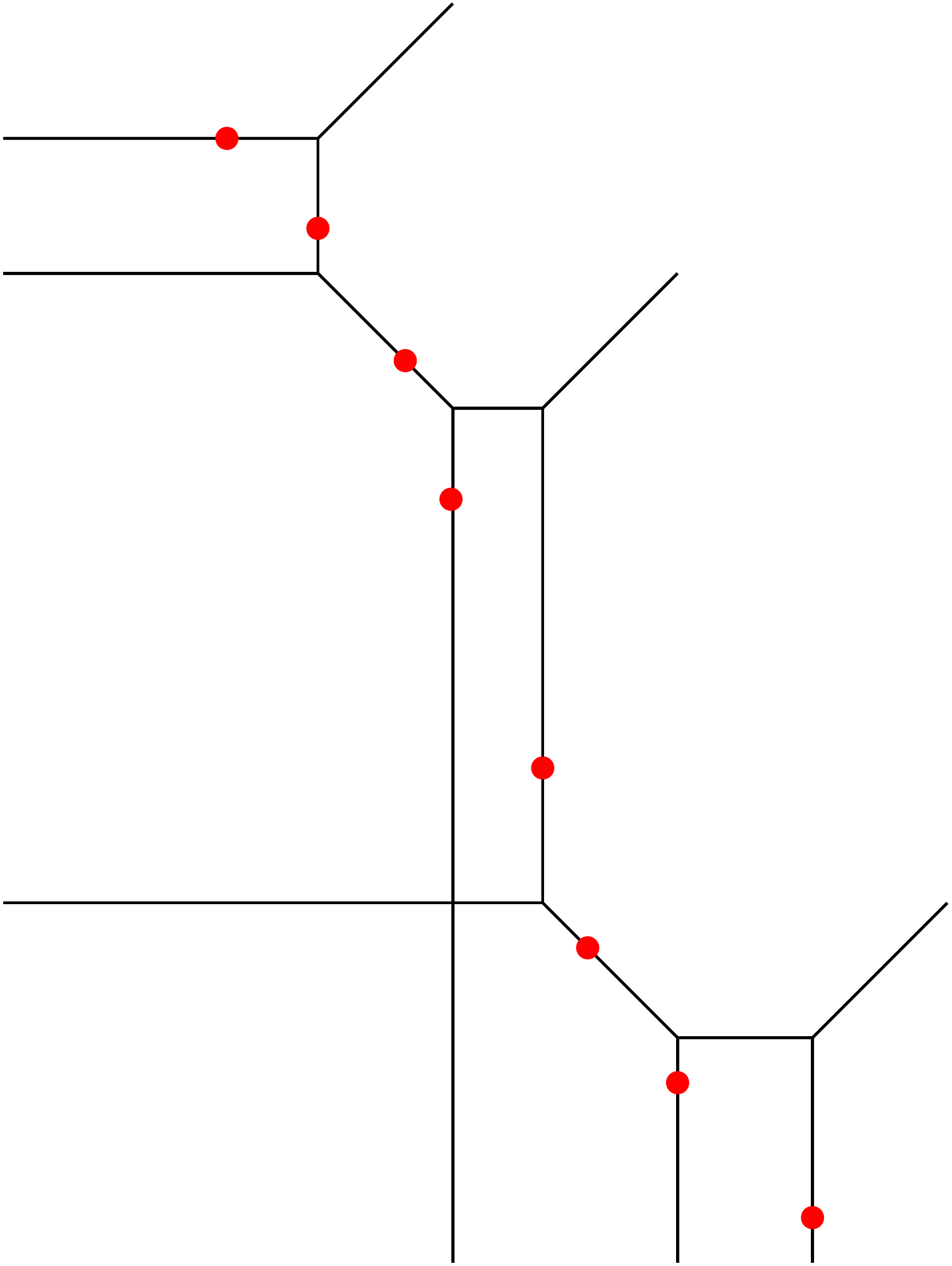}&&
 \includegraphics[width=2cm, angle=0]{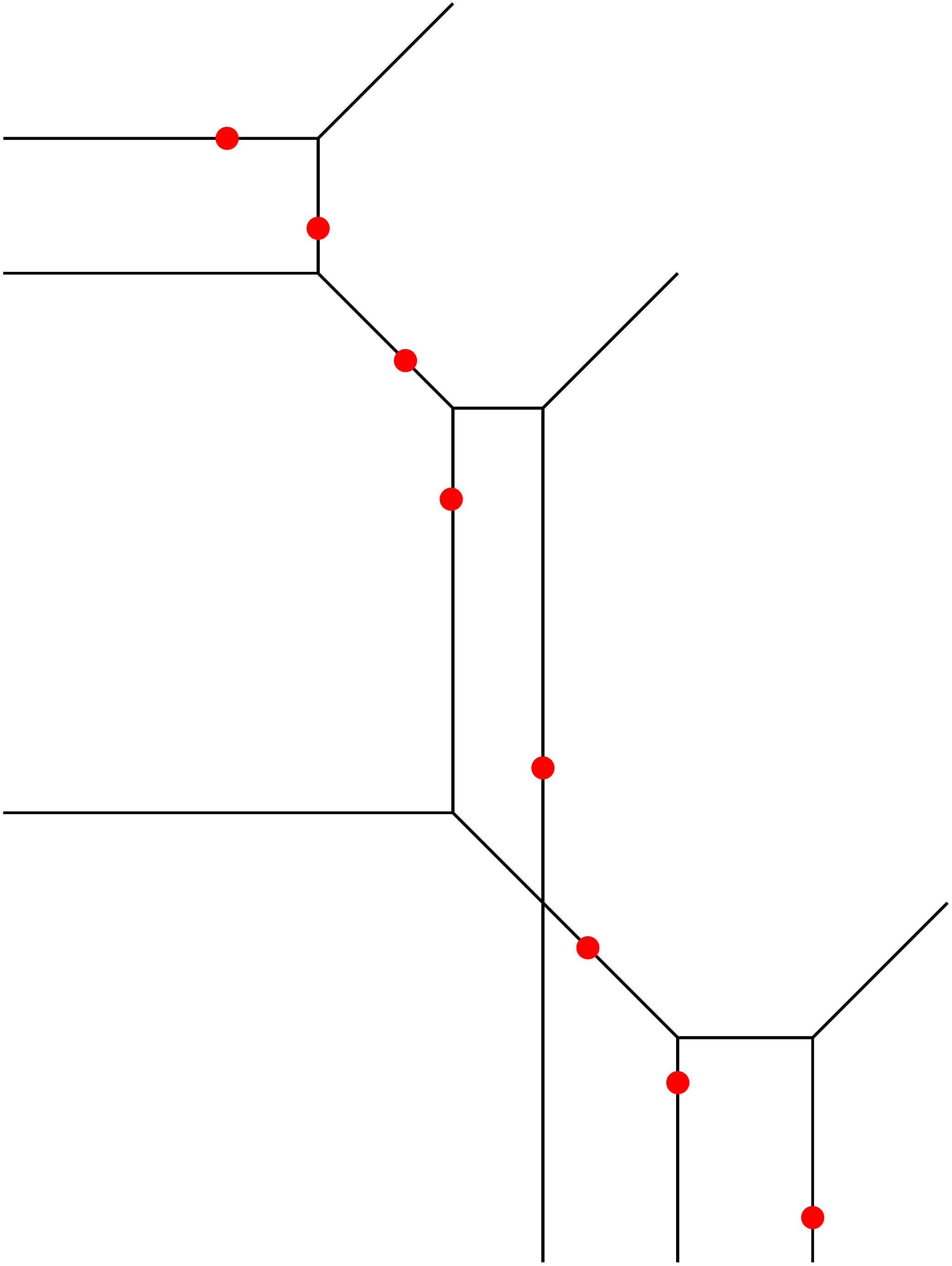} 
\\ $m_{\CC} = m_{\RR}=G= 1$ && $m_{\CC} = m_{\RR}=G= 1$&&$m_{\CC} = m_{\RR}=G= 1$

\\
\\
\includegraphics[width=2cm, angle=0]{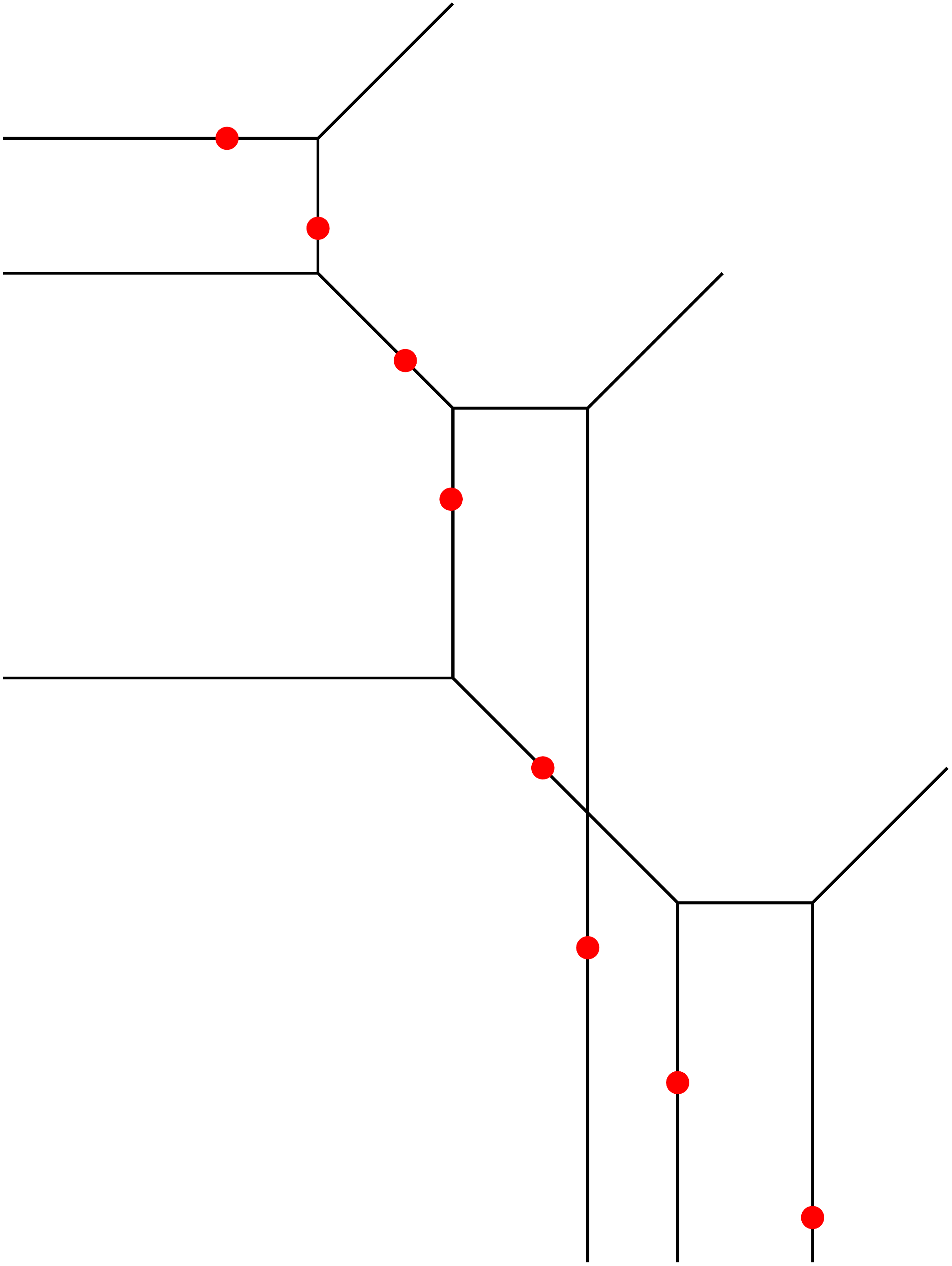}& &
\includegraphics[width=2cm, angle=0]{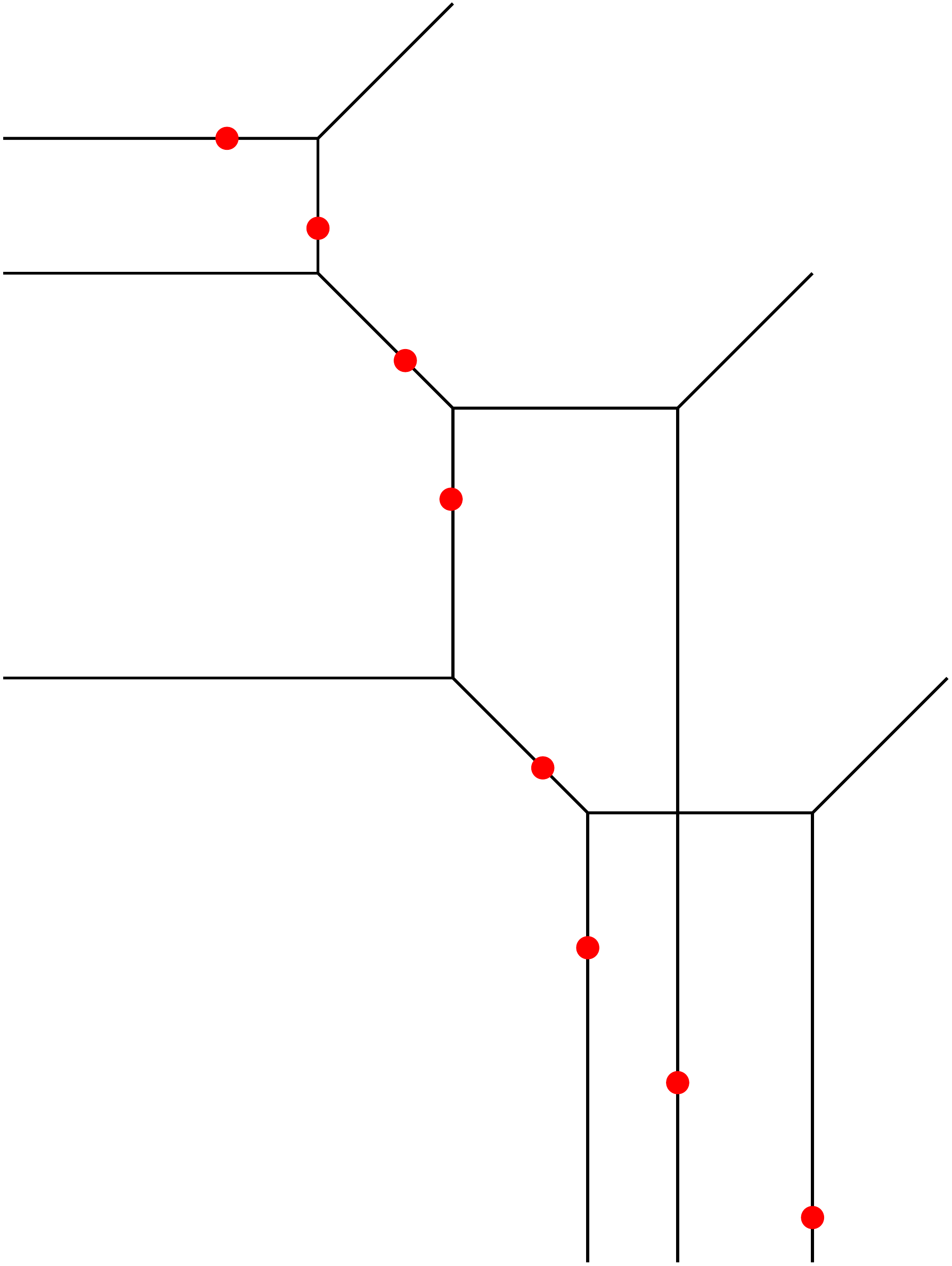}&&
\includegraphics[width=2cm, angle=0]{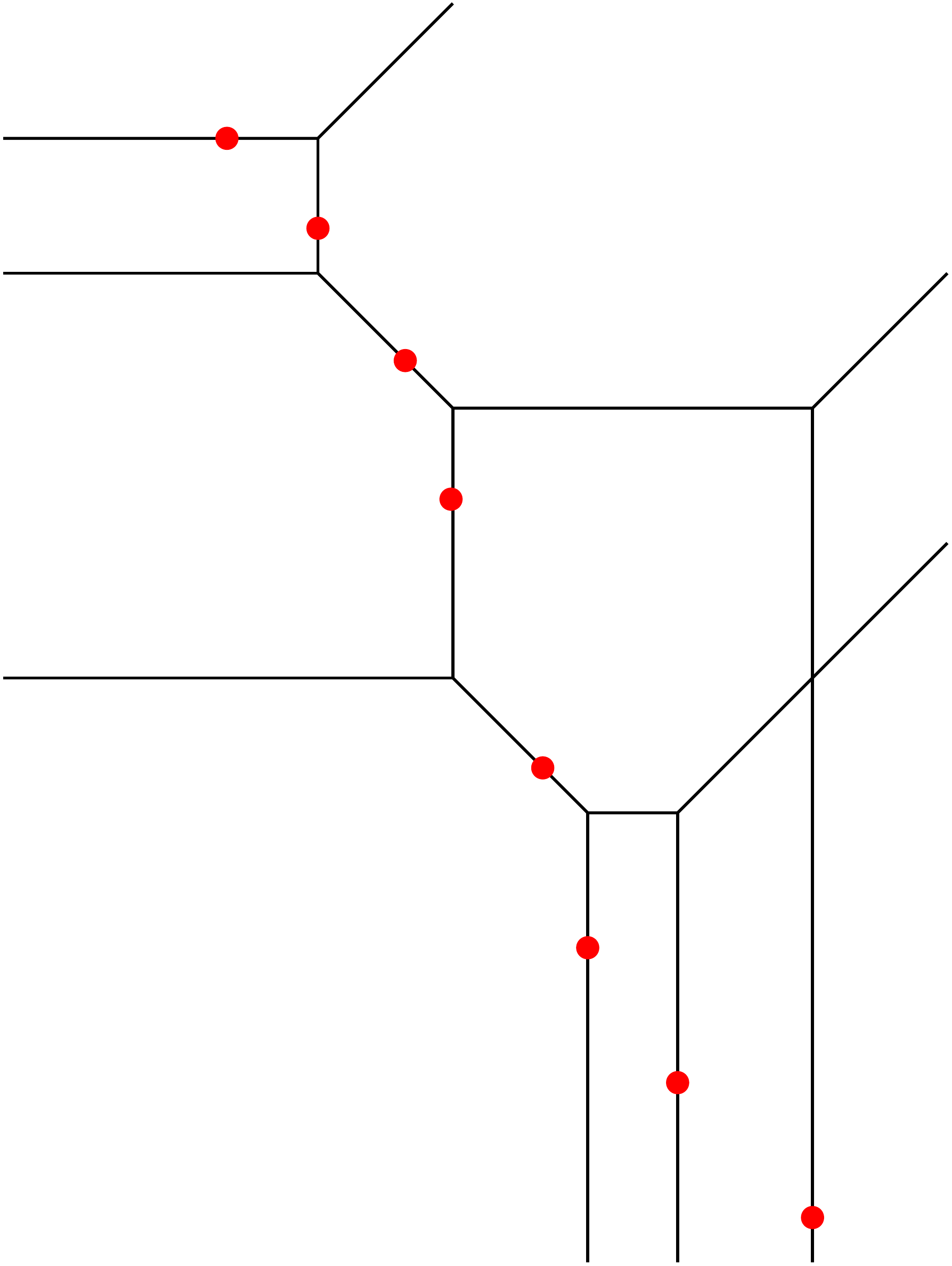} 
\\ $m_{\CC} = m_{\RR}=G= 1$&&$m_{\CC} = m_{\RR}=G= 1$ && $m_{\CC} = m_{\RR}=G= 1$
\end{tabular}
\caption{}
\label{fig:enum 30}
\end{figure}
\end{exa}

\subsection{Quantum enumeration of tropical curves}\label{quantum} 
As we saw in the previous section, the numbers $N_{d, g}$ and $W_d$ can be expressed 
in terms of certain tropical invariants. It turned out that there are many other enumerative tropical invariants
(and complex/real analogs for most of them are unknown).
For example, in the tropical world Welschinger invariants exist in arbitrary genus.

\begin{theorem}[\cite{IKS2}]\label{tropical_invariance} 
Let $\omega$ be a collection of $3d - 1 + g$ points in general position in $\RR^2$. Then, the number
of irreducible nodal tropical curves $C$ of degree $d$ and genus $g$ in $\RR^2$,
counted with multiplicities $m_{\RR}(C)$, which pass through the points of $\omega$
does not depend on the choice of a {\rm (}generic{\rm )\/} collection $\omega$.
\end{theorem} 

Denote by $W^{trop}_{d, g}$ the invariant provided by Theorem \ref{tropical_invariance}.
This 
theorem
can be generalized in the following way. 

F.~Block and L.~G\"ottsche \cite{BlGo14} proposed a new type of multiplicities for tropical curves 
(a motivation for these multiplicities is provided by a Caporaso-Harris type calculation of the refined Severi degrees;
the latter degrees were introduced by G\"ottsche in connection with \cite{KoShTh11}). 
Consider again a collection $\omega$ of $3d-1+g$ points in general position in $\RR^2$,
and let $C$ be an irreducible nodal tropical curve of degree $d$ and genus $g$ which pass through the points of $\omega$.
To each trivalent vertex $v$ of $C$ we associate
$$
G(v) = \frac{q^{m_{\CC}(v)/2} - q^{-m_{\CC}(v)/2}}{q^{1/2} - q^{-1/2}}. 
$$ 
Put
$$
G(C) = \prod_v G(v).
$$
where the product is taken over all trivalent vertices of $C$.
The value of the Block-G\"ottsche multiplicity $G(C)$ at $q = 1$ is $m_{\CC}(C)$.
It is not difficult to check that the value of $G(C)$ at $q = -1$ is equal to $m_{\RR}(C)$,
so the Block-G\"ottsche multiplicities interpolate between the complex and real multiplicities for tropical curves. 

\begin{theorem}[\cite{IteMik13}]\label{refined_invariance}   
Let $\omega$ be a collection of $3d - 1 + g$ points in general position in $\RR^2$.  
Then, the sum of the  Block-G\"ottsche multiplicities $G(C)$
over the irreducible nodal tropical curves $C$ of degree $d$ and genus $g$ in $\RR^2$ 
which pass through the points of $\omega$
does not depend on the choice of a {\rm (}generic{\rm )\/} collection $\omega$.
\end{theorem} 

Denote by $G_{d, g}$ the invariant provided by Theorem \ref{refined_invariance}.
Again, we have $G_{d, g} = 0$ whenever $g > \frac{(d - 1)(d - 2)}{2}$. 
\begin{exa}\label{ex:enum smooth}
We have $G_{1, 0} = G_{2, 0}= G_{3, 1} = 1$ (see Figure \ref{fig:enum smooth}). 
More generally, one 
easily 
shows 
that  
$$G_{d, \frac{(d - 1)(d - 2)}{2}} =1.$$ 
\end{exa}
\begin{exa}\label{ex:enum 30}
We have (see Figure \ref{fig:enum 30}) 
$$G_{3, 0} = q^{-1} + 10 + q.$$
\end{exa} 
\begin{exa}\label{ex:enum deg 4}
Using the technique of floor diagrams, presented in the next section,
one can compute the invariants $G_{4,g}$:
$$
\displaylines{
G_{4, 2} = 3q^{-1} + 21 + 3q, \cr 
G_{4, 1} = 3q^{-2} + 33q^{-1} +153 + 33q + 3q^2, \cr
G_{4, 0} = q^{-3} +13q^{-2} + 94q^{-1} + 404 + 94q + 13q^2 + q^3. 
} 
$$
\end{exa} 
Each coefficient of $G_{d, g}$ is an integer valued tropical invariant.
The sum of these coefficients is equal to $N_{d, g}$,
and the
alternating sum of the coefficients of $G_{d, g}$ is equal to $W^{trop}_{d, g}$.
It is not clear what is a complex enumerative interpretation of individual coefficients
of $G_{d, g}$. 

Theorem \ref{refined_invariance} has the following 
corollary. 
\begin{cor}[cf. \cite{IteMik13}]
Fix a non-negative integer $g$ and a positive integer $k$. 
Then, for any sufficiently large integer $d$ and any 
generic collection 
$\omega$ of $3d - 1 + g$ points in
$\RR^2$, there exists an
irreducible nodal tropical curve $C$ of degree $d$ and genus $g$ in $\RR^2$
such that $C$ passes through the points of $\omega$ and 
$m_{\CC}(C)\ge k$. 
\end{cor} 

\subsection{Floor diagrams}\label{sec:floor_diagrams}
Theorem \ref{correspondence} reduces the problem of 
enumeration of complex (or real) curves 
to calculation of the corresponding tropical invariants.  
One of the most efficient techniques for computation of enumerative tropical invariants
(and, in particular, Gromov-Witten invariants $N_{d, g}$ and Welschinger invariants $W_d$)
is based on so-called {\it floor diagrams}. 
Floor diagrams are  related to
the Caporaso-Harris approach \cite{CapHar1};
we refer to \cite{Bru14}
for more details. 

\begin{definition}\label{floor_diagrams}
A {\rm (}plane{\rm )\/} floor diagram of degree $d$ and genus $g$ is the data of 
a connected 
oriented graph $\mathcal D$ 
{\rm (}considered as a topological object; edges of $\mathcal D$ are not necessarily compact{\rm )\/}   
which satisfy the following conditions: 
\begin{itemize}
\item the oriented graph $\mathcal D$ is acyclic;
\item $\mathcal D$ has exactly $d$ vertices; 
\item the first Betti number $b_1({\mathcal D})$ of $\mathcal D$ is equal to $g$;
\item 
 each edge  has a weight which is a positive integer number; 
\item there are exactly $d$ non-compact edges of $\mathcal D$; 
all of them are of weight $1$ and are oriented towards their unique adjacent vertex;
\item for each vertex of $\mathcal D$, the sum of weights of incoming edges
is greater by $1$ than the sum of weights of outgoing edges. 
\end{itemize} 
\end{definition}

A floor diagram inherits a partial ordering from its
orientation. A map $m$ between two partially ordered sets is said to be
{\it increasing} if 
$$m(i) > m(j) \Rightarrow i > j.$$ 
\begin{definition}\label{marking} 
A marking of a floor diagram $\mathcal D$ of degree $d$ and genus $g$ 
is an increasing map $m: \{1, \ldots, 3d - 1 + g\} \to {\mathcal D}$ 
such that for any edge or vertex $x$ of $\mathcal D$,
the set $m^{-1}(x)$ consists of exactly one element.
A floor diagram enhanced with a marking is called a marked floor diagram. 
\end{definition}

We consider the floor diagrams up to a natural equivalence: 
two floor diagrams $\mathcal D$ and $\mathcal D'$ are equivalent
if there exists a homeomorphism of oriented graphs $\mathcal D$ and $\mathcal D'$
which respects the weights of all edges. 
Similarly, two marked floor diagrams $({\mathcal D}, m)$ and $({\mathcal D}', m')$  
are equivalent if there exists a homeomorphism of oriented graphs $\varphi: {\mathcal D} \to {\mathcal D}'$ 
which respects the weights of all edges and such that 
$m' = \varphi \circ m$. 

\begin{exa}
Figure \ref{comp cplx1}
shows all floor diagrams (up to equivalence) of degree at most $3$ 
and indicates for each of them the number of possible
markings. Similarly, Figures \ref{degree 4 g=3,2}, \ref{degree 4 g=1}, and \ref{degree 4 g=0}
show all floor diagrams  of degree  $4$.

We use the following convention to depict floor diagrams: 
vertices of $\mathcal D$
are represented by white ellipses, 
edges
are represented by vertical lines, and the orientation is implicitly
from down to up. We specify the 
weight of an edge only if this 
weight is at least $2$. 

\begin{figure}[h]
\centering
\begin{tabular}{cccccc}
\includegraphics[height=1.3cm, angle=0]{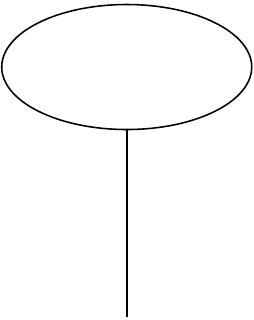}&
\includegraphics[height=2.3cm, angle=0]{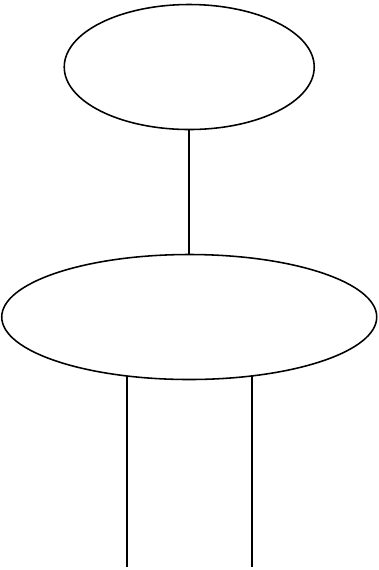}&
\includegraphics[height=3cm, angle=0]{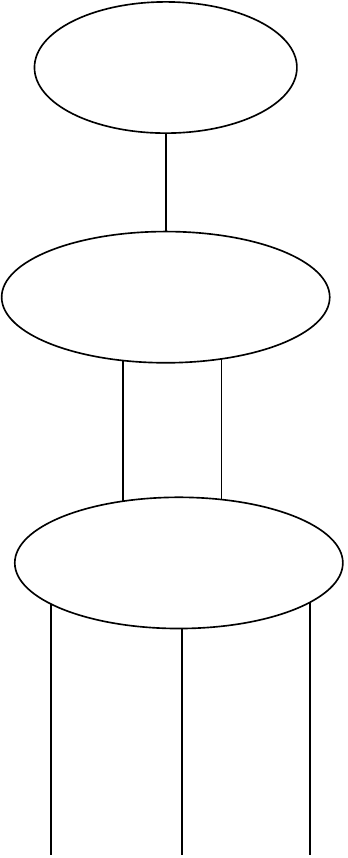}&
\includegraphics[height=3cm, angle=0]{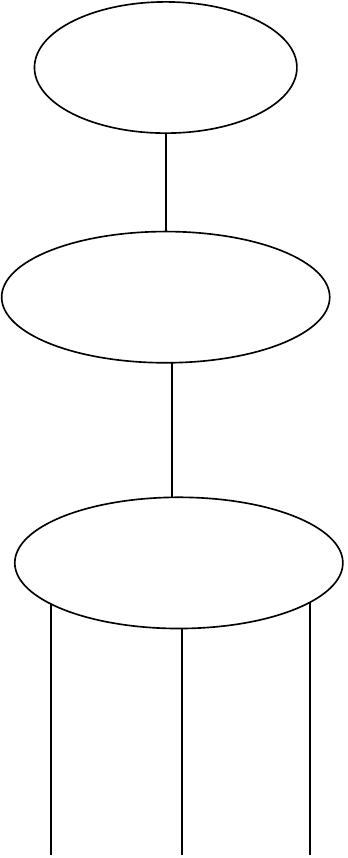}
\put(-15,40){\small{$2$}}&
\includegraphics[height=3cm, angle=0]{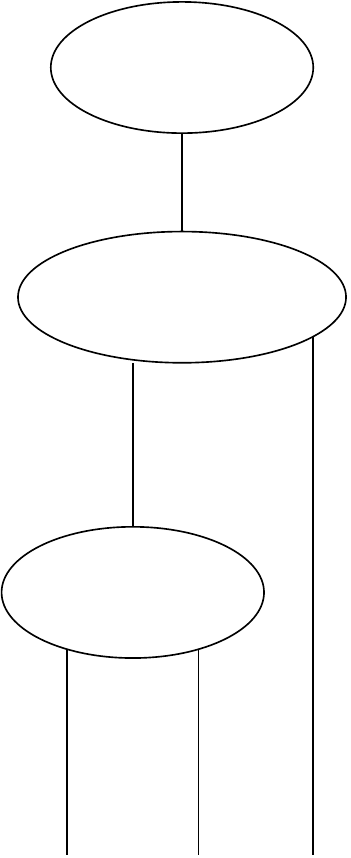}&
\includegraphics[height=3cm, angle=0]{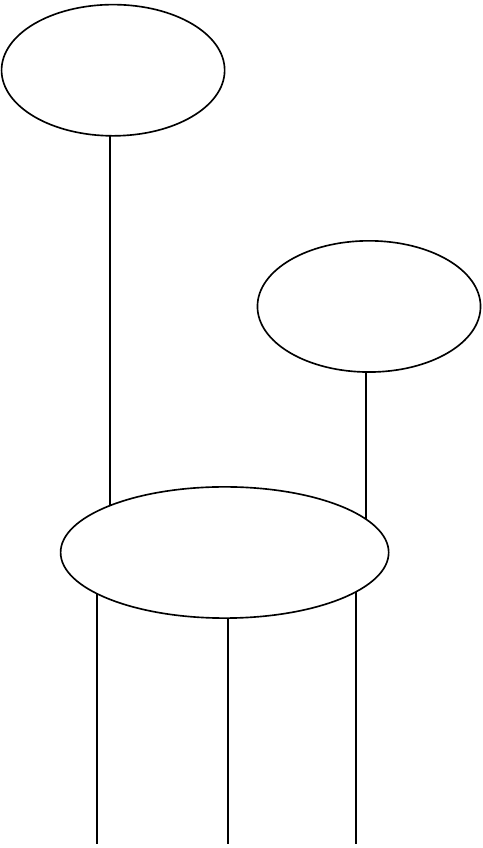}
\\
\\  1 marking &   1 marking &   1 marking &   1 marking & 
5   markings &   3 markings
\end{tabular}
\caption{Floor diagrams of degree $\leq 3$}
\label{comp cplx1}
\end{figure}

\begin{figure}[h]
\centering
\begin{tabular}{cccccc}
\includegraphics[height=3.5cm, angle=0]{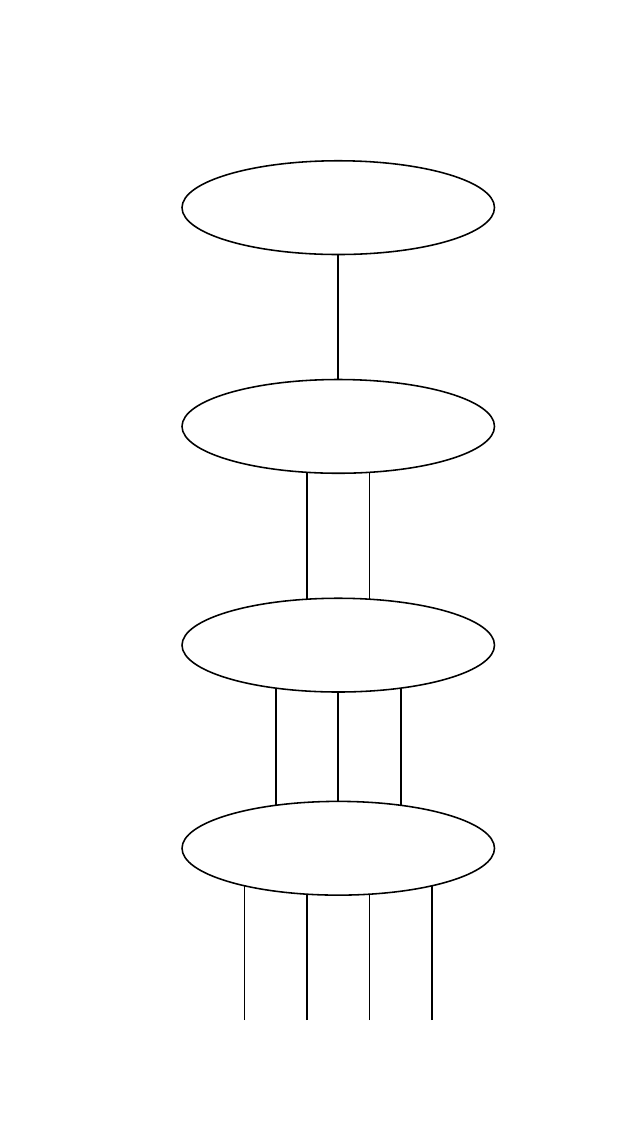}&
\includegraphics[height=3.5cm, angle=0]{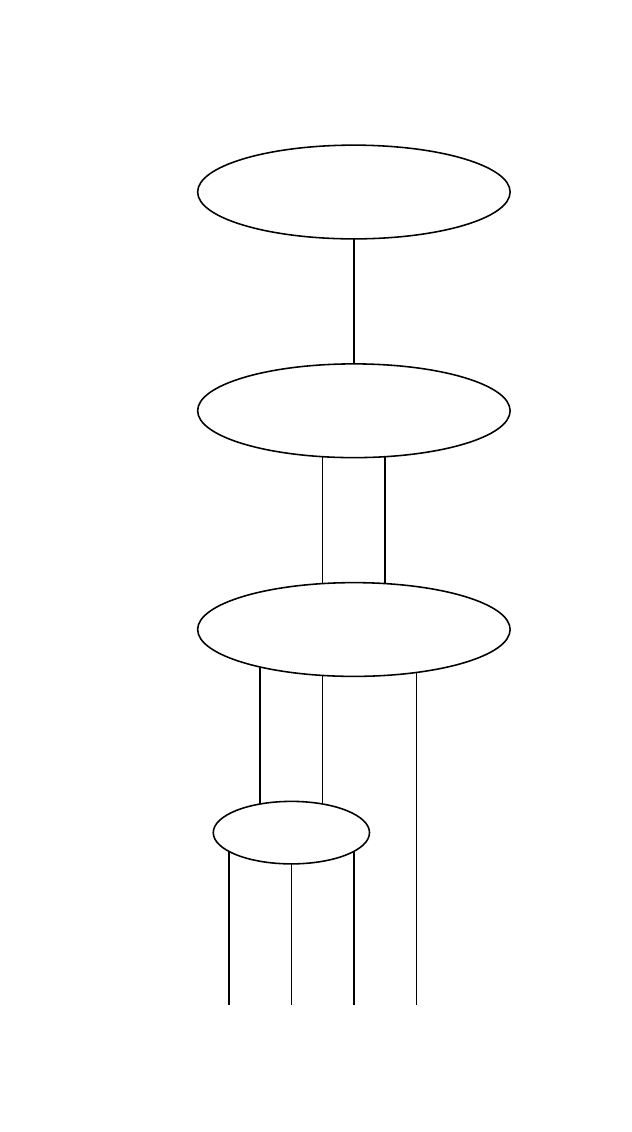}&
\includegraphics[height=3.5cm, angle=0]{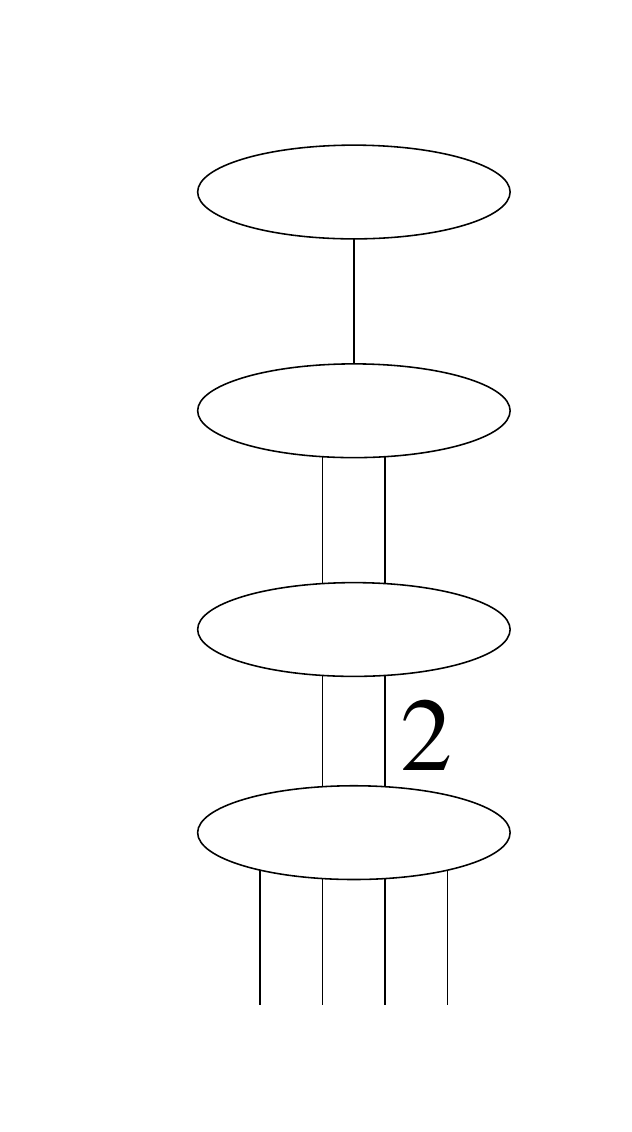}&
\includegraphics[height=3.5cm, angle=0]{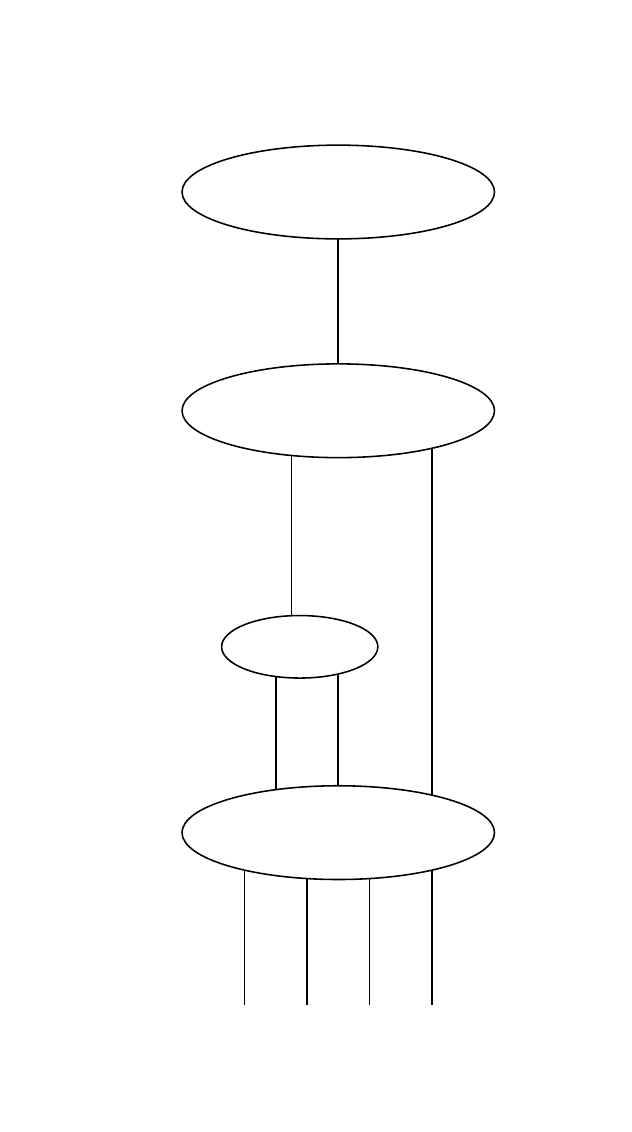}&
\includegraphics[height=3.5cm, angle=0]{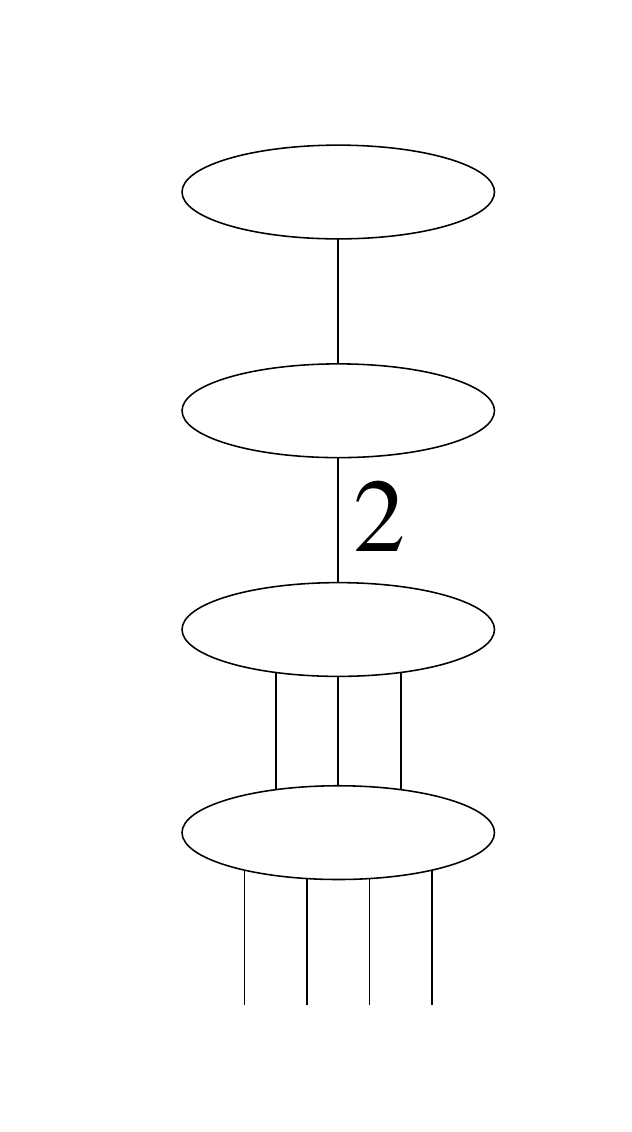}&
\includegraphics[height=3.5cm, angle=0]{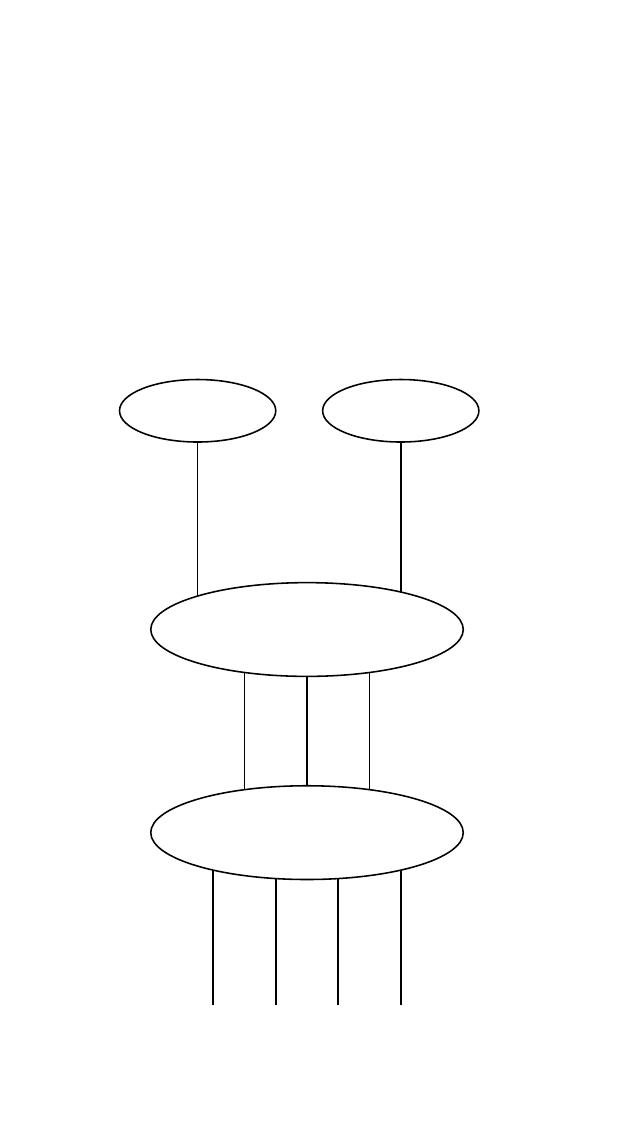}
\\  1 marking &   7 markings &   2 markings &   5 markings & 
1   marking &   3 markings
\end{tabular}
\caption{Floor diagrams of degree $4$ and genus 3 or 2}
\label{degree 4 g=3,2}
\end{figure}

\begin{figure}[h]
\centering
\begin{tabular}{cccccc}
\includegraphics[height=3.5cm, angle=0]{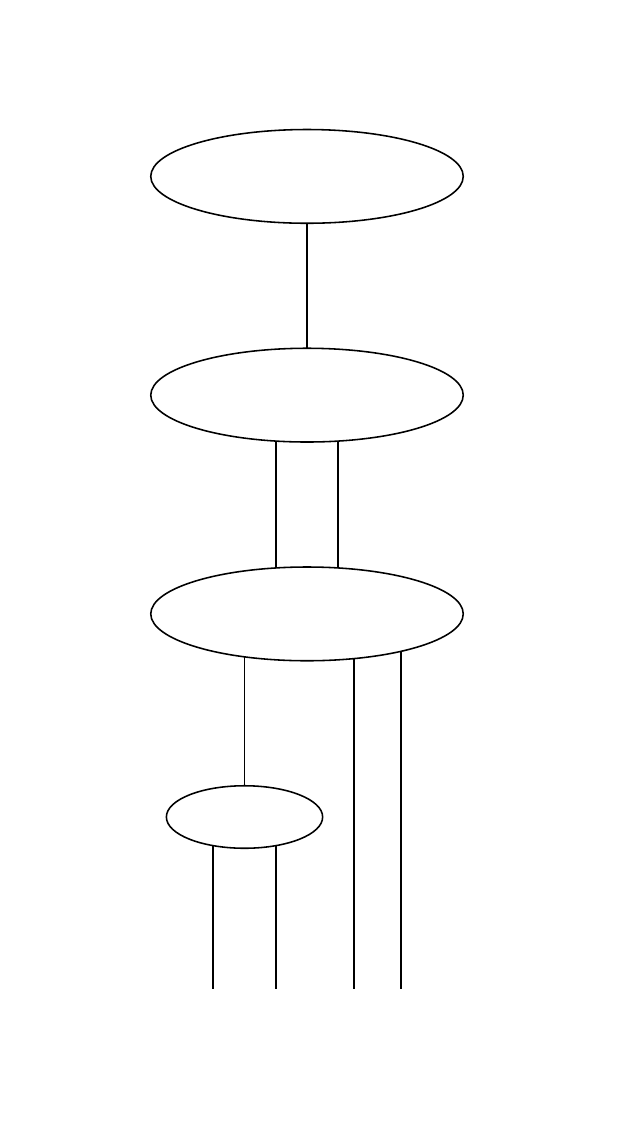}&
\includegraphics[height=3.5cm, angle=0]{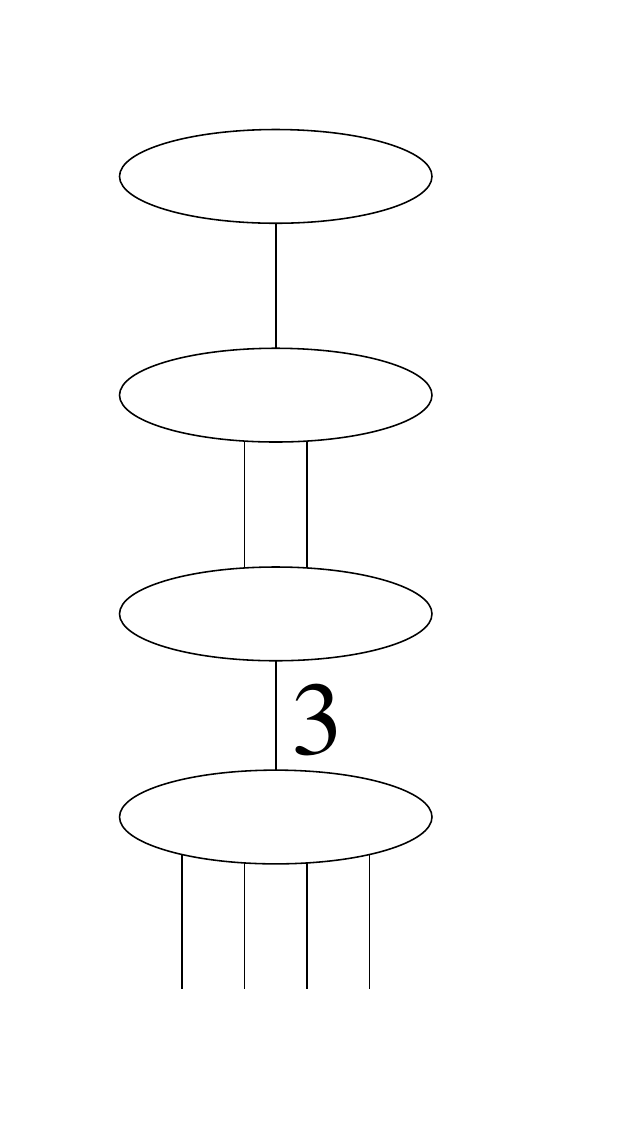}&
\includegraphics[height=3.5cm, angle=0]{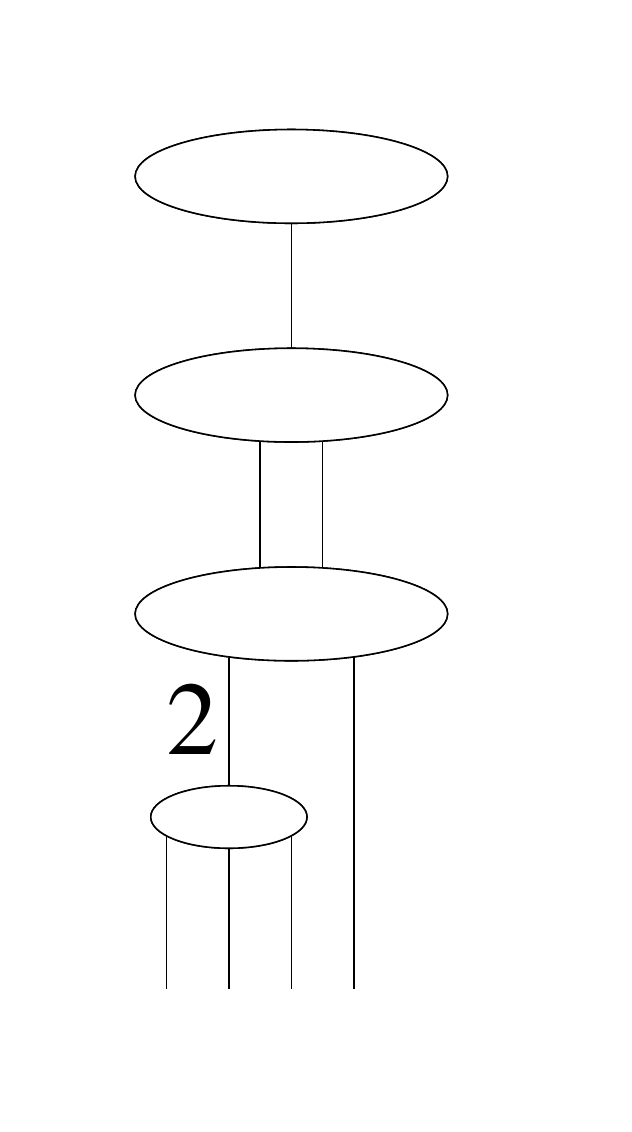}&
\includegraphics[height=3.5cm, angle=0]{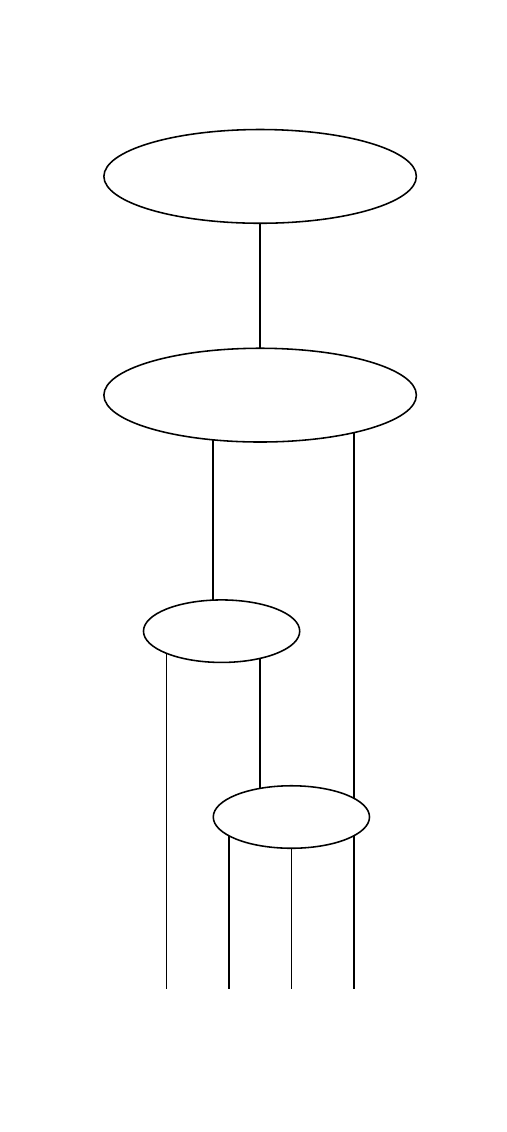}&
\includegraphics[height=3.5cm, angle=0]{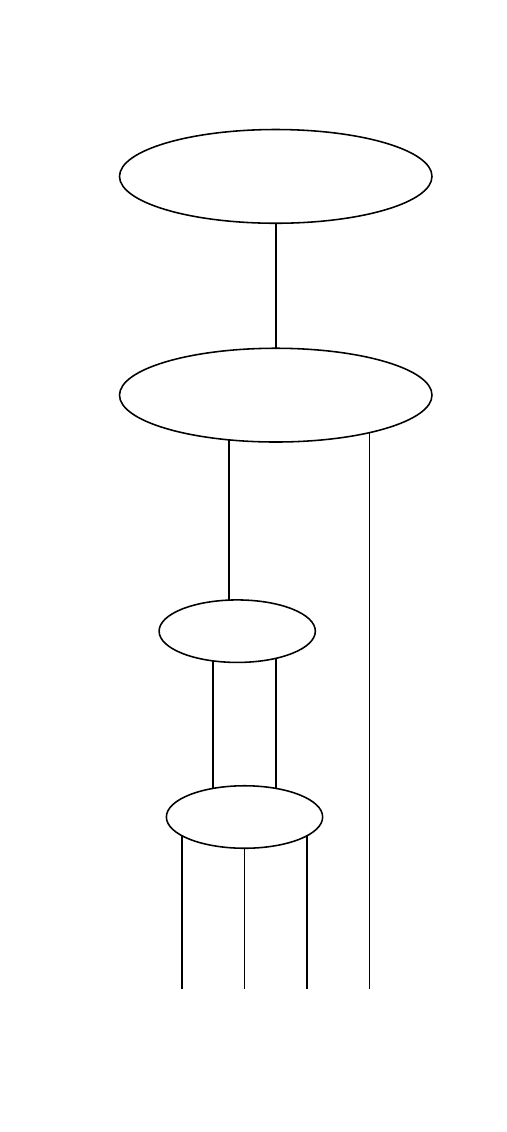}&
\includegraphics[height=3.5cm, angle=0]{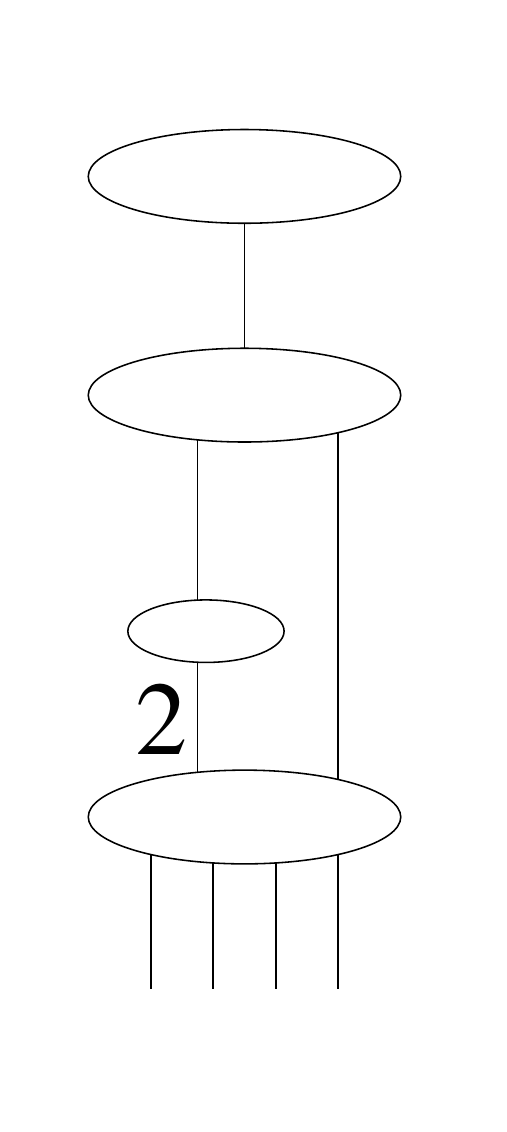}
\\  15 markings &   1 marking &   6 markings &   26 markings & 
9   markings &   4 markings
\\
\includegraphics[height=3.5cm, angle=0]{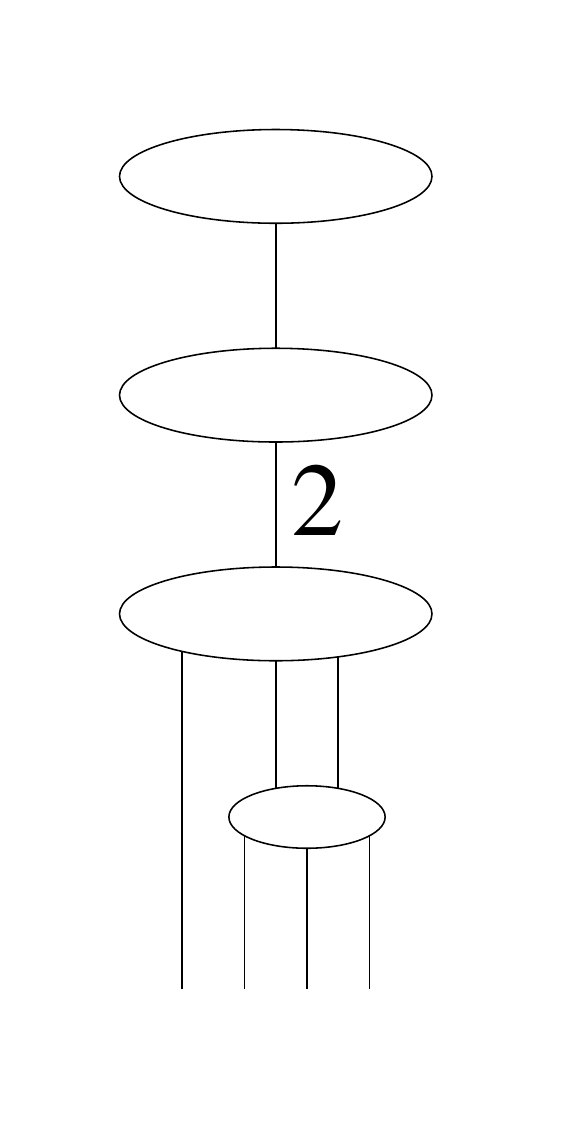}&
\includegraphics[height=3.5cm, angle=0]{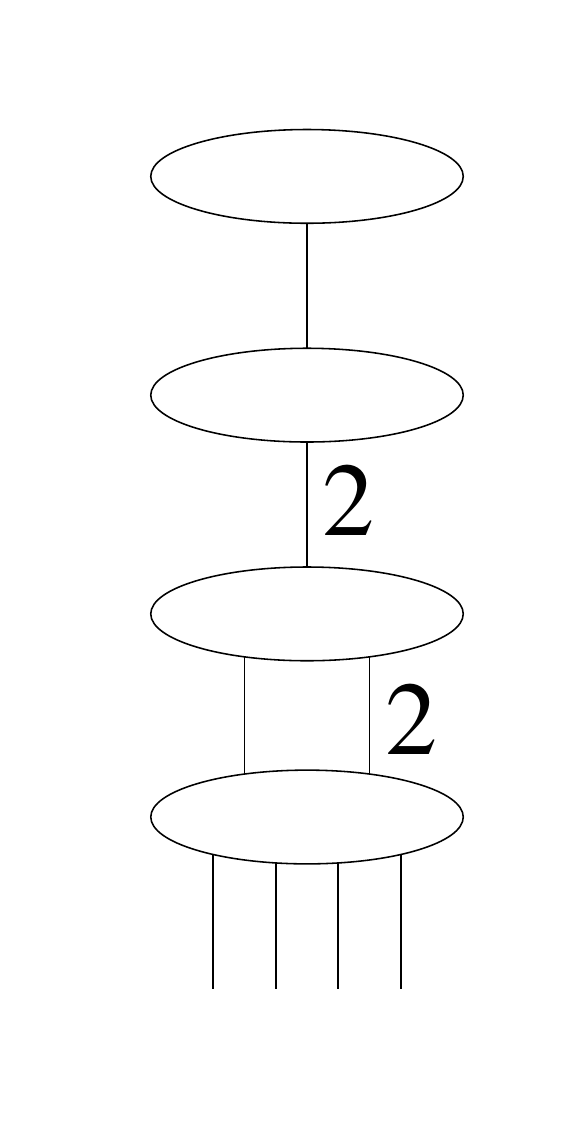}&
\includegraphics[height=3.5cm, angle=0]{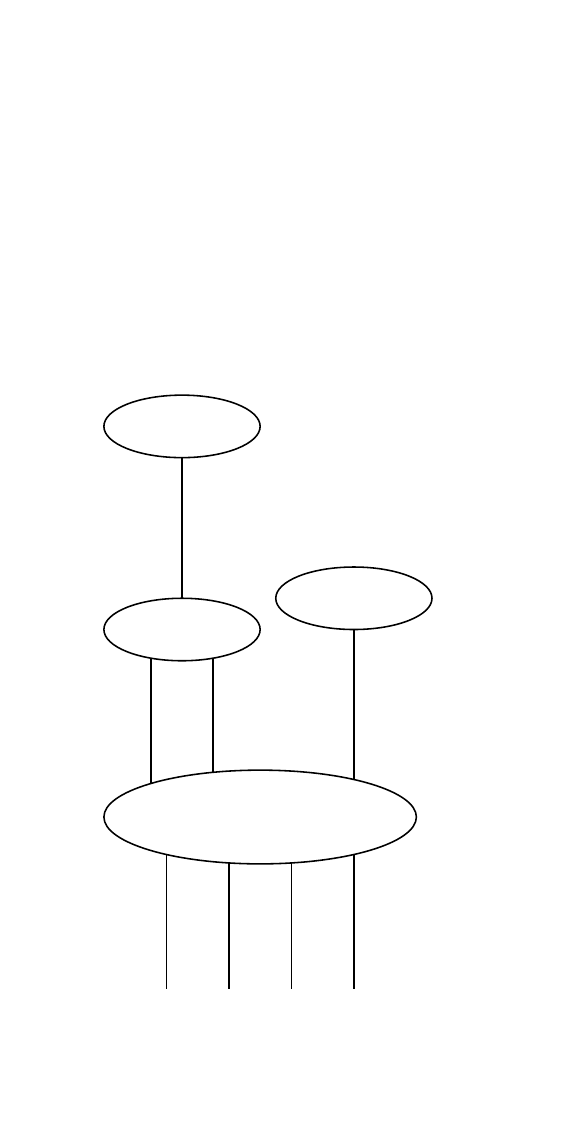}&
\includegraphics[height=3.5cm, angle=0]{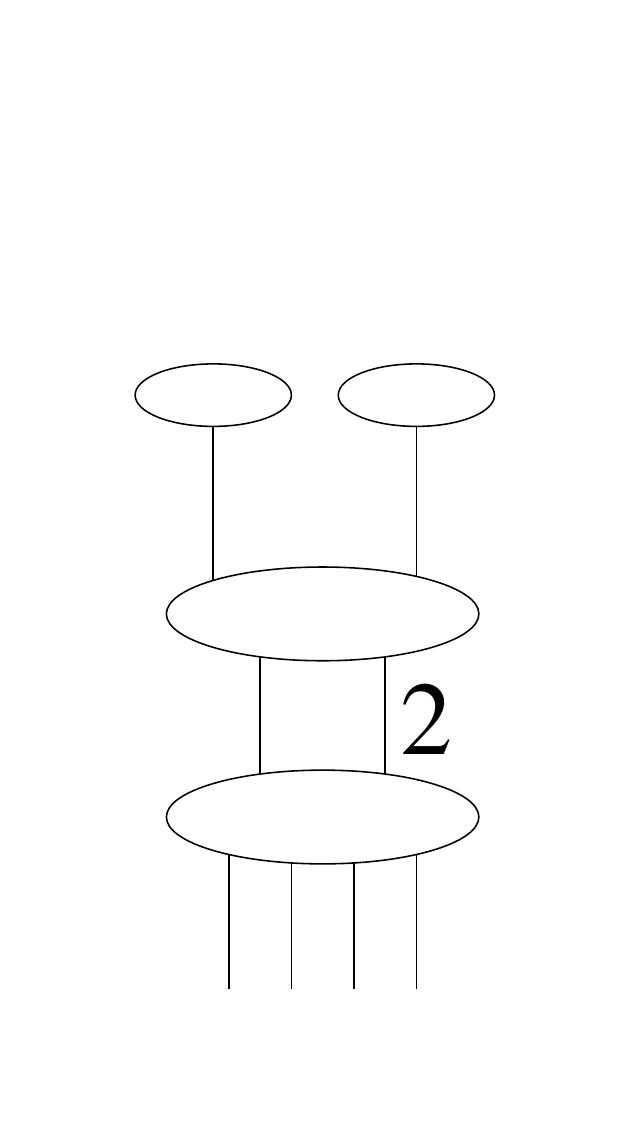}&
\includegraphics[height=3.5cm, angle=0]{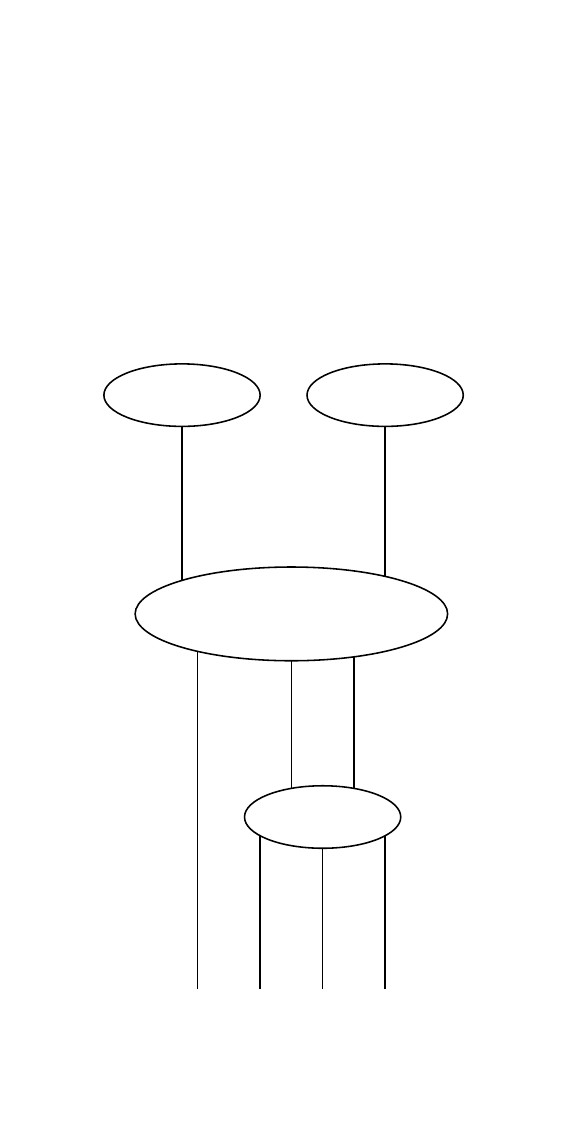}&
\\  7 marking &   2 markings &   21 markings &   6 markings & 
21   markings &  

\end{tabular}
\caption{Floor diagrams of degree $4$ and genus 1}
\label{degree 4 g=1}
\end{figure}

\begin{figure}[h]
\centering
\begin{tabular}{cccccc}
\includegraphics[height=3.5cm, angle=0]{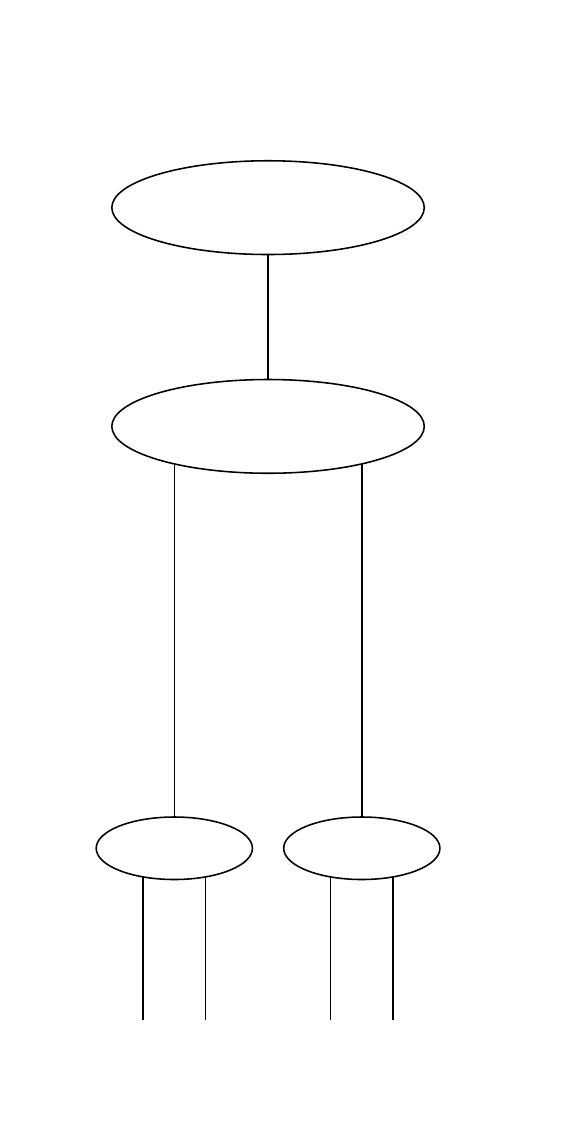}&
\includegraphics[height=3.5cm, angle=0]{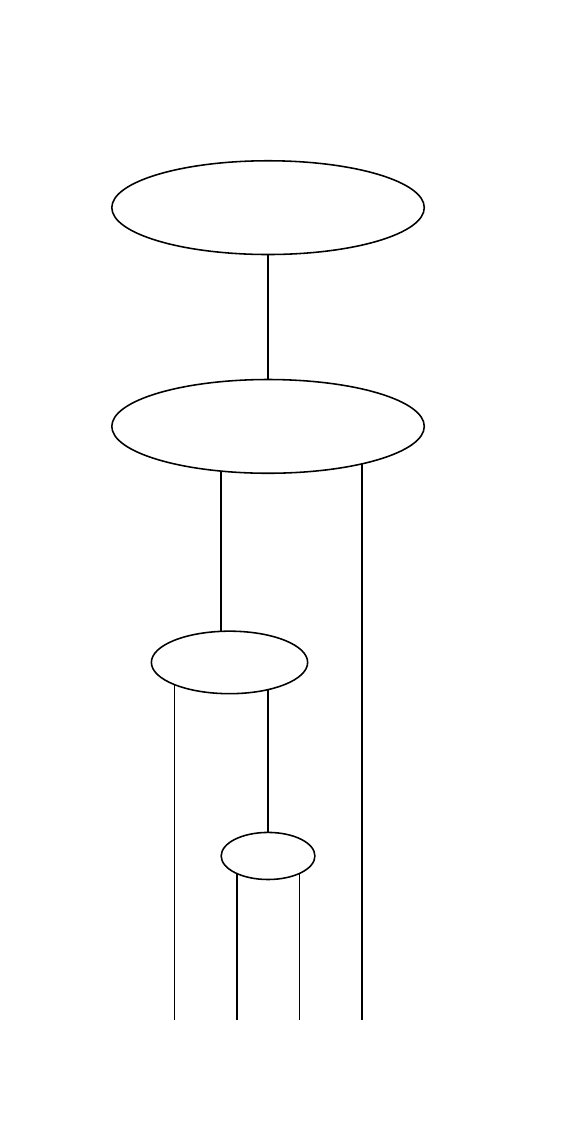}&
\includegraphics[height=3.5cm, angle=0]{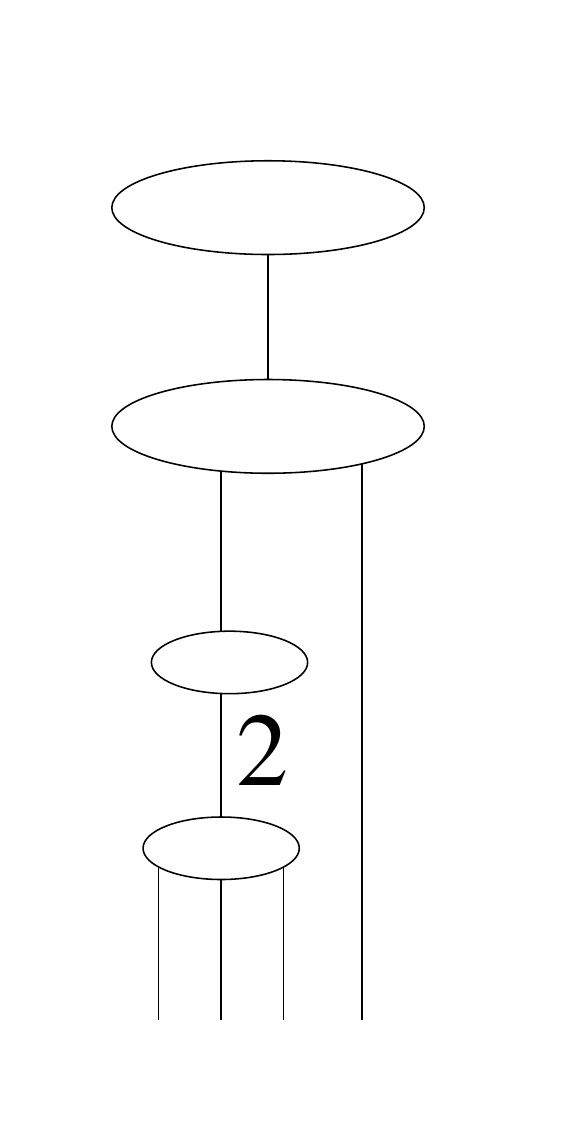}&
\includegraphics[height=3.5cm, angle=0]{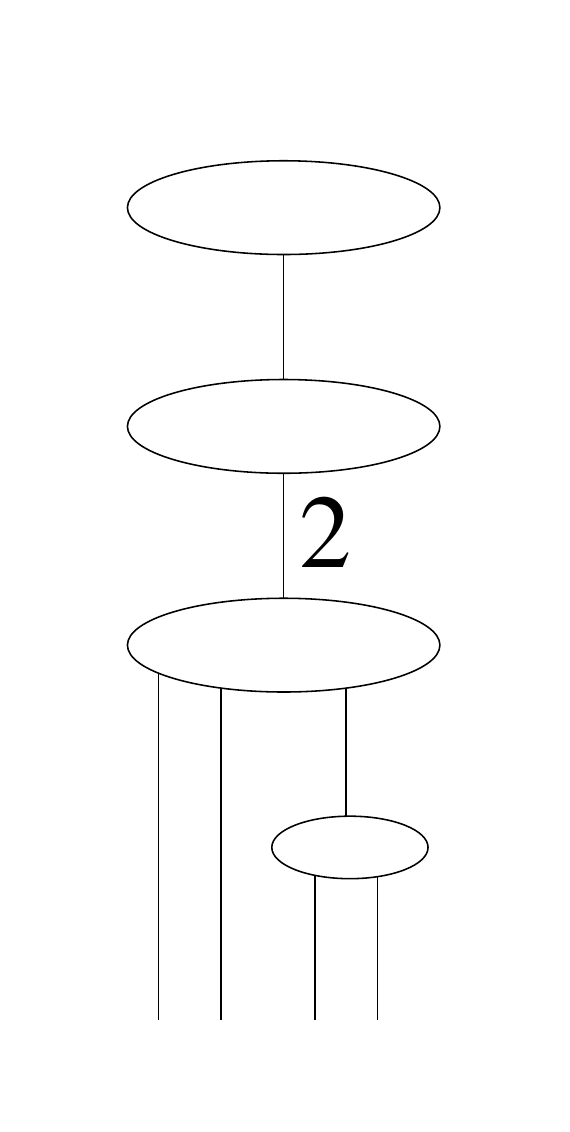}&
\includegraphics[height=3.5cm, angle=0]{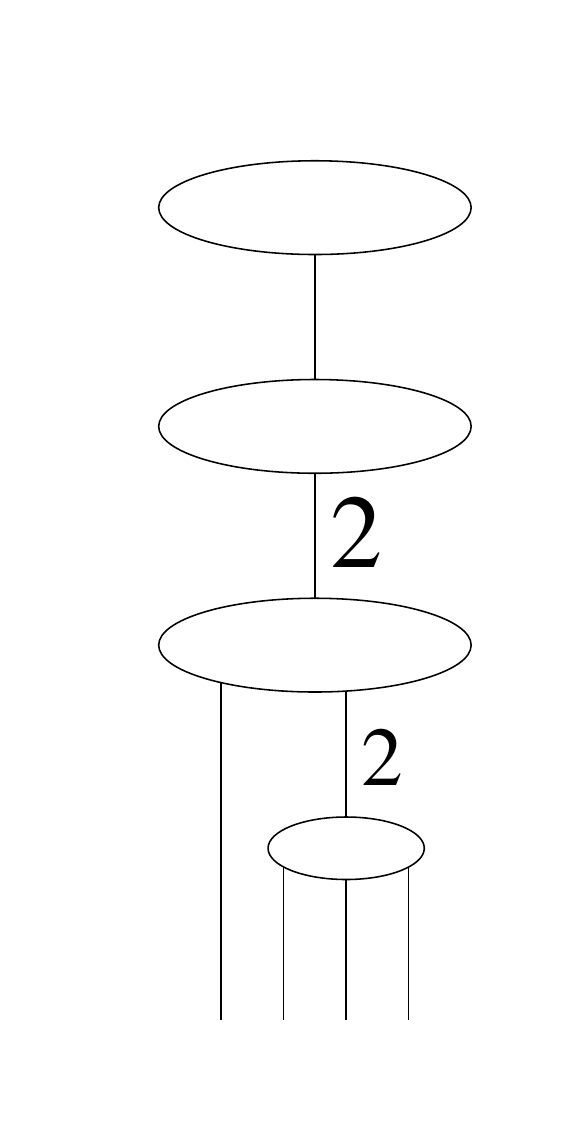}&
\includegraphics[height=3.5cm, angle=0]{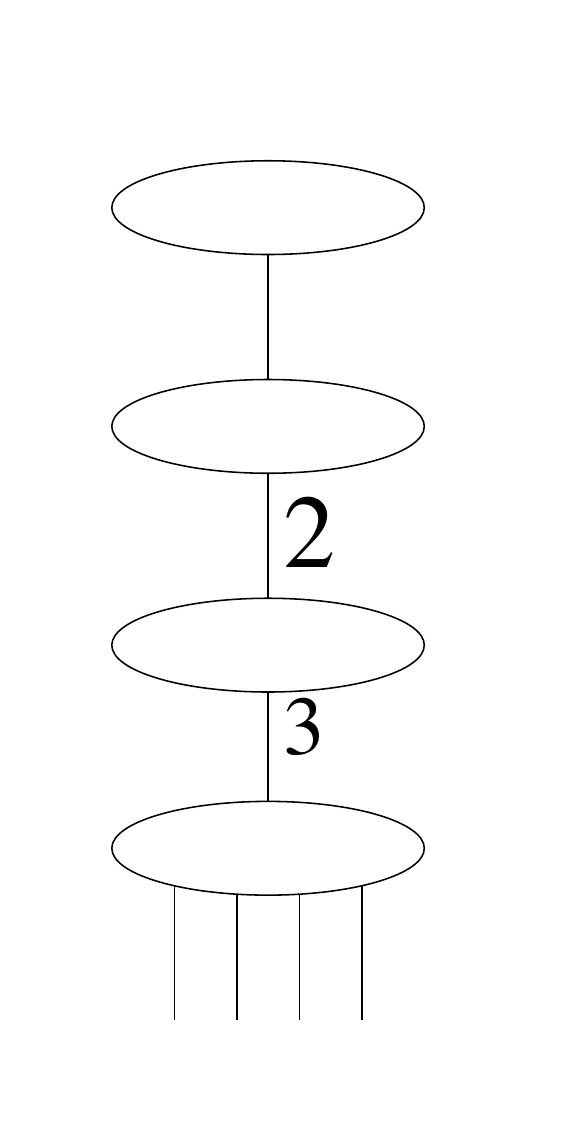}
\\  35 markings &   40 markings &   8 markings &   15 markings & 
6   markings &   1 marking
\\
\includegraphics[height=3.5cm, angle=0]{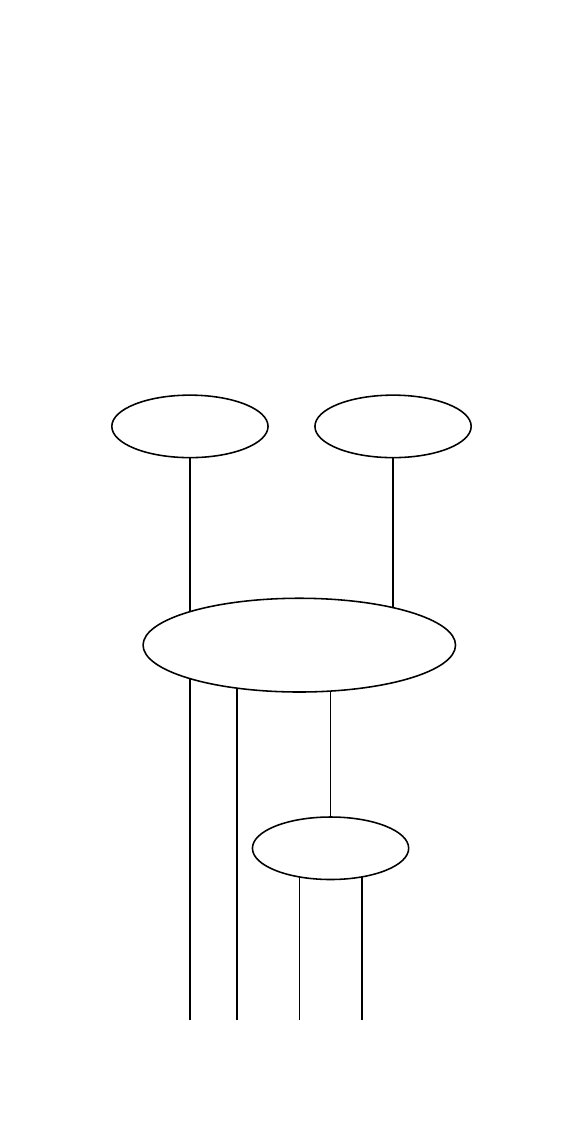}&
\includegraphics[height=3.5cm, angle=0]{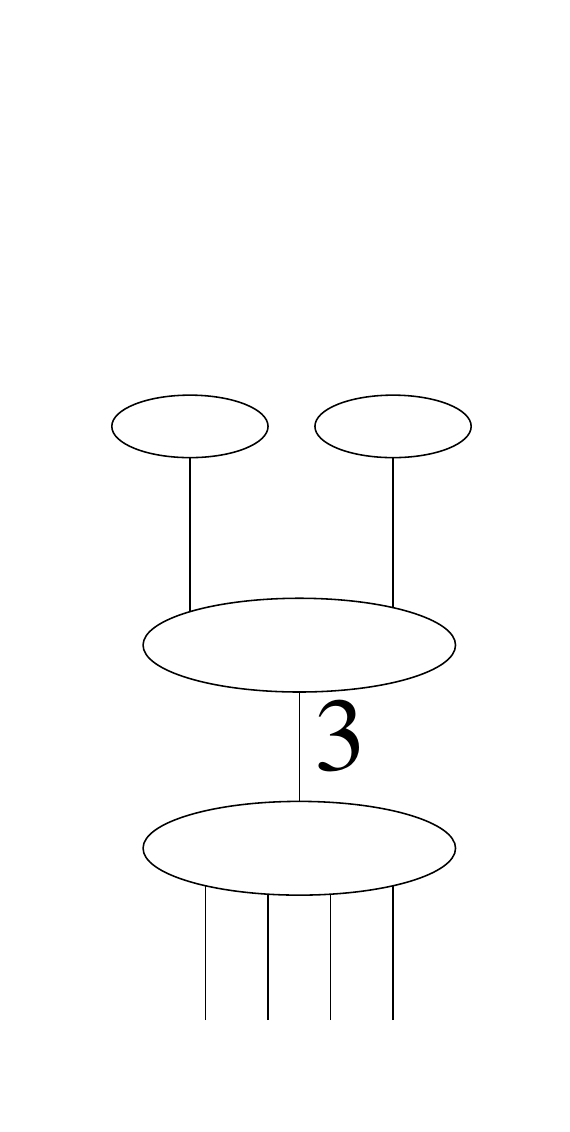}&
\includegraphics[height=3.5cm, angle=0]{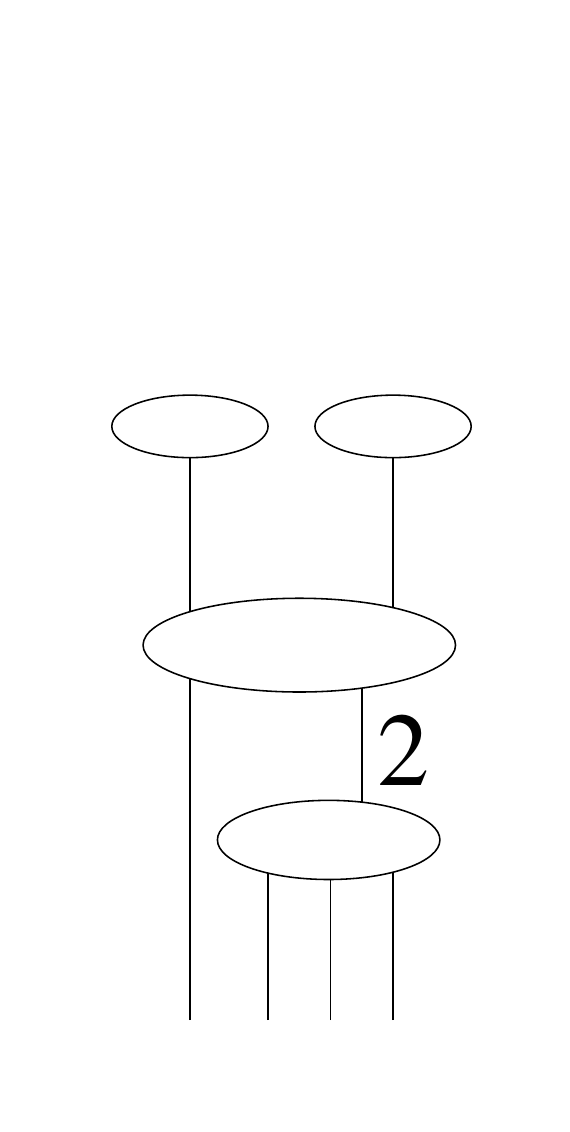}&
\includegraphics[height=3.5cm, angle=0]{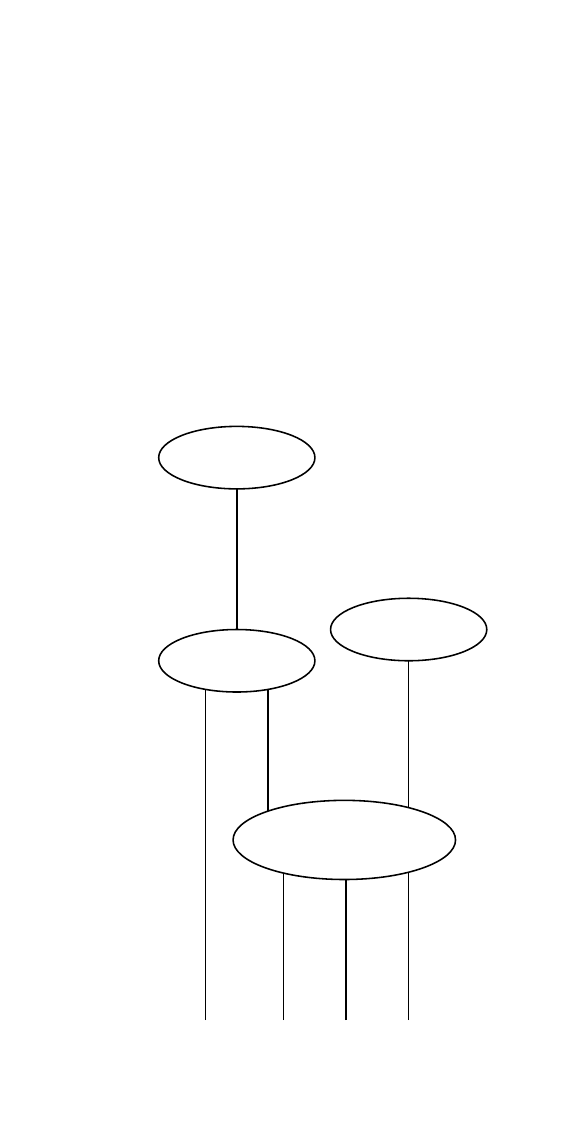}&
\includegraphics[height=3.5cm, angle=0]{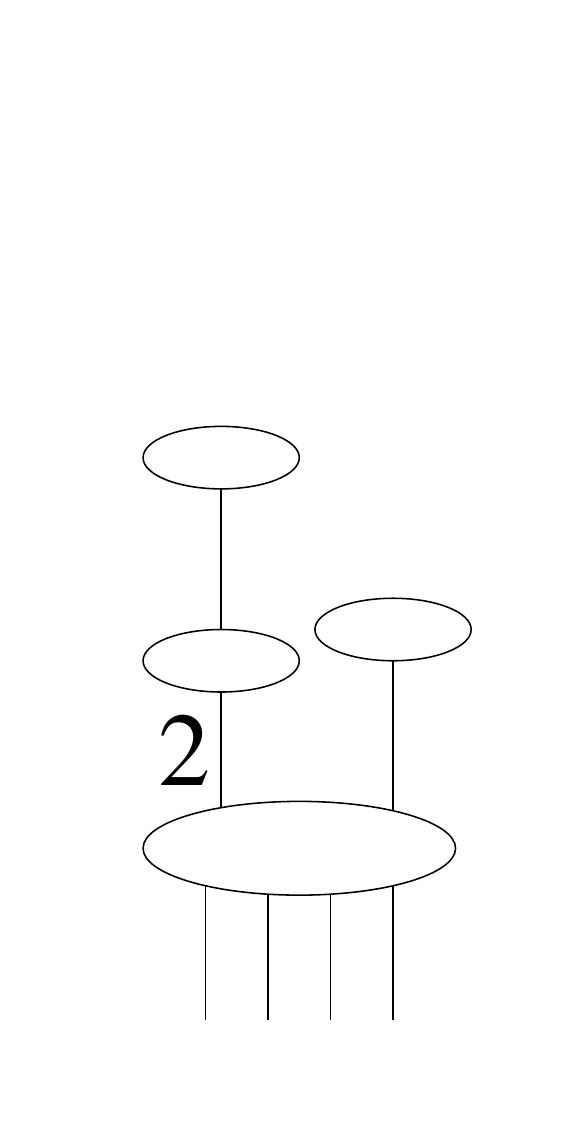}&
\includegraphics[height=3.5cm, angle=0]{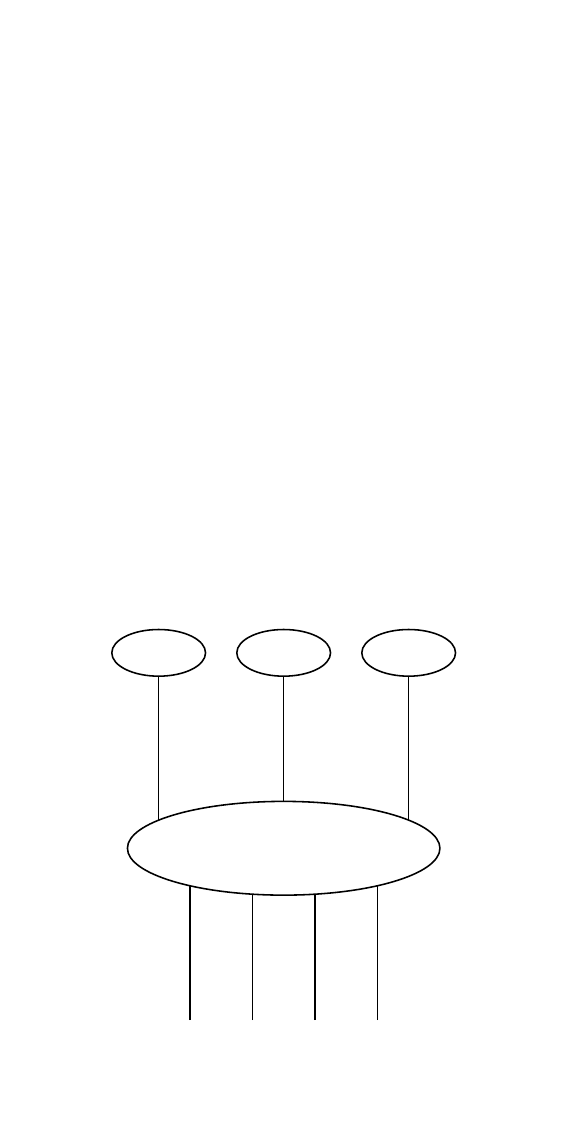}
\\  45 marking &   3 markings &   18 markings &   102 markings & 
15   markings &  15  markings

\end{tabular}
\caption{Floor diagrams of degree $4$ and genus 0}
\label{degree 4 g=0}
\end{figure}
\end{exa}

To each floor diagram $\mathcal D$ of degree $d$ and genus $g$ we can associate
its complex, real and Block-G\"ottsche multiplicities
putting, respectively, 
$$
\displaylines{
m_{\CC}(\mathcal D) = \prod_e (w(e))^2, \cr
m_{\RR}(\mathcal D) = \prod_e r(e), \cr
G({\mathcal D}) = \prod_e \big(\frac{q^{w(e)/2} - q^{-w(e)/2}}{q^{1/2} - q^{-1/2}}\big)^2,
}
$$
where the products are taken over all edges $e$ of $\mathcal D$, the number $w(e)$ 
is the weight of $e$, and $r(e)$ is equal to $0$ if $w(e)$ is even and
is equal to $1$ otherwise. 
Note that $m_{\CC}(\mathcal D)$ is the value of $G({\mathcal D})$ at $q = 1$, and
$m_{\RR}(\mathcal D)$ is the value of $G({\mathcal D})$ at $q = - 1$.   

\begin{figure}[!h]
\centering
\begin{tabular}{ccccc}
\includegraphics[width=3cm, angle=0]{Figures/Cubic7.pdf}&
\hspace{5ex} &
\includegraphics[width=3cm, angle=0]{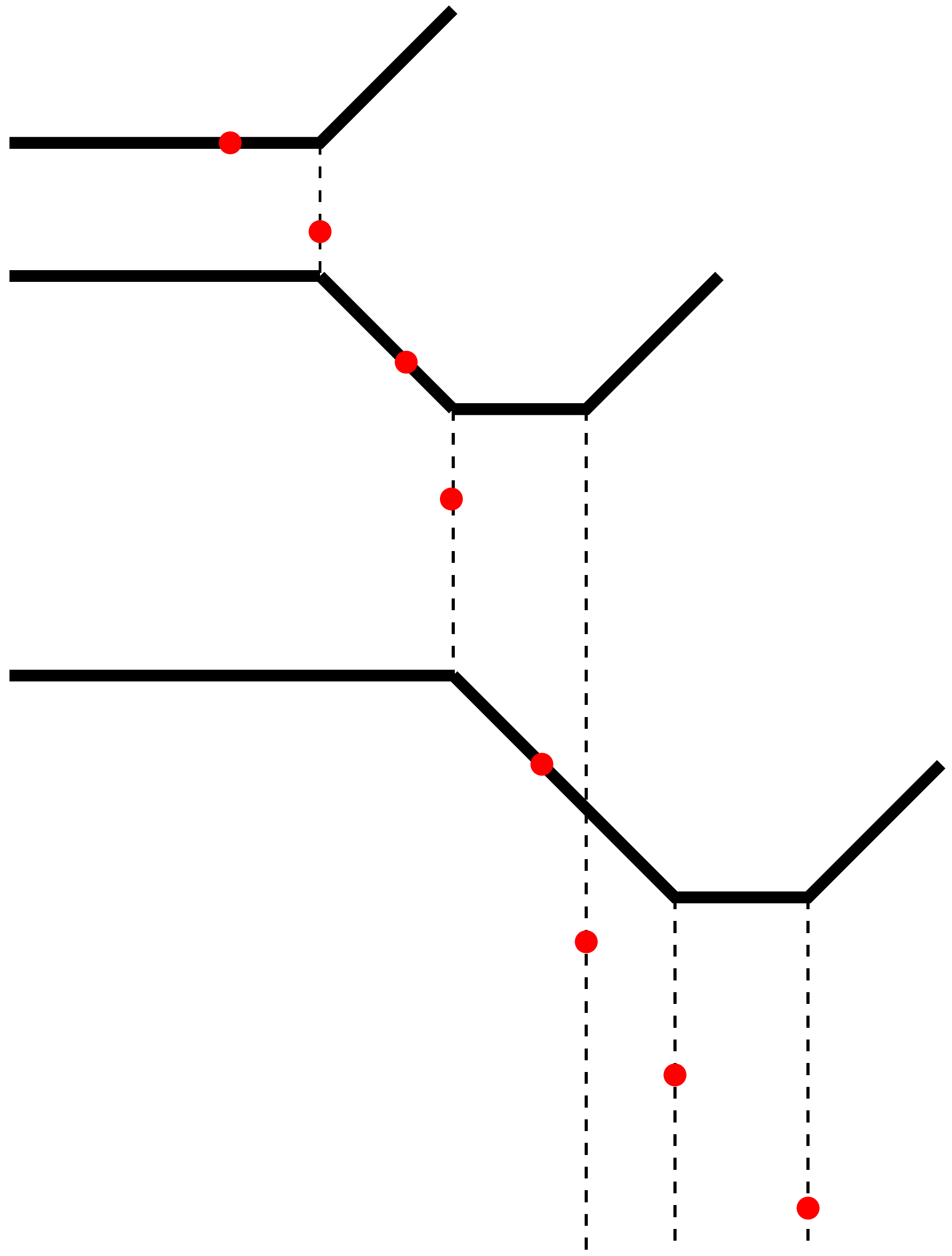}&
\hspace{5ex} &
\includegraphics[width=1.5cm, angle=0]{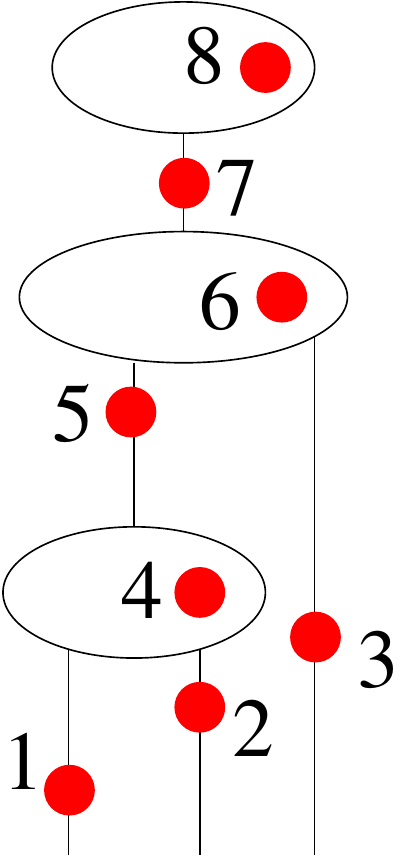}
\end{tabular}
\caption{From tropical curves to floor diagrams}
\label{fig:trop to floor}
\end{figure}

\medskip
Let us relate floor diagrams to enumeration of tropical curves. 
A collection $\omega$ of $3d-1+g$ points in $\RR^2 $
is said to be 
 {\it vertically stretched} if
 the absolute value of the difference between the second coordinates
 of any two of these points 
is much larger than the absolute value of the difference between the first coordinates of any two
of the points considered 
(in other words, 
all points of $\omega$ are in
a very narrow strip $[a,a+\varepsilon]\times\RR$). 
Fix 
a vertically stretched collection $\omega$ of $3 d - 1 + g$ points in $\RR^2$, and 
associate to the points of $\omega$ the numbers $1$, $\ldots$, $3d - 1 + g$ in such a way
that higher point always has a 
larger 
number. 
Given 
an irreducible nodal tropical curve of degree $d$ and genus $g$ in $\RR^2$ 
such that $C$ passes 
through the points of $\omega$, 
one can show that 
each vertical edge (an {\it elevator}) of $C$, as well as 
each connected component (a {\it floor}) of the complement in $C$ of the union 
of interiors of elevators,  
contains exactly one point of $\omega$
(see 
Figures \ref{fig:enum smooth}, \ref{fig:enum 30}, and 
\ref{fig:trop to floor}). 

Contracting each floor of $C$,  
we obtain a weighted graph whose edges correspond to vertical edges of $C$;
orient these edges in the direction of increasing of the second coordinate.
As it was shown by E.~Brugall\'e and 
Mikhalkin \cite{Br6b}, the result is a floor diagram 
of degree $d$ and genus $g$, and the $3d - 1 + d$ points in $\omega$ provide a marking
of this floor diagram (see Figure 
\ref{fig:trop to floor}). 
Conversely, any marked floor diagram of degree $d$ and genus $g$
corresponds to exactly one irreducible tropical curve of degree $d$ and genus $g$
 passing through the points of $\omega$.
This leads to the following statement (which is an immediate generalization
of a theorem proved in \cite{Br6b}).

\begin{theorem}[cf. \cite{Br6b}]\label{floor_diagrams_enumeration}
One has
$$
G_{d, g} = \sum_{\mathcal D} G({\mathcal D}), 
$$
where the sum is taken over all marked floor diagrams $\mathcal D$ of degree $d$ and genus $g$.
\end{theorem}

\begin{exa}
Combining Theorem \ref{floor_diagrams_enumeration} 
with the lists of Figures 
\ref{comp cplx1},\ref{degree 4 g=3,2}, \ref{degree 4 g=1}, and \ref{degree 4 g=0},
we obtain the values of $G_{d,g}$ given in Examples \ref{ex:enum
smooth}, \ref{ex:enum 30}, and \ref{ex:enum deg 4}.
\end{exa} 

Theorem \ref{floor_diagrams_enumeration} immediately implies the following formulas
for $N_{d, g}$ and $W^{trop}_{d, g}$:
$$
N_{d, g} = \sum_{\mathcal D} m_{\CC}({\mathcal D}), \;\;\; W^{trop}_{d, g} = \sum_{\mathcal D} m_{\RR}({\mathcal D}), 
$$ 
where the sums are again taken over all marked floor diagrams $\mathcal D$ of degree $d$ and genus $g$.

Beyond providing an efficient tool for explicit computations of
enumerative invariants, 
floor diagrams  also turned out to be a powerful technique in the study of
 (piecewise-)polynomial behaviour of 
Gromov-Witten 
invariants of complex surfaces, see for example
 \cite{FM,Blo11,ArdBlo,LiuOsse14,BA13}.

\begin{exo}
\

\begin{enumerate}
\item Fix a positive integer $d$, and put $g=\frac{(d-1)(d-2)}{2}-1$.
By adapting the 
approach used in Examples \ref{exa:rat cubic} and
  \ref{exa:rat real cubic}, prove that 
$$N_{d,g}=3(d-1)^2 ,$$
and that the quantity
$$\sum_C (-1)^{s(C)}$$
from Theorem \ref{real_enumeration} depends on $\omega$ 
for $d\ge 4$.

\item Fix a positive integer $d$. It is clear from the definition that the numbers $W_d$ and $N_d$ are equal
modulo $2$. Prove that 
$$W_d=N_d \ \mbox{mod}\ 4.$$ 

\item Show that  either  9 or 10 distinct rational tropical cubics pass
 through a given generic configuration of 8
  points in $\RR^2$.

\item With the help of Figures \ref{degree 4 g=3,2}, \ref{degree 4
  g=1}, and \ref{degree 4 g=0}, work out 
the computations of Example \ref{ex:enum deg 4}. 

\item Fix a positive integer $d$.
Using floor diagrams, prove that
$$G_{d,\frac{(d-1)(d-2)}{2}}=1\quad\mbox{and} \quad
  G_{d,\frac{(d-1)(d-2)}{2}-1}= (d-1)\cdot\left[\frac{d-2}{2}\cdot 
q^{-1}+2d-1+\frac{d-2}{2}\cdot q\right].$$

\item Fix a positive integer $d$ and an integer $0 \leq g \leq \frac{(d-1)(d-2)}{2}$. 
Prove that the highest power of $q$ which appears in $G_{d, g}$ is
$\frac{(d-1)(d-2)}{2}-g$, and that the
coefficient of the corresponding monomial of $G_{d, g}$ is equal to 
$$\left(\begin{array}{c}\frac{(d-1)(d-2)}{2} 
\\ g
\end{array}\right).$$ 

\end{enumerate}
\end{exo}

%% file: Kristin.tex
\section{Tropical subvarieties of $\RR^n$ and $\TT^n$}\label{sec:generalization} 
We focused so far  on  plane tropical curves. 
In this section we define tropical subvarieties of higher dimension and
codimension 
in $\R^n$. 
We have seen three equivalent definitions of a tropical curve in $\RR^2$:
\begin{enumerate}
\item an algebraic one via tropical polynomials;
\item a combinatorial one via balanced graphs;
\item a geometric one via limits of amoebas.
\end{enumerate}
All these three definitions can be generalized to arbitrary
dimensions.
In the case of  tropical hypersurfaces of $\RR^n$, all these three
definitions remain equivalent. Moreover, the proof that these yield equivalent definitions
is exactly the same as for tropical curves in $\R^2$.  
However,  these three
definitions produce different objects in higher codimension.

\subsection{Tropical hypersurfaces of $\R^n$}\label{sec:hyp}
As stated above, the situation for tropical hypersurfaces is entirely
similar to the case of tropical curves in $\RR^2$.
Let $n$ be a positive integer number. 
A tropical hypersurface in $\R^n$ is defined by a tropical polynomial
in $n$ variables
$$P(x_1, \dots  , x_n) = \tg\sum_{i \in A} a_{i}x^{i}\td =
 \max_{i \in A}  \{ a_{i} + \langle x, {i} \rangle \},$$ 
where $A \subset (\ZZ_{\geq 0})^n$ is a finite subset, 
$x^i=x_1^{i_1}x_2^{i_2}\ldots x_n^{i_n}$, and
$\langle x, {i} \rangle$ denotes the standard inner product 
in $\R^n$. 
As usual we are interested only in the function defined by $P$, and not in the 
polynomial expression.
Once again $P$ is a piecewise integer affine function on $\R^n$, and 
the definition  from Section \ref{sec:deftropcurves} of the set 
$\widetilde V(P)$ for tropical curves generalizes directly to higher dimensions: 
$$
\widetilde V(P)=\left\{ x\in\RR^n\ | \ 
\exists i \ne j ,\quad   P(x)= \tg a_{i}x^{i} \td
=\tg a_{j}x^{j} \td\right\}.
$$
The set  $\widetilde V(P)$ is a
finite union of convex polyhedral domains of
dimension $n-1$ forming a polyhedral complex 
(see Section \ref{sec:subvar}),
 and its facets
are equipped with a weight
in the same way as in  Definition \ref{def:trop curve}. 
The \emph{tropical hypersurface defined by $P$},
denoted by $V(P)$, is the set $\widetilde V(P)$ equipped with this
weight function on the facets.

\begin{figure}[h]
\centering
\begin{tabular}{ccc}
\includegraphics[scale=1]{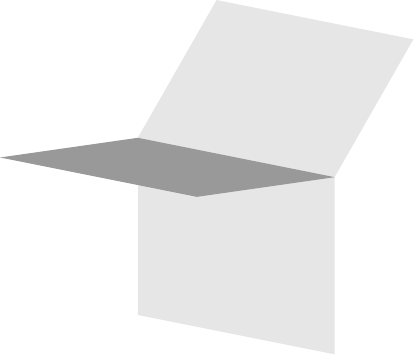}
&\hspace{10ex} &
\includegraphics[scale=0.3]{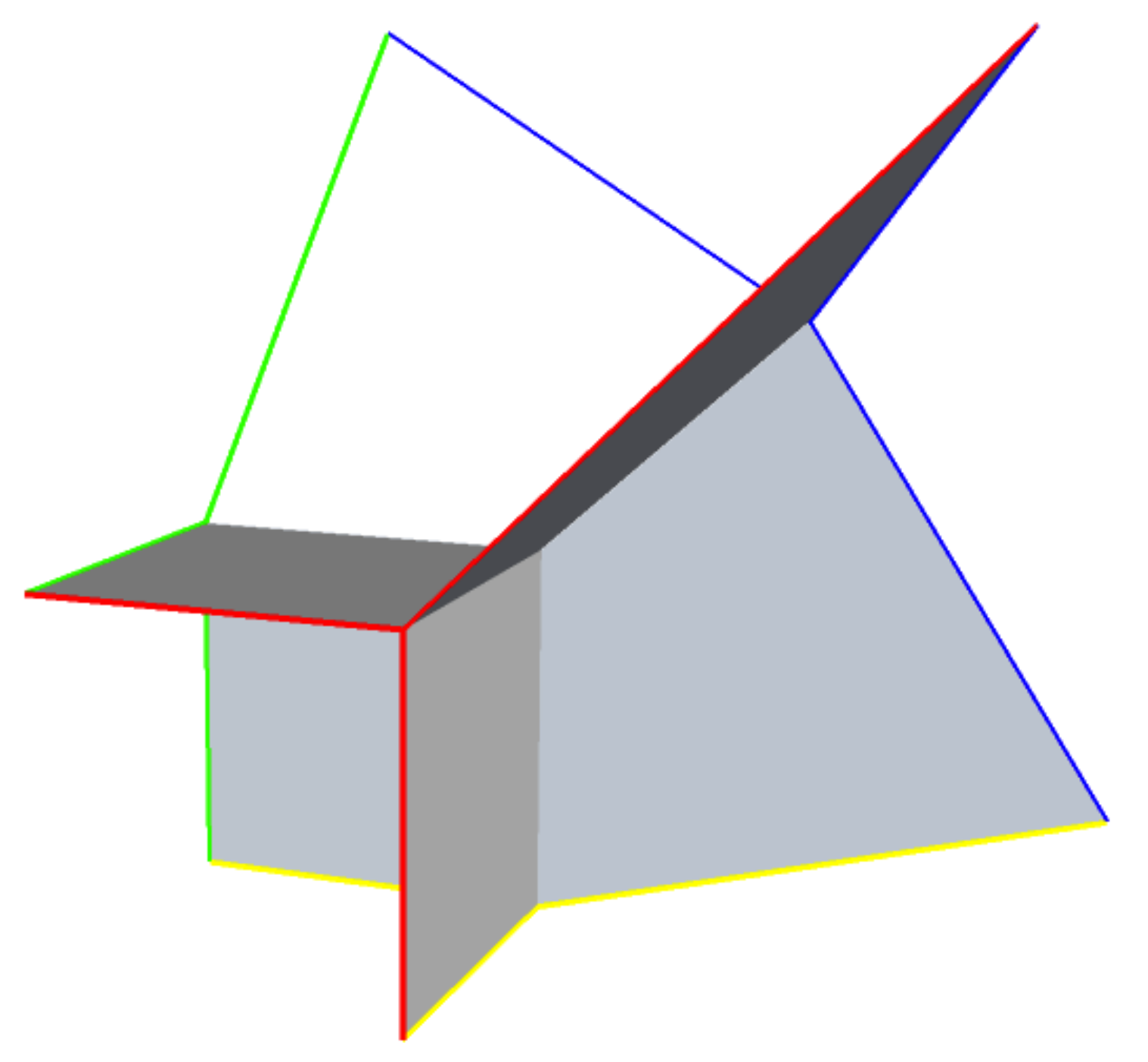}
\\ a) &&b)
\end{tabular}
\bigskip

\begin{tabular}{c}
\includegraphics[scale=0.02]{Figures/Quadricfloor.pdf}
\\ c)
\end{tabular}
\caption{Two
tropical planes in a) and b) and a  tropical quadric 
in $\T^3$ in c). }
\label{fig:surfaces}
\end{figure}

\begin{exa}\label{ex:trophyper degen}
The tropical hypersurface of $\RR^3$ defined by the tropical polynomial 
$$P (x,y,z) =``x+z+0"$$
is a cylinder in the $y$-direction.
It is 
formed
by the three facets
$$\begin{array}{c}
F_1: \ z=0 \ \mbox{and}\ x\le 0
\\F_2:  \ x=0 \ \mbox{and}\ z\le 0
\\ F_3:\ x=z \ \mbox{and}\ x\ge 0
\end{array}$$
that meet along the line $E$ with equation $x=z=0$, and equipped with the constant
weight function equal to 1 (see Figure \ref{fig:surfaces}a).
\end{exa}

\begin{example}\label{ex:trophyper}
Consider the linear tropical polynomial 
$$P (x_1,\ldots,x_n) =``a_1x_1 + \dots + a_nx_n + a_{0}".$$
The 
tropical hyperplane   
$V(P)$
is a fan of dimension $n-1$ equipped with weight one on all of its
facets. 
Any such facet has 
the
point
$(a_0-a_1,\ldots, a_0-a_n)$ as vertex,
and is generated by $n-1$ of the vectors
$$(-1,0,\ldots,0), \ (0,-1,0,\ldots,0),\ldots, \  (0,\ldots,
0,-1),\ (1,\ldots, 1). $$
Conversely, any $n-1$ elements of this set of $n+1$ vectors define a
facet of $V(P)$.
If $n=2$,
we have again 
a tropical line in the plane 
(cf. Example \ref{ref:line}). 
A tropical plane in $\RR^3$ is 
depicted 
in
Figure \ref{fig:surfaces}b. 
Such a tropical plane has $4$ rays, in the directions 
$$(-1,0,0), \quad (0,-1,0), \quad (0,0,-1), \quad (1,1,1),$$
and $6$ top dimensional faces, one spanned by each pair of 
rays. 
\end{example}

\begin{exa}
A tropical quadric surface in $\RR^3$ is depicted in Figure \ref{fig:surfaces}c. 
\end{exa}

Just as 
described in Section
\ref{ref:dualSub}
for tropical curves in $\RR^2$, any tropical 
polynomial $P$ induces a  subdivision of its
Newton polygon 
$$\Delta(P)=Conv\left(\left\{i\in(\ZZ_{\geq 0})^n \ | \ a_i\ne-\infty
\right\}\right)\subset \RR^n.$$ 
The tropical hypersurface 
$V(P)$ is
dual to this subdivision
in the sense of Proposition \ref{prop:dual subd}.
By this duality, a  top
dimensional face $F$  of $ V(P)$ is dual to 
an edge $e$
of the subdivision of 
$\Delta(P)$,
and  the weight of $F$ is equal to the lattice length of 
$e$. As a result,
the tropical hypersurface satisfies a generalization of the 
balancing condition from Section \ref{sec:balancing} along 
the faces of 
dimension $n-2$.
Let $E$ be a face of dimension $n-2$  of  $V(P)$, and let 
$F_1, \dots, F_k$ be the faces of dimensions $n-1$ 
of $V(P)$ that are adjacent
to $E$. 
Denote by 
$w_1, \dots , w_k$
the weights of $F_1, \dots, F_k$.
Let $v_i$, $i = 1$, $\ldots$, $k$, 
be the primitive integer vector orthogonal to $E$ and such that
$x+\epsilon v_i\in F_i$ if $x\in E$ and $1>>\epsilon >0$. 

\begin{prop}[Balancing condition for tropical hypersurfaces]
One has
$$
\sum_{i=1}^k w_i v_i=0. 
$$
\end{prop}

The converse holds, {\it i.e.} Theorem \ref{prop:trop balanced}
generalizes to  tropical hypersurfaces of $\RR^n$.

\begin{exa}
For the tropical surfaces of Example \ref{ex:trophyper degen} we have
$$v_1=(-1,0,0), \quad v_2=(0,0,-1), \quad v_3=(1, 0,1), $$
and $v_1+v_2+v_3=0. $
\end{exa}

\medskip
As in the case of curves, any tropical hypersurface of $\RR^n$ is the
limit of amoebas of algebraic hypersurfaces of $(\CC^\times)^n$.
The next theorem generalizes Theorem \ref{approx}.

\begin{thm}\label{approx any dim} 
{\rm (}cf. \cite{Kap1,Mik12}{\rm )} 
Let $P_t(z_1,\ldots,z_n)=\sum_{i}\alpha_{i}(t)z^i$ be a polynomial whose
coefficients are functions $\alpha_{i}:\RR \to \CC$, and suppose
that $\alpha_{i}(t)\sim \gamma_{i}t^{a_{i}}$ when
$t$ goes to $+\infty$ with $\gamma_{i}\in\CC^\times$ and
$a_{i}\in\TT$. 
If $\V_t$ denotes 
the hypersurface in $(\CC^\times)^n$ defined by the polynomial $P_t(z)$, 
then  
the amoeba $\text{Log}_t(\V_t)$ converges to the tropical hypersurface
of $\RR^n$ defined
by the tropical polynomial $P_{trop}(x)=\tg \sum_{i}a_{i}x^i \td$.
\end{thm}

\medskip

A tropical hypersurface is said to be 
\emph{non-singular} if any top dimensional cell of its
 dual subdivision has Euclidean volume $\frac{1}{n!}$, 
 as the standard simplex in $\RR^n$. Equivalently, a tropical
 hypersurface is 
 non-singular 
 if each vertex has a  
 neighborhood which is the same, up to
 translations and 
 the action of $GL_n(\Z)$, as
 a neighborhood of the  tropical hyperplane from Example \ref{ex:trophyper}.
More generally, a non-singular tropical subvariety of $\RR^n$ 
locally looks like a tropical linear space. 
We treat these objects in Section \ref{sec:fanlinear}.

\subsection{Tropical subvarieties of $\RR^n$}\label{sec:subvar}

A $k$-dimensional \emph{finite polyhedral complex}  in $\RR^n$ is the union of 
 finitely many $k$-dimensional
convex polyhedral domains $F_1,\ldots,F_s\subset \RR^n$ (each $F_j$
is the intersection of some half-spaces in $\RR^n$) such that the
intersection of two polyhedra $F_i$ and $F_j$ is either empty or
 a face of both $F_i$ and $F_j$.
We define tropical subvarieties of $\RR^n$ by generalizing the above
balancing condition to 
finite
polyhedral complexes of any dimension in $\RR^n$.
A finite polyhedral complex $V\subset\RR^n$ is said to be 
\begin{itemize}
\item 
\emph{rational}
if 
for each face $F\subset V$
the integer vectors tangent to $F$ 
form a lattice of rank equal to $\dim F$; 
\item \emph{weighted} if each of 
the top dimensional
faces of $V$ is equipped with a \emph{weight} that is an
integer. 
\end{itemize}

\begin{defi}[General balancing condition in $\RR^n$]\label{def:bal}
Let $V$ be a 
weighted 
rational finite polyhedral complex 
of 
dimension $k$ in
$\R^n$, and let $E \subset V$ be a codimension one face of
$V$. Let  $F_1,
\dots, F_s$ be the facets adjacent to $E$,
and let $\Lambda_{F_i} \subset \ZZ^n$ denote the lattice parallel to $F_i$, (analogously for $\Lambda_E$). 
Let $v_i$ be a primitive
integer vector such that, together $v_i$ and $\Lambda_E$ generate $\Lambda_{F_i}$, 
and 
 for any $x \in E$, 
one has
$x + \epsilon v_i \in F_i$
for $1>>\epsilon >0$.
We say that $V$ is \emph{balanced} along $E$
if  the vector
$$ \sum_{i = 1}^s w_{F_i}v_i$$
is in $\Lambda_E$, 
 where $w_{F_i}$ is the weight of the facet $F_i$. 
\end{defi}

Using this we introduce the definition of a tropical subvariety of $\R^n$. 

\begin{defi}
\label{def:tropvar}
A \emph{tropical cycle} $V$ of $\R^n$ of dimension $k$ is a 
$k$-dimensional 
weighted rational finite polyhedral complex 
which
is balanced along every
codimension 
one face. 

A tropical cycle  of $\R^n$ equipped with non-negative integer
weights is called 
\emph{effective}, 
or also a 
 \emph{tropical subvariety} of $\RR^n$.
\end{defi}

\begin{exa}
A tropical curve in $\RR^2$ is precisely a balanced
graph as defined in Section \ref{sec:balancing}.
Note that this definition from Section
\ref{sec:balancing}
generalizes immediately to  any weighted
rectilinear graph in $\RR^n$ with rational slopes, and this
generalization is again
equivalent to Definition \ref{def:tropvar}.
\end{exa} 

In the case of tropical subvarieties of $\RR^n$ of dimension $n-1$,
Definition \ref{def:tropvar} is equivalent to the definition of tropical
hypersurfaces of $\RR^n$ we gave in Section \ref{sec:hyp}.
\begin{exa}
Consider again  the tropical surface in $\RR^3$ from Example \ref{ex:trophyper degen}.
This is indeed a tropical surface in the sense of Definition \ref{def:bal}, since we have 
$$v_1=(-1,a,0), \quad v_2=(0,b,-1), \quad v_3=(1, c,1), $$
and so
$$v_1+v_2+v_3=(0,a+b+c,0) $$
which is 
parallel to $E$.
\end{exa}

\newcommand{\ignore}[1]{\relax}
\ignore{
\subsection{Tropical limits of  algebraic varieties}
As in the case of 
hypersurfaces,
families
of algebraic subvarieties of $(\CC^\times)^n$ are related to tropical
subvarieties of $\RR^n$. 
We briefly indicate below this relation
and refer to \cite{IKMZ} for
more details. 

Let $(\V_t)_{t\in\RR_{>1}}$ be a 
family of complex algebraic subvarieties  of $(\CC^\times)^n$ of dimension
$k$, and
of bounded degree in some
compactification of $(\CC^\times)^n$ to $\CC P^n$. Assume furthermore that
the amoebas $\Log_t(\V_t)$ converge to a set $V$, with respect to the Hausdorff metric on
compact sets of $\RR^n$. It turns out that $V$
admits the structure of
a rational finite polyhedral complex of  dimension $k$, 
cf. \cite{BiGr} for a relevant statement for field valuations.

\begin{exa}
An integer matrix $M\in\mathcal M_{n,m}(\ZZ)$
 defines a multiplicative map
$\Phi_{M}:(\CC^\times)^m\to (\CC^\times)^n$. 
 In this case the set $V \subset \RR^n$  for the
 constant family $\V = \Phi_{M}((\CC^\times)^m)$ 
is
 the vector subspace
of $\RR^n$ spanned by the columns of $M$.
\end{exa}

 Given a facet $F$ 
 of $V$, we associate  to it a positive integer,  $w(F)$, known as 
 the \emph{weight} of $F$, as 
follows: pick a point $(p_1,\ldots,p_n)$ 
in the relative interior of $F$, and choose
a basis $(e_1,\ldots,e_n)$ of $\ZZ^n\subset\RR^n$ such that
$(e_1,\ldots,e_m)$ is a basis of the subspace parallel to  $F$. Denote by
$\mathcal Y_{F,t}\subset(\CC^\times)^n$ 
the 
multiplicative translation of $\Phi_{(e_{m+1},\ldots,e_{n})}((\CC^\times)^{n-m})$ by
$(t^{p_1},\ldots,t^{p_n})$, and define
$w(F)$ as the number (counted with multiplicity) of intersection
points of $\V_t$ and $\mathcal Y_{F,t}$ whose amoeba in base $t$
converges to $(p_1,\ldots,p_n)$, 
when $t\to +\infty$.
It turns out
that $w(F)$ does not depend on the 
choices made. 

\begin{exa}\label{ex:matrix}
A matrix $M\in\mathcal M_{n,m}(\ZZ)$ with $Ker\ M=\{0\}$ maps the
lattice $\ZZ^m$ to a sub-lattice $\Lambda'$ of $\Lambda=\ZZ^n\cap Im \ M$.
The weight 
of 
$\Phi_{M}((\CC^\times)^m)$ is the index of $\Lambda'$ in $\Lambda$.
\end{exa}

\begin{defi}\label{def:trop limit}
The \emph{tropical limit} of the family $(\V_t)_{t\in\RR_{>1}}$,
denoted by $Trop(\V_t)$, is the
finite polyhedral complex $V$ equipped with the weight function 
described above
on
its facets.
\end{defi}

The following theorem generalizes Theorems \ref{approx} and
\ref{approx any dim}.
\begin{thm}[\cite{IKMZ}, cf also \cite{BiGr,Spe1} in the
    non-Archimedean setting]
The tropical limit $\lt(\V_t)$ is a tropical subvariety of $\RR^n$.
\end{thm}

\begin{defi}\label{def:approx}
A tropical subvariety $V$ of $\R^n$ is said to be \emph{approximable} if there
exists a family of complex algebraic subvarieties $\V_t$ of $(\CC^\times)^n$ such%
that $\lt(\V_t)=V$.  
\end{defi}

 \begin{rem}\label{rem:strong approx}
As we work here with subvarieties of 
$(\CC^\times)^n$ (which is non-compact), 
there are several version of
tropical approximability. 
The definition above refers to tropical approximability in the {\em
  weak sense}, and we may define a \emph{stronger} notion. 
Recall that the subvariety $\V_t$ defines a class  $[\V_t]$ in 
$\tilde{H}_{2k}((\CC^\times)^n)$, the
homology group of the universal toric compactification of
$(\CC^\times)^n$. 
A tropical subvariety $V$ of $\R^n$
of dimension $k$ also defines a class $[V] \in
\tilde{H}_{2k}((\CC^\times)^n)$. Indeed,  it is the balancing condition
which ensures that the chain coming from $V$ is closed. The recipe for
this class is given in \cite{FS} when $V$ is a weighted balanced fan.

A stronger definition of approximability 
 includes the additional condition that the degree of the approximating
complex subvarieties is the same as the degree of the tropical limit.
Definition \ref{def:approx} is weaker since it
 does not include this condition on homology classes.
As a result the degree of the approximating complex subvarieties might be higher, but the excess
degree escapes to infinity. 
\end{rem}

\begin{exa}
We illustrate the difference between the strong and weak notions of
approximation of tropical varieties in the simple 
case of a tropical hypersurface $V$ of $\RR^n$. 
The
strong approximation mentioned in
Remark \ref{rem:strong approx} requires the existence
of a  family $(\V_t)_{t\in\R_{>1}}$ such
that $\lt(\V_t)=V$ and $\Delta(V)=\Delta(\V_t)$, whereas the
weak notion from Definition \ref{def:approx} only implies that
$\Delta(V)\subset \Delta(\V_t)$.

For example, the tropical curve in $\RR^2$ defined by the tropical
polynomial $\tg 0+x\td $ is the tropical limit of
the family $(\C_t)_{t\in\R_{>1}}$ of curves in $(\CC^\times)^2$ with
equation
$$1+z+ t^{-t}w=0 $$
in the sense of Definition \ref{def:approx}, but not in the strong sense of Remark 
\ref{rem:strong approx}.
\end{exa}

\begin{rem}
Specifying only the weights of the facets in the tropical limit as in
Definition \ref{def:trop limit} may also be refined
to specifying a finite covering of degree $w(F)$ for the $k$-torus $(S^1)^k$ corresponding to a facet $F$.
Note that the sublattice $\Lambda'\subset\Lambda$ defines such a
covering in Example \ref{ex:matrix}. 
This refinement may subdivide a facet of the tropical limit to smaller subfacets with different refinements.
\end{rem}

Not every tropical subvariety of $\RR^n$ is approximable, 
as illustrated by the 
following example.
Determining which tropical subvarieties 
of codimension different from one are approximable is quite
difficult, even in the case of curves.

\begin{exa}\label{exo:spatialcubic}
Consider
the
plane tropical cubic $C$ of genus 1 depicted in Figure \ref{equil}a, and
draw it in the affine plane with equation $z=0$ in $\RR^3$. On each of the three unbounded edges in 
the direction $(1, 1, 0)$, choose a point in such a way that these three
points are not contained in a tropical line in $z=0$. Now at  these  three points,
replace the unbounded part of $C$ in the direction $(1,1,0)$ by two
unbounded edges, one in the direction $(0,0,-1)$ and one in the
direction $(1,1,1)$ to obtain a spatial 
tropical cubic $\tilde{C}$ (see Figure \ref{fig:spatial cubic}).

\begin{figure}[h]
\begin{center}
\begin{tabular}{c}
\includegraphics[width=7cm, angle=0]{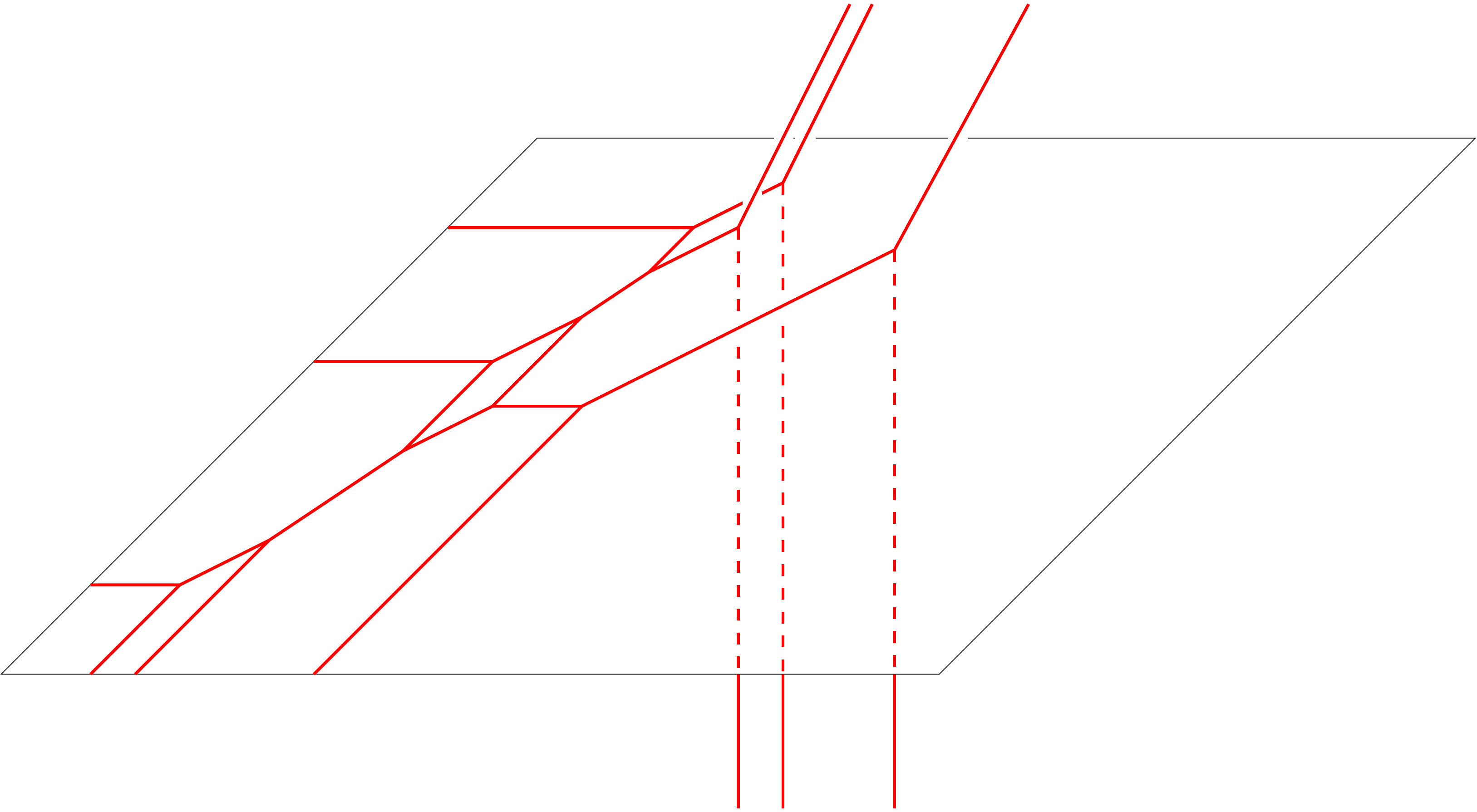}
\end{tabular}
\end{center}
\caption{A non-approximable tropical spatial cubic curve.}
\label{fig:spatial cubic}
\end{figure}
The choice we made of the three points ensure that 
$\tilde{C}$
is not contained in any tropical hyperplane; since as
we 
saw in Example \ref{ex:trophyper}, the unbounded rays in the
$-e_3$ direction contained in 
a tropical 
plane in $\RR^3$ 
project to a tropical line on the plane $z=0$.
 For the sake of simplicity, we show that $\tilde{C}$ is not
strongly  approximable in the sense of Remark \ref{rem:strong approx},
{\it i.e.} there does not exist a family $(\C_t)_{t\in\R_{>1}}$ of  spatial
elliptic 
cubics with $\lt(\C_t)=\tilde{C}$. 
The tropical curve $\tilde{C}$ is not approximable in the sense of
Definition \ref{def:approx} either 
({\it i.e.} even if the curves $\C_t$ are of  degree greater than $3$, they 
cannot approximate $\tilde{C}$), 
however the proof is more involved.

Any  spatial algebraic elliptic cubic curve over $\CC$
is planar.
Hence if 
an appropriate family  $(\C_t)_{t\in\R_{>1}}$  of elliptic cubics exists, 
each member  $\C_t$ would be necessarily contained in 
 a plane $\P_t$.
 By passing to a subsequence $\{t_s\}$ if necessary, one can prove that
the tropical limit $\lt(\P_{t_s})$ exists. 
 Thus
 $\lt(\C_t)$
would be contained in 
$\lt(\P_{t_s})$, which is a tropical plane in $\RR^3$.

 This proves that the curve $\tilde{C} \subset \RR^3$ 
is not
 approximable in the sense of 
Remark \ref{rem:strong approx}.
\end{exa}

We refer to \cite{Spe2,Br12} for generalizations of Example \ref{exo:spatialcubic}.
}

\input{section53.tex}

\subsection{Linear spaces and hyperplane arrangements}\label{sec:fanlinear}
Now we turn to tropical limits of constant families of linear spaces,
and more generally to \emph{fan tropical linear spaces}. 
These objects will make a reappearance in Section \ref{sec:manifold} where they will play the role of the local models for tropical manifolds. 
Regarding more general tropical linear spaces which are not
necessarily fans, we refer to \cite{Spe3}.

 Let us start with an affine linear  space $\tilde{\L} \subset \CC^n$ of dimension $k$, 
 and assume
$\L  = \tilde{\L} \cap (\CC^\times)^n\subset (\CC^\times)^n$ is non-empty.
We obtain a natural hyperplane
 arrangement $\A$ in $\CC P^k$ by compactifying 
$\L$
to $\bar{\L} \cong \CC P^k \subset \CC  P^n$. This arrangement 
  consists  of 
the intersection of $\bar{\L}$ with all coordinates hyperplanes in $\CC P^n$. 
  Here we  
give a description of 
the tropical variety $\lt \L$.
In this case, it turns out that the weights 
of the top
dimensional faces 
of 
$\lt \L$ are all 
equal to $1$.

We already mentioned
that given 
a 
complex curve $\C$, the tropical limit 
$\lt \C$
has only 
one vertex and directions corresponding to the asymptotic directions
of the amoeba.  The same 
is true
for $\lt \L$. 
It can be equipped with the structure of a finite polyhedral
fan, and the faces record the asymptotic directions.  
 This type of limit for any $\V \subset (\CC^\times)^n$ is also known
 as the
Bergman fan of the variety. These objects were considered by G.M.~Bergman before 
the birth of tropical geometry \cite{Berg71}.

\begin{prop}\cite{SturmPoly} 
Let $\L \subset (\CC^\times)^n$ be a linear space,  then $\lt \L$
depends only on the intersection properties of the hyperplanes in
$\mathcal{A}$.  
\end{prop}

\begin{figure}[h]
\includegraphics[scale = 0.7]{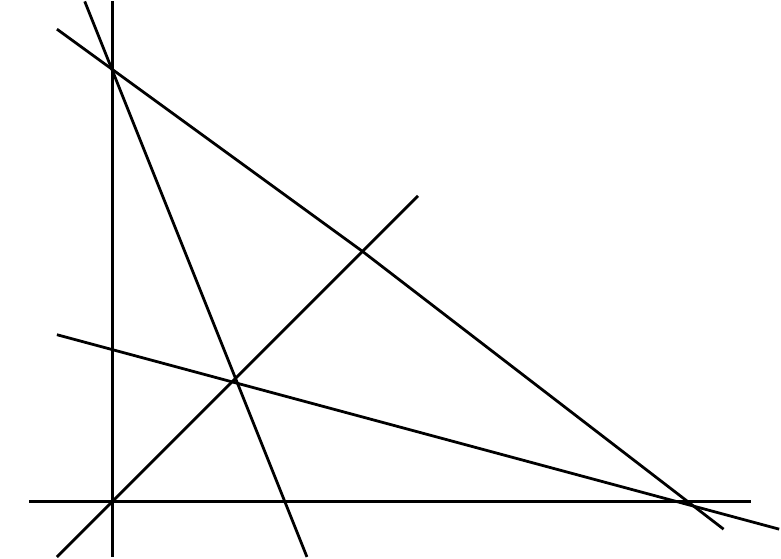} 
\put(-113, 55){$4$}
\put(-60, 50){$0$}
\put(-60, 0){$2$}
\put(-145, 65){$1$}
\put(-155, 45){$5$}
\put(-70, 75){$3$}
\hspace{1.5cm}
\includegraphics[scale = 0.2]{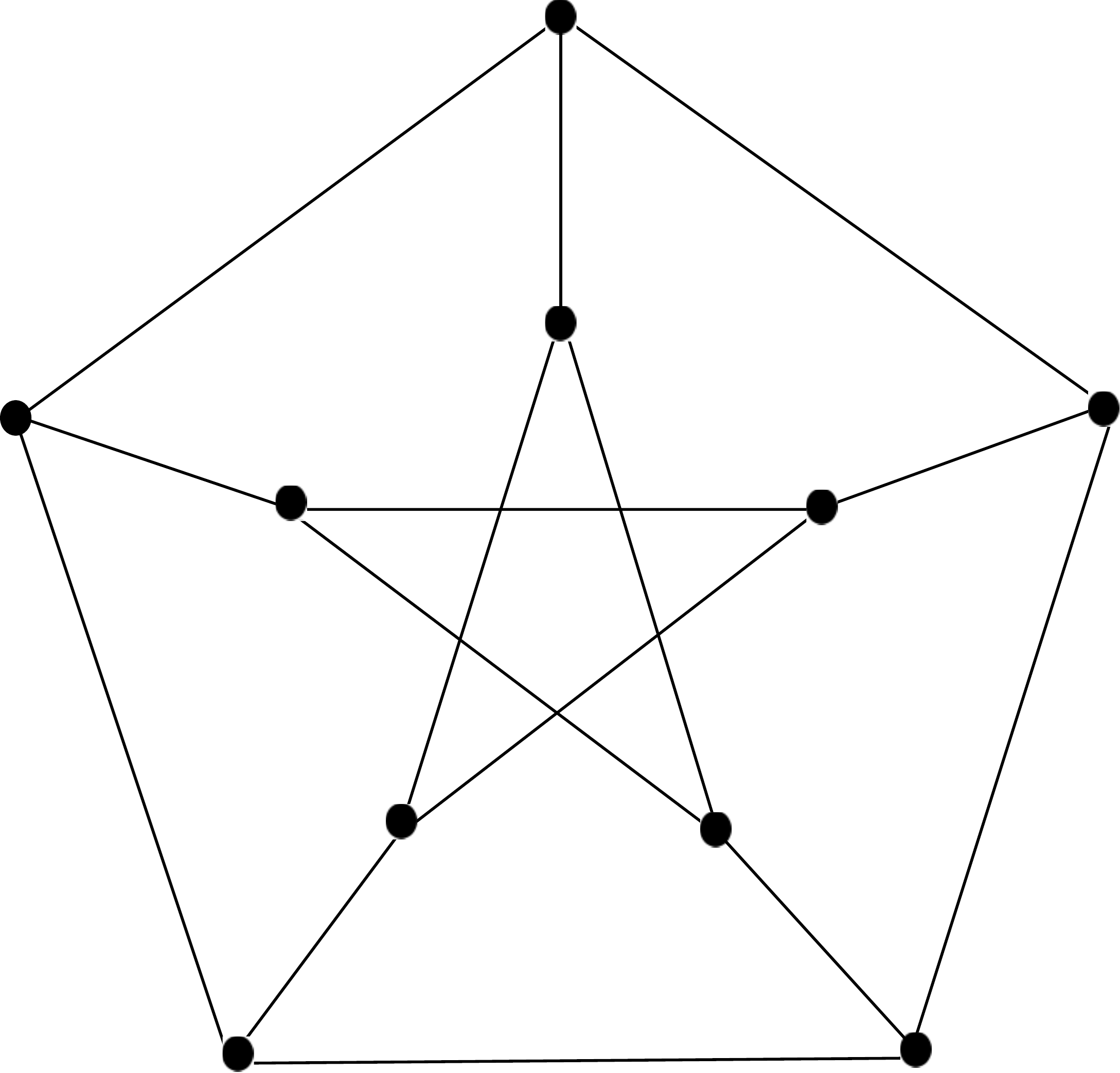} 
\put(-55, 85){$4$}
\put(5, 70){$0$}
\put(-15, 0){$2$}
\put(-40, 30){$3$}
\put(-135, 70){$1$}
\put(-90, 70){$5$}
\caption{The braid arrangement of lines and the link of the vertex of its  tropical limit.}
\label{fig:M_05}
\end{figure}
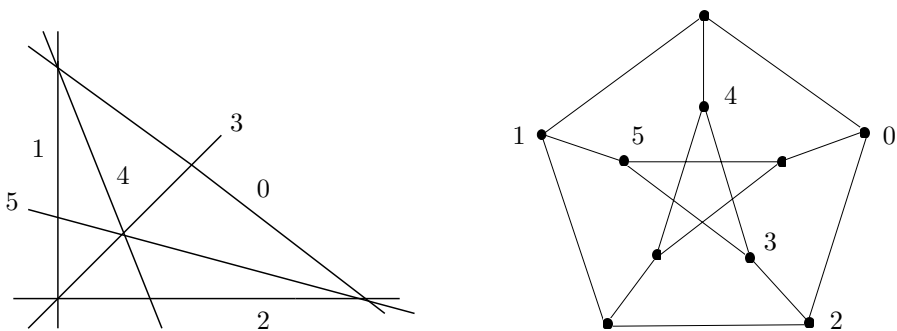

Let us explain what do we mean by ``intersection properties". A hyperplane arrangement $\A = \{H_0, \dots , H_n\}$ in $\CC P^k$
is a  stratified space, here we  
refer to the strata as \emph{flats}. 
Each flat can be indexed by the maximum subset $I \subset \A$ of hyperplanes which contains it; label such a flat $F_I$. See Figure \ref{fig:M_05} for a line arrangement where some of the flats are labeled.  
The flats form a  partially ordered set, the  order being  given  by inclusion and is known as the \emph{lattice of flats} of the arrangement $\A$. 
 Now a more precise statement of the above proposition is that for a linear space, the set $\lt \L$ can be 
determined 
from the lattice of flats of $\mathcal{A}$ and does not depend on the position of the hyperplanes.

Now we 
show how to construct 
the fan $\lt \L$ explicitly.

\begin{construction}[\cite{Ard}]\label{construction:matroidfan}
Set $v_i = -e_i$ for $i = 1, \dots , n$ and $e_0 = \sum_{i=1}^n e_i$, where $e_i$ are the standard basis 
of $\RR^n$. For any $I \subset \{0, \dots , n\}$ let $v_I = \sum_{i\in I} v_i$. 
A chain of flats is a collection of flats satisfying  $$F_{I_1} \subset F_{I_2} \subset  \dots \subset F_{I_{l-1}} \subset F_{I_l},$$ 
such that $\dim(F_{I_i})  = \dim(F_{I_{i-1}}) - 1$. For every chain of flats, there is a 
cone in $\lt \L$ spanned by the vectors $ v_{I_1}, \dots , v_{I_l}.$
Then, $\lt \L$ is the underlying set which is the union of all such cones, 
every top dimensional face being equipped with weight $1$. 
\end{construction} 

\begin{exa}\label{ex:uniformlin}
A linear space $\L \subset (\CC^\times)^n$ is called 
$\partial$-transversal if 
the corresponding hyperplane arrangement satisfies $\codim ( \cap_{i  \in I} H_i  ) = |I|$ for any subset $I \subset \{0, \dots , n\}$. 
Such a hyperplane arrangement is also known as \emph{uniform}. 
The condition for a linear space to be 
$\partial$-transversal implies that any subset $I \subset \{0, \dots , n\}$ gives a flat $F_I$, 
just as any chain of subsets of $\{0, \dots , n\}$ gives a cone in $\lt \L$.

To understand $\lt \L$ as a set, notice that there are $n+1$
one-dimensional rays in directions $v_0, v_1, \dots , v_n$. For every 
$k$-tuple $\{i_1, \dots , i_k\} \subset \{0, \dots , n\}$,
there is a cone of dimension $k$ in $\lt \L$ spanned by the vectors $v_{i_1}, \dots , v_{i_k}$.  Any other cone from the construction above simply subdivides one of the above cones, thus does not add anything to the set $\lt \L$. 

In particular,   
the tropical limit of $\partial$-transversal  
hyperplane is exactly
the standard hyperplane from Example \ref{ex:trophyper}. The tropical limit of a 
$\partial$-transversal 
linear space of dimension $k$ in 
$(\CC^\times)^n$ 
is the
$k$-skeleton of the standard  tropical hyperplane, meaning it consists of all
faces of dimension less than or equal to $k$
of the standard  tropical hyperplane.

The special case $k=1$ is particularly easy to describe: a 
$\partial$-transversal fan tropical
line in $\RR^n$ is made of $n+1$ rays emanating from the origin
and going to infinity in the directions
$$(-1,0,\ldots,0), \ (0,-1,0,\ldots,0),\ldots, \  (0,\ldots,
0,-1),\ \mbox{and}\ (1,\ldots, 1). $$
\end{exa}

\begin{exa}\label{ex:M05}
Take the arrangement of six hyperplanes defined in $\CC P^2$ by 
$$
\A = \{x=0, y=0, z=0, x=y, x=z, y=z\}.
$$ 
This is the well known braid arrangement in dimension $2$. 
The line arrangement is shown on the left of Figure
\ref{fig:M_05}. The complement $\CC P^2 \backslash \A$ can be
identified with the moduli space $\mathcal{M}_{0,5}$ of $5$-marked
rational curves up to automorphism. Moreover, the complement has a
linear embedding to $(\CC^\times)^5$.  As for the flats of this
arrangement, there are the lines themselves, 
there are
$4$
points which are the intersection of $3$ lines, and $3$ points which
are the intersection of just $2$ lines.  
The link of the singularity of the tropical limit  
({\it i.e.} the intersection of $\lt \L$ with a small sphere centered at
the vertex of $\lt \L$)
is the Petersen graph drawn on the right-hand side of Figure
\ref{fig:M_05}.
This tropical space is also the moduli space of $5$-marked tropical
rational curves
 (see \cite{Mik7}).
\end{exa}

There is a combinatorial object, introduced by H. Whitney, which seeks to
capture and generalize the notion of independence coming from linear
algebra and also graph theory. This object is 
called a
\emph{matroid}. 
We refer,
for example, to \cite{Ox} for an introduction to the rich theory of matroids. 

 The lattice of flats of a hyperplane  arrangement  in $\CC P^k$ encodes a matroid. Yet there exists many matroids which are not realizable by the lattice of flats of a hyperplane arrangement over  a fixed field, and some which are not realizable over any field at all. As examples, the Fano plane arrangement is only realizable in characteristic two, and the ``non-Pappus" matroid (which violates Pappus' hexagon theorem) is not realizable over any field.

The construction \ref{construction:matroidfan}
above utilizes only the information from the
intersection lattice of the hyperplane arrangement. Any matroid comes
with this data, thus we can construct a finite polyhedral fan as above and
moreover the resulting space is a  tropical subvariety of $\RR^n$ \cite{SturmPoly}. So in a sense, any matroid has a geometric representation as a tropical object.

\begin{definition}\label{def-lfan}
A fan tropical linear space $L \subset \RR^n$ is a
 tropical 
subvariety of $\RR^n$
obtained from the Construction \ref{construction:matroidfan} applied
to 
a matroid.  
\end{definition}

Fan tropical linear spaces which come from a matroid 
non-realizable 
over $\CC$  provide a
class of examples  of tropical subvarieties of $\RR^n$ 
 which are not uniformly approximable
by linear spaces in $(\CC^{\times})^n$.

\subsection{Tropical subvarieties of $\TT^n$}
The tropical affine space
$\TT^n$ 
is naturally a stratified space.  
For any $I \subset [n]:= \{0, \dots , n\}$ let 
$$
\T^{J}:= \{x \in \T^n \ | \ x_i = -\infty \text{ for all } i \notin J \}.
$$
Then $\partial \T^n$ is the union of 
$\T^{[n] \backslash i}$
for all $1 \leq i \leq n$, and 
$\T^{[n] \backslash i}$
is of codimension one 
in 
$\T^n$.
Notice that $\partial \T^n$ looks like $n$  simple normal crossing divisors. 
The \emph{sedentarity} $I(x)$ of a point $x \in \T^n$ is the set of
coordinates of $x$ which are equal to $-\infty$. 
For  $J = [n] \backslash I$ 
denote by $\RR^J \subset \TT^J$
the subset of points of sedentarity 
$I$. 

The coordinate-wise logarithm map $(\CC^{\ast})^n \to \R^n$ has fibers $(S^1)^n$. 
We extend this map to $\T^n$, by setting $\log(0) = -\infty$. 
Notice that different strata have different fibers. 
Over a point in the interior of the  boundary strata  
$\T^J$ the fiber is 
 $(S^1)^{|J|}$.
Over the tropical origin $(-\infty, \dots , -\infty)$
there is only the point $(0, \dots , 0)$. 

We can also extend our definitions of tropical subvarieties and
tropical limits  to $\TT^n$ and $\CC^n$.

\begin{defi}
A {\it tropical subvariety $V$ of $\T^n$ of sedentarity $I$} is the closure
in $\T^n$ of a tropical subvariety $V^o$ of 
$\R^{J}$ for  $J = [n] \backslash I$. 
A {\it tropical subvariety 
of $ \T^n$} is a union of tropical subvarieties
of $\T^n$  of possibly different sedentarities.  
\end{defi} 

As algebraic curves in $\CC^n$, 
 tropical curves in $\TT^n$ also have a degree. 
\begin{defi}\label{def:degree}
The \emph{degree} of a
 tropical curve  $C\subset \TT^n$ 
 is defined as
$$
\text{deg}(C)=\sum\limits_{e}w_{e}\ {\max}_{j=1}^n\{0,s_j(e)\},
$$
where the sum is taken over all unbounded edges of $C$, 
the vector $(s_1(e),\ldots, s_n(e))$ is the primitive integer vector of $e$
 in the
outgoing direction, and $w_e$ is the weight of the edge $e$. 
(By convention $s_j(e)=0$ if $j\in I$, where $I$ is the sedentarity of $C$.)
\end{defi}

The tropical and classical notions of degree are related by the
following proposition.

\begin{prop}\label{lem:degree}
Let $\C \subset \CC^N$ be a complex algebraic curve, and 
put $C = \lt \C$.  
Then $\deg(\overline{\C}) = \deg(C)$, where
$\overline{\C}$ is the closure of $\C$ in the compactification 
of $\CC^N$ to $\CC P^N$. 
\end{prop}

\subsection{Tropical modifications}\label{sec:mod}
From all the examples and pictures of tropical subvarieties of $\RR^n$
we provided
so far, one can observe a particular feature of tropical geometry:
different tropical limits of the same classical variety can
have different topologies depending on
the embedding of
the classical variety. 
There is a way to understand how to relate these different tropical
models, it is called \emph{tropical modification} and was introduced
in  \cite{Mik3}.  

Given a tropical subvariety $V$ of $\T^n$ and a tropical polynomial
$P:\T^n \to \T$, a 
tropical modification is a map $\pi_P: \tilde{V} \to V$, where
$\tilde{V}$  is a 
tropical subvariety of $\T^{n+1}$ defined below 
 and $\pi_P$ is simply the linear projection of
$\T^{n+1} \to \T^n$ with kernel $e_{n+1}$.  
The tropical modification is a space $\tilde{V}$ along with the map
$\pi_P$, so our  notation is similar to that for birational
modifications in classical geometry.  

Now we describe precisely how to obtain the space
$\tilde{V}$. Firstly, consider the graph $\Gamma_P(V)$ of the
piecewise integer affine function $P$ restricted to $V$.
The graph $\Gamma_P(V)$ is 
equipped with the
weight function inherited from the weight function on $V$. 
In general, 
this graph is not a tropical subvariety 
of $\T^{n+1}$
as it may not be balanced in the $e_{n+1}$ direction. Recall our first considerations of graphs of tropical polynomials:
the graphs shown in Figure \ref{graphes} 
are 
not tropical curves.  At every codimension 
one face of $\Gamma_P(V)$ which fails to satisfy the balancing 
condition, add a new top dimensional face in the $-e_{n+1}$ direction
equipped with the unique integer weight so that balancing condition is
now satisfied. This is what the dashed vertical line segments
represent in Figure \ref{graphes}.  

\begin{defi}\label{def:trop mod}
Let $V$ be a tropical subvariety of 
$\T^n$ 
with 
sedentarity $\emptyset$,
and $P: \T^n \to \T$ a
tropical polynomial. 
The map $\pi_P: \tilde{V} \to V$ described above  is the \emph{tropical modification} 
 of $V$ along the function $P$.

The \emph{divisor} $\text{div}_{V}(P) \subset V$ of $P$ restricted to $V$ is the
union of the set of the points $x \in V$ such that 
$P(x) = -\infty$ and of the projection of the corner locus of the graph of $\Gamma_P(V)$. 
The weight of a top dimensional face $F$  of $\text{div}_{V}(P)$ is  described by the 
following two cases.
\begin{itemize}
\item  
If $F$ is a face 
with empty
sedentarity, 
the weight  of $F$  is the unique integer weight on 
$$\tilde{F} = \{ (x, P(x)) -te_{n+1} \ | \ x \in F \text{ and } t > 0 \}   $$ 
required to make $\Gamma_P(V) \cup \tilde{F}$ balanced 
along 
$\Gamma_P(F)$. 
\item  
If $F$ is a face 
with non-empty
sedentarity 
({\it i.e.}~if $P(F) = \{-\infty\}$), 
let $F_1, \dots F_k$ 
denote the facets 
of $V$
adjacent to $F$, and let
 $w_i$, $i = 1$, $\ldots$, $k$, 
 denote the weight of the  facet $F_i$ of $V$.
 For each $i = 1$, $\ldots$, $k$, 
there exists a unique primitive integer vector 
$v_i \in \ZZ^n$ such that 
$\lim_{t \to +\infty}  x_i - tv_i $ is contained in the relative
interior of $F$ if
 $x_i \in F_i$.
In a neighborhood of $F$, 
the value of
$P$ restricted to $F_i$ is 
given by
a single monomial 
$``a_{\alpha_i}x^{\alpha_i}"$,
and the weight of $F$ in $\text{div}_{V}(P)$ is 
defined as
$$\sum_{i=1}^k w_i \langle  \alpha _i,  v_i \rangle.$$
\end{itemize}
\end{defi}

The divisor  $\text{div}_{V}(P) $ is a union of tropical subvarieties of $\TT^n$ of possibly 
different sedentarities.  
Definition of divisors 
can 
be extended to tropical rational functions $R = ``\frac{P}{Q}" = P-Q$
which are tropical
quotients of polynomials. 
In this case, a divisor is a tropical cycle  in the sense of Definition \ref{def:tropvar}.  When the divisor 
is effective, {\it i.e.}~is a tropical 
subvariety of $\TT^n$, 
the definition of a tropical modification
can also be extended. 

\begin{exa}\label{exa:modline}
Consider $V = \T^2$ and $P(x, y) = ``x + y +0"$. The function $P$ has three domains of linearity in $\T^2$, 
$$
P(x, y) = \begin{cases} x &\mbox{if } x \geq \max\{y, 0\} \\
y & \mbox{if } y \geq \max\{x, 0\} \\
0 & \mbox{if } 0 \geq \max\{x, y\} \end{cases}.$$ 
Therefore, the graph of $\Gamma_P$ 
has three one-dimensional faces, 
and $\text{div}(P)$ is precisely the tropical line defined by $P$.
In addition, the surface $\tilde{V}$ is a tropical plane in $\TT^3$, and
the tropical modification $\pi_P: \tilde{V} \to V$ is the vertical projection,
see Figure \ref{Plane}.
\begin{figure}[h]
\begin{center}
\includegraphics[scale=1]{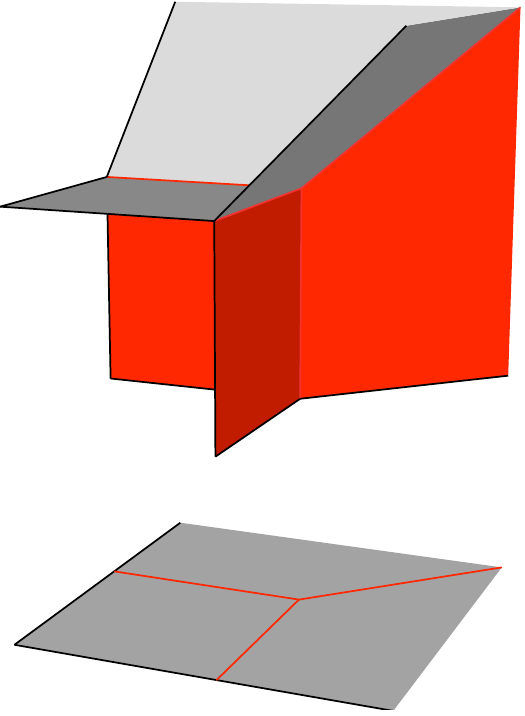}
\put(-110,20){0}
\put(-90, 45){$y$}
\put(-50, 15){$x$}
\put(-15, 15){$\T^2$}
\put(-150, 160){\tiny{$(-t,0,0)$}}
\put(0, 205){\tiny{$(t,t,t)$}}
\put(-122, 135){\tiny{$(0,-t,0)$}}
\put(-60, 110){\tiny{$(0,0, -t)$}}
\end{center}
\vspace{0.5cm}
\caption{Tropical modification  of  the tropical affine plane $\T^2$
  along $\tg x+y+0\td$.}
\label{Plane}
\end{figure}
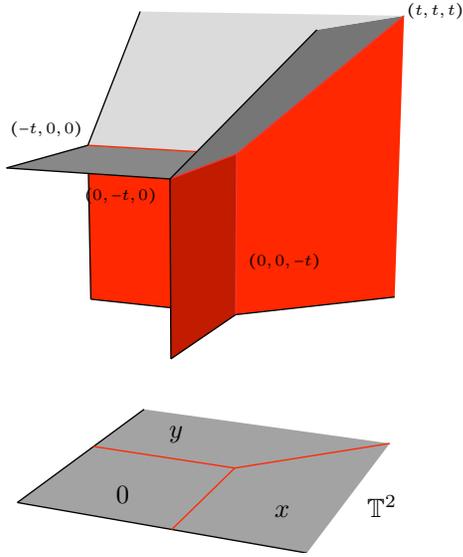
\end{exa}
\begin{exa}
More generally, if $V = \T^n$ and $P(x_1,\ldots,x_n)$  is a tropical polynomial, the divisor
$\text{div}(P)\subset \TT^n$ is
the tropical hypersurface 
defined by $P(x_1,\ldots,x_n)$, and $\tilde{V}\subset\TT^{n+1}$ is the tropical
hypersurface defined by the tropical polynomial $\tg x_{n+1}+P(x_1,\ldots,x_n)\td$.
\end{exa}

The next proposition relates tropical modifications and tropical limits of
subvarieties of $\CC^n$.

\begin{prop}\label{prop:tropgraphmod}
Let $V \subset \T^n$ be  the tropical limit  of a family of complex
algebraic
subvarieties $(\V_t)_{t\in\R_{>1}}$ of $\CC^n$, and $P: \T^n \to \T$ be a  tropical
polynomial.  Denote by $\pi_P: \tilde{V} \to V$ the tropical 
modification of $V$ along
$P$. 
Choose a family of complex
polynomials $P_t$ 
such that the tropical limit of this family is $P$, and 
denote by $\Gamma_{P_t}(\V_t)
\subset \CC^{n+1}$
the graph of $\V_t$ along the function
$\P_t$.  
Then, for a generic choice of the family $P_t$,  we have 
$$\lt \Gamma_{P_t}(\V_t) = \tilde{V}.$$ 
\end{prop}

\begin{figure}[h]
$$
\begin{array}{ccc}
\includegraphics[width=3cm, angle=0]{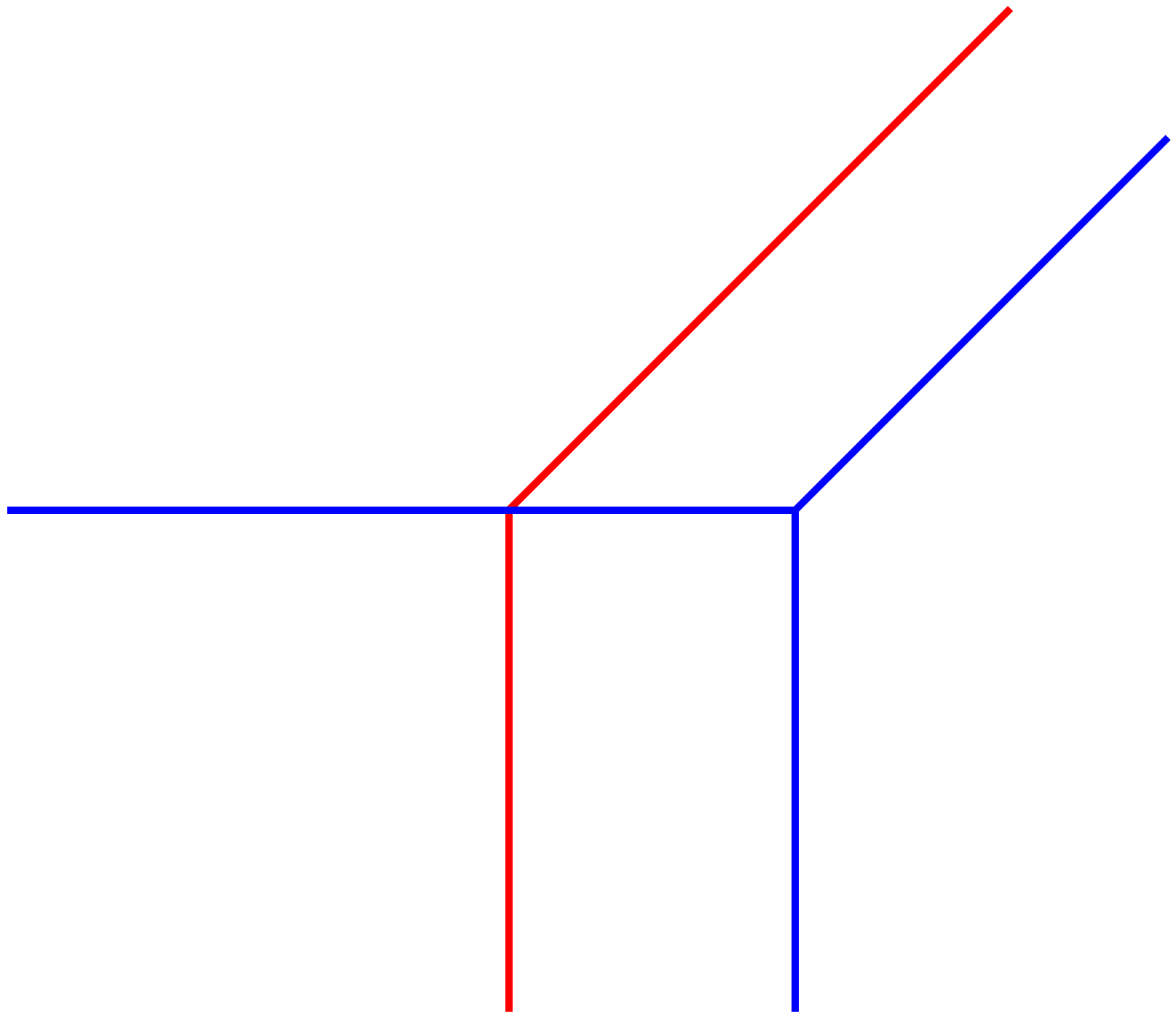}&
\hspace{10ex}&
\includegraphics[width=4cm, angle=0]{Figures/planelines.pdf}
\\ a) && b)
\end{array}
$$
\caption{The two lines from Example \ref{exa:modline} and the tropical limit of the configuration in the modified tropical plane}\label{fig:modintersection}
\end{figure}

Upon taking the tropical limit, properties of varieties which may be of interest
sometimes are no longer 
visible, and 
the tropical limit of a different embedding can reveal new features. 
We make 
this vague
remark 
explicit
with two examples, returning to curves in the plane. 

\begin{exa}
For a pair of  planar tropical curves $C$ and $C^{\prime}$
intersecting transversally, Proposition \ref{limit inter} 
relates 
the
tropical limit of the intersection points of the families $(\C_t)_{t\in\R_{>1}}$ and
$(\C^{\prime}_t)_{t\in\R_{>1}}$  counted with multiplicities to the area of polygons
in the subdivision dual to $C \cup C^{\prime}$. When $C$ and
$C^{\prime}$ do not intersect transversally, we cannot determine  the
precise location of the limit of the intersection points $\C_t \cap
\C^{\prime}_t$. 
By 
embedding suitably $\CC^2$ to $\CC^3$, we may
reveal the location of the tropical limit  of $\C_t
\cap \C^{\prime}_t$. 

Consider the constant family $\C$ defined by $P(z, w) = z+w+1$ and let $\C^{\prime}_t$ be defined by 
 $$Q_t(z, w) =  (t+1)z+w+(1-t^{b+1})=0\} \textrm{ with } b\le -1.$$ 
 A simple substitution verifies that $\C$ and $\C^{\prime}_t$ intersect at ${\bf{p}}_t  = (t^b, -1 -t^b)$ whose tropical limit is  to $(b, 0) \in \R^2$.
In the tropical  limit  we obtain two tropical lines which intersect along a real half line $\{ (t, 0) \ | \ t \leq -1\}$, thus the limit of the point 
${ \bf p}_t$ is not visible in this tropical limit 
(see Figure \ref{fig:modintersection}a). 
Consider the reembedding given by taking the graph along $P$, {\it i.e.}~
the graph $\Gamma_P: 
\CC^2 \to \CC^3$ 
of the function $(z, w) \mapsto (z, w, P(z, w))$.
This graph is simply a hyperplane in $\CC^3$ and  $\lt \Gamma_P(\CC^2)$ 
is a
tropical plane $\Pi$
from Example \ref{ex:trophyper}. Moreover 
the projection $ \Pi \to \T^2$  
is a tropical modification along the tropical function $``x+y +0"$. 

The reembedding of the  family of lines $\C^{\prime}_t$ is 
$$
\Gamma_P(\C^{\prime}_t) = \{ (z, t^{b+1} - 1 -(t+1)z,  t^{b+1} -tz ) \ | \ z \in \CC\} \subset \CC^3. 
$$
Denote the tropical limit  of  $\Gamma_P(\C^{\prime}_t)$  by $\tilde{C}^{\prime}$. Then $\tilde{C}^{\prime}$ is a tropical curve 
which must be contained in $\Pi$ and satisfy $\pi(\tilde{C}^{\prime}) = C^{\prime}$. Moreover, the intersection point $\Gamma_P({\bf p}_t) = (t^b, -1 -t^b, 0)$ 
is sent to $(b, 0, -\infty)$ in the tropical limit. This implies that
the tropical curve $\tilde{C}^{\prime}$ must have an unbounded ray of
the form $(b, 0, -s)$ for $s>>0$.    In the 
modified
picture the position of the intersection point of $\C$ and $\C^{\prime}_t$ is revealed
(see Figure \ref{fig:modintersection}b).
The use of tropical modifications to study intersection points in the non-transverse case is a technique used in \cite{BrLop}.
\end{exa} 

\begin{figure}[h]
\begin{tabular}{ccc}
\includegraphics[height=2.5cm, angle=0]{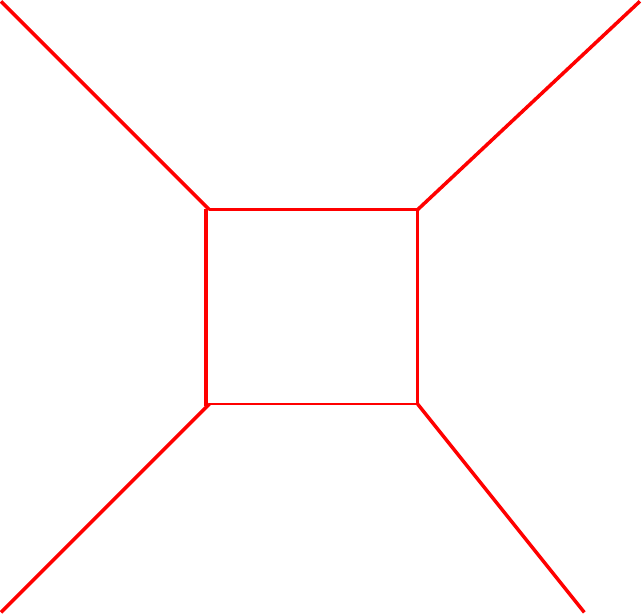}
&\hspace{1ex} &
\includegraphics[height=2.5cm, angle=0]{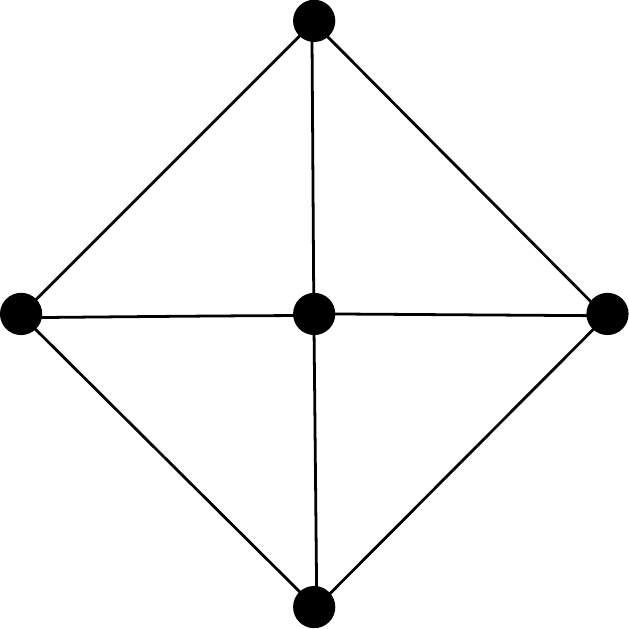}
\\a)  & &b) In coordinates $(x_1, x_2)$
\\ \\ \\ 
\includegraphics[height=2.5cm, angle=0]{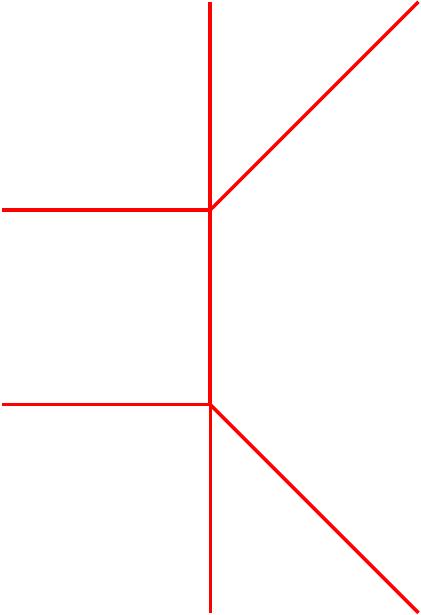}
\put(-20, 30){$2$}
&\hspace{1ex} &
\includegraphics[height=2.5cm]{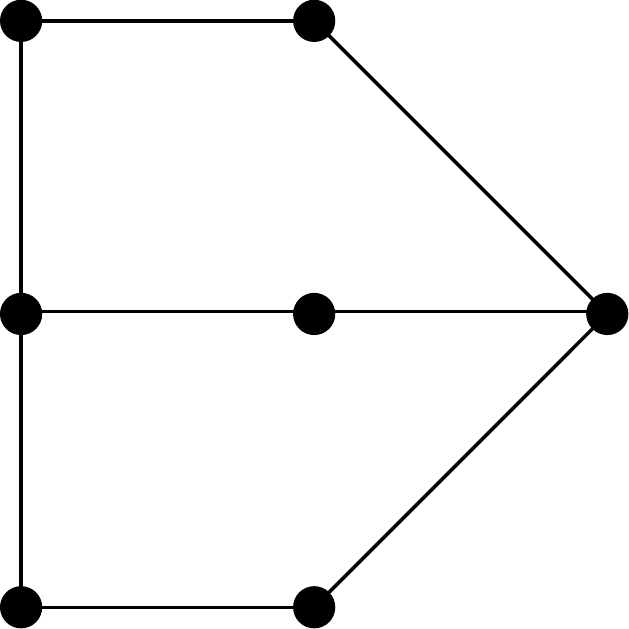} \\
c) & & d) In coordinates $(x_3, x_2)$
\end{tabular}
\caption{The two tropical limits  of the curve in Example \ref{ex:cycmod} with respect to the two coordinate systems. Beside them are the respective subdivisions of their Newton polygons. \label{fig:cyclemod}}
\end{figure}

\begin{exa}\label{ex:cycmod}
Consider the family of complex  curves in $\CC^2$ defined by the
polynomials 
$P_t(z_1, z_2) = z_1 + z_2 + t^az_1z_2+ z_1^2z_2 + z_1z_2^2$
for $a > 0$. It can be checked that 
these curves are 
smooth 
and of genus one. Then, 
$\lt \C_t$
is 
the tropical curve defined by
the tropical polynomial $``x_1 + x_2 + ax_1x_2 + x_1^2x_2 + x_1x_2^2"$.
This tropical curve is non-singular and  has a cycle of length $4a$.
The tropical curve and its dual subdivision are depicted in Figures \ref{fig:cyclemod}a,b. 

Now perform a linear 
coordinate change 
$z_3= z_1 - \alpha t^a$  for some complex number $\alpha$. Then, we have
$$
P_t(z_3, z_2) = \alpha t^a + z_3 + (1 + \alpha t^{2a} + \alpha^2 t^{2a})z_2 + (t^a+ 2\alpha t^a)z_2z_3 + z_3^2z_2 +z_2^2z_3 +\alpha t^az_2^2.
$$ 

As long as $\alpha \neq -1$ or $-1/2$ the tropical limit of this
family is given by the tropical polynomial $$P_{trop}(x_3, x_2) =  ``a
+x_3+2ax_2 + ax_2x_3 +x_3^2x_2 + x_2^2x_3 + ax_2^2".$$ 
The tropical limit is dual to the subdivision of the Newton polygon depicted 
in
Figures \ref{fig:cyclemod}c,d. Therefore, the tropical curve defined by this change of 
coordinates contains no cycle. 

Of the two tropical limits  above the one we started with is arguably better. First of all, the  first Betti number of $C$ is equal to the genus of a member of $\C_t$. Moreover, $C$ is non-singular and this implies a relation between the integer  length of the cycle of $C$  and the limit as $t \to +\infty$  under $\log_t$  of the $j$-invariants of $\C_t$ \cite{MarkKatz:jinvar}. 
There is also a tropical picture, which captures both of these
tropical limits, it is obtained by a tropical modification as follows. 
Consider the family of embeddings  $i_t: \CC^2 \to \CC^3$ given by
taking the graph along the function $z_1 - \alpha t^a$. 
Taking the tropical limit,
we obtain the tropical modification of $\T^2$ along the vertical line
$x_1= a$. It consists of three $2$-dimensional faces intersecting in a
line shown in Figure \ref{fig:cyclemod2}. Taking the tropical limit of
$i_t(\C_t) \subset \CC^3$
we obtain the tropical curve contained in the union of these three faces also depicted in the same figure. 

Both tropical limits above can be seen from this picture. Indeed we
obtain the first curve if we project onto the $(x_1, x_2)$ coordinates
and we obtain the second picture if we project onto the $(x_3, x_2)$
coordinates.  We could imagine that had we started with a curve
defined by a polynomial in  the second set of coordinates, we would
wish to find the change of coordinates  which produces a smooth
tropical curve.  ``Repairing" tropical limit of curves using
tropical modifications is the subject of \cite{CuetoMarkwig}. 
\begin{figure}[h]
\includegraphics[scale=1]{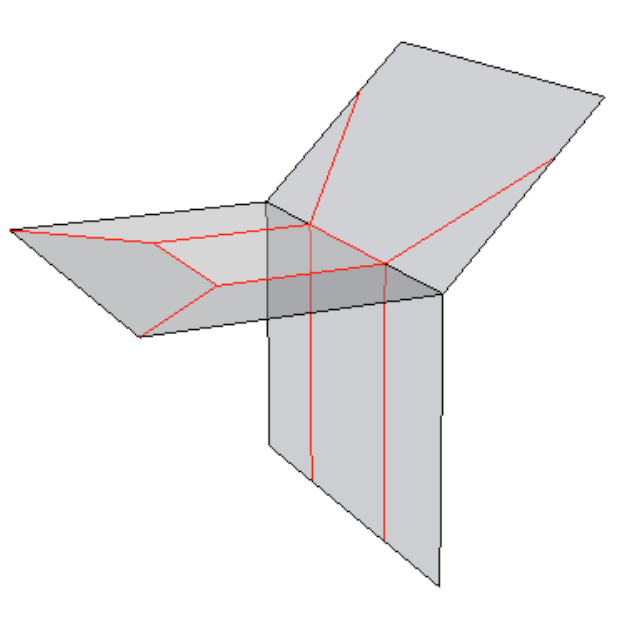}
\put(-190, 90){$\nwarrow x_2$}
\put(-135, 70){$\rightarrow x_1$}
\put(-50, 30){$\uparrow x_3$}
\caption{The tropical modification of $\T^2$ which captures the tropical limit of the family of curves $(\C_t)_{t\in\R_{>1}}$ from Example \ref{ex:cycmod} with respect to both systems of coordinates. \label{fig:cyclemod2}}
\end{figure}
\end{exa}

\begin{exo}
\

\begin{enumerate}
\item Prove that a  tropical hypersurface $V$ in 
$\RR^n$ 
with Newton polytope
the standard $n$-simplex of size $d$
has
  at most $d^n$ vertices, and that $V$ is non-singular if 
and only if
equality
  holds (compare with Exercise $2(2)$).
\item 
Let $S$ be a 
non-singular tropical 
surface
in $\RR^3$ 
with Newton polytope
the standard $3$-dimensional simplex 
of size $2$.
Show that $S$
 has a unique
  compact facet. (Hint:
  one has to prove that the dual subdivision of the simplex has
  a unique edge not contained on the boundary; use the fact that the 
  tetrahedron has Euler characteristic $1$).

\item \label{exo:skeleton} For a 
$\partial$-transversal linear space 
$\L \subset (\CC^*)^n$ 
of dimension $k$, show that $\lt \L$ is the $k$-skeleton of the hyperplane in $\T^n$ 
from Example \ref{ex:trophyper}. 

\item Verify that the braid arrangement from \ref{ex:M05} produces the cone over the Petersen graph from Figure \ref{fig:M_05}.

\item Let 
$\L \subset (\CC^*)^n$ 
be a 2-dimensional linear space. Find a formula for the Euler characteristic of $\L$
in terms of 
$\lt \L$.

\item  Let 
$\L \subset (\CC^*)^n$ 
be a linear space. 
Using Construction \ref{construction:matroidfan}, 
show that  $\lt \L$ satisfies the balancing condition.

\item Let $\Pi \subset \T^3$ be the 
tropical
plane 
from Example \ref{ex:trophyper} with vertex at the
origin, and let $L = \{ (x, x, 0) \ | \ x \in \T \} \subset 
\Pi$, 
equipped with weight one on each of its two edges. Find a 
function $P$ such that 
$div_\Pi(P) = L$. 
Can $P$ be a tropical
polynomial? 
\item Find the image of the tropical limit of the curve from Example \ref{ex:cycmod} when $\alpha = -1$ and $-\frac{1}{2}$.
\end{enumerate}
\end{exo}

%% file: section53.tex
\subsection{Tropical limits of  algebraic varieties}
As in the case of 
hypersurfaces,
families
of algebraic subvarieties of $(\CC^\times)^n$ are related to tropical
subvarieties of $\RR^n$. 
We briefly indicate below this relation
and refer to \cite{IKMZ} for
more details. 

Let 
$(\V_t)_{t\in A}$
be a 
family of proper complex analytic subvarieties
of dimension $k$
in $(\CC^\times)^n$. 
Here, $A\subset (1,+\infty)$ is any subset not bounded from above
(our main examples are $A= (1,+\infty)$ and $A=\ZZ\cap (1,+\infty)$).
%
The amoebas 
$\Log_t(\V_t)\subset\R^n$ form a family of closed subsets.
\begin{defi}
We say that $\Log_t(\V_t)$ {\em uniformly converges}
 to a closed
subset $V\subset\R^N$ if $V$ is the limit of $\Log_t(\V_t)$ with respect
to the Hausdorff metric on closed subsets of the metric space $\R^n$.

We say that $\Log_t(\V_t)$  {\em converges} 
 to a closed
subset $V\subset\R^N$ if it uniformly converges on compacts in $\R^n$,
{\it i.e.} for any compact $K\subset\R^n$
the family $K\cap\Log_t(\V_t)$ uniformly converges to $K\cap V$.
\end{defi}
We write $V=\lt \V_t$ for any of these convergences.
It turns out that $\lt \V_t$ admits the structure of
a rational finite polyhedral complex of  dimension $k$. 

\begin{exa}\label{exa-M}
Any linear map $M:\Z^n\to\Z^m$ 
is the multiplicative character map for a homomorphism
$\Phi_{M}:(\CC^\times)^m\to (\CC^\times)^n$.
Also it is (additively) dual to a map $\phi_M:\R^m\to\R^n$.
Amoebas $\Log_t(\Phi_M((\CC^\times)^m))\subset\R^n$ of the constant
(independent of $t$) family $\Phi_M$ coincide with $\phi_M(\R^m)$ for any $t$,
so we have
$$\phi_M(\R^m)=\lt 
(\Phi_M((\CC^\times)^m)).$$
In this case the tropical limit is a $k$-dimensional linear subspace of $\R^n$.
\end{exa}

Let $V=\lt \V_t$ for some family $\V_t\subset (\CC^\times)^n$ of $k$-dimensional
varieties.
Consider a facet $F\subset V$.
Since $F$ is a $k$-dimensional affine subspace of $\R^n$ we may find
an $(n-k)$-dimensional affine subspace $Y\subset\R^n$ intersecting $F$
transversely in a single generic point $p$ in the relative interior of
$F$ and such 
that the integer vectors tangent to $Y$ and those tangent to $F$ generate
the entire lattice $\Z^n$. We may find a linear map $M:\Z^n\to\Z^{n-k}$
such that $Y=p+\phi_M(\R^{n-k})$ so that $Y=\lt
t^p\Phi_M((\CC^\times)^{n-k})$, 
where $t^p(z_1,\ldots,z_n)=(t^{p_1}z_1,\ldots,t^{p_n}z_n)$.

\begin{defi}\label{def:trop limit}
We say that the face $F\subset \lt (\V_t)$ is {\em of weight $w$} if for 
any sufficiently small open set $U\ni p$ in $\R^n$ the set
$$t^p\Phi_M((\CC^\times)^{n-k})\cap \V_t\cap\Log_t^{-1}(U)$$
consists of $w$ points (counted with the intersection multiplicities
of cycles of complimentary dimension)
whenever $t$ is sufficiently large.  
It is a fact that
this definition of weight
does not depend on the choice of a generic point $p\in F$ 
or the subspace $Y$ with properties as above.

We say that $V=\lt \V_t$ is the {\em tropical limit} of the family $\V_t$ if
all the faces of $V$ acquire a well-defined weight.
If the convergence of $\Log_t(V_t)$
is uniform then we say that the tropical limit $V=\lt \V_t$ is uniform. 
\end{defi}

\begin{exa}\label{ex:matrix}
Suppose that $M:\Z^n\to\Z^k$ is a homomorphism such that the
cokernel group $\Z^k/\Im(M)$ is finite.
Then the weight of the affine $k$-space $\lt \Phi_M((\CC^\times)^k)$
from Example \ref{exa-M} coincides with the cardinality of $\Z^k/\Im(M)$.
\end{exa}

The following theorem generalizes Theorems \ref{approx} and
\ref{approx any dim}.
\begin{thm}[\cite{IKMZ}]
The tropical limit $\lt \V_t$ is a tropical subvariety of $\RR^n$.
\end{thm}

\begin{defi}\label{def:approx}
A tropical subvariety $V\subset\R^n$ is called \emph{approximable} if there
exists a family of complex algebraic subvarieties $\V_t$ of $(\CC^\times)^n$ such
that $\lt \V_t=V$.  It is called {\em uniformly approximable} if the convergence can be made uniform.
\end{defi}

The following example illustrates the difference between uniform
and conventional (that is uniform only on compacts) tropical convergence.  
If a family of hypersurfaces $\V_t\subset(\CC^\times)^n$
with the Newton polyhedron $\Delta$ has the tropical 
hypersurface $V\subset\R^N$ as its uniform tropical limit, then $\Delta(V)=\Delta$
(up to translation). 
However, if the convergence is uniform only on compacts, then 
it might happen that $\Delta(V)$ is strictly smaller than $\Delta$ as 
a part of the hypersurface may escape to infinity.
\begin{exa}
The tropical curve in $\RR^2$ defined by the tropical
polynomial $\tg 0+x\td $ is the tropical limit of
the family $(\C_t)_{t\in\R_{>1}}$ of curves in $(\CC^\times)^2$ with
equation
$$1+z+ t^{-t}w=0,$$
but is not a uniform tropical limit of this family.
\end{exa}

\begin{rem}
Specifying the weights of the facets in the tropical limit as in
Definition \ref{def:trop limit} may be further refined.
Namely, we can specify a finite covering of degree $w$ for the $k$-torus $(S^1)^k$ corresponding to
each facet $F$ of weight $w$. This refined weight is essential for consideration of tropical subvarieties,
but for the sake of keeping it simple we ignore it in this survey (replacing refined
tropical subvarieties with effective tropical cycles). Note though that there is no
difference if the weights of all facets are 1.

The weight refinement can be incorporated to the notion of the tropical limit.
This refinement may force subdividing
facets of the tropical limit to smaller subfacets with different pattern
of refined weights.
There is a tropical compactness theorem ensuring that any family of algebraic varieties
admits a tropically converging subfamily.
By passing to a subfamily
we can also ensure that the tropical limit in the refined sense exists.
\end{rem}

Not every tropical subvariety of $\RR^n$ is approximable, 
as illustrated by the 
following example.
Determining which tropical subvarieties 
of codimension different from one are approximable is quite
difficult, even in the case of curves.

\begin{exa}\label{exo:spatialcubic}
Consider
the
plane tropical cubic $C$ of genus 1 depicted in Figure \ref{equil}a, and
draw it in the affine plane with equation $z=0$ in $\RR^3$.
On each of the three unbounded edges in 
the direction $(1, 1, 0)$, choose a point in such a way that these three
points are not contained in a tropical line in $z=0$. Now at  these  three points,
replace the unbounded part of $C$ in the direction $(1,1,0)$ by two
unbounded edges, one in the direction $(0,0,-1)$ and one in the
direction $(1,1,1)$ to obtain a spatial 
tropical cubic $\tilde{C}$ (see Figure \ref{fig:spatial cubic}).

\begin{figure}[h]
\begin{center}
\begin{tabular}{c}
\includegraphics[width=7cm, angle=0]{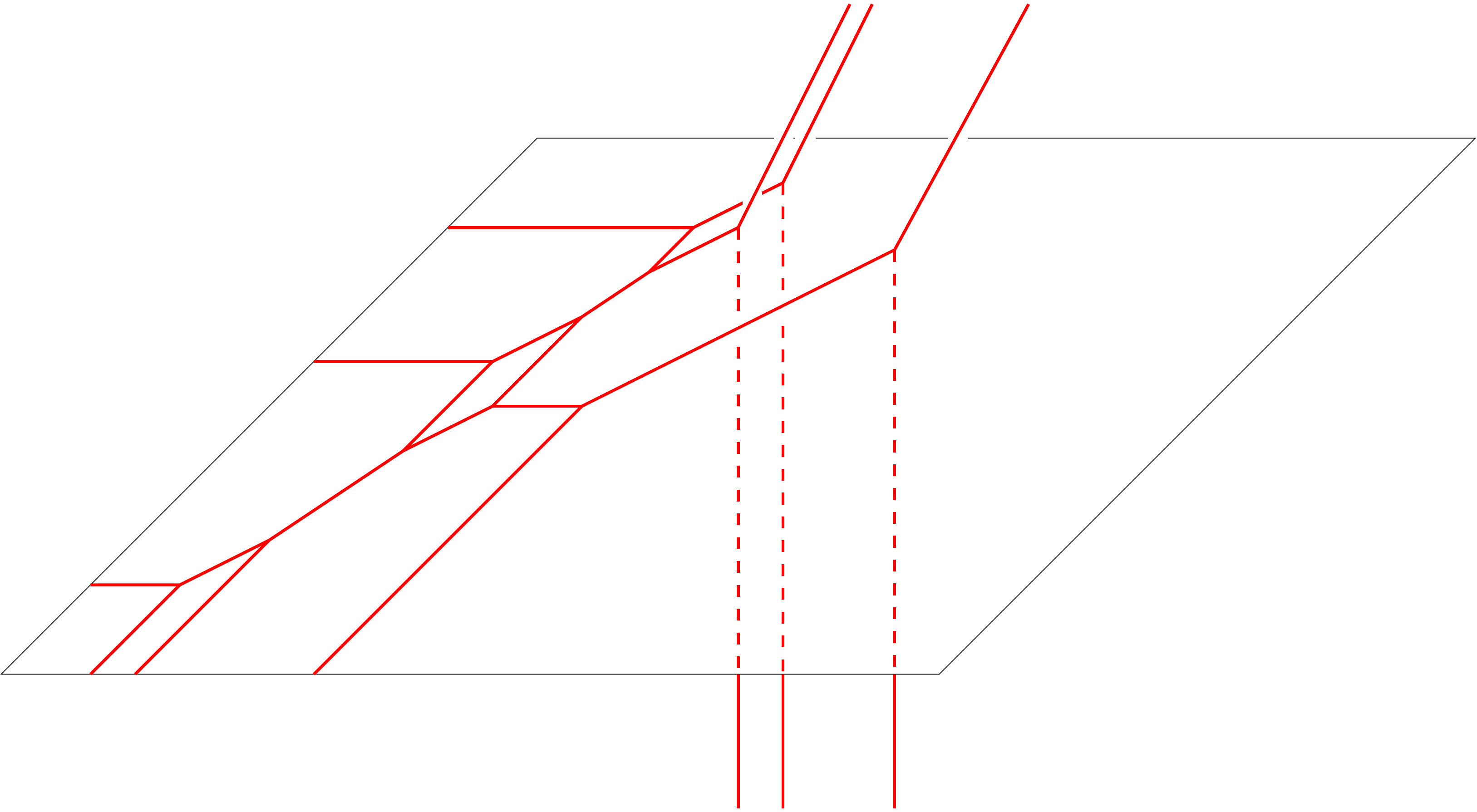}
\end{tabular}
\end{center}
\caption{A non-approximable tropical spatial cubic curve.}
\label{fig:spatial cubic}
\end{figure}
The choice we made of the three points ensure that 
$\tilde{C}$
is not contained in any tropical hyperplane; since as
we 
saw in Example \ref{ex:trophyper}, the unbounded rays in the
$-e_3$ direction contained in 
a tropical 
plane in $\RR^3$ 
project to a tropical line on the plane $z=0$.
 For the sake of simplicity, we show that $\tilde{C}$ is not
uniformly approximable 
{\it i.e.} there does not exist a family $(\C_t)_{t\in\R_{>1}}$ of  spatial
elliptic cubics with $\lt(\C_t)=\tilde{C}$. 
The tropical curve $\tilde{C}$ is not approximable even without
requiring uniform convergence,
however the proof is more involved.

Any  spatial algebraic elliptic cubic curve 
is planar.
Hence if 
an appropriate family  $(\C_t)_{t\in\R_{>1}}$  of elliptic cubics exists, 
each member  $\C_t$ would be necessarily contained in 
 a plane $\P_t$.
 After passing to a subfamily we may assume that $\lt \P_t$
 exists.
 It is a tropical plane  that contains
$\lt \C_t$. However our initial tropical curve 
is not contained in any tropical plane.
\end{exa}

Example \ref{exo:spatialcubic} admits generalisations, see
 \cite{Spe2,Nishinou,Br12}.

%% file: Kristin2.tex
\section{Tropical curves in tropical surfaces}

We have had quite an in-depth look at the applications of tropical
geometry 
in several questions related to 
curves in $\CC P^2$, or more generally in
toric surfaces.  
The patchworking construction
allows  
one
to construct real algebraic curves with prescribed topology. 
The correspondence theorem \ref{correspondence}
tells us
how to
count tropical curves with multiplicities in order to
obtain both Gromov-Witten and Welschinger invariants. 

Here we take a look at what happens with tropical curves in more
general 
tropical
surfaces in $\RR^n$ or $\TT^n$, and 
simultaneous approximation of curves and surfaces. 
We  mainly 
restrict to 
approximations by \emph{constant families} of pairs of 
 \emph{fan tropical curves} ({\it i.e.} tropical
curves with at most one vertex, at the origin) in tropical
planes.
This is not only an easier  particular case of the 
 approximation problem, 
instead it constitutes also a \emph{local} approximation problem needed to further study of any
 \emph{global} approximation. We provide more details in Section \ref{sec:global}. 

\subsection{Approximation of pairs}
 We are
interested in the following problem: given $C\subset S\subset \TT^n$ 
(or $C\subset S\subset \RR^n$)
a
tropical curve $C$ contained in a tropical surface $S$, does there exists
two families $\C_t\subset \S_t\subset \CC^n$ 
(or $\C_t\subset \S_t\subset (\CC^\times)^n$)
of complex algebraic
curves and surfaces such that $\lt \C_t=C$
and $\lt \S_t=S$?

It turns out that even if $C$ and $S$ are both approximable by
families $(\C_t)_{t\in\R_{>1}}$ and $(\S_t)_{t\in\R_{>1}}$, 
it might not be possible to find families
satisfying the extra condition  that $\C_t\subset \S_t$.

In the next section we will focus only on approximation 
of pairs by constant  families $\C \subset \P \subset (\CC^*)^n$,
where $\C$ is an algebraic curves and $\P$ is a plane. Recall that 
the tropical limits of $\C$ and $\P$ must both be fans
in this case.

\begin{exa}\label{ex:twocusps}
Consider the 
tropical plane 
$\tP \subset \T^3$ centered at the origin,
and the tropical curve $C$ made of
three rays in the directions
$$(-2,-3,0), \qquad (0,1,1), \qquad (2,2,-1),$$
each ray being equipped with  weight $1$. 
We already saw that $\tP$ is approximable.
The tropical curve $C$ is also approximable, 
for example one can check that $C=\lt \C$ where
$$\C=\{(\frac{u^2}{(u-1)^2}, \frac{u^3}{(u-1)^2}, u-1), \ u\in\CC\}.$$ 
However, we claim that there are  no
pairs $\C \subset \P\subset \CC^3$ for which 
$$C=\lt \C  \qquad \text{ and } \qquad  \tP=\lt \P.$$ 
We give several (related) proofs of this fact further in the text, nevertheless
we can already explain the reason why such a pair $\C \subset
\P\subset \CC^3$ does not exist:
if it did, the curve
$\C$ should have 
degree $3$, and have two cusps
singularities. This contradicts the fact that a plane cubic curve
cannot have more than one cusp singularity.
\end{exa}

\subsection{Intersection of fan tropical curves in fan tropical
  planes}\label{sec:fan curves}
As mentionned above, we
restrict in this section to 
approximations by \emph{constant families} of pairs of 
fan tropical curves in fan
tropical
planes. 

Intersection theory can be used to detect tropical curves in tropical
surfaces which are not approximable as pairs.
For the sake of simplicity, we restrict to the intersection of fan
tropical curves in fan tropical planes. We refer to \cite{Shaw1,AlRa1}
for more elaborate tropical intersection
theories.
Also we will restrict to planes 
$\P\subset(\CC^\times)^n$ 
that are 
$\partial$-transversal, 
{\it i.e.} no three lines in the corresponding line arrangement are 
concurrent. We refer to \cite{Br17} for the general case.

Recall that in 
Proposition  \ref{limit inter} we used the Euclidean area of
subdivisions dual to $C \cup C^{\prime}$ to determine the 
intersection multiplicities of curves $\C_t$ and $\C_t^{\prime}$. This
tool is no longer available to us when the tropical curves are in
general
fan
linear 
planes.
However, 
what is  visible from a fan tropical curve $\lt \C \subset \lt \P$
are the \emph{Newton diagrams of the curve} in different coordinates given by the
intersection points of the lines in the arrangement  
$\A$ corresponding to $\P$. 

Let $P(x,y)$ be a 
polynomial with Newton
polygon $\Delta(P)$, and let
$$
\bar{\Delta}(P) = Conv\{\Delta(P) \cup (0,0)\} \quad \mbox{and}\quad
\Delta(P)^c = \bar{\Delta}(P) \backslash \Delta(P).
$$

The polygon $\Delta(P)^c$ is the Newton 
diagram of the singularity   in coordinates $(x,
y)$ at the origin
of the curve defined by $P$. 
For a
different
system of coordinates the curve  has different Newton 
diagrams.
For tropical curves in $\T^2$
containing none of the coordinate 
axes,  
by duality there is a correspondence between the unbounded rays of the tropical curve heading toward $(-\infty, -\infty)$ 
and the non vertical/horizontal edges of $\Delta^c$. The unbounded
rays of the 
tropical limit
$\lt \C$ of a curve 
$\C\subset (\CC^\times)^n$ 
contained in a
$\partial$-transversal plane $\P$
provide an 
information 
on the  
Newton 
diagrams
of the curve 
in $\CC P^2\simeq\overline\P\subset \CC P^n$ 
in the systems of coordinates coming from pairs of lines in the arrangement $\A$ defined by $\P$. 

\begin{exa}\label{ex:cusp}
Let $\P=\CC^2$ equipped with some coordinate system, and let $\C\subset\P$ 
be a complex algebraic curve. 
Then, $\Delta(\C)^c$ is the triangle with vertices $(0,0)$, $(0,2)$,
and $(3,0)$ if and only if 
$\lt \C$ has a unique ray passing through $(-\infty,-\infty)$, and
this latter is of weight 1 with direction  $(-2,-3)$.
Note that this is also equivalent to the fact that the curve $\C$ has
a cusp at the origin, with the tangent 
at the origin being the abscissa axis. 
\end{exa}

\begin{exa}
Returning  to the curve from Example \ref{ex:twocusps}, each ray is
contained in a different face of the standard tropical plane $\tP = \lt \P$,
where $\P$ is defined by the equation $z_1 + z_2 + z_3 = 0$. The ray
in direction $(-2, -3,0) $ 
is
contained in the face spanned by $-e_1$ and $-e_2$.
It
 tells us that in coordinates given by lines $\L_1 = \{z_1 = 0\}$ and $\L_2
 = \{z_2 = 0\}$,
 an approximating curve $\C$ must have the triangle with vertices $(0,0)$, $(0,2)$,
and $(3,0)$ as Newton 
diagram.
According to Example \ref{ex:cusp}, 
the curve $\C$  
must have
a cusp at $\L_1 \cap \L_2$. Similarly, the ray in direction $(2, 2,
-1)$ is
contained in the face spanned by $e_0$ and $-e_3$ since
as this direction is $(2, 2,-1)=2e_0 - 3e_3$.
Again, this indicates 
that $\C$ has a cusp at the intersection of
$\L_3 = \{z_3 = 0\}$  with the
line at infinity. 
The last
direction $(0,1,1)$ tells us that  the curve is of multiplicity one at
the point $\L_0 \cap \L_1$.   

This along with the fact that  an approximating curve must have degree
three, tells us that the tropical curve from Example \ref{ex:twocusps}
can not be approximated by a
complex cubic curve
in $\P$: 
each cusp singularity decreases the genus of a curve by one, and a degree three curve in a plane has genus at most one.
\end{exa}

Now 
we
define the intersection of two fan tropical curves in a fan
tropical plane.
Recall that if $\P\subset \CC^n$ is a 
$\partial$-transversal plane,
each two-dimensional face of $\lt \P$ is the cone generated by two vectors $v_i, v_j$ corresponding to a pair of lines of 
the line arrangement $\A$ associated with $\P$.
Thus, the two-dimensional faces of $\lt \P$ are 
in bijection with pairs of lines $\L_i, \L_j$
of $\A$.

\begin{defi}\label{def:cornerint}
Let 
$\P  \subset (\CC^\times)^n$ 
be a 
$\partial$-transversal plane.
Given two fan tropical curves 
$C_1, C_2 \subset \lt \P\subset\RR^n$, suppose that for a  
two-dimensional 
face $F_{ij}$ of $\lt \P$ spanned by $v_i, v_j$, the curves $C_1, C_2$ have each exactly one ray in its interior. 
Suppose the ray of $C_1$ has weight $w_1$ and is in direction $p_1v_i + q_1v_j$
and the ray of $C_2$ has weight $w_2$ and is in direction $p_2v_i + q_2v_j$. 
Define the \emph{corner intersection multiplicity} 
of $C_1$ and $C_2$ in $F_{ij}$ as 
$$(C_1\cdot C_2)_{ij} = w_1w_2\min \{p_1q_2, q_1p_2\}.$$ 
 
\end{defi}

When the curves $C_1, C_2$ have 
several 
rays in the interior of a face, the above definition is extended by distributivity.

\begin{definition}\label{def:localInt}
Let
$\P  \subset (\CC^\times)^n$ 
 be a 
$\partial$-transversal plane. Given  two fan tropical curves
$C_1$ and $C_2$ in $\lt \P$, their  \emph{tropical intersection multiplicity}
at the origin  of $ \lt \P$ is defined as
$$(C_1\cdot C_2)_0= \deg(C_1)\cdot\deg(C_2) - \sum_{F_{ij}}
(C_1.C_2)_{ij},$$
where the sum is taken over all two-dimensional faces of $\lt \P$. 
\end{definition}

Here $\deg(C)$ is the degree of the tropical curve considered in
$\TT^n$ 
(see Definition \ref{def:degree}). 
Notice that the above definition also applies to tropical curves 
which have edges in common, 
and even to self-intersections.

\begin{exa}
Suppose that $\P=(\CC^\times)^2$, and that $C_1$ and $C_2$ are two fan
tropical curves in $ \lt \P=\RR^2$ centered at the origin. Then 
$(C_1\cdot C_2)_0$ is equal to the mixed volume of $\Delta(C_1)$ and
$\Delta(C_2)$, {\it i.e.}
$$(C_1\cdot C_2)_0= Area(\Delta(C_1\cup C_2)) - Area(\Delta(C_1)) - Area(\Delta(C_2)),$$
where $Area(\Delta(C))$ is the Euclidean area of the Newton polygon
of $C$ (compare with Exercise $2(4)$). 
In particular $(C_1\cdot C_2)_0\ge 0$.
\end{exa}

\begin{exa}\label{ex:trop inters}
Let $C$ be the 
 degree $3$ tropical curve from Example \ref{ex:twocusps}, and let 
$$L = \{ (x, x, 0) \ | \ x \in \RR \}$$
be  the degree 1 tropical curve 
equipped with weight one on its edges. Both tropical curves are contained in the
 tropical plane $\Pi$ centered at the origin in $\RR^3$, and 
 we have the following intersection numbers in $\Pi$
$$(L^2)_0=-1,\quad (C^2)_0=-4,\quad\mbox{and}\quad (C\cdot L)_0=-1.$$
\end{exa}

Classical and tropical intersections are related by the following
theorem.
Recall that two complex algebraic curves intersect \emph{properly} if they
intersect in finitely many points.
\begin{thm}[\cite{Br17}]\label{thm:relation inter}
Let 
$\C_1, \C_2\subset (\CC^\times)^n$ 
be two complex algebraic curves in a 
$\partial$-transversal
plane 
$\P\subset (\CC^\times)^n$.  
We denote  respectively by $C_1$ and $C_2$ the tropical limits
of $\C_1$ and $\C_2$, and
by 
$m(\C_1 \cdot \C_2)$ 
their number of intersection points
in $\P$ 
counted
with multiplicity.  
 If $\C_1$ and $\C_2$ intersect properly,
then we have
$$m(\C_1 \cdot \C_2) \leq (C_1\cdot C_2)_0.$$
If in addition 
the intersection of $C_1$ and $C_2$ is reduced to a point, 
then this inequality is an equality.
\end{thm}
\begin{rem}
A similar result holds without the assumption that $\P$ is 
$\partial$-transversal.
\end{rem}

Recall that intersection points of two complex algebraic curves are
always positive. As a consequence,
we deduce two immediate corollaries
from Theorem \ref{thm:relation inter}. 

\begin{corollary}[\cite{Br17}]\label{cor:CD<0}
Let  
$\P \subset (\CC^\times)^n$ 
be a 
$\partial$-transversal
plane.
Suppose there exists 
an irreducible and reduced 
 complex algebraic curve $\C \subset \P$ such that $\lt \C = C$. 
If $D \subset \lt \P$ is 
another
 fan tropical curve such that 
 $D \neq C$ and $(C.D)_ 0 <0$, then 
$D$ is not the tropical limit of any 
 irreducible  
 complex algebraic curve $\mathcal D \subset \P$.
\end{corollary}

\begin{exa}
It follows from Corollary \ref{cor:CD<0} that at most one of the two
tropical curves from Example \ref{ex:trop inters} can be 
the tropical limit of a
 complex algebraic curve contained in the plane $\P$ with equation
$z_1+z_2+z_3+1=0$.
Since one easily sees that $L$ is 
the tropical limit of
a line  
$\L\subset\P$,
 this proves that
the tropical curve $C$ is not the tropical limit of any irreducible
 complex algebraic curve $\mathcal C \subset \P$.
\end{exa}

\begin{corollary}[\cite{Br17}]\label{cor:CC<0}
Let 
$\P \subset (\CC^\times)^n$ 
be a non-degenerate plane, and suppose that 
there exists 
a reduced
irreducible curve $\C \subset \P$ such that $\lt \C = C$. 
If 
$(C^2)_0 <0$, then $\C$ is the unique complex
algebraic curve in $\P$ 
whose tropical limit is $C$.
\end{corollary}

\begin{exa}
The tropical line $L$ from Example \ref{ex:trop inters} is
approximable by a unique line contained in the plane $\P$ with equation
$z_1+z_2+z_3+1=0$.
\end{exa}

Combining Theorem \ref{thm:relation inter} with classical results from
algebraic geometry, one may obtain further obstructions to the
approximability of pairs. As an example, 
the following 
theorem can be deduced
from Theorem \ref{thm:relation inter} together
with the adjunction formula. 

\begin{thm}[\cite{Br17}]\label{thm:localadj}
Let 
$\P \subset (\CC^\times)^n$
be a 
$\partial$-transversal 
plane, 
and let $C \subset \lt \P \subset \RR^n$ be a fan tropical curve of 
degree $d$.
If there exists an irreducible and reduced complex algebraic curve $\C \subset \P $
such that $\lt \C=C$, then 
$$(C^2)_0 + (n-2)d - \sum_{e} w_{e} + 2 \geq 2g(\C),$$
where the sum goes over all edges $e$ of $C$,  $w_e$ is the weight
of $e$, and $g(\C)$ is the geometric genus of $\C$.
In particular, if the left hand side is negative, then $C$ is not
approximable by a reduced and irreducible 
complex algebraic
curve in $\P$.  
\end{thm} 

\begin{exa}
We provide 
another 
proof of the fact that the tropical curve $C$ from Example \ref{ex:twocusps}
does not form an approximable pair with the tropical plane $\Pi$. 
Indeed, since $(C^2)_0=-4$, the left hand side is equal to $-2$.
\end{exa}

Another example of application of 
the techniques introduced in this
section is the classification of approximable trivalent fan tropical
curves in fan tropical planes. The next statement is a special case of
this classification.
\begin{thm}[\cite{Br17}]\label{thm:trivalent}
Let 
$\P \subset (\CC^\times)^3$ 
be a 
$\partial$-transversal 
plane, 
and $C \subset \lt \P$ be a fan tropical curve 
made of at most three rays.
Then, there exists an irreducible and reduced complex curve $\C \subset \P $
such that $\lt \C=C$ if and only if $(C^2)_0=0$ or $(C^2)_0=-1$.
\end{thm}

\subsection{From local to global}\label{sec:global}
So far 
we have mostly restricted
to considering 
tropical limits of constant families of planes and
curves defined over $\CC$. This restriction always  
produces tropical spaces which are fans. 
Considering these cases is 
still useful when we pass to 
tropical limits of families 
of varieties 
due to a \emph{localization procedure}. Suppose that 
a family  
$(\V_t)_{t\in A}$ 
in $(\CC^\times)^n$ 
has $V$ as tropical limit. 
For 
any 
point 
$x$ of the tropical variety $V$, 
we denote by $V(x)$
the  
fan composed 
of all vectors $v\in\RR^n$
such that $x+\varepsilon v$ is contained in $V$ for 
a sufficiently small 
positive real number $\varepsilon$;
this fan is 
equipped with the weight function on its facets which is inherited from $V$. 
The localization 
procedure provides,
for each point $x \in V$,
a 
complex algebraic
variety  $\V(x) \subset (\CC^\times)^n$
such that $\lt \V(x) = V(x)$. 
This localization procedure extends to  approximations of 
a pair (or any 
tuple) of tropical varieties.  
Therefore, if 
a pair consisting of
a tropical curve
in a tropical surface 
is 
not locally approximable 
at some point,
this pair is not 
globally approximable. 
However, there are also global obstructions to approximating tropical
curves in surfaces.
\begin{exa}
Let $C$ be the tropical curve whose directions to infinity
are
$$ (-2,1,1), \quad (1,-2,1),  \quad (0,0,-1), 
\quad \text{and} \quad (1,1,-1),$$
and which contains a bounded edge in direction $(1, 1, -2)$ (see Figure \ref{fig:global2}).
The curve $C$ is contained in the tropical plane $\Pi$ centered at the
origin, and the
pair is locally approximable (see Theorem \ref{thm:trivalent}).
However global approximability 
by planar cubics
would
imply the existence of an algebraic cubic in the projective plane
together with a line passing through exactly two of its inflection
points. This contradicts
the fact that a line passing
through two inflection points of a cubic actually intersects this cubic
in three inflection points.
\begin{figure}[h]
\includegraphics[scale=0.015]{Figures/plane2.pdf}
\put(-70, 50){\small{$(1, 1, -2)$}}
\put(-175, 95){\small{$(1, -2, 1)$}}
\put(-107, 100){\small{$(-2, 1, 1)$}}
\caption{}
\label{fig:global2}
\end{figure}
\end{exa}

Nevertheless, 
 some situations only require  the   consideration of
local obstructions. 
An example is given by
 the generalization 
of the 
Correspondence Theorem \ref{correspondence} to the
enumeration of tropical curves in the tropical surface from Example
\ref{ex:trophyper degen} (see \cite{Br18}).
Another example is provided by the study of tropical lines in tropical surfaces. 

In \cite{Vig1} and \cite{Vig2}, M. Vigeland exhibited generic
non-singular tropical surfaces of degree $d\ge 4$ containing tropical
lines, and  generic
non-singular tropical surfaces of degree $d=3$ containing infinitely
many tropical
lines. The 
following 
theorem shows that when we restrict our attention to
the tropical lines which are approximable in the surface,  the situation turns
out to be analogous to the case of complex algebraic surfaces.
\begin{thm}[\cite{Br17}]\label{prohib Vigeland intro}
Let $S$ be a
generic
non-singular 
tropical surface in 
$\TT^3$
of degree $d$.
If $d=3$,  
then there exist only 
finitely many tropical lines $L\subset S$ such that $L$ and $S$ 
form an approximable  pair. 

If $d\ge 4$, 
then there exist
no tropical lines $L\subset S$ such that $L$ and $S$ 
form an approximable pair. 
\end{thm}

\begin{exo}
\

\begin{enumerate}
\item Show that the tropical curve from Example \ref{ex:twocusps} is contained
  in $\tP$, and show that it has degree 3.

\item  Recheck the computation of $(L \cdot C)_0$ of Example \ref{ex:trop inters}.

\item Show that the tropical curve $L$ of Example \ref{ex:trop inters}
together with $\Pi$  
form an approximable pair. 

\item Let 
$\tP \subset \RR^n$ be the 
$\partial$-transversal tropical 
plane
centered at the origin. Show that the intersection number of the 
$\partial$-transversal tropical line in $\tP$ 
with any fan tropical curve $C$ in $\Pi$ is equal to $\deg(C)$. 
\end{enumerate} 
\end{exo}

%% file: Grisha.tex
\newcommand{\dd}{\partial}
\newcommand{\cp}{{\mathbb C}{ P}}
\newcommand{\tp}{{\mathbb T}{ P}}
\newcommand{\Hom}{\operatorname{Hom}}
\renewcommand{\FF}{\mathcal F}
\newcommand{\WW}{\mathcal W}

\section{Tropical manifolds and their homology groups}\label{sec:abstract}
\subsection{Abstract tropical manifolds}\label{sec:manifold}
So far we have seen examples of tropical 
subvarieties
in $\R^n$ and $\T^n$.
In this section,
we introduce the notion of {\it tropical manifold}.
The notion of {\it abstract tropical variety} was first introduced in \cite{Mik3}. 
Tropical manifolds have the restriction that they are locally modeled 
on fan tropical linear spaces 
(defined in Section \ref{sec:fanlinear}). Namely, we say that $M\subset\R^N\times\T^s\subset\T^{N+s}$ is 
a {\em tropical smooth local model} if $M=L\times\T^s$, where $L\subset\R^N$ is a fan tropical linear space.
The dimension of such a tropical smooth local model is $\dim L + s$. 

\begin{definition}\label{def:manifold}

An {\it $n$-dimensional tropical manifold} $X$ is a Hausdorff topological space equipped 
with an atlas of charts $\{ (U_\alpha , \Phi_\alpha) \}$, 
with
$\Phi_\alpha : U_\alpha \rightarrow X_{\alpha} \subset 
\T^{N_{\alpha}}$, such that the following hold: 

\begin{enumerate}

\item for every $\alpha$, 
the map $\Phi_\alpha : U_\alpha \rightarrow X_{\alpha} \subset \T^{N_\alpha}$
is such that 
$X_\alpha$ is an $n$-dimensional 
tropical smooth model, 
and $\Phi_{\alpha}$ is 
an open embedding of $U_\alpha$ in $X_\alpha$;

\item for every $\alpha_1, \alpha_2$,
the overlapping map $\Phi_{\alpha_1} \circ \Phi_{\alpha_2}^{-1}$, defined on $\Phi_{\alpha_2}(U_{\alpha_1} \cap U_{\alpha_2})$, 
is the restriction of an 
integer affine linear 
map 
$\T^{N_{\alpha_2}}  \rightarrow \T^{N_{\alpha_1}}$
{\rm (}{\it i.e.}, of the continuous extension of an integer affine linear map 
$\R^{N_{\alpha_2}} \longrightarrow \R^{N_{\alpha_1}}${\rm )\/}; 

\item $X$ is of finite type, {\it i.e.}, there is a finite collection of open sets $\{W_i\}_{i=1}^m$ 
such that 
$\bigcup_{i=1}^m W_i = X$ 
and, for each $i$, there exists $\alpha$ satisfying the conditions $W_i \subset U_{\alpha}$ 
and $\overline{\Phi_{\alpha}(W_i)} \subset \Phi_{\alpha}(U_{\alpha}) \subset \T^{N_{\alpha}}$. 

\end{enumerate}

\end{definition}

As usual, two atlases on $X$ are called \emph{equivalent} if their
union is again an atlas on $X$; 
any equivalence class of atlases on $X$ has a unique saturated (or
maximal) representative. We always implicitly consider a tropical manifold
equipped with its maximal atlas, even when  defining its tropical
structure using  a non-maximal one.
\begin{exa}
The set of tropical numbers $\TT$ equipped with the identity chart
$Id:\TT\to\TT$ is a tropical manifold of dimension 1.

The tropical torus $\TT^\times=\RR$ equipped with the unique chart 
$Id:\RR\to\RR$ is not a tropical manifold since it does not satisfy
the finite type condition. 
Nevertheless one can enlarge this atlas with the two charts
$$\begin{array}{ccc}
(0,+\infty) & \longrightarrow & \RR
\\x&\longmapsto & x
\end{array} \quad \mbox{and}\quad
\begin{array}{ccc}
(-\infty,1) & \longrightarrow & \RR
\\x&\longmapsto & -x
\end{array}$$
which turn $\RR$ into a tropical manifold of dimension 1.

Analogously, the set $\RR_{>0}$ equipped with the 
inclusion chart $\RR_{>0}\hookrightarrow \RR$ 
is not a tropical manifold. However there is no way to
complete this atlas to turn $\RR_{>0}$ into a tropical manifold.
\end{exa}
\begin{exa}
The product of two tropical manifolds, equipped with the product atlas,
is a tropical manifold. In particular, the affine tropical spaces
$\TT^n$ and
the tropical torus $\RR^n=(\TT^\times)^n$
for any $n$ are tropical manifolds. 
\end{exa} 

\begin{exa}
Consider a tropical subvariety $V$ of $\RR^n$, equipped
with the atlas induced by the one on $\RR^n$.
Then $V$ is a $k$-dimensional tropical
manifold if and only if the fan $V(x)$ defined in Section \ref{sec:global} is a
tropical linear space for all points $x$ of $V$.
\end{exa}
\begin{exa}
Consider a lattice $\Lambda$ of rank $k$ in $\RR^n$, and the atlas on $\RR^n$
given by the identity map $Id:U\to U$ on open sets $U$ satisfying
$U\cap (U+\Lambda)=\emptyset$. This atlas induces a structure of
tropical manifold of dimension $n$ on the quotient space
$\RR^n/\Lambda$. Note that this
atlas also turns  $\RR^n/\Lambda$ into a differentiable manifold
diffeomorphic to $(S^1)^{k}\times\RR^{n-k}$.
\end{exa}

It is also possible to think of a tropical manifold as a locally
ringed space $(X, \mathcal{O}_X)$. There is a sheaf of \emph{regular
  functions} on $\T^n$, coming from the pre-sheaf of tropical
polynomials. By restricting the sheaf of regular functions on
$\T^{N_{\alpha}}$  to $X_{\alpha}$ for each chart we obtain a sheaf
$\mathcal{O}_{U_{\alpha}}$. The condition on the overlapping maps  
in Definition \ref{def:manifold} 
ensures that the local sheaves are compatible, 
that is,
the restrictions of $\mathcal{O}_{U_{\alpha_1}}$ and $\mathcal{O}_{U_{\alpha_2}}$ to $U_{\alpha_1} \cap U_{\alpha_2}$ agree.
We direct the reader to   \cite{Mik3} or \cite{MikZhar:Eigenwave} for more details.

If a tropical manifold $X$ is compact, 
then $X$ can be enhanced with a structure of 
a finite polyhedral complex of 
pure
dimension $n$. 

\subsection{Abstract tropical curves}\label{sec:abstract curves}
Since $GL_1(\ZZ) = O_1(\RR) = \{\pm 1\}$, a compact smooth tropical
curve 
gives rise 
to a finite graph equipped with a complete metric on the complement of
the 
set of
$1$-valent vertices. 
Conversely, each finite graph (without isolated vertices)
 equipped with a complete inner metric 
on the complement of the 
set of
$1$-valent vertices 
can be seen as a compact smooth tropical curve:
the interiors of edges can be identified by isometries with open intervals in $\RR$ 
(these intervals are unbounded for the edges adjacent to $1$-valent vertices), and 
for each vertex $v$ 
with valence $N_{v} + 1\ge 3$,
 we choose a chart
$\Phi_v : U_v \rightarrow X_{v} \subset \RR^{N_v} \subset 
\T^{N_{v}}$, where $ U_{v}$ is a neighborhood of $v$, 
and $X_v$ 
is the 
$\partial$-transversal fan tropical line in $\RR^{N_v}$
(see Example \ref{ex:uniformlin}). 

Recall tropical modifications of subvarieties of 
$\TT^n$ from Section \ref{sec:mod}. 
 Let $\Gamma$ and $\Gamma'$ be two 
 abstract
 tropical curves,
and let $p$ be a point in the complement of
the 
set of
$1$-valent vertices of $\Gamma$. 
We say that  $\Gamma'$ is  the \emph{elementary tropical
  modification of $\Gamma$ at $p$} if
$\Gamma'$ is obtained
 by gluing 
$\Gamma$ and $[-\infty,0]\subset\TT$ at $p$ and $0$ (the metric
considered on $(-\infty,0]$ is the standard Euclidean metric).
Notice that this is equivalent to performing,  in a single chart, a tropical modification 
in the sense of Section \ref{sec:mod}. 
We say that  $\Gamma'$ is  a \emph{tropical
  modification of $\Gamma$} if
$\Gamma'$ is obtained by a finite sequence of elementary tropical modifications
of $\Gamma$.
Tropical modifications can also be performed on tropical manifolds of arbitrary dimension
(see \cite{Mik3}). However, unlike for tropical curves,  it can be quite difficult to determine if two 
tropical manifolds of higher dimension are related by this operation.

One can introduce the notion of {\it tropical morphism}
between tropical 
manifolds.
For the sake of brevity  
once again, 
we do it in this text only in the case
of morphisms from a tropical curve to $\RR^n$.
Let $\Gamma$ be a compact connected smooth tropical curve.
As we already noticed, the curve $\Gamma$ can be seen as a graph
equipped with a 
complete inner
metric on the complement 
of 
the set of
$1$-valent vertices of $\Gamma$.
\begin{defi}\label{def:trop morphism}
Let $\Gamma^0 \subset \Gamma$ be  the complement of 
a subset of
$1$-valent vertices of $\Gamma$ (we do not have to remove all 1-valent vertices
from $\Gamma$, but only those sent to infinity by the morphism).
A {\it tropical morphism} from $\Gamma^0$ to $\RR^n$
is a 
proper
continuous map $f: \Gamma^0 \to \RR^n$ subject to the following two properties.
\begin{description}
\item[Integrality] the restriction of $f$ to each edge of $\Gamma^0$
is integer affine linear. Equivalently, the image of any unit tangent vector
to $\Gamma^0$ under the differential $df$ is integer, i.e. an element of $\ZZ^n\subset\RR^n$.
\item[Balancing] for each vertex $v$ of $\Gamma^0$, we have
$$
\sum_e u(e) = 0,
$$
where the sum is taken over all edges adjacent to $v$,
and $u(e)$ is the image under $df$ 
of the unit tangent vector to $e$ such that
this vector points outward of $v$. 
\end{description}
\end{defi}
Notice that an edge $e$ of $\Gamma^0$ is mapped to a point if and only
if $u(e)=0$.
Since in Definition \ref{def:trop morphism} we required $f$ to be proper
we have the following property:
an edge of $\Gamma$ 
adjacent to a 1-valent vertex $v$ 
is mapped to a point by $f$
 if and only if $v\in  \Gamma^0$. 

If $f: \Gamma^0 \to \RR^n$ is a tropical morphism, 
then $f(\Gamma^0)$ gives rise to a tropical curve $C$ in $\R^n$.
For each edge $e$ of $\Gamma^0$
with $u(e)\ne 0$ (we may take $u(e)$ with respect to any vertex adjacent to $e$),
 the weight $w(e)$ of $f(e)$
is the positive integer such that $\frac{u(e)}{w(e)}$ is a primitive
integer vector. 
If for several edges $e_1$, $\ldots$, $e_r$ of $\Gamma^0$ the intersection $\cap_{i = 1}^r f(e_i)$
is a segment, we put the weight of this segment to be equal to
$\sum_{i = 1}^r w(e_i)$.
We say that $\Gamma$ is a {\it parameterization} of $C$.

\begin{defi}\label{def:genus}
The  \emph{genus} of 
a tropical curve 
$\Gamma$ is 
defined as
the 
first Betti number $b_1(\Gamma)$
 of $\Gamma$. 
The \emph{genus} of an irreducible tropical curve $C$ in $\RR^n$ is
the minimal genus of parameterizations of
$C$.
\end{defi}

One can easily check that in the case of irreducible nodal 
tropical
curves in $\RR^2$,
this definition of genus coincides with the definition given in Section \ref{tropical_enumeration}.
There exist tropical curves whose only tropical morphisms to $\RR^n$
are the constant maps. Nevertheless, up to tropical modifications and
after removing 
 1-valent vertices,
 any tropical curve 
admits a non-constant tropical
morphism to some $\RR^n$.

\begin{exa}\label{ex:parameterization}
Figure \ref{parametrization}a shows an example of a tropical curve 
$\Gamma$
of genus 1. Clearly, the curve $\Gamma$ does not admit any non-constant 
tropical morphism to $\RR^n$.
By a sequence of nine elementary tropical modifications, one obtains the
 tropical curve $\Gamma'$ 
depicted in Figure \ref{parametrization}b, which parameterizes (after
removing its nine 1-valent vertices) the tropical cubic in $\R^2$
depicted in
Figure \ref{parametrization}c. 

\begin{figure}[h]
\centering
\begin{tabular}{ccccc}
\includegraphics[height=4cm, angle=0]{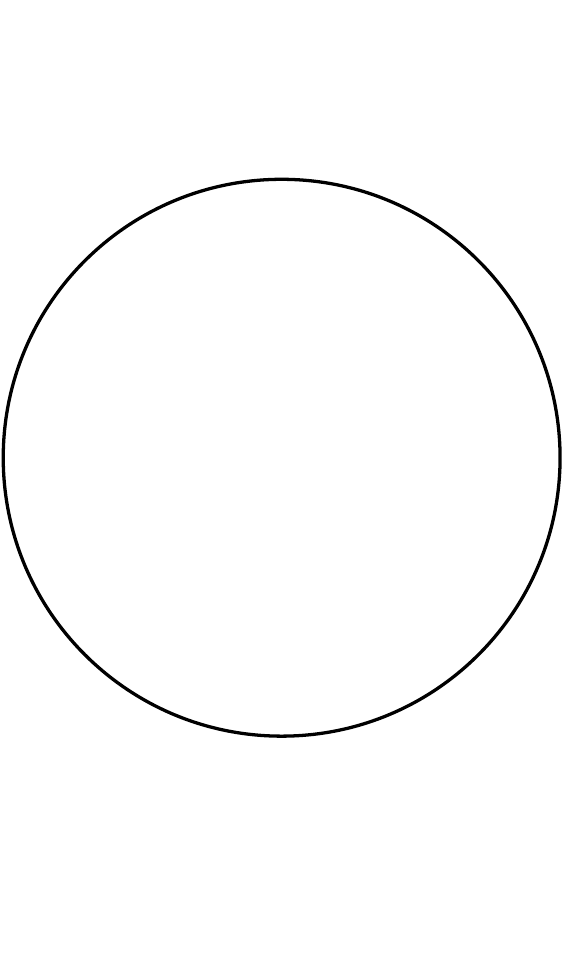}& 
\hspace{3ex} &
\includegraphics[height=4cm, angle=0]{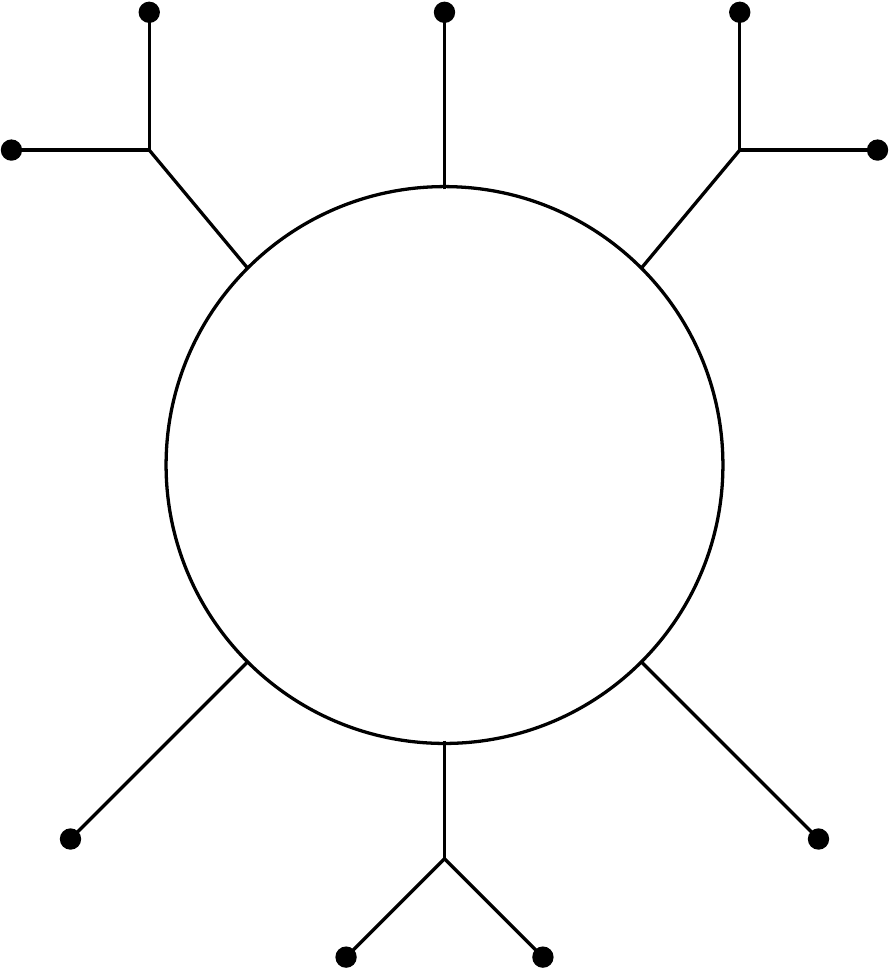}& 
\hspace{3ex} &
\includegraphics[height=5cm, angle=0]{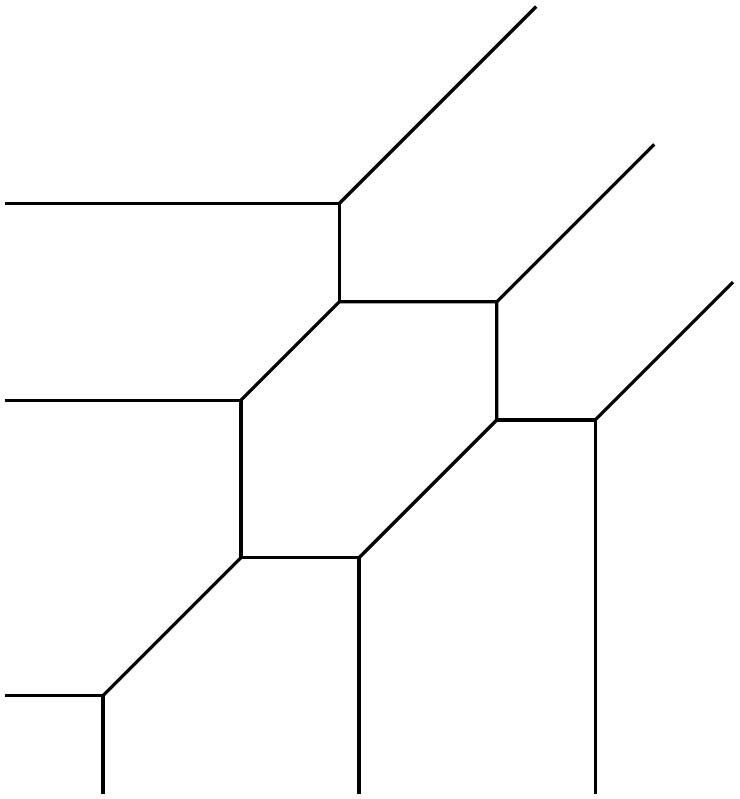} 
\\ \\ a) && b) && c)
\end{tabular}
\caption{
Tropical curve, tropical modification, and parameterization from Example \ref{ex:parameterization}
} 
\label{parametrization} 
\end{figure} 
\end{exa} 

\subsection{Tropical toric varieties}\label{proj spaces} 
The logarithm  transforms
multiplication to addition.  
As a consequence, any operation performed in complex algebraic
geometry  using only monomial maps translates {\it mutatis mutandis} in the 
tropical setting. In other words, 
non-singular
tropical toric
varieties 
are constructed analogously as in complex geometry. 
Let us illustrate this 
in the case of 
projective spaces. 

The projective line $\CC P^1$ may be obtained by taking two copies of
$\CC$, with coordinates $z_1$ and $z_2$,
 and  gluing 
 these copies 
 along $\CC^\times$ via the identification
 $z_2=z_1^{-1}$. Similarly, the projective plane $\CC P^2$
can be constructed by
taking three copies of $\CC^2$, with coordinates $(z_1,w_1)$,
$(z_2,w_2)$, and $(z_3,w_3)$, and  gluing them along
 $(\CC^\times)^2$ via the identifications
$$(z_2,w_2)=(z_1^{-1},w_1z_1^{-1})\quad \text{and}\quad  (z_3,w_3)=(z_1w_1^{-1},w_1^{-1}).$$

\begin{figure}[h]
\centering
\begin{tabular}{ccc}
\includegraphics[width=4cm, angle=0]{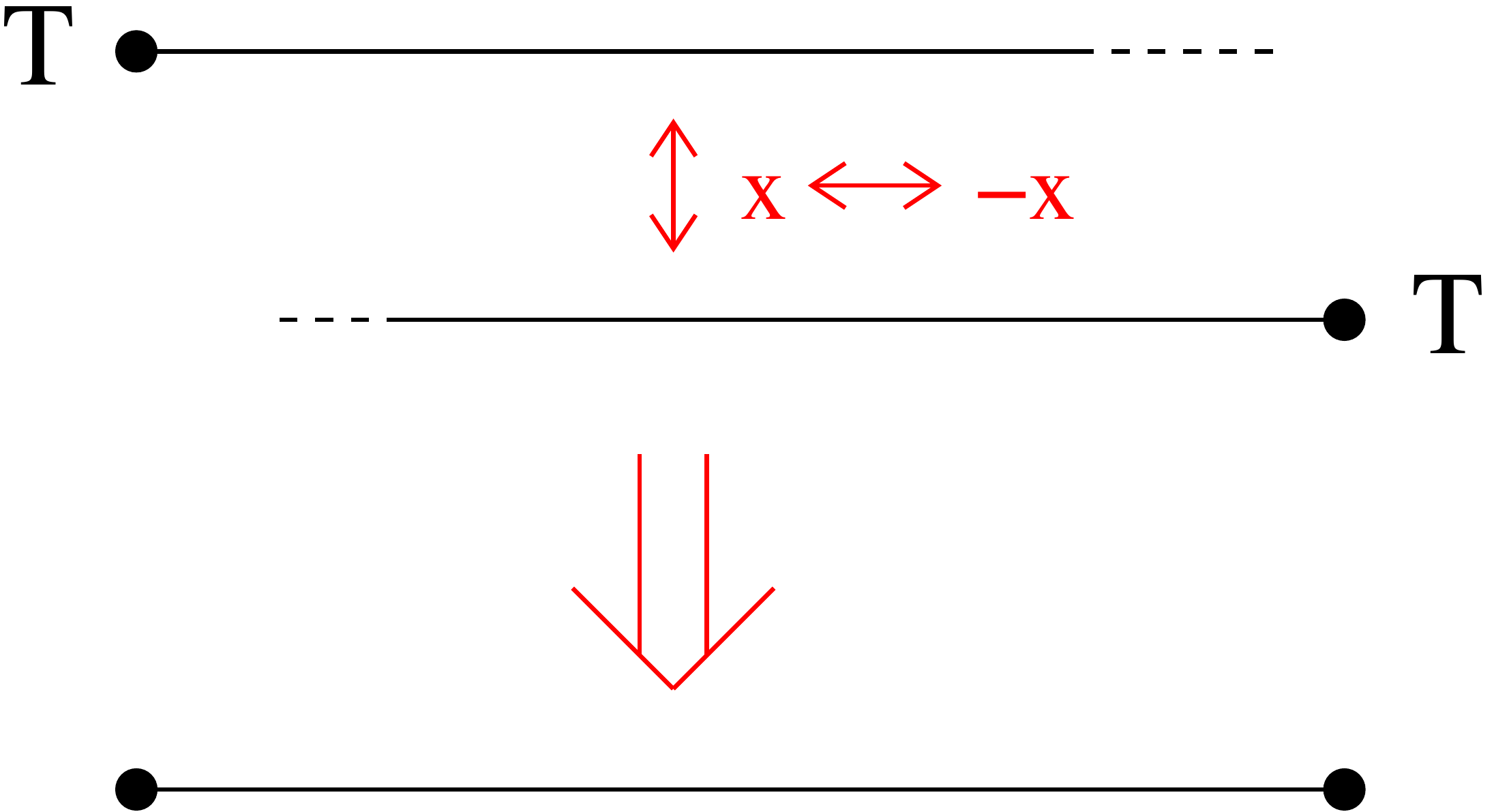}& \hspace{5ex} &
\includegraphics[width=6cm, angle=0]{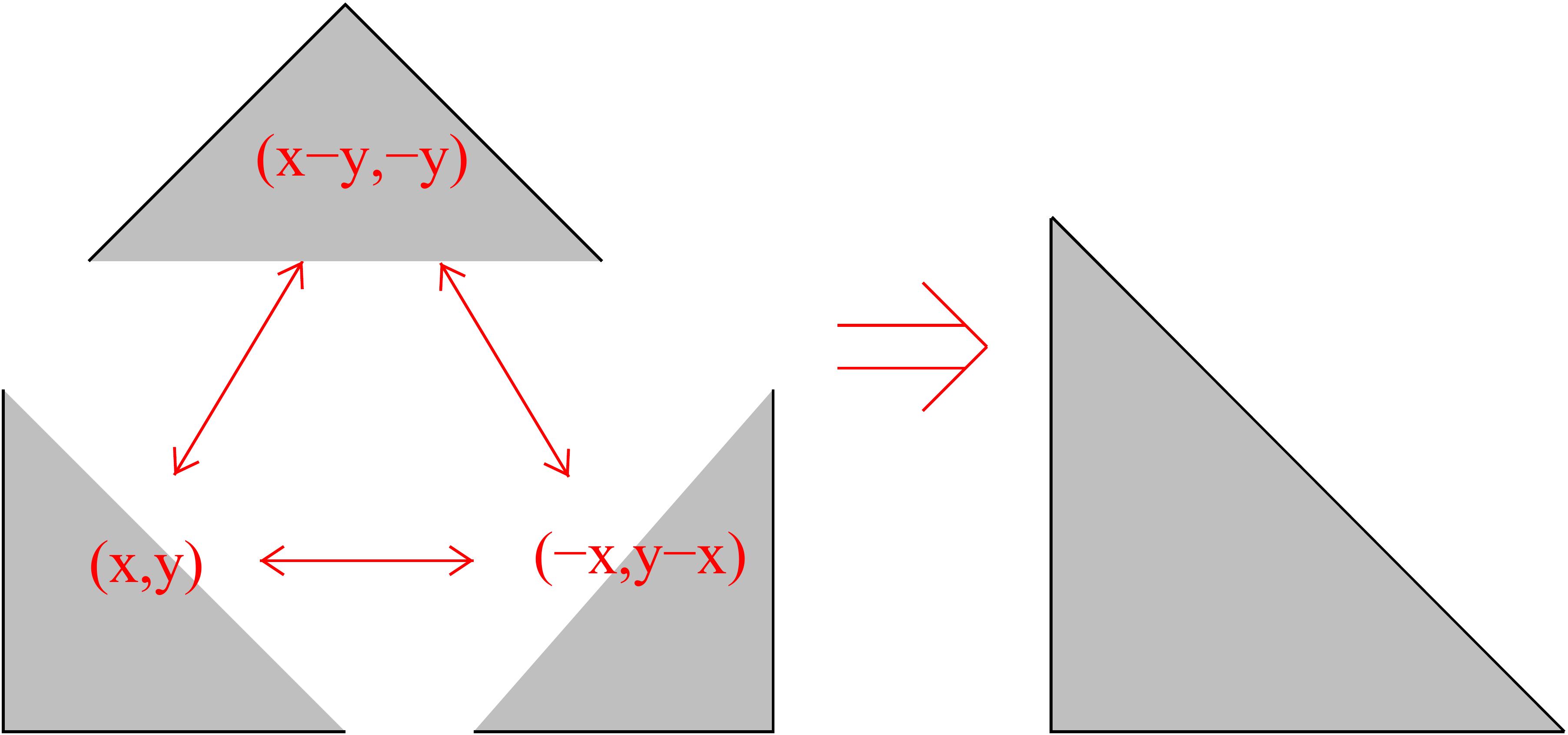}
\\a) $\TT P^1$&&
b)   $\TT P^2$
\\ \\
\includegraphics[width=3cm, angle=0]{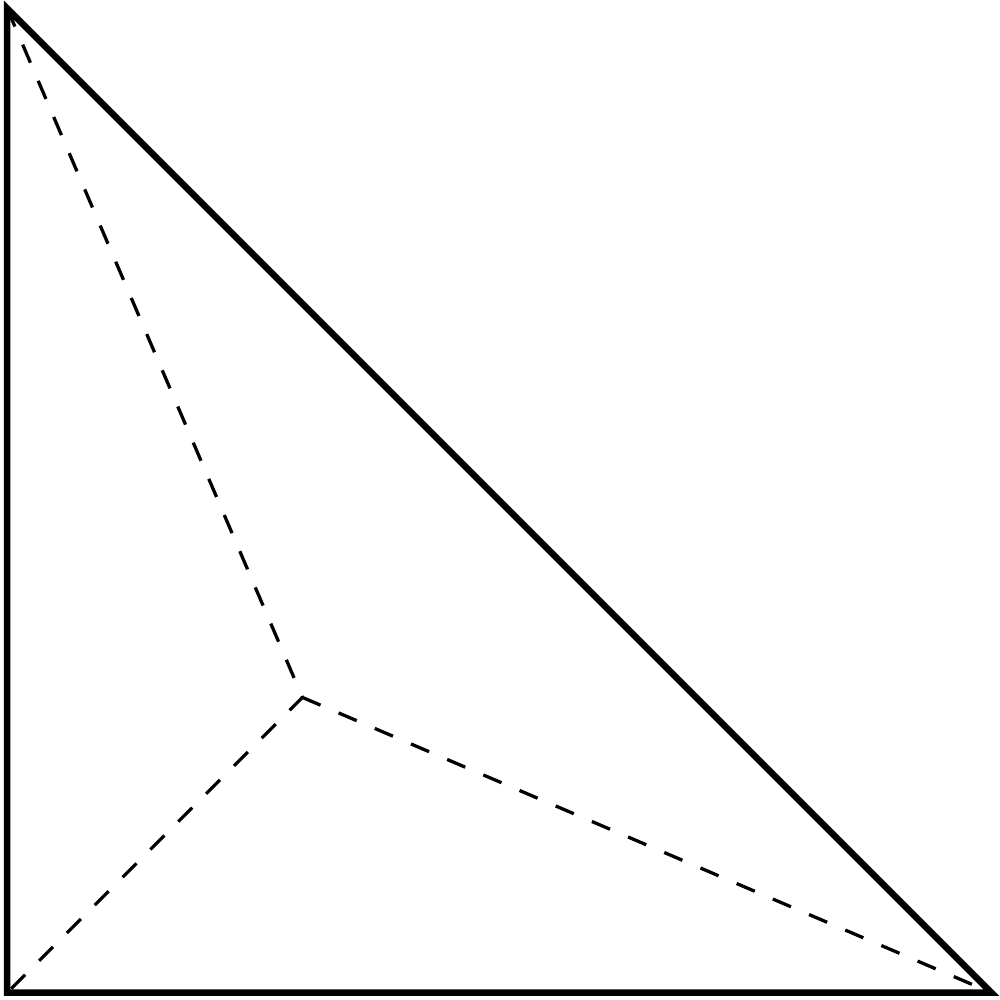}& &
\includegraphics[width=2.5cm, angle=0]{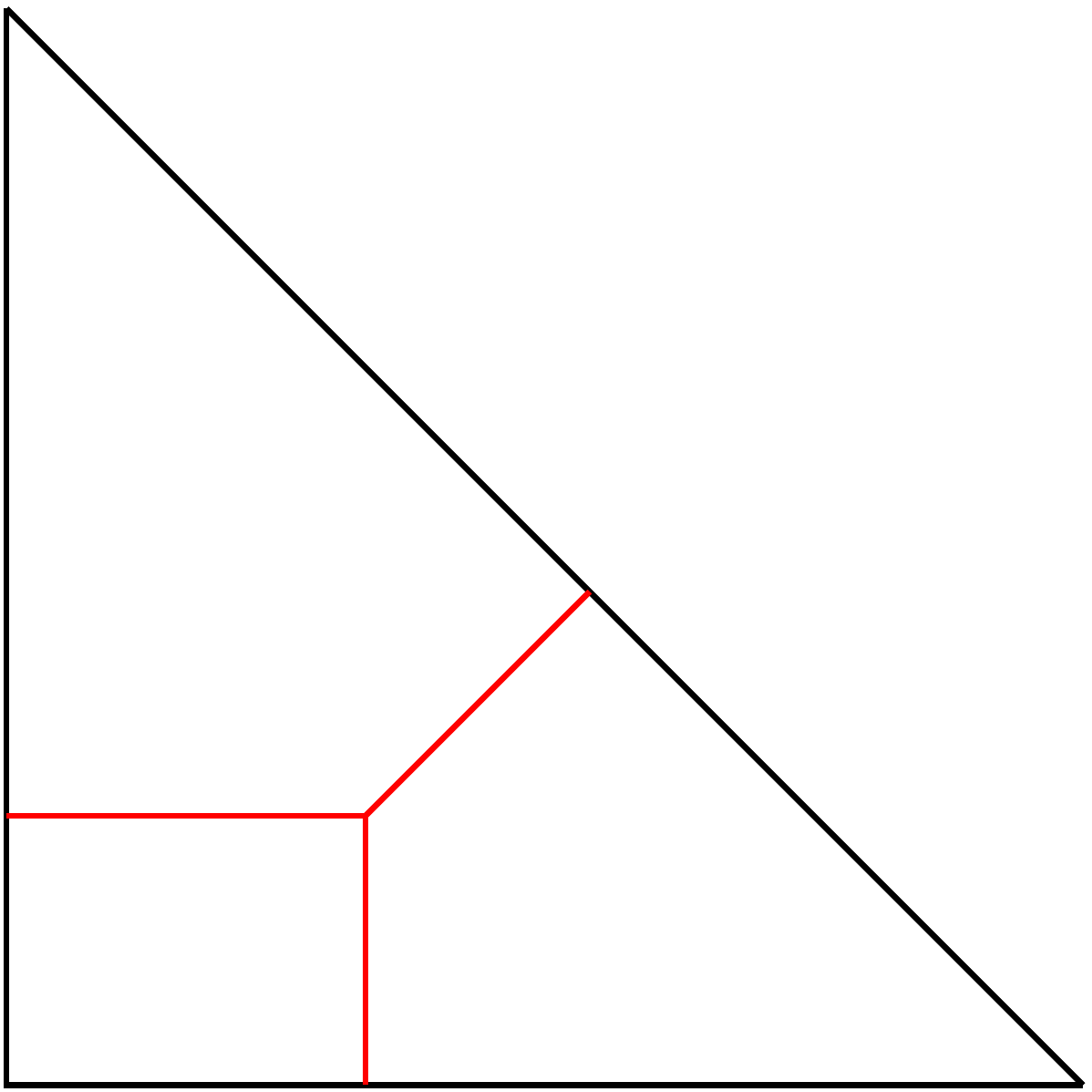}
\\
c) $\TT P^3$&&
d) a line in $\TT P^2$
\end{tabular}
\caption{Tropical projective spaces}\label{proj}
\end{figure}

Taking into account that $\tg x^{-1}\td =-x$, the above constructions 
over $\TT$
also yield
 the projective tropical line
 $\TT P^1$ and plane $\TT P^2$.  
In particular, we see that $\TT P^1$ is a segment
(Figure \ref{proj}a), 
and $\TT P^2$ is a triangle (Figure \ref{proj}b). 
More generally, the projective space $\TT
P^n$ is a simplex of dimension $n$, each of its faces corresponding to
a coordinate hyperplane. 
For example, the tropical 3-space $\TT P^3$ is a tetrahedron (Figure \ref{proj}c). 
Note that tropical toric varieties have more structure than just a bare topological space. 
Since all gluing maps are 
classical linear maps with integer coefficients,
each open face of dimension $q$ can be
identified with $\RR^q$ together with the lattice $\ZZ^q$ inside. 
As usual, the affine space 
$\TT^n$ 
embeds naturally into 
$\TT P^n$,
and any 
tropical 
subvariety of $\TT^n$
has a closure in 
$\TT P^n$. 
For example, we
depicted in Figure \ref{proj}d the closure in $\TT P^2$ of a tropical
line in 
$\TT^2$.

\vspace{0.5cm} 

\subsection{Tropical homology}\label{tropical_homology} 

This section delivers an introduction to tropical homology and cohomology. For a more in depth 
look the reader is referred to  \cite{IKMZ} and \cite{MikZhar:Eigenwave}. 
One of the main interests of tropical homology comes from its connection to the 
Hodge theory  of complex algebraic varieties. 
As in the case of tropical subvarieties of $\RR^n$,
some 
abstract tropical manifolds appear as a tropical limit 
of a family 
$({\mathcal X}_t)$ 
of
complex algebraic varieties of 
the same dimension.
It turns out that if the varieties ${\mathcal X}_t$ are projective and under certain additional conditions,
the Hodge numbers of 
${\mathcal X}_t$ are encoded in the tropical limit 
as ranks of the tropical $(p, q)$-homology groups, 
(see 
Theorem \ref{thm:homology} or
\cite{IKMZ}).

Let $M\subset\R^N\times\T^s\subset\T^{N+s}$ be
a tropical smooth local model,
and let  
$x \in M$ be a point. 
Recall the sedentarity $I(x)$ of $x$ from Section \ref{sec:mod}, and set $J(x) = \{1, 2, \dots, N+s\} \setminus I(x)$. 
Consider the stratum $\T^{J(x)} \subset \T^{N+s}$ and its relative interior $\R^{J(x)}.$
Let
$T_x(\R^{J(x)})$
be the vector space tangent to $\R^{J(x)}$ 
at $x$. 
For a face $E\subset M\cap\R^{J(x)}$ adjacent to 
$x$, 
denote by
 $T_x(E)\subset T_x(\R^{J(x)})$ 
the cone formed by the tangent vectors to $E$ 
that are directed towards $E$ from $x$.

\begin{definition}
The {\em tropical tangent space} 
$\FF_1(x)$ is the linear subspace of 
$T_x(\R^{J(x)})$
generated by $T_x(E)\subset T_x(\R^{J(x)})$
for all faces $E$ adjacent to $x$.
\end{definition} 

Let $X$ be a tropical 
manifold,
and let $x\in X$ be a point. 
Different charts at $x$ exhibit neighborhoods of $x$
as tropical smooth local models in various tropical
spaces $\T^N$ (perhaps even of different dimensions $N$). 
However, the differentials of overlapping maps at $x$ 
provide  canonical isomorphisms among the corresponding
tangent spaces $\FF_1(x)$. Thus, 
the tangent space $\FF_1(x)$ 
of $X$ at $x$ is well defined since it
does not depend on the choice
of a chart.
In addition to the tangent space we have
multitangent spaces $\FF_p(x)$ that are
spanned by multivectors tangent to the same face
adjacent to~$x$.
\begin{definition}
For any integer $p \geq 0$, 
the multitangent space $\FF_p(x) \subset \Lambda^p(T_x(\R^{J(x)}))$ 
of $X$ at $x$
is the linear subspace 
generated by the $p$-vectors 
of the type $\lambda_1\wedge\dots\wedge\lambda_p$,
where $\lambda_1,\dots,\lambda_p\in T_x(E)$ for a
face 
$E \subset X$
{\rm (}in a tropical smooth local model{\rm )} adjacent to $x$.
\end{definition} 

In particular, $\FF_0(x)=\R$. 
\begin{exa}
Let $X=\TT^n$, 
and let 
$x \in \TT^n$ be a point of sedentarity $I \subset \{1, 2, \ldots, n\}$.  
Then, $\FF_p(x)$ is isomorphic to $\Lambda^p(\R^{J(x)})$. 
\end{exa}

\begin{exa}\label{ex:hom line} 
Consider the tropical line in $\TT^2$ depicted in Figure \ref{fig:hom
  line}. The tangent spaces $\FF_1(x_5)$ and $\FF_1(x_6)$ are null.
The tangent space  $\FF_1(x_2)\subset\RR^2$ is generated by the
vector $(1,0)$, and so is 
isomorphic to $\RR$.
Analogously, the  tangent spaces  $\FF_1(x_3)$ and $\FF_1(x_4)$ are
isomorphic to $\RR$.
The tangent space $\FF_1(x_1)\subset\RR^2$ is generated by the vectors 
$$(-1,0),\quad (0,-1),\quad \mbox{and}\quad (1,1), $$
and so is isomorphic to $\RR^2$.

The multitangent space $\FF_2(x_2)\subset \Lambda^2(\RR^2)$ is 
equal to $\Lambda^2(\RR(1,0))=\{0\}$. 
Analogously, the  multitangent spaces  $\FF_2(x_3)$ and $\FF_2(x_4)$
are null. Since the multitangent space $\FF_2(x_1)\subset
\Lambda^2(\RR^2)$ 
is generated by $\FF_2(x_2)$,  $\FF_2(x_3)$, and $\FF_2(x_4)$, we
obtain that  $\FF_2(x_1)$ is null as well.
\begin{figure}[h]
\begin{center}
\includegraphics[scale=0.3]{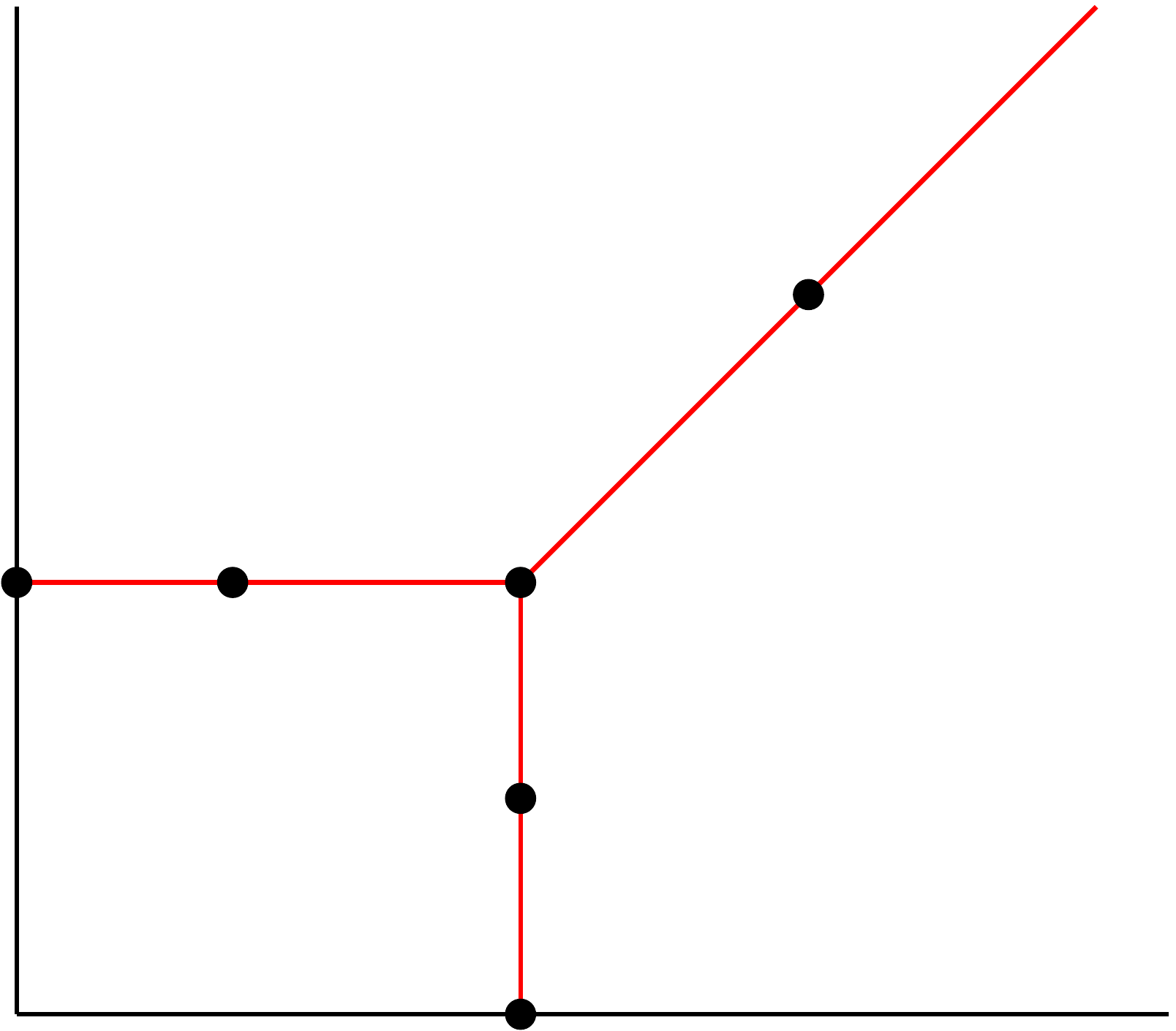}
\put(-72,25){$x_2$}
\put(-90, 45){$x_1$}
\put(-115, 60){$x_3$}
\put(-40, 80){$x_4$}
\put(-80, -10){$x_5$}
\put(-155, 50){$x_6$}
\end{center}
\caption{Multitangent spaces of an affine tropical line.}
\label{fig:hom line}
\end{figure}
\end{exa}

Note that a tropical smooth local model
$M\subset\T^N$ is
naturally stratified by maximal
linear subspaces of $\R^{J}\subset\T^J\subset\T^N$
locally contained in $M$.
Namely, for  $x\in M$ 
the union of the cones $T_x(E)$, for $E$ adjacent to $x$, 
is the fan $V(x) \subset T_x(\R^{J(x)})$ described in Section \ref{sec:global}.
{\em A priori}, there might be several maximal linear subspaces contained inside the fan $V(x)$.
Let $W_x\subset V(x)$ be the intersection of all such maximal subspaces.
The intersection $W_x$ itself is a linear subspace of $T_x(\R^{J(x)})$.
Note that $x+W_x$ is locally (near $x$) contained in $M$
by the balancing condition.
\begin{defn}[cf. \cite{MikZhar:Eigenwave}]
A {\em tropical stratum} is a subset $E\subset X$
that locally looks like $x+W_x$ in charts.
\end{defn}
Any tropical stratum $E$ is 
a differentiable manifold by construction. 
Tropical strata of $X$ are partially ordered:
we say that $E$ is {\em dominated} by $E'$
if the topological closure $\overline{E'}$
of $E'$
in $X$ contains $E$. 
Since
an $n$-dimensional tropical manifold
 $X$ 
admits a structure of
 (a pure) $n$-dimensional
polyhedral complex, any 
tropical 
stratum
is dominated by a $n$-dimensional 
tropical 
stratum.
These strata are called {\em open facets} of $X$.
\begin{exa}
Let $X$ be a non-singular tropical hypersurface of $\RR^n$. 
Tropical duality provides a correspondence between tropical strata of
dimension $k$ of $X$ and cells of dimension $n-k$ of the subdivision of
$\Delta(X)$ dual to $X$. In particular, open facets of $X$ correspond
to edges of this subdivision. 
\end{exa}
\begin{exa}
If $x$ is a point in an open facet 
of $X$, 
then 
$\FF_p(x)$ is isomorphic to 
$\Lambda^p(\R^{\dim(X)})$.
Indeed, in 
this case we have $W_x=T_x(\R^{J(x)})$
 which is isomorphic to 
$\R^{\dim(X)}$.
\end{exa}

Notice that if 
$E$ is a tropical stratum of $X$,
any path $\gamma: [0, 1] \to E$ gives rise to an identification
(by parallel translation) 
of the multitangent spaces $\FF_p(\gamma(t))$, $t \in [0, 1]$, 
for each $p$. 
Let $E$ be a tropical stratum of $X$, and let $E'$ be a tropical stratum
such that $E$ is dominated by $E'$.
Consider a path $\gamma: [0, 1] \to (E' \cup E)$ such that $\gamma([0, 1)) \subset E'$
and $\gamma(1) \in E$.
Put $x = \gamma(0)$ and $y = \gamma(1)$. 
For each $p$, the identifications by parallel translations of 
 $\FF_p(\gamma(t))$,  for $t \in [0, 1)$,
extend to a map 
\begin{equation}
\label{rho-gamma}
\rho_p(\gamma): \FF_p(x) \to \FF_p(y).  
\end{equation} 
If the points $x$ and $y$ have the same sedentarity, then
$\rho_p(\gamma): \FF_p(x) \to \FF_p(y)$ is a monomorphism,
since 
by definition the tangent space $\FF_1(\gamma(t))$
 is a subspace of $\FF_1(y)$ when $\gamma(t)$ is contained in a
 tropical smooth local model of $X$ at $y$.
If $I(x) \ne I(y)$, then 
since $E'$ dominates $E$ we have $I(x) \subset I(y)$. 
Then $\rho_p(\gamma): \FF_p(x) \to \FF_p(y)$
is given by 
the composition
of the projection along the directions indexed by $I(y) \setminus
I(x)$, 
with 
the same 
monomorphism as in the case of equal 
sedentarity. 
We obtain immediately the following statement. 

\begin{prop}
\label{prop-cosheaf}
Let 
$E'\ni x$ and $E\ni y$ 
be tropical strata of $X$ such that $E$ is dominated
by $E'$. 
Consider two paths $\gamma_1$ and $\gamma_2$ 
to $E\cup E'$
such that
\begin{itemize}
\item $\gamma_1(0) = \gamma_2(0) = x$ and $\gamma_1(1) = \gamma_2(1) = y$,
\item $\gamma_1([0,1)) \subset E'$ and $\gamma_2([0,1)) \subset E'$,
\item $\gamma_1$ and $\gamma_2$ are homotopic
among paths satisfying the two above conditions.
\end{itemize}
Then, for each $p$, the maps $\rho_p(\gamma_1)$ and $\rho_p(\gamma_2)$ coincide.
\end{prop} 

\begin{exa}
Let $I_1\subset I_2\subset\{1,\ldots,n\}$. Then, the tropical 
stratum $E=\RR^{I_2}$ is dominated by the 
tropical 
stratum $E'=\RR^{I_1}$ inside the tropical manifold $\T^n$. The map 
$$\Lambda^p(\R^{n-|I_1|}) \to\Lambda^p(\R^{n-|I_2|}) $$
induced by any path as above is 
given 
by the projection 
$\R^{n-|I_1| }\to\R^{n-|I_2|}$ 
 along the directions indexed by $I_2\setminus I_1$.
\end{exa}
\begin{exa}\label{ex:hom line2}
For the tropical line of  Example \ref{ex:hom line}, the two morphisms
$\FF_1(x_2) \to \FF_1(x_5)$ and $\FF_1(x_3) \to \FF_1(x_6)$ are null.
The tangent space $\FF_1(x_1)$ is naturally identified with $\RR^2$.
We identify the tangent spaces $\FF_1(x_2)$, $\FF_1(x_3)$, and 
$\FF_1(x_4)$
with 
$\RR$ by respectively choosing
$(0,1)$, $(1,0)$, and $(-1,-1)$ as directing vectors. Then, 
paths from $x_2$, 
(respectively, $x_3$ or $x_4$) to $x_1$ induce the morphisms
$$\begin{array}{cccc}
\rho_2: &\FF_1(x_2) & \longrightarrow & \FF_1(x_1)
\\ &x &\longmapsto & (0,x)
\end{array}\qquad \begin{array}{cccc}
\rho_3:&\FF_1(x_3) & \longrightarrow & \FF_1(x_1)
\\ &x &\longmapsto & (x,0)
\end{array}$$
and
$$ \begin{array}{cccc}
\rho_4:&\FF_1(x_4) & \longrightarrow & \FF_1(x_1)
\\& x &\longmapsto & (-x,-x)
\end{array}.$$
Notice that 
$$\rho_2(1)+ \rho_3(1)+\rho_4(1)=0 \in\FF_1(x_1).$$
\end{exa} 

\medskip
The multitangent 
spaces
$\FF_p(x)$ can be
thought of as a kind of coefficient system
parameterized by points of $X$.
 However, this coefficient system is not locally constant
as its value may jump at smaller-dimensional
strata.
But, as we saw, 
the coefficient groups $\FF_p(x)$ at different points $x$ 
are related 
via  
the
morphisms
$\rho_p$ described above.
Such coefficient systems are known as
{\em constructible cosheaves}, and can
be used as coefficients for
homology groups of $X$. 

Let us review, for example, the construction 
from \cite{IKMZ,MikZhar:Eigenwave}
of singular homology, 
with coefficients
in $\FF_p$, 
of a tropical manifold $X$. 
Given 
a closed standard $q$-dimensional simplex $\sigma_q$, 
we say that a singular simplex
$f:\sigma_q\to X$
is {\em compatible with tropical stratification}
if for each open face $\sigma'\subset\sigma_q$
there exists a tropical stratum $E_{\sigma'}$
such that $f(\sigma')\subset E_{\sigma'}$. 
The singular simplex 
$f:\sigma_q \to X$ 
is then equipped with a coefficient 
$\phi(f)$ 
in
$\mathcal{F}_p(f(c))$, where $c$ is the 
barycenter of $\sigma_q$. 
The term 
$\phi(f) \cdot f$ 
is called a \emph{$(p,q)$-cell}, 
and the 
coefficient 
$\phi(f)$ 
is 
often referred to as the {\it framing}
of the $(p,q)$-cell
$\phi(f) \cdot f$. 

We define the {\em tropical chain group}
$C_{p,q}(X)$ 
 to be the direct sum
of 
$\FF_p(f(c))$ 
over all $q$-dimensional
singular simplicies $f$ which are compatible with the tropical
stratification
of $X$. 
Note that  
due to Proposition 
\ref{prop-cosheaf}
 we have
 a well-defined boundary map
$$\dd:C_{p,q}(X)\to C_{p,q-1}(X).$$
This map is the usual simplicial boundary map 
along with the restrictions of the coefficients 
given by the maps $\rho$ described above.
It is well defined since the homotopy class of the path
from the barycenter of a simplex 
to the 
barycenter
of any of its faces
is unique.
Furthermore, we have $\dd\circ\dd=0$
by the same argument as in the case
of singular homology groups with 
constant coefficients.

\begin{definition}[Tropical homology, 
cf. \cite{IKMZ}]
The {\it tropical homology 
group}   
$H_{p,q}(X)$ of $X$ is the $q$-th homology group 
of the complex
$$
\ldots \to C_{p, q + 1}(X) \to C_{p, q}(X) \to C_{p, q - 1}(X) \to \ldots 
$$  
\end{definition} 
\begin{rem}\label{rem:Z hom}
The differentials of overlapping maps in the atlas of a tropical
manifold are integer linear maps. In particular, each tangent space 
$\FF_1(x)$ contains a full rank lattice, 
and one could 
consider integer multitangent spaces $\ZZ\FF_p(x)\subset
\Lambda^p(T_x(\ZZ^{J(x)}))$ instead of the  multitangent spaces
$\FF_p(x)$ as above. The corresponding homology groups of $X$ are
called \emph{tropical integer homology groups} of $X$.
\end{rem}

\subsection{Some examples of homology computations}
\begin{exa}\label{ex:hom q=0}
Since $\FF_0$ is locally constant with stalk $\RR$,
  any tropical manifold $X$ satisfies
$$H_{0,q}(X)=H_q(X;\R).$$
\end{exa}
\begin{exa}\label{ex:hom torus}
The contraction $\phi(f)\cdot tf$ 
by a factor $t\in [0,1]$ of a $(p,q)$-cell $\phi(f)\cdot f$ in $\RR^n$ 
is again a $(p,q)$-cell in $\RR^n$. Moreover, the boundary map clearly
commutes with the contraction.
As a consequence, any
$(p,q)$-cycle in $\RR^n$ is homologous to a $(p,q)$-cycle supported at
the origin of $\RR^n$, and so
$$H_{p,0}(\RR^n)=\Lambda^p(\RR^n)\quad
\mbox{and}\quad
H_{p,q}(\RR^n)=0 \
\mbox{for any} \ q\ge 1.$$ 
\end{exa}
\begin{exa}\label{ex:hom affine}
Given $t\in \TT$, consider the map
$$\begin{array}{ccccl}
\tau_t:&\TT^n&\longrightarrow&\TT^n&
\\ & (x_1,\ldots,x_n)&\longmapsto & \tg t(x_1,\ldots,x_n)\td &=(x_1+t,\ldots,x_n+t)
\end{array} .$$
The boundary map commutes with $\tau_t$, and 
$\phi(f)\cdot (\tau_t\circ f)$ is
a $(p,q)$-cell in $\TT^n$ for any $t\in\TT$ and any $(p,q)$-cell
$\phi(f)\cdot f$ in $\TT^n$. Hence, any
$(p,q)$-cycle in $\TT^n$ is homologous to a $(p,q)$-cycle supported at
$(-\infty,\ldots,-\infty)$, and so
$$H_{0,0}(\TT^n)=\RR\quad
\mbox{and}\quad
H_{p,q}(\TT^n)=0 \quad \mbox{if } p+q\ge 1.$$ 
\end{exa}
\begin{exa}\label{ex:hom line3}
We compute the tropical homology of the affine tropical line $L$ of
Example \ref{ex:hom line}. By Example \ref{ex:hom q=0} we have
$$H_{0,0}(L)=\RR \quad \mbox{and}\quad H_{0,1}(L)=0.$$
We use the same identifications of $\FF_1(x_2)$, $\FF_1(x_3)$, and
$\FF_1(x_4)$ with $\RR$ as in Example \ref{ex:hom line2}.
The group $H_{1,0}(L)$ is clearly generated by $1\cdot x_2$, $1\cdot x_3$,
and $1\cdot x_4$. The obvious path from $x_2$ to $x_5$ equipped with the
framing $1$ gives $1\cdot x_2=0$ in $H_{1,0}(L)$. Analogously, we have 
$1\cdot x_3=0$ in $H_{1,0}(L)$. It follows from Example \ref{ex:hom
  line2} that
the obvious path from $x_4$ to $x_1$
equipped with the 
framing $1$ gives $1\cdot x_4=-1\cdot x_2-1\cdot x_3$ in
$H_{1,0}(L)$. Altogether we obtain
$$H_{1,0}(L)=0. $$
Analogously, one sees that any $(1,1)$-cycle is homologous to a cycle
with support disjoint from  the open edge containing $x_4$, implying
that
$$H_{1,1}(L)=0. $$

\end{exa}
\begin{exa}\label{ex:puncturedLineHomo}
Let us consider the  tropical line $L'$ in $\RR^2$ obtained by
removing the vertices $x_5$ and $x_6$ from the affine line of
Example \ref{ex:hom line}. The same computations as in Example
\ref{ex:hom line3} 
give 
$$H_{0,0}(L')=\RR, \qquad H_{0,1}(L')=H_{1,1}(L')=0,$$
 and
$$H_{1,0}(L')=\RR(1\cdot x_2)\oplus \RR(1\cdot x_3)\oplus
\RR(1\cdot x_4)/\RR
 v\simeq\RR^2, $$
where $v=1\cdot x_2+1\cdot x_3+1\cdot x_4$.
\end{exa}
\begin{exa}\label{ex:hom TP1}
Since $\TT P^1$ is contractible as a topological space, we have
$$H_{0,0}(\TT P^1)=\RR \quad \mbox{and}\quad H_{0,1}(\TT P^1)=0. $$
Moreover as in Example \ref{ex:hom line3}, any $(1,0)$-cell is a
boundary, and so
$$H_{1,0}(\TT P^1)=0. $$
Let $\phi$ be a non-zero element of $\FF_1(\RR)$.
The $(1,1)$-cell $\phi\cdot \TT P^1$ is clearly 
 a $(1,1)$-cycle, and any  $(1,1)$-cycle is a multiple of 
 $\phi\cdot \TT P^1$. Thus we have
$$H_{1,1}(\TT P^1)=\RR(\phi\cdot \TT P^1)\simeq\RR.$$
\end{exa}
\begin{exa}\label{ex:hom TP2}
Consider a point $x\in\RR^2\subset\TT P^2$
and 
a simplicial subdivision 
of $\TT P^2$ into three triangles $T_1$, $T_2$, and $T_3$, induced by
the point $x$ and the three vertices of $\TT P^2$.
Suppose that $\phi_1$, $\phi_2$, and $\phi_3$ are three
elements of $\FF_p(\RR^2)$ such that at least one of them is not zero.
Note that
$$\partial\left(\phi_1\cdot T_1+\phi_2\cdot T_2+ \phi_3\cdot T_3\right)=0
\Leftrightarrow \phi_1=\phi_2=\phi_3\mbox{ and }p=2,$$
from which we deduce that
$$H_{2,2}(\TT P^2)=\Lambda^2(\RR^2)\simeq \RR \quad \mbox{and}\quad
H_{0,2}(\TT P^2)=H_{1,2}(\TT P^2)=0. $$ 
As in the classical case, there are homogeneous coordinates $[x:y:z]$
on the tropical projective plane $\TT P^2$.
Any $(p,q)$-cycle in $\TT P^2$ with $p<2$ is homologous to a
$(p,q)$-cycle in $\TT P^2$ whose support does not contain the point
$[-\infty:-\infty:0]$. As in 
Examples \ref{ex:hom torus} and 
\ref{ex:hom affine}, the maps 
$$\begin{array}{ccc}
\TT P^2\setminus\{[-\infty:-\infty:0]\}&\longrightarrow&\TT
P^2\setminus\{[-\infty:-\infty:0]\} 
\\ \ [x:y:z] &\longmapsto & [x:y:z+t]
\end{array}  $$
can be used to show that any $(p,q)$-cycle in $\TT P^2\setminus\{[-\infty:-\infty:0]\}$  is
homologous to a  $(p,q)$-cycle with support contained in $\{z=-\infty\}$. Since
this latter is a tropical projective line, by Example \ref{ex:hom TP1}
we have
$$H_{0,0}(\TT P^2)\simeq H_{1,1}(\TT P^2)\simeq \RR $$
(one sees easily
that the generators of these groups remain non-homologous to zero
when we pass to $\TT P^2$)
and
$$ H_{0,1}(\TT P^2)=H_{1,0}(\TT P^2)=H_{2,0}(\TT P^2)=H_{2,1}(\TT P^2)=0. $$
\end{exa}

\subsection{Straight cycles}\label{exa:fundamental}

Let $X$ be a tropical manifold, and suppose that $Z\subset X$ is 
a finite {\em weighted balanced polyhedral subcomplex
of 
dimension $p$}
in $X$. 
This means that $Z$ is a
weighted balanced polyhedral subcomplex 
of pure dimension $p$
of $X$ in each chart of $X$, and the weight 
of a point in a facet does not depend
on a choice of the chart.
A choice of an orientation on 
a facet $F$ of $Z$ 
produces a canonical framing in 
$\ZZ \FF_p(X)$,
equal to the weight of 
$F$
multiplied by the primitive
integer element of 
$\Lambda^p(\FF_1(F))$
that agrees
with the chosen orientation of $F$. (Recall that
$\FF_1(F)$ comes with a lattice of full rank, see Remark \ref{rem:Z hom}.)

We form the {\em fundamental 
class}
 $[Z]$ of $Z$
as follows. Choose an orientation for each facet of $Z$,
subdivide each facet into
        singular $p$-simplices, and enhance each $p$-simplex
with the corresponding canonical framing in 
$\ZZ \FF_p(X)$. 
The sum of all these $(p,p)$-cells
 is an integer 
$(p,p)$-chain $[Z]$.
The balancing condition is equivalent
to the condition $\dd [Z]=0$,
   {\it i.e.} 
   $[Z]$ is a cycle (compare with Example \ref{ex:hom
     line2}). 
Note that the choice of orientation 
of 
facets of $Z$, 
and their subdivision into singular $p$-simplices  is irrelevant for the
resulting cycle class $[Z]$ in $H_{p,p}(X)$. 
 This is because the orientation
of each facet is used in the construction
of $[Z]$ twice: once in the orientation
of the underlying facet, and
once in the choice of the $p$-framing.

\begin{defi}[cf. \cite{MikZhar:Eigenwave}]
The $(p, p)$-cycles 
that can be presented
as fundamental class $[Z]$ for some weighted balanced subcomplex
$Z\subset X$ are called {\em straight
cycles}.
\end{defi}
Straight cycles
generate a subspace in $H_{p,p}(X)$
that behaves in a semicontinuous way with respect
to deformation of the tropical structure of $X$. 
The 
problem 
of 
detecting whether a cycle in $H_{p,p}(X)$
is straight 
is a very interesting question
related to the 
famous Hodge conjecture.

\subsection{Tropical cohomology}
There is also a dual theory
of tropical cohomology. We may consider
the spaces $$\FF^p(x)=\Hom(\FF_p(x);\R)$$
together with the morphisms 
\begin{equation}
\label{co-rho-gamma}
\FF^p(y)\to\
\FF^p(x) 
\end{equation} 
dual to morphisms (\ref{rho-gamma}).  
The coefficient system $\FF^p$ forms 
a {\em constructible sheaf} and we may
take cohomology groups with coefficients in $\FF^p$.
Namely, we form the cochain groups 
$$C^{p,q}(X)=\Hom(C_{p,q}(X),\R).$$ 
An element of this group can be interpreted
as a functional associating an element of
$\FF^{p}(f)=\Hom(\FF_p(f(c)),\R)$
to each compatible singular $q$-simplex 
$f:\sigma_q\to X$, where $c$ is the barycenter of the $q$-dimensional simplex $\sigma_q$.
Accordingly, we have a coboundary map
$$
\delta:C^{p,q}(X)\to C^{p,q+1}(X)
$$
with $\delta\circ\delta=0$.
\begin{defi}[Tropical cohomology, cf. \cite{IKMZ}]
The {\it tropical cohomology group} 
$H^{p,q}(X)$ is the $q$-th cohomology group 
of the complex 
$$
\ldots \to C^{p, q - 1}(X) \to C^{p, q}(X) \to C^{p, q + 1} (X) \to \ldots 
$$ 
 
\end{defi} 
As usual, cohomology groups admit the cup product. 
It is based on the following 
observation. Suppose that $f_1$ and $f_2$ 
are two
singular 
simplices  of dimensions $q_1$ and $q_2$, respectively, 
which are faces of a 
compatible
singular simplex $f$. Then, we have 
the composed homomorphism 
$$\FF^{p_1}(f_1)\otimes\FF^{p_2}(f_2)\to\FF^{p_1}(f)\otimes\FF^{p_2}(f)
\to\FF^{p_1+p_2}(f).$$
The first homomorphism in this composition
is obtained with the help of \eqref{co-rho-gamma}, 
while the second homomorphism is given
by the exterior product.
This product in coefficients descends to the cup product
$$H^{p_1,q_1}(X)\otimes H^{p_2,q_2}(X)\to H^{p_1+p_2,q_1+q_2}(X)$$
with the usual super-commutativity property
$$\alpha\cup\beta= (-1)^{p_1p_2+q_1q_2}\beta\cup\alpha$$
for $\alpha\in H^{p_1,q_1}(X)$ and $\beta\in H^{p_2,q_2}(X)$. 

\begin{exa}\label{ex:cohom q=0}
As for tropical homology, we have
$$H^{0,q}=H^q(X;\RR). $$
\end{exa}
\begin{exa}\label{ex:cohom contr}
As in Examples \ref{ex:hom torus} and \ref{ex:hom affine} we compute
$$H^{p,0}(\RR^n)=\Lambda^p(\RR^n) \quad\mbox{and}\quad
H^{p,q}(\RR^n)=0 
\ \mbox{for any} \ q\ge 1,$$ 
and
$$H^{0,0}(\TT^n)=\RR \quad\mbox{and}\quad
H^{p,q}(\TT^n)=0 \quad \mbox{if } p+q\ge 1 .$$
\end{exa}
\begin{exa}
Let $L$ be the tropical line from Example \ref{ex:hom line}.
Denote by $e_i$, $i = 1$, $\ldots$, $5$, the edge of $L$ with vertices $x_i$ 
and $x_1$ oriented toward $x_1$.  We use the same identifications of
$\FF_1(x_2)$, $\FF_1(x_3)$, and 
$\FF_1(x_4)$ with $\RR$ as in Example \ref{ex:hom line2}.
Given 
$\Phi_p \in C^{p, 0}(X)$, 
we have
$$\delta\Phi_0(e_i)(1)= \Phi_0(x_1)(1)-\Phi_0(x_i)(1),  $$
and
$$\delta\Phi_1(e_5)(1)= \Phi_1(x_1)(0,1),\quad
 \delta\Phi_1(e_6)(1)= \Phi_1(x_1)(1,0),$$
$$\delta\Phi_1(e_4)(1)= \Phi_1(x_1)(-1,-1) -\Phi(x_4)(1). $$
 From this we deduce 
$$H^{0,0}(L)=\RR \quad \mbox{and}\quad H^{0,1}(L)=H^{1,0}(L)=H^{1,1}(L)=0 .$$
Whereas the cohomology groups of  $L' = L \cap \R^2$ are
$$H^{0,0}(L')=\RR,\quad H^{1,0}(L')=\RR^2, \quad \mbox{and} \quad
H^{0,1}(L')=H^{1,1}(L')=0.$$
The homology groups of $L$ and $L'$ are calculated 
in Examples \ref{ex:hom line3} and \ref{ex:puncturedLineHomo}, respectively.  
Next example generalizes the computation of the groups $H^{p,q}(L')$
for any tropical limit  of  a linear space.

\end{exa}
\begin{exa}
Let 
a fan  $L \subset \R^n$ be the tropical limit of a linear space
$\mathcal{\L} \subset (\CC^\times)^n$ 
(see Section \ref{sec:fanlinear}).
As in Example \ref{ex:cohom contr}, we have
$$H^{p, 0}(L) = \mathcal{F}^p(x) \quad\mbox{and}\quad H^{p, q}(L) =0
\ \mbox{for any} \ q\ge 1,$$ 
where $x$ is the vertex of the fan 
$L$. 
Along with the cup product, 
$H^\bullet(L)$
is isomorphic to the cohomology ring  of $\mathcal{\L}$  which is known as an Orlik-Solomon algebra, cf.  \cite{Zharkov:Orlik}. 
\end{exa}

As seen in the  last example, 
it turns out that tropical cohomology
groups capture the cohomology groups of complex varieties 
$\X_t \subset\cp^n$ in the case when 
our tropical 
manifold $X\subset\tp^n$
comes as the tropical limit of 
a family 
$(\X_t)_{t \in U}$. 
Here for a parameterizing set $U$
we take a subset $U\subset \CC$ such that $\CC\setminus U$ is
bounded, {\it i.e.} $U$ is
        a punctured neighborhood of $\infty$
        in $\cp^1=\CC\cup\{\infty\}$.
The definition of tropical limit that we use here
is a projective version of Definition \ref{def:approx}.
This means that
$X$ is the limit of the subsets
$\Log_{|t|} \X_t\subset\tp^n$
(called {\em amoebas} of $\X_t$),
    where $\Log(z_0:\dots:z_N)=
    (\log|z_0|:\dots:\log|z_N|)$.
As in Definition \ref{def:approx} the tropical limit
should be enhanced with weights. 
Given a tropical manifold $X$, we define
$$h^{p,q}(X)=\dim H^{p,q}(X). $$
\begin{thm}[\cite{IKMZ}]\label{thm:homology}
Let $X$ be a tropical submanifold of $\TT P^n$ {\rm (}in particular, 
the weights 
of facets of $X$ are all equal to $1${\rm )}. 
If $X$ 
is the tropical limit of 
a complex analytic family 
$(\X_t)_{t \in U}$
of projective varieties, then for sufficiently
large $|t|$ the complex variety $\X_t\subset\cp^n$ is smooth.
Furthermore, we have 
the following
relation between the Hodge numbers of $\X_t$
{\rm (}for sufficiently large $|t|${\rm )}  
and the 
dimensions of tropical 
cohomology groups 
of $X$:
$$h^{p,q}(\X_t)=h^{p,q}(X).$$ 
\end{thm} 
As a consequence, if an $n$-dimensional tropical submanifold $X$ of $\TT P^n$ is
 the tropical limit of 
a complex analytic family $\X_t$ of projective varieties, 
then $X$ must satisfy
$$ h^{p,q}(X)=h^{q,p}(X)=h^{n-p,n-q}(X) 
\ \mbox{for any} \ p,q \ge 0.$$ 

\subsection{Cohomology of tropical curves}
Let $\Gamma$ be a connected compact tropical curve 
(see Section \ref{sec:abstract curves}).   
According to Example \ref{ex:cohom q=0}, 
 we have 
 $H^{0,0}(\Gamma) =  \R$ and
$H^{0,1}(\Gamma)\simeq\R^g$,
where $g=b_1(\Gamma)$ is the 
genus of 
$\Gamma$. 

Consider the space $H^{1,0}(\Gamma)=H^0(\Gamma;\FF^1)$. 
There is a map 
$C^{1,0}(\Gamma) \to C_{0, 1}(\Gamma)$ which descends to an isomorphism between 
the cohomology and homology groups. 
To describe this map, first fix an orientation of every edge of $\Gamma$. 
For an oriented edge  
$E$ of $\Gamma$, 
choose a point $x$ inside $E$ 
and enhance $x$ 
with the unit tangent vector in the direction 
of $E$ to produce a  $0$-chain with coefficients in $\FF_1(x)$. Denote this chain by $\tau_E \in C_{1, 0}(\Gamma)$. 
A 
cochain  
$\beta\in C^{1,0}(\Gamma)$ can be evaluated at $\tau_E$, 
and we set
$$Z_\beta=\sum\limits_{E \subset \text{Edge}(\Gamma)} \beta(\tau_E) E.$$
The condition that $\beta$ is a cocycle implies that $Z_\beta$
is a cycle, 
and thus $\dim H^{1,0}(\Gamma)=\dim H_{0,1}(\Gamma)=g$.

Finally, we compute the group $H^{1, 1}(\Gamma) = \R$ directly from the sequence,
$$C^{1,0}(\Gamma) \to C^{1, 1}(\Gamma) \to 0,$$
along with the observation that it  suffices to take $C^{1, 0}(\Gamma)$ and $C^{1, 1}(\Gamma)$ as simplicial 
cochains. 
Therefore, these groups are finite dimensional and we have explicitly, 
$$C^{1, 0}(\Gamma) = \bigoplus_{v \in \mbox{Vert}(\Gamma)} \RR^{val(v) -1}
= \RR^{2E - V}
\quad \mbox{and} \quad C^{1, 1}(\Gamma) = \RR^E,
$$
where $\mbox{Vert}(\Gamma)$ is the set of vertices of $\Gamma$, and $V$ (respectively, $E$) 
denotes the number of vertices (respectively, edges) of $\Gamma$.  
The kernel of the map $C^{1,0}(\Gamma) \to C^{1, 1}(\Gamma)$ is $H^{1, 0}(\Gamma) \simeq \RR^g$ by our previous computation. Thus 
$H^{1, 1}(\Gamma) \simeq \RR$. 
Notice also that  a $(1, 1)$-cochain can be evaluated on
the fundamental $(1, 1)$-cycle $[\Gamma]$ described in 
Section
\ref{exa:fundamental}. 
One sees easily 
that the image of $\delta$ is contained in subspace $\{\phi \ | \  \phi([\Gamma]) = 0 \} \subseteq C^{1, 1}(\Gamma)$. Moreover, 
both subspaces are of codimension one in $C^{1, 1}(\Gamma)$, so they are equal. 
Therefore, we obtain a non-degenerate pairing between the cohomology and homology groups 
$H^{1, 1}(\Gamma) \times H_{1, 1}(\Gamma) \to \RR$. 

We have determined the tropical cohomology groups of the curve $\Gamma$. 
Their 
dimensions 
can be arranged in a diamond shape which is exactly the same as the Hodge diamond of a
Riemann surface of genus $g$:

$$\left.\begin{array}{ccc} & 1 & \\g &  & g \\  & 1 & \end{array}\right. $$

\begin{rem}[Electric networks interpretation]
We may think of 
the tropical curve $\Gamma$ as an electric network, 
where each edge has resistance equal to its length. Note, 
in particular, that resistance 
is additive,
so it agrees with the length interpretation: if an edge is subdivided
into two smaller edges by a two-valent vertex, then its resistance
is the sum of two smaller resistances. 
Usual cohomology and homology of graphs have interpretations in terms of electrical circuits,
and tropical 
cohomology groups 
provide an even better framework for such interpretation.
The measure of the magnitude and direction of a stationary electrical current flowing through 
$\Gamma$
can be viewed as  
 a $(1,0)$-cochain  
 $I \in C^{1,0}$.
Indeed, for a point  $x\in\Gamma$ that is not a vertex
and a unit tangent vector $u$ at $x$ (an element of $\FF_1(x)$)
we may insert the ammeter at $x$ in the direction of $u$ and measure the current.
It follows from Kirchhoff's current law (divergence-free current)  that 
such a cochain 
$I$ is a cocycle.

Similarly, with a voltmeter we can  measure voltage
between any two points of $\Gamma$. Measuring it at the endpoints
of an oriented edge 
gives us a 1-cochain 
$V_0\in C^{0, 1}(\Gamma) = C^1(\Gamma;\R)$.
Note that $V_0$ must 
be a  coboundary 
by Kirchhoff's voltage law, 
which  also implies it is a 
cocycle.
This reflects the fact that no stationary electric current can be present
in such networks.  Since energy dissipates in resistance, to support
a stationary current we need power sources in our network. 
We may think of these power elements to be localized in some edges.
For an edge $E\subset\Gamma$ we define $V(E)$ to be $V_0(E)$ minus the voltage
of the power element (taken with sign).

Let $E\subset\Gamma$ be an oriented edge disjoint from power elements. Suppose that
$R$ is the resistance of $E$.
Let $I$ be a current through $E$
and $V$ be the voltage
drop at the endpoints of $E$.
Recall Ohm's law: $V=RI$. 
This may be interpreted through the so-called
{\em eigenwave} of the tropical variety $\Gamma$,
    see \cite{MikZhar:Eigenwave}.
This is a particular
1-cochain $\Phi$ with the coefficients
in $\WW_1$, the constructible sheaf with 
$\WW_1(x)=T_x(E)$ if $x$ is a point inside
an edge $E$ and $\WW_1(x)=0$ if $x$ is a vertex
of $\Gamma$ (which we assume to have valence
        greater than 2).
The eigenwave $\Phi$ is determined by 
\begin{equation}\label{phi-eigen}
\Phi(E)=\int\limits_{E}\omega,
\end{equation}
where $E\subset\Gamma$ is an edge 
and $\omega\in \operatorname{Hom}(\WW_1(E),\RR)=T_x^*(E)$
is the coefficient at the 1-simplex $E$.
In \eqref{phi-eigen} we interpret $\omega$
as the (constant) 1-form on the edge $E$.    
Recall that the length of $E$ represents the resistance.

In this way Ohm's law comes as
taking the cup-product with $\Phi$ producing
the homomorphism
$$H^{1,0}(\Gamma)\to H^{0,1}(\Gamma)$$
responsible for the definition of 
the Jacobian of $\Gamma$, see \cite{Mik6}.

\end{rem}

\begin{exo}
\

\begin{enumerate}
\item Consider  $\RR_{>0}$ equipped with the 
inclusion chart 
$\RR_{>0}\hookrightarrow\RR$. Prove that there is no way to
complete this atlas to turn $\RR_{>0}$ into a tropical manifold.

\item Revise Example \ref{ex:parameterization} by showing that three tropical modifications 
of the circle from Figure \ref{parametrization}a suffice to find a morphism
to a cubic curve in $\R^2$.

\item For each example of
computation of tropical homology and cohomology
performed in Section \ref{tropical_homology}, 
compute the corresponding tropical integral homology and cohomology (see Remark \ref{rem:Z hom}). 

\item Let $X$ be the Klein bottle obtained by
quotienting the unit square $[0,1]^2\subset\RR^2$ by
  the subgroup of $O_2(\RR)$ generated by $(x,y)\mapsto (x,y+1)$ and
  $(x,y)\mapsto (x+1,-y+1)$. Show that $X$ inherits a tropical structure
  from the tropical structure of $\RR^2$, and compute tropical (real
  and integral) homology and
  cohomology of $X$.

\item  Compute tropical homology and cohomology of $\tp^n$. In
  particular, show that the tropical cohomology ring of  $\tp^n$ is
  isomorphic to the cohomology ring of $\CC P^n$.

\item Find explicit generators of the homology and cohomology groups of
  a compact tropical curve.

\item Consider two  tropical curves $\Gamma$ and $\Gamma'$ such that
  $\Gamma'$ is a tropical modification of $\Gamma$. Show that the
  groups $H_{p,q}(\Gamma)$ and $H_{p,q}(\Gamma')$ are canonically
  isomorphic, as well as the groups $H^{p,q}(\Gamma)$ and $H^{p,q}(\Gamma')$.

\item Express the energy loss (the Joule-Lenz law) in terms of the cup
product in tropical cohomology (followed by evaluation on the fundamental
class) on an electric network considered as a tropical curve.

\item Let  $X$ be the product of two compact tropical curves of genus 1
  without 1-valent vertices ({\it i.e.} as in Figure \ref{parametrization}a). In
  particular $X$ is diffeomorphic to $S^1\times S^1$.
 Show that all
$(1,1)$-cycles in $X$ are straight. Construct a deformation
of $X$ ({\it i.e.}, a tropical 3-fold $Y$ and a tropical map $h: Y \to \R$ 
with $X=h^{-1}(0)$) such that the general fiber $Z$ does not have
a single straight $(1,1)$-cycle. 
Deduce that if $Z'$ is a connected tropical manifold that admits a
non-constant map $Z'\to Z$ which is given by affine maps in charts,
then $Z'$ is not projective ({\it i.e.} is not embeddable to $\tp^n$).

\end{enumerate}

\end{exo}